\documentclass{amsart}
\usepackage{amssymb}
\usepackage{mathrsfs}
\usepackage{stmaryrd}
\usepackage{bbm}
\usepackage{oldgerm}
\usepackage[frenchb,english]{babel}
\usepackage[utf8]{inputenc}
\usepackage[T1]{fontenc}
\usepackage{lmodern}
\usepackage[all]{xy}
\usepackage{hyperref}
\usepackage{upref}

\def\oo{\mathfrak{o}}
\def\mm{\mathfrak{m}}
\def\BB{\mathfrak{V}}
\def\XX{\mathfrak{X}}

\def\Mod{\mathbf{Mod}}
\def\Ens{\mathbf{Ens}}
\def\Rep{\mathbf{Rep}}
\def\Ob{\mathbf{Ob}}

\def\HB{\mathbf{B}}
\def\FHom{\mathscr{H}om}
\def\FExt{\mathscr{E}xt}
\def\FTor{\mathscr{T}or}
\def\rT{\mathsf{T}}

\def\bfD{\mathbf{D}}
\def\Ab{\mathbf{Ab}}
\def\Pt{\mathbf{Pt}}
\def\Sysl{\mathbf{LL}^{\textnormal{tf}}}
\def\SyWT{\Rep^{\WT}_{\mathfrak{C}}(\pi_1(C,\overline{x}))}
\def\LPtf{\mathbf{LP}^{\textnormal{tf}}}
\def\Alg{\mathbf{Alg}}
\def\RHom{\mathbf{R}\textnormal{Hom}}

\makeatletter
\def\rightharpoonfill@{\arrowfill@\relbar\relbar\rightharpoonup}
\DeclareRobustCommand{\overrightharpoon}{\mathpalette{\overarrow@\rightharpoonfill@}}
\makeatother

\DeclareMathOperator{\rmp}{p}
\DeclareMathOperator{\ltf}{ltf}
\DeclareMathOperator{\ptf}{ptf}
\DeclareMathOperator{\aptf}{\alpha ptf}
\DeclareMathOperator{\atf}{atf}
\DeclareMathOperator{\sat}{sat}
\DeclareMathOperator{\tf}{tf}
\DeclareMathOperator{\sol}{sol}
\DeclareMathOperator{\Dolb}{Dolb}
\DeclareMathOperator{\lt}{pltf}
\DeclareMathOperator{\op}{op}
\DeclareMathOperator{\coh}{coh}

\DeclareMathOperator{\scoh}{scoh}
\DeclareMathOperator{\id}{id}
\DeclareMathOperator{\Et}{\textnormal{\textbf{\'{E}t}}}
\DeclareMathOperator{\et}{\textnormal{\'{e}t}}
\DeclareMathOperator{\fet}{\textnormal{f\'{e}t}}
\DeclareMathOperator{\zar}{zar}

\DeclareMathOperator{\Aut}{Aut}
\DeclareMathOperator{\Pic}{Pic}

\DeclareMathOperator{\Ker}{Ker}
\DeclareMathOperator{\Coker}{Coker}
\DeclareMathOperator{\Image}{Im}
\DeclareMathOperator{\Coim}{Coim}
\DeclareMathOperator{\Hom}{Hom}

\DeclareMathOperator{\End}{End}
\DeclareMathOperator{\Ext}{Ext}
\DeclareMathOperator{\Tor}{Tor}

\DeclareMathOperator{\Spec}{Spec}
\DeclareMathOperator{\Spf}{Spf}
\DeclareMathOperator{\Vect}{Vect}
\DeclareMathOperator{\cont}{cont}
\DeclareMathOperator{\rH}{H}
\DeclareMathOperator{\rB}{B}
\DeclareMathOperator{\rZ}{Z}

\DeclareMathOperator{\rE}{E}
\DeclareMathOperator{\rR}{R}

\DeclareMathOperator{\rL}{L}
\DeclareMathOperator{\rI}{I}
\DeclareMathOperator{\rII}{II}
\DeclareMathOperator{\Fr}{Fr}

\DeclareMathOperator{\WT}{WT}
\DeclareMathOperator{\DW}{DW}
\DeclareMathOperator{\rig}{rig}
\DeclareMathOperator{\an}{an}
\DeclareMathOperator{\ad}{ad}
\DeclareMathOperator{\tor}{tor}

\newtheorem{theorem}{Théorème}[section]
\newtheorem{prop}[theorem]{Proposition}
\newtheorem{lemma}[theorem]{Lemme}

\newtheorem{coro}[theorem]{Corollaire}

\theoremstyle{definition}
\newtheorem{rem}[theorem]{Remarque}
\newtheorem{definition}[theorem]{Définition}
\newtheorem{nothing}[theorem]{}

\numberwithin{equation}{section}
\numberwithin{equation}{theorem}

\title{Transport parallèle et correspondance de Simpson $p$-adique}
\author{Daxin Xu}
\date{}
\AtEndDocument{\bigskip{\footnotesize%
  \textsc{Daxin Xu, Laboratoire Alexander Grothendieck, ERL 9216 du CNRS, Institut des Hautes Études Scientifiques, 35 route de Chartres, 91440 Bures-sur-Yvette, France.} \par
  \textit{E-mail address}: \texttt{daxinxu@ihes.fr} \par
}}

\begin{document}
\maketitle
\selectlanguage{english}
\begin{abstract}
	Deninger and Werner developed an analogue for $p$-adic curves of the classical correspondence of Narasimhan and Seshadri between stable bundles of degree zero and unitary representations of the topological fundamental group for a complex smooth proper curve. Using parallel transport, they associated functorially to every vector bundle on a $p$-adic curve whose reduction is strongly semi-stable of degree $0$ a $p$-adic representation of the fundamental group of the curve. They asked several questions: whether their functor is fully faithful; whether the cohomology of the local systems produced by this functor admits a Hodge-Tate filtration; and whether their construction is compatible with the $p$-adic Simpson correspondence developed by Faltings. We answer these questions in this article.	
\end{abstract}
\selectlanguage{french}
\begin{abstract}
	Deninger et Werner ont développé un analogue pour les courbes $p$-adiques de la correspondance classique de Narasimhan et Seshadri entre les fibrés vectoriels stables de degré zéro et les représentations unitaires du groupe fondamental topologique pour une courbe complexe propre et lisse. Par transport parallèle, ils ont associé fonctoriellement à chaque fibré vectoriel sur une courbe $p$-adique dont la réduction est fortement semi-stable de degré $0$ une représentation $p$-adique du groupe fondamental de la courbe. Ils se sont posés quelques questions: si leur foncteur est pleinement fidèle; si la cohomologie des systèmes locaux fournis par leur foncteur admet une filtration de Hodge-Tate; et si leur construction est compatible avec la correspondance de Simpson $p$-adique développée par Faltings. Nous répondons à ces questions dans cet article.
\end{abstract}
\section{Introduction} \label{Introduction}
\begin{nothing}
	Narasimhan et Seshadri \cite{NS65} ont établi en 1965 une correspondance bijective entre l'ensemble des classes d'équivalence de représentations unitaires irréductibles du groupe fondamental topologique d'une surface de Riemann compacte $X$ de genre $\ge 2$ et l'ensemble des classes d'isomorphisme de fibrés vectoriels stables de degré $0$ sur $X$. Deninger et Werner ont développé récemment un analogue partiel pour les courbes $p$-adiques \cite{DW05}. Leur construction est basée sur la notion de fibré \textit{fortement semi-stable} sur une courbe en caractéristique $p>0$. Rappelons qu'un fibré vectoriel sur une courbe $C$ propre et lisse sur un corps algébriquement clos de caractéristique $p$ est dit fortement semi-stable si ses images inverses par toutes les puissances entières du Frobenius absolu de $C$ sont semi-stables.
\end{nothing}
\begin{nothing}
	Soient $K$ un corps de valuation discrète complet de caractéristique $0$ de corps résiduel une clôture algébrique d'un corps fini $\mathbb{F}_{p}$, $\overline{K}$ une clôture algébrique de $K$, $\mathfrak{C}$ le complété $p$-adique de $\overline{K}$. On note $\mathcal{O}_{K}$ (resp. $\mathcal{O}_{\overline{K}}$, resp. $\oo$) l'anneau de valuation de $K$ (resp. $\overline{K}$, resp. $\mathfrak{C}$). On pose $S=\Spec(\mathcal{O}_{K})$ et on note $s$ (resp. $\eta$) le point fermé (resp. générique) de $S$ et $\overline{\eta}$ le point géométrique de $S$ associé à $\overline{K}$.
	
	Soient $X$ une $S$-courbe plate et propre à fibre générique lisse et géométriquement connexe et $\overline{x}$ un point géométrique de $X_{\overline{\eta}}$. On dit qu'un fibré $\mathcal{F}$ sur $X\otimes_{\mathcal{O}_{K}}\oo$ est de \textit{Deninger-Werner} si l'image inverse de $\mathcal{F}$ sur la normalisation de chaque composante irréductible de $X_{s}$ est fortement semi-stable de degré $0$. \`{A} un tel fibré $\mathcal{F}$, Deninger et Werner associent une représentation du groupe fondamental $\pi_{1}(X_{\overline{\eta}},\overline{x})$ définie comme suit.
	
	Pour tout entier $n\ge 1$, quitte à remplacer $K$ par une extension finie, il existe un morphisme propre $\varphi:X'\to X$ tel que $\varphi_{\eta}$ soit étale fini, que $X'$ soit une $S$-courbe semi-stable et que l'image réciproque de $\mathcal{F}_{n}=\mathcal{F}/p^{n}\mathcal{F}$ par la réduction modulo $p^{n}$ de $\varphi$ soit triviale \eqref{DW fortement semistable}. Le ``transport parallèle'' permet alors de construire une représentation de $\pi_{1}(X_{\overline{\eta}},\overline{x})$ sur la fibre de $\mathcal{F}_{n}$ en $\overline{x}$ (cf. \ref{definition of functor rho}). Par passage à la limite projective, on obtient une $\oo$-représentation continue $p$-adique de $\pi_{1}(X_{\overline{\eta}},\overline{x})$.
\end{nothing}
\begin{nothing}
	Soient $C$ une $\overline{K}$-courbe propre, lisse et connexe, $\check{C}=C\otimes_{\overline{K}}\mathfrak{C}$ et $\overline{x}$ un point géométrique de $C$. On dit qu'un fibré $F$ sur $\check{C}$ est de \textit{Deninger-Werner} si, quitte à remplacer $K$ par une extension finie, il existe un $S$-modèle propre et plat $X$ de $C$ et un fibré vectoriel de Deninger-Werner $\mathcal{F}$ sur $X\otimes_{\mathcal{O}_{K}}\oo$ de fibre générique $F$. On note $\BB_{\check{C}}^{\DW}$ la catégorie de tels fibrés et $\Rep_{\mathfrak{C}}^{\cont}(\pi_{1}(C,\overline{x}))$ la catégorie des $\mathfrak{C}$-représentations continues $p$-adiques de $\pi_{1}(C,\overline{x})$ sur des $\mathfrak{C}$-espaces vectoriels de dimension finie (cf. \ref{oo oon reprentations general}). La construction de Deninger-Werner induit alors un foncteur \eqref{rho C check}
	\begin{equation}
		\mathbb{V}:\BB_{\check{C}}^{\DW}\to \Rep_{\mathfrak{C}}^{\cont}(\pi_{1}(C,\overline{x})).
		\label{DW foncteur intro}
	\end{equation}

	Inspirés par le cas complexe, Deninger et Werner ont posé quelques questions concernant ce foncteur. Ils se sont demandés s'il est pleinement fidèle; si la cohomologie des systèmes locaux fournis par leur foncteur admet une filtration de Hodge-Tate; et si leur construction est compatible avec la correspondance de Simpson $p$-adique introduite par Faltings \cite{Fal05}. Le but de cet article est de répondre à ces questions.
\end{nothing}

\begin{nothing}
	Simultanément, Faltings a développé une correspondance pour les systèmes locaux $p$-adiques sur les variétés sur des corps $p$-adiques, inspiré par les travaux de Simpson dans le cas complexe. Rappelons que ce dernier a étendu le résultat de Narasimhan et Seshadri aux représentations linéaires quelconques du groupe fondamental topologique d'une variété complexe projective et lisse. Pour ce faire, on a besoin de la notion de \textit{fibré de Higgs}: si $X$ est un schéma lisse de type fini sur un corps $F$, un fibré de Higgs sur $X$ est un couple $(M, \theta)$ formé d'un fibré vectoriel $M$ sur $X$ et d'un morphisme $\mathscr{O}_{X}$-linéaire $\theta: M\to M\otimes_{\mathscr{O}_{X}}\Omega_{X/F}^{1}$ tel que $\theta \wedge \theta =0$ (cf. \ref{topos module higgs}). Le résultat principal de Simpson \cite{Sim92} établit une équivalence de catégories entre la catégorie des représentations linéaires de dimension finie (à valeurs complexes) du groupe fondamental topologique d'une variété complexe projective et lisse $X$ et celle des fibrés de Higgs semi-stables de classes de Chern nulles sur $X$.
\end{nothing}
\begin{nothing}
	La construction de Faltings, appelée classiquement \textit{correspondance de Simpson $p$-adique}, utilise sa théorie des extensions presque-étales et prolonge ses travaux en théorie de Hodge $p$-adique, en particulier ceux qui concernent la décomposition de Hodge-Tate pour une variété propre et lisse sur un corps $p$-adique. L'objet principal d'étude est la notion de \textit{représentation généralisée}, qui étend celle de représentation $p$-adique du groupe fondamental géométrique. Ce sont, en termes simplifiés, des représentations semi-linéaires $p$-adiques continues du groupe fondamental géométrique dans des modules sur un certain anneau $p$-adique muni d'une action continue du groupe fondamental géométrique. Faltings construit un foncteur de la catégorie de ces représentations dans la catégories des fibrés de Higgs \cite{Fal05}. Nous utiliserons la variante développée par Abbes et Gros \cite{AGT}, qui s'applique à une classe de représentations généralisées vérifiant une condition d'admissibilité à la Fontaine, dites de \textit{Dolbeault}. Nous introduisons dans le même esprit une autre condition d'admissibilité, plus forte que celle de Dolbeault, que nous qualifions de \textit{Weil-Tate}. Elle correspond aux fibrés de Higgs à champ de Higgs nul dans la correspondance de Simpson $p$-adique.
	
	Par ailleurs, en s'inspirant de la construction de Deninger-Werner, nous associons à certaines représentations généralisées, qualifiées de \textit{potentiellement libres de type fini}, des représentations continues $p$-adiques du groupe fondamental géométrique. Les fibrés de Deninger-Werner définissent naturellement des représentations généralisées potentiellement libres de type fini, et on retrouve ainsi le foncteur de Deninger-Werner. Les représentations généralisées de Weil-Tate sont aussi potentiellement libres de type fini. La principale question est de comparer ces deux sous-catégories. Notre principal résultat est que les représentations généralisées fournies par les fibrés de Deninger-Werner sont de Weil-Tate. Ce résultat nous permet de répondre aux trois questions de Deninger et Werner.
\end{nothing}
\begin{nothing}	\label{Faltings topos intro}
	Pour définir la notion de représentation généralisée, nous avons besoin du \textit{topos de Faltings}. Soit $X$ un $S$-schéma de type fini à réduction semi-stable tel que $X_{\overline{\eta}}$ soit connexe. On désigne par $E$ la catégorie des morphismes de schémas $V\to U$ au-dessus du morphisme canonique $X_{\overline{\eta}}\to X$ tels que le morphisme $U\to X$ soit étale et que le morphisme $V\to U_{\overline{\eta}}$ soit fini étale \eqref{basic topos de Faltings}. On équipe $E$ de \textit{la topologie co-évanescente} engendrée par les recouvrements $\{(V_{i}\to U_{i})\to (V\to U)\}_{i\in I}$ des deux types suivants:
\begin{itemize}
	\item[(v)] $U_i = U$ pour tout $i\in I$, et $(V_{i}\to V)_{i\in I}$ est un recouvrement.
	\item[(c)] $(U_{i}\to U)_{i\in I}$ est un recouvrement et $V_i = U_i \times_{U} V$ pour tout $i\in I$.
\end{itemize}
Le site $E$ ainsi défini est appelé \textit{site de Faltings de $X$}. On désigne par $\widetilde{E}$ et l'on appelle \textit{topos de Faltings de $X$}, le topos des faisceaux d'ensembles sur $E$.

Pour tout schéma $Y$, on note $\Et_{/Y}$ (resp. $\Et_{\textnormal{f}/Y}$) la catégorie des schémas étales (resp. finis et étales) au-dessus de $Y$ munie de la topologie étale, et $Y_{\et}$ (resp. $Y_{\fet}$) le topos des faisceaux d'ensembles sur $\Et_{/Y}$ (resp. $\Et_{\textnormal{f}/Y}$). Les foncteurs
\begin{eqnarray}
	\Et_{\textnormal{f}/X_{\overline{\eta}}}\to E \qquad V\mapsto (V\to X)\\
	\Et_{/X}\to E\qquad U\mapsto (U_{\overline{\eta}}\to U)
\end{eqnarray}
sont continus et exacts à gauches. Ils définissent donc deux morphismes de topos
\begin{eqnarray}
	\beta:\widetilde{E}&\to& X_{\overline{\eta},\fet}, \label{beta intro}\\
	\sigma: \widetilde{E}&\to& X_{\et}. \label{sigma intro}
\end{eqnarray}

Pour chaque objet $(V\to U)$ de $E$, on note $\overline{U}^{V}$ la clôture intégrale de $\overline{U}=U\otimes_{\mathcal{O}_{K}}\mathcal{O}_{\overline{K}}$ dans $V$. On désigne par $\overline{\mathscr{B}}$ le préfaisceau d'anneaux sur $E$ défini pour tout $(V\to U)\in \Ob(E)$ par
\begin{equation}
	\overline{\mathscr{B}}(V\to U)=\Gamma(\overline{U}^{V},\mathscr{O}_{\overline{U}^{V}}).
\end{equation}
En fait, celui-ci est un faisceau sur $E$ (\cite{AGT} III.8.16). Les représentations généralisées sont essentiellement les $\overline{\mathscr{B}}$-modules de $\widetilde{E}$. Cependant, pour tenir compte de la topologie $p$-adique, nous travaillons avec le système projectif des anneaux $\overline{\mathscr{B}}_{n}=\overline{\mathscr{B}}/p^{n}\overline{\mathscr{B}}$, $n\ge 1$.

Pour tout entier $n\ge 1$, $\overline{\mathscr{B}}_{n}$ est un objet de la fibre spéciale $\widetilde{E}_{s}$ de $\widetilde{E}$, c'est-à-dire du sous-topos fermé de $\widetilde{E}$ complémentaire de l'ouvert $\sigma^{*}(X_{\eta})$ \eqref{sous-topos ferme}. Celui-ci s'insère dans un diagramme commutatif à isomorphisme canonique près
\begin{equation}
	\xymatrix{
		\widetilde{E}_{s}\ar[r]^{\delta} \ar[d]_{\sigma_{s}} & \widetilde{E} \ar[d]^{\sigma} \\
		X_{s,\et} \ar[r]^{a} & X_{\et}
	}
\end{equation}
où $a$ est l'injection canonique, $\delta$ est le plongement canonique et $\sigma_{s}$ est induit par $\sigma$.

Les systèmes projectifs d'objets de $\widetilde{E}_{s}$ indexés par l'ensemble ordonné des entiers naturels $\mathbb{N}$, forment un topos que l'on note $\widetilde{E}_{s}^{\mathbb{N}^{\circ}}$. On le munit de l'anneau $\breve{\overline{\mathscr{B}}}=(\overline{\mathscr{B}}_{n})_{n\ge 1}$ de $\widetilde{E}_{s}^{\mathbb{N}^{\circ}}$. On dit qu'un $\breve{\overline{\mathscr{B}}}$-module $(M_{n})_{n\ge 1}$ est \textit{adique} si, pour tout entier $i\ge 1$, le morphisme $M_{i+1}\otimes_{\overline{\mathscr{B}}_{i+1}}\overline{\mathscr{B}}_{i}\to M_{i}$, déduit du morphisme de transition $M_{i+1}\to M_{i}$, est un isomorphisme.


On note $\Mod_{\mathbb{Q}}(\breve{\overline{\mathscr{B}}})$ la catégorie des $\breve{\overline{\mathscr{B}}}$-modules à isogénie près \eqref{categorie a isogenies}. Si $X$ est propre sur $S$, d'après le principal théorème de comparaison de Faltings \eqref{thm iso de Faltings} (cf. \cite{Fal02} Theorem 8, \cite{AG15} 2.4.16), l'image réciproque associée au morphisme $\beta$ induit un foncteur pleinement fidèle \eqref{beta pullback pleinement fidele}
	\begin{equation}
		\Rep_{\mathfrak{C}}^{\cont}(\pi_{1}(X_{\overline{\eta}},\overline{x}))\to \Mod_{\mathbb{Q}}(\breve{\overline{\mathscr{B}}}).
		\label{C Rep vers B modules}
	\end{equation}
\end{nothing}


\begin{nothing}
	Supposons que $X$ soit une $S$-courbe semi-stable et posons $C=X_{\overline{\eta}}$. Soient $\overline{x}$ un point géométrique de $C$ et $n\ge 1$ un entier. Si $X'$ est une $S$-courbe semi-stable, $(\widetilde{E}',\overline{\mathscr{B}}')$ le topos annelé de Faltings associé à $X'$ et $\varphi:X'\to X$ un $S$-morphisme propre, alors $\varphi$ induit par fonctorialité un morphisme de topos annelés $\Phi_{n}:(\widetilde{E}_{s}',\overline{\mathscr{B}}'_{n})\to (\widetilde{E}_{s},\overline{\mathscr{B}}_{n})$. On dit qu'un $\overline{\mathscr{B}}_{n}$-module $M_{n}$ est \textit{potentiellement libre de type fini} s'il est de type fini et si, quitte à remplacer $K$ par une extension finie, il existe un $S$-morphisme propre $\varphi:X'\to X$ tel que $X'$ soit une $S$-courbe semi-stable, $\varphi_{\eta}$ soit fini étale et qu'avec les notations précédentes $\Phi_{n}^{*}(M_{n})$ soit un $\overline{\mathscr{B}}_{n}'$-module libre. Par ``transport parallèle'', on associe à tout $\overline{\mathscr{B}}_{n}$-module potentiellement libre de type fini une représentation de $\pi_{1}(C,\overline{x})$ \eqref{construction of varrho}.
	
	On dit qu'un $\breve{\overline{\mathscr{B}}}$-module $M=(M_{n})_{n\ge 1}$ est \textit{potentiellement libre de type fini} s'il est adique de type fini et si, pour tout entier $n\ge 1$, $M_{n}$ est potentiellement libre de type fini. On désigne par $\Mod^{\lt}_{\mathbb{Q}}(\breve{\overline{\mathscr{B}}})$ la catégorie des $\breve{\overline{\mathscr{B}}}$-modules potentiellement libres de type fini à isogénie près. La construction précédente induit un foncteur \eqref{varrho entier}
	\begin{equation}
		\mathscr{W}: \Mod^{\lt}_{\mathbb{Q}}(\breve{\overline{\mathscr{B}}})\to \Rep^{\cont}_{\mathfrak{C}}(\pi_{1}(C,\overline{x})).
		\label{foncteur W intro}
	\end{equation}
	Celui-ci généralise le foncteur $\mathbb{V}$ \eqref{DW foncteur intro}. En effet, notons $\Vect_{\check{C}}$ la catégorie des fibrés vectoriels sur $\check{C}=C\otimes_{\overline{K}}\mathfrak{C}$. L'image réciproque associée au morphisme $\sigma$ induit un foncteur (cf. \eqref{jmathXX coh} et \eqref{morphisme de topos T})
	\begin{eqnarray}
		\Vect_{\check{C}} \to  \Mod_{\mathbb{Q}}(\breve{\overline{\mathscr{B}}}).
		\label{Vect bundle vers B modules}
	\end{eqnarray}
	La restriction de ce dernier à $\BB_{\check{C}}^{\DW}$ se factorise à travers $\Mod^{\lt}_{\mathbb{Q}}(\breve{\overline{\mathscr{B}}})$. Le foncteur $\mathscr{W}$ induit alors le foncteur $\mathbb{V}$ \eqref{compatible of 2 rho}.
	
	Par ailleurs, le foncteur \eqref{C Rep vers B modules} se factorise à travers la sous-catégorie $\Mod^{\lt}_{\mathbb{Q}}(\breve{\overline{\mathscr{B}}})$ de $\Mod_{\mathbb{Q}}(\breve{\overline{\mathscr{B}}})$. Le composé de ce dernier et de $\mathscr{W}$ est isomorphe au foncteur identique (\ref{rep generalise rep vrai}(ii)).

	On dit qu'un fibré vectoriel sur $\check{C}$ et une $\mathfrak{C}$-représentation continue $V$ de $\pi_{1}(C,\overline{x})$ sont \textit{$\breve{\overline{\mathscr{B}}}_{\mathbb{Q}}$-associés} si leurs images par les foncteurs \eqref{Vect bundle vers B modules} et \eqref{C Rep vers B modules} dans $\Mod_{\mathbb{Q}}(\breve{\overline{\mathscr{B}}})$ sont isomorphes (cf. \ref{B associe}).
	\label{intro DW via Faltings}
\end{nothing}
\begin{nothing}
	Soient $C$ une courbe propre et lisse sur $\overline{K}$, $\overline{x}$ un point géométrique de $C$ et $\check{C}=C\otimes_{\overline{K}}\mathfrak{C}$. On dit qu'un fibré vectoriel $F$ sur $\check{C}$ est de \textit{Weil-Tate} si, quitte à remplacer $K$ par une extension finie, il existe un $S$-modèle semi-stable $X$ de $C$ et une $\mathfrak{C}$-représentation continue $V$ de $\pi_{1}(C,\overline{x})$ tels que $F$ et $V$ soient $\breve{\overline{\mathscr{B}}}_{\mathbb{Q}}$-associés dans le topos annelé de Faltings relatif à $X$. On définit symétriquement la notion de \textit{$\mathfrak{C}$-représentation continue de Weil-Tate de $\pi_{1}(C,\overline{x})$}. On note $\BB_{\check{C}}^{\WT}$ (resp. $\Rep_{\mathfrak{C}}^{\WT}(\pi_{1}(C,\overline{x}))$) la catégorie de tels fibrés (resp. telles $\mathfrak{C}$-représentations). On montre le résultat suivant:
\end{nothing}
\begin{theorem}[cf. \ref{theorem WT}] \label{thm WT equi intro}
	Il existe des équivalences de catégories quasi-inverses l'une de l'autre
	\begin{equation}
		\mathscr{V}: \BB_{\check{C}}^{\WT}\to \Rep^{\WT}_{\mathfrak{C}}(\pi_{1}(C,\overline{x}))\qquad\textnormal{et}\qquad \mathscr{T}:\Rep^{\WT}_{\mathfrak{C}}(\pi_{1}(C,\overline{x}))\to \BB_{\check{C}}^{\WT}.
		\label{V T intro}
	\end{equation}
\end{theorem}

Notre principal résultat est le suivant.

\begin{theorem}[cf. \ref{DW implique WT coro}]
	Tout fibré vectoriel de Deninger-Werner sur $\check{C}$ est de Weil-Tate.
	\label{thm DW implique WT intro}
\end{theorem}

La condition de Deninger-Werner  peut se voir comme une condition d'admissibilité à la Fontaine. Elle présente toutefois des différences notables avec la condition d'admissibilité de Weil-Tate. D'une part, le fibré est trivialisé par une courbe semi-stable au-dessus de $X$ qui est finie et étale au-dessus de $C$, mais qui n'est a priori pas finie au dessus de $X$. D'autre part, on a besoin d'une infinité de tels revêtements, un pour chaque réduction modulo $p^{n}$, $n\ge 1$. La première difficulté peut être surmontée grâce à un résultat de Raynaud sur le quotient d'une courbe semi-stable par un groupe fini (\ref{domine par bon revêtement}(iii) cf. aussi \cite{Ray90}). La seconde difficulté est plus sérieuse. On la surmonte grâce à la théorie des déformations des presque-modules dans le sens de Faltings \eqref{theorem principal deformation alpha modules}.

\begin{prop}[cf. \ref{foncteur DW WT pf}]
	La restriction du foncteur $\mathscr{V}$ à $\BB_{\check{C}}^{\DW}$ s'identifie au foncteur $\mathbb{V}$.
\end{prop}

	On en déduit la pleine fidélité du foncteur de Deninger-Werner \eqref{DW foncteur intro}.

\begin{nothing}
	On démontre que pour qu'une $\mathfrak{C}$-représentation continue $V$ de $\pi_{1}(C,\overline{x})$ soit de Weil-Tate, il faut et il suffit que son image, par le foncteur \eqref{C Rep vers B modules}, soit de Dolbeault et que le fibré de Higgs associé par la correspondance de Simpson $p$-adique	soit muni du champ de Higgs nul (cf. \ref{Higgs nul et Weil-Tate}). De plus, le fibré vectoriel sous-jacent à ce fibré de Higgs est donné par $\mathscr{T}(V)$ \eqref{V T intro}. Cela signifie que l'équivalence de catégories $\mathscr{T}$ \eqref{V T intro} n'est autre que la correspondance de Simpson $p$-adique pour les fibrés de Higgs à champ de Higgs nul. Dans \ref{HT ss prop}, on démontre l'existence d'une filtration de Hodge-Tate pour la cohomologie des représentations de Weil-Tate.
\end{nothing}

\begin{nothing}
	Après avoir fixé les notations générales dans la section \ref{Not et pre}, nous développons dans les sections \ref{Alg homologique} à \ref{courbes rel} quelques préliminaires utiles pour la suite de cet article. Dans la section \ref{Alg homologique}, nous rappelons la structure des $\oo$-modules de type $\alpha$-fini et nous étudions les représentations à coefficients dans ces modules. Dans la section \ref{Geo rigide}, nous rappelons la fidélité du foncteur GAGA pour les modules cohérents sur un $\mathfrak{C}$-schéma de type fini. La section \ref{courbes rel} contient des rappels et des compléments sur les courbes relatives. Dans la section \ref{Cor DW}, nous rappelons les résultats de Deninger et Werner. Suivant (\cite{AGT} III et VI), nous présentons le topos annelé de Faltings dans la section \ref{Topos de Faltings}. Ensuite, nous rappelons le principal théorème de comparaison de Faltings et nous donnons quelques corollaires. Nous finissons cette section par la description du topos annelé de Faltings associé à un trait \eqref{S bar et et coevanescent}. Dans la section \ref{DW via Faltings}, nous présentons la construction du foncteur $\mathscr{W}$ \eqref{foncteur W intro} et ses propriétés. Nous établissons un résultat de descente galoisienne pour le topos de Faltings dans la section \ref{Descente Gal}. Grâce à ce résultat, nous démontrons une variante modulo $p^{n}$ du théorème \ref{thm DW implique WT intro} (cf. \ref{Mod Weil-Tate mod pn}). Dans la section \ref{FV Rep WT}, nous introduisons les notions de fibré vectoriel de Weil-Tate et de $\mathfrak{C}$-représentation continue de Weil-Tate. Nous construisons les foncteurs $\mathscr{V}$ et $\mathscr{T}$ \eqref{V T intro} et nous démontrons le théorème \ref{thm WT equi intro}. Dans la section \ref{Simpson p adique}, nous comparons la correspondance de Simpson $p$-adique et la correspondance $\mathscr{T}$ \eqref{V T intro}. Dans la section \ref{Faisceaux alpha}, nous rappelons la notion de \textit{faisceaux de $\alpha$-modules} suivant (\cite{AG15} 1.4). Nous développons ensuite une théorie de déformations des faisceaux de $\alpha$-modules suivant l'approche d'Illusie \cite{Il71} dans la section \ref{Deformation alpha}. Nous étudions la relation entre les déformations des représentations du groupe fondamental et les déformations des $\alpha$-$\overline{\mathscr{B}}$-modules et nous démontrons le théorème \ref{thm DW implique WT intro} dans la section \ref{Dem finale}.
\end{nothing}

\textbf{Remerciement.}
Ce travail fait partie de ma thèse préparée à l'université de Paris-Sud. Je remercie profondément mon directeur de thèse Ahmed Abbes pour les nombreuses discussions, ainsi que pour ses remarques et ses critiques sur ce texte. Je remercie Bhargav Bhatt de m'avoir signalé une erreur dans une version préliminaire de cet article. Je remercie également Michel Gros pour ses commentaires et suggestions. Ce travail a bénéficié du soutien du programme ANR Théorie de Hodge $p$-adique et Développements (ThéHopadD) ANR-11-BS01-005.

\tableofcontents

\section{Notations et préliminaires} \label{Not et pre}
\begin{nothing}
Dans cet article, $p$ désigne un nombre premier, $K$ un corps de valuation discrète complet de caractéristique $0$ dont le corps résiduel $k$ est \textit{une clôture algébrique de $\mathbb{F}_{p}$} et $\overline{K}$ une clôture algébrique de $K$. On note $\mathcal{O}_K$ l'anneau de valuation de $K$ et $\mathcal{O}_{\overline{K}}$ la clôture intégrale de $\mathcal{O}_K$ dans $\overline{K}$. On désigne par $\oo$ le séparé complété $p$-adique de $\mathcal{O}_{\overline{K}}$, par $\mathfrak{C}$ son corps des fractions, par $\mm$ son idéal maximal et par $v$ la valuation de $\oo$ normalisée par $v(p)=1$. Pour tout entier $n\ge 1$, on pose $\oo_n=\mathcal{O}_{\overline{K}}/p^n\mathcal{O}_{\overline{K}}$.

	On pose $S=\Spec(\mathcal{O}_K)$, $\overline{S}=\Spec(\mathcal{O}_{\overline{K}})$, $\check{\overline{S}}=\Spec(\oo)$, $\mathscr{S}=\Spf(\oo)$ et pour tout entier $n\ge 1$, $S_n=\Spec(\mathcal{O}_K/p^n\mathcal{O}_K)$. On note $\eta$ (resp. $s$) le point générique (resp. fermé) de $S$ et $\overline{\eta}$ (resp. $\check{\overline{\eta}}$) le point géométrique générique correspondant à $\overline{K}$ (resp. $\mathfrak{C}$). Pour tous $S$-schémas $X$ et $S'$, on pose $X_{S'}=X\times_{S}S'$. Pour tout entier $n\ge 1$, on pose
\begin{equation}
	X_{n}=X\times_{S}S_{n},\qquad \overline{X}=X\times_{S}\overline{S},\qquad \check{\overline{X}}=X\times_{S}\check{\overline{S}}.
	\label{basic notation}
\end{equation}
Pour tout morphisme de $S$-schémas $\pi:Y\to X$, on pose $\pi_n=\pi\times_{S}S_n$. Pour tout $\mathscr{O}_{X}$-module $M$ sur un $S$-schéma $X$, on note $M_n$ le $\mathscr{O}_{X_n}$-module $M\otimes_{\mathscr{O}_{S}}\mathscr{O}_{S_n}$ sur $X_n$.

Si $T$ est un trait (i.e. le spectre d'un anneau de valuation discrète) et $\tau$ son point générique, on dit que $(T,\tau)$ est un \textit{trait génériquement pointé} (ou simplement \textit{trait} lorsqu'il n'y a aucun risque d'ambiguïté). Si $(T,\tau)$ et $(T',\tau')$ sont deux traits génériquement pointés, un morphisme $(T',\tau')\to (T,\tau)$ est un morphisme dominant $T'\to T$. On dit qu'un tel morphisme est fini si le morphisme sous-jacent $T'\to T$ est fini.
\label{notations 11}
\end{nothing}
$\hspace*{-1.2em}\bf{\arabic{section}.\stepcounter{theorem}\arabic{theorem}.}$
On désigne par $\Ens$ la catégorie des ensembles que l'on considère ainsi comme un topos (ponctuel) et on le note aussi $\Pt$ (\cite{SGAIV} IV 2.2). On désigne par $\Ab$ la catégorie des groupes abéliens.

\begin{nothing}	\label{def Ext cat derivee}
	Soit $\mathscr{A}$ une catégorie abélienne et $M$, $N$ deux objets de $\mathscr{A}$. On désigne par $\mathbf{D}(\mathscr{A})$ sa catégorie dérivée. On note abusivement $M$ le complexe concentré en degré $0$ de valeur $M$. Pour tout entier $i$, on pose (\cite{Weib} 10.7.1)
	\begin{equation}
		\Ext^{i}_{\mathscr{A}}(M,N)=\Hom_{\mathbf{D}(\mathscr{A})}(M,N[i]),
		\label{def Ext general}
	\end{equation}
	où $[i]$ désigne le décalage de degré $i$.

	Soit $F:\mathscr{A}\to \mathscr{A}'$ un foncteur exact entre catégories abéliennes. Il s'étend en un foncteur de $\mathbf{D}(\mathscr{A})$ dans $\mathbf{D}(\mathscr{A}')$. En vertu de \eqref{def Ext general}, il induit, pour tous objets $M$ et $N$ de $\mathscr{A}$ et tout entier $i$, un morphisme canonique
	\begin{equation}
		\Ext^{i}_{\mathscr{A}}(M,N)\to \Ext^{i}_{\mathscr{A}'}(F(M),F(N)).
		\label{exact foncteur induit Ext morphisme}
	\end{equation}
	Pour tous objets $M$ et $N$ de $\mathscr{A}$, le groupe $\Ext^{1}_{\mathscr{A}}(M,N)$ classifie les extensions de $M$ par $N$ dans $\mathscr{A}$ (cf. \cite{Stacks} 13.27.6).
\end{nothing}

\begin{nothing}
	Soient $\mathscr{A}$ une catégorie abélienne ayant suffisamment d'injectifs et $M$ un objet de $\mathscr{A}$. Le foncteur $\Hom_{\mathscr{A}}(M,-):\mathscr{A}\to \Ab$ est exact à gauche. On désigne par
	\begin{equation}
		\Ext^{i}_{\mathscr{A}}(M,-):\mathscr{A} \to \Ab \qquad \forall i\ge 0
	\end{equation}
	les foncteurs dérivées à droite de $\Hom_{\mathscr{A}}(M,-)$. On pose $\Ext^{i}_{\mathscr{A}}(M,-)=0$ pour $i<0$. Cette définition est compatible avec \eqref{def Ext general} (cf. \cite{Weib} 10.7.4 et 10.7.5).
	\label{definition Ext using Hom derive}
\end{nothing}

\begin{nothing}	\label{partie suite spectrale}
	Soient $\mathscr{A}$ une catégorie abélienne et $\rE$ une suite spectrale dans $\mathscr{A}$
	\begin{equation}
		\rE^{i,j}_{2}\Rightarrow \rE^{i+j}
		\label{suite spectrale generale}
	\end{equation}
	telle que $\rE_{2}^{i,j}=0$ si $i<0$ ou $j<0$. Les termes de bas degré fournissent une suite exacte (\cite{Milne} Appendix B)
	\begin{equation}
		0\to \rE_{2}^{1,0}\to \rE^{1} \to \rE_{2}^{0,1}\to \rE_{2}^{2,0}\to \rE_{1}^{2} \to \rE_{2}^{1,1},
		\label{suite exacte terme petite}
	\end{equation}
	où $\rE_{1}^{2}=\Ker(\rE^{2}\to \rE^{0,2}_{2})$. Soient $\rE'$ une autre suite spectrale 
	\begin{equation}
		\rE'^{i,j}_{2}\Rightarrow \rE'^{i+j}	
	\end{equation}
	telle que $\rE'^{i,j}_{2}=0$ si $i<0$ ou $j<0$, et $u:\rE\to \rE'$ un morphisme de suites spectrales (\cite{EGA III} 0.11.1.2). On a alors un diagramme commutatif 
	\begin{equation}
		\xymatrix{
			0\ar[r]& \rE_{2}^{1,0}\ar[r] \ar[d]_{u^{1,0}_{2}} & \rE^{1} \ar[r] \ar[d]^{u^{1}}& \rE_{2}^{0,1}\ar[r] \ar[d]^{u^{0,1}_{2}}& \rE_{2}^{2,0}\ar[r] \ar[d]^{u_{2}^{2,0}} & \rE_{1}^{2} \ar[r]\ar[d]^{u_{1}^{2}}& \rE_{2}^{1,1} \ar[d]^{u_{2}^{1,1}}\\
			0\ar[r]& \rE'^{1,0}_{2}\ar[r]  & \rE'^{1} \ar[r]& \rE'^{0,1}_{2}\ar[r] & \rE'^{2,0}_{2}\ar[r] & \rE'^{2}_{1} \ar[r]& \rE'^{1,1}_{2}
		}
		\label{morphsme de suite spectrale petit}
	\end{equation}
\end{nothing}
\begin{lemma} \label{lemma suite spectrale}
	Conservons les notations de \ref{partie suite spectrale} et supposons que les morphismes $u^{i,j}_{2}$ sont isomorphes pour $i+j \le 1$ et sont injectifs pour $i+j=2$. Alors, le morphisme $u^{n}:\rE^{n}\to\rE'^{n}$ est un isomorphisme si $n=0,1$, et est un monomorphisme si $n=2$. 
\end{lemma}
\textit{Preuve}. Il est clair que $u^{0}$ est un isomorphisme. D'après \eqref{morphsme de suite spectrale petit}, le lemme des cinq et le lemme du serpent, on déduit que le morphisme $u^{1}$ est un isomorphisme et que le morphisme $u_{1}^{2}$ est un monomorphisme. Par suite, on en déduit un diagramme commutatif
\begin{equation}
	\xymatrix{
		0\ar[r]& \rE^{2}_{1} \ar[r]\ar@{^{(}->}[d]_{u_{1}^{2}}& \rE^{2} \ar[r] \ar[d]_{u^{2}}& \rE^{0,2}_{2} \ar@{^{(}->}[d]^{u_{2}^{0,2}} \\
		0\ar[r]& \rE'^{2}_{1} \ar[r]& \rE'^{2} \ar[r] & \rE'^{0,2}_{2},
	}
\end{equation}
d'où l'injectivité de $u^{2}$.

\begin{nothing}
	Pour tout schéma $X$, on note $\Pic(X)=\rH^{1}(X,\mathscr{O}_{X}^{\times})$ le groupe des classes d'isomorphismes de $\mathscr{O}_{X}$-modules inversibles. Soit $f:X\to T$ un morphisme propre de schémas. On désigne par $\Pic_{X/T}$ \textit{le foncteur de Picard relatif de $X$ au-dessus de $T$}, défini pour tout $T$-schéma $T'$ par
	\begin{equation}
		\Pic_{X/T}(T')=\rH^{0}(T',\rR^{1}_{\textnormal{fppf}}f'_{*}(\mathbb{G}_{m})),
		\label{prefaisceaux Pic}
	\end{equation}
	où $f':X\times_{T}T'\to T'$ est la projection canonique et $\mathbb{G}_{m}$ est le faisceau fppf qui à chaque schéma $Y$ associe le groupe $\Gamma(Y,\mathscr{O}_{Y}^{*})$ (cf. \cite{Ray70} 1.2). Soit $t$ un point de $T$. D'après \eqref{prefaisceaux Pic}, on a $\Pic_{X/T}(t)=\Pic_{X_{t}/t}(t)$. Le foncteur $\Pic_{X_{t}/t}$ est représentable par un schéma en groupes localement de type fini (cf. \cite{Ray70} 1.5.2). Il possède un plus petit sous-groupe ouvert connexe, sa composante neutre $\Pic_{X_{t}/t}^{0}$ (\cite{SGAIII} $\textnormal{VI}_{\textnormal{A}}$, 2). On définit le sous-foncteur $\Pic_{X/T}^{0}$ de $\Pic_{X/T}$ comme suit: pour tout $T$-schéma $T'$, $\Pic_{X/T}^{0}(T')$ est le sous-groupe de $\Pic_{X/T}(T')$ formé des éléments qui pour chaque point $t$ de $T'$ induisent par fonctorialité un élément de $\Pic_{X_{t}/t}^{0}(t)$.
\label{not of Picard}
\end{nothing}
\begin{nothing}
Soit $X$ un schéma. On dit qu'un $\mathscr{O}_{X}$-module quasi-cohérent $F$ est \textit{un fibré vectoriel} sur $X$ s'il est localement libre de type fini. On désigne par $\Vect_{X}$ la sous-catégorie pleine de la catégorie des $\mathscr{O}_{X}$-modules formée des fibrés vectoriels.
\label{notations vect bundle}
\end{nothing}


\begin{nothing} \label{LPtf}
	Soit $(\mathcal{T},A)$ un topos annelé. On désigne par $\Mod(\mathcal{T},A)$ la catégorie des $A$-modules de $\mathcal{T}$. On dit qu'un $A$-module $M$ de $\mathcal{T}$ est localement projectif de type fini si les conditions équivalentes suivantes sont satisfaites (cf. \cite{AGT} III.2.8):

(i) $M$ est de type fini et le foncteur $\FHom_{A}(M,\cdot)$ est exact;

(ii) $M$ est de type fini et tout épimorphisme de $A$-modules $N\to M$ admet localement une section;

(iii) $M$ est localement facteur direct d'un $A$-module libre de type fini.

On désigne par $\Sysl(\mathcal{T},A)$ (resp. $\LPtf(\mathcal{T},A)$) la sous-catégorie pleine de $\Mod(\mathcal{T},A)$ formée des $A$-modules localement libres (resp. localement projectifs) de type fini de $\mathcal{T}$.

Un $A$-module localement projectif de type fini est plat en vertu de (iii). Soient $M'$ et $M''$ deux $A$-modules localement projectifs de type fini. D'après (iii), le $A$-module $\FHom_{A}(M',M'')$ est aussi localement projectif de type fini. Soit $0\to M'\to M\to M''\to 0$ une suite exacte de $\Mod(\mathcal{T},A)$. On en déduit que, pour tout $i\ge 1$, le foncteur $\FExt^{i}_{A}(M,-)$ est nul. Donc, $M$ est aussi localement projectif de type fini.

\'{E}tant donné un $A$-module localement projectif de type fini $M$, la suite spectrale qui relie les Ext locaux et globaux (\cite{SGAIV} V 6.1) induit, pour tout $A$-module $N$ et tout entier $i\ge 0$, un isomorphisme
\begin{equation}
	\Ext^{i}_{A}(M,N)\xrightarrow{\sim}\rH^{i}(\mathcal{T},\FHom_{A}(M,N)).
	\label{Ext 1 iso to H1}
\end{equation}
\end{nothing}
\begin{nothing}\label{topos module higgs}
	Soient $(\mathcal{T},A)$ un topos annelé et $E$ un $A$-module de $\mathcal{T}$. Un \textit{$A$-module de Higgs à coefficients dans $E$} est un couple $(M, \theta)$ formé d'un $A$-module $M$ de $\mathcal{T}$ et d'un morphisme $A$-linéaire $\theta:M\to M\otimes_{A}E$ tel que $\theta\wedge \theta=0$.
\end{nothing}
%
%
%
%
\begin{nothing}
Soit $\mathcal{T}$ un topos. Les systèmes projectifs d'objets de $\mathcal{T}$ indexés par l'ensemble ordonné des entiers naturels $\mathbb{N}$, forment un topos que l'on note $\mathcal{T}^{\mathbb{N}^{\circ}}$ (cf. \cite{AGT} VI.7.1). Le foncteur
\begin{equation}
	\lambda^{*}:\mathcal{T}\to \mathcal{T}^{\mathbb{N}^{\circ}}
	\label{lambda pushforward}
\end{equation}
qui à un objet $F$ associe le foncteur constant $\mathbb{N}^{\circ}\to \mathcal{T}$ de valeur $F$ est exact à gauche. Il admet pour adjoint à droite le foncteur
\begin{equation}
	\lambda_{*}:\mathcal{T}^{\mathbb{N}^{\circ}}\to \mathcal{T}
	\label{lambda pullback}
\end{equation}
qui à un foncteur $\mathbb{N}^{\circ}\to \mathcal{T}$ associe sa limite projective (\cite{AGT} III.7.4). Le couple $(\lambda_{*}, \lambda^{*})$ définit donc un morphisme de topos
\begin{equation}
	\lambda: \mathcal{T}^{\mathbb{N}^{\circ}}\to \mathcal{T}.
	\label{morphisme limite proj}
\end{equation}

\'{E}tant donné un anneau $\breve{A}=(A_{n})_{n\ge 1}$ de $\mathcal{T}^{\mathbb{N}^{\circ}}$, on dit qu'un $\breve{A}$-module $M=(M_{n})_{n\ge 1}$ de $\mathcal{T}^{\mathbb{N}^{\circ}}$ est \textit{adique} si pour tous entiers $i$ et $j$ tels que $1\le i\le j$, le morphisme $M_j\otimes_{A_j}A_i\to M_i$ déduit du morphisme de transition $M_j\to M_i$ est un isomorphisme.
\label{limite projective de topos}
\end{nothing}
\begin{nothing}
	\'{E}tant donnés un morphisme de topos $\gamma:\mathcal{T}_1\to \mathcal{T}_3$ et un groupe abélien $\mathscr{F}$ de $\mathcal{T}_3$, on définit, pour tout entier $i\ge 0$, un morphisme
\begin{equation}
	\gamma^{*}: \rH^{i}(\mathcal{T}_{3},\mathscr{F}) \to \rH^i(\mathcal{T}_1,\gamma^*(\mathscr{F}))
	\label{cohomologie pullback par gamma}
\end{equation}
comme le composé
\begin{equation}
	\rH^i(\mathcal{T}_3,\mathscr{F})\to \rH^{i}(\mathcal{T}_{3},\gamma_{*}\gamma^{*}(\mathscr{F})) \to \rH^i(\mathcal{T}_1,\gamma^*(\mathscr{F})),
	\label{cohomologie pullback par gamma 2}
\end{equation}
où la première flèche est induite par le morphisme d'adjonction $\id \to \gamma_{*}\gamma^{*}$ et la seconde flèche est induite par la suite spectrale de Cartan-Leray.

Supposons que $\gamma$ soit le composé de deux morphismes de topos $\mathcal{T}_1\xrightarrow{\alpha}\mathcal{T}_2\xrightarrow{\beta}\mathcal{T}_3$. Alors, pour tout entier $i\ge 0$ et tout faisceau abélien $\mathscr{F}$ de $\mathcal{T}_3$, le morphisme composé
\begin{equation}
	\rH^i(\mathcal{T}_3,\mathscr{F})\xrightarrow{\beta^*}\rH^i(\mathcal{T}_2,\beta^*(\mathscr{F}))\to \rH^i(\mathcal{T}_2,\alpha_*\gamma^*(\mathscr{F}))\to \rH^i(\mathcal{T}_1,\gamma^*(\mathscr{F})),
	\label{pullback and forward}
\end{equation}
où la deuxième flèche est induite par le morphisme d'adjonction $\beta^*\to \alpha_*\alpha^*\beta^*$ et la troisième flèche est induite par la suite spectrale de Cartan-Leray, s'identifie au morphisme \eqref{cohomologie pullback par gamma}. En effet, si on note, pour $i=1,2,3$, $\Gamma_{i}: \mathcal{T}_{i}\to \Ab$ le foncteur ``sections globales'', le morphisme composé \eqref{pullback and forward} est induit par le morphisme composé de la catégorie dérivée $\bfD(\Ab)$:
\begin{equation}
	\rR\Gamma_{3} \mathscr{F}\to\rR\Gamma_{3}\rR\beta_*\beta^*\mathscr{F} \to \rR\Gamma_{2}\rR\alpha_{*}\alpha^{*}\beta^{*}\mathscr{F}= \rR\Gamma_{1}\alpha^*\beta^*\mathscr{F}.
	\label{cohomlogie adjunction 1}
\end{equation}
Celui-ci s'identifie au morphisme induit par adjonction
\begin{equation}
	\rR\Gamma_{3}\mathscr{F}\to \rR\Gamma_{3} \rR \gamma_* \gamma^*\mathscr{F}=\rR\Gamma_{1}\alpha^*\beta^*\mathscr{F},
	\label{cohomlogie adjunction 2}
\end{equation}
qui fournit le morphisme \eqref{cohomologie pullback par gamma}.
\label{general morphism of topos}
\end{nothing}

\begin{nothing} \label{fet et}
	Pour tout schéma $X$, on note $\Et_{/X}$ (resp. $X_{\et}$) le site (resp. topos) étale de $X$. On désigne par $\Et_{\coh/X}$ (resp. $\Et_{\scoh/X}$) la sous-catégorie pleine de $\Et_{/X}$ formée des schémas étales de présentation finie sur $X$ (resp. étales, séparés et de présentation finie sur $X$), munie de la topologie induite par celle de $\Et_{/X}$. Si $X$ est quasi-séparé, le foncteur de restriction de $X_{\et}$ dans le topos des faisceaux d'ensembles sur $\Et_{\coh/X}$ (resp. $\Et_{\scoh/X}$) est une équivalence de catégories (\cite{SGAIV} VII 3.1 et 3.2).

	On désigne par $X_{\zar}$ le topos de Zariski de $X$ et par
	\begin{equation}
		u_{X}:X_{\et}\to X_{\zar}
		\label{uX et to zar}
	\end{equation}
	le morphisme canonique (\cite{SGAIV} VII 4.2.2). Si $F$ est un $\mathscr{O}_{X}$-module quasi-cohérent de $X_{\zar}$, on note encore $F$ le faisceau de $X_{\et}$ défini pour tout $X$-schéma étale $X'$ par (\cite{SGAIV} VII 2 c), cf. aussi \cite{AGT} III.2.9)
	\begin{equation}
		F(X')=\Gamma(X',F\otimes_{\mathscr{O}_{X}}\mathscr{O}_{X'}).
		\label{F zar to et}
	\end{equation}

On désigne par $\Et_{\textnormal{f}/X}$ le site fini étale de $X$, c'est-à-dire la sous-catégorie pleine de $\Et_{/X}$ formée des schémas étales et finis sur $X$, munie de la topologie induite par celle de $\Et_{/X}$. On note $X_{\fet}$ le topos fini étale de $X$, c'est-à-dire le topos des faisceaux d'ensembles sur $\Et_{\textnormal{f}/X}$. Le foncteur d'injection canonique $\Et_{\textnormal{f}/X} \to \Et_{/X}$ est continu. Il induit donc un morphisme de topos que l'on note
\begin{equation}
	\rho_{X}: X_{\et} \to X_{\fet}.
	\label{topos varsigma X}
\end{equation}
\end{nothing}
\begin{nothing}
Soient $X$ un schéma connexe, $\overline{x}$ un point géométrique de $X$. On désigne par $\pi_1(X,\overline{x})$ le groupe fondamental de $X$ en $\overline{x}$ et par $\HB_{\pi_1(X,\overline{x})}$ le topos classifiant du groupe pro-fini $\pi_1(X,\overline{x})$, c'est-à-dire la catégorie des ensembles discrets munis d'une action continue à gauche de $\pi_1(X,\overline{x})$ (\cite{SGAIV} IV 2.7). Le foncteur fibre en $\overline{x}$
\begin{equation}
	\omega_{\overline{x}}:\Et_{\textnormal{f}/X}\to \HB_{\pi_1(X,\overline{x})}
	\label{foncteur fibre en x}
\end{equation}
qui à tout revêtement étale $Y$ de $X$ associe l'ensemble des points géométriques de $Y$ au-dessus de $\overline{x}$, induit une équivalence de catégories (\cite{AGT} VI.9.8)
\begin{equation}
	\mu_{\overline{x}}: \HB_{\pi_1(X,\overline{x})}\xrightarrow{\sim} X_{\fet}.
	\label{foncteur mu}
\end{equation}

Le foncteur $\omega_{\overline{x}}$ est pro-représentable par un système projectif $(X_i,\overline{x}_{i},\phi_{ij})_{i\in I}$ de revêtements étales galoisiens pointés de $X$, indexé par un ensemble ordonné filtrant $I$, où $\overline{x}_{i}$ est un point géométrique de $X_{i}$ au-dessus de $\overline{x}$ et $\phi_{ij}:X_{i}\to X_{j}$ ($i\ge j$) est un épimorphisme qui envoie $\overline{x}_{i}$ sur $\overline{x}_{j}$. Il est appelé \textit{pro-objet universel galoisien} de $X$ en $\overline{x}$. Le foncteur
\begin{equation}
	\nu_{\overline{x}}:X_{\fet}\to \HB_{\pi_1(X,\overline{x})} \qquad F\mapsto \varinjlim_{i\in I} F(X_i),
	\label{foncteur nu}
\end{equation}
est une équivalence de catégories, quasi-inverse de $\mu_{\overline{x}}$ (\cite{AGT} VI.9.8).
\label{rep and sysl}
\end{nothing}

\begin{lemma}
	Soient $X$ un schéma quasi-compact, quasi-séparé et connexe et $\mathbb{L}$ un faisceau abélien de torsion de $X_{\fet}$. Le morphisme canonique
\begin{equation}
	\rho_{X}^{*}: \rH^i(X_{\fet},\mathbb{L})\to \rH^i(X_{\et},\rho^*_{X}(\mathbb{L})).
	\label{varsigma pullback}
\end{equation}
est un isomorphisme si $i=0,1$, et est un monomorphisme si $i=2$.
	\label{H1 et fet}
\end{lemma}
\textit{Preuve}. Le morphisme d'adjonction $\id \to \rho_{X*}\rho_{X}^{*}$ est un isomorphisme (\cite{AGT} VI.9.18), d'où la proposition pour $i=0$. Montrons-là pour $i=1$. Les topos $X_{\et}$ et $X_{\fet}$ sont cohérents (\cite{AGT} VI.9.12). Par (\cite{AGT} VI.9.20) et (\cite{SGAIV} VI 5.1), on peut se borner au cas où $\rho_{X}^{*}(\mathbb{L})$ est un faisceau abélien localement constant et constructible, i.e. représentable par un schéma en groupe fini et étale $G$ sur $X$. Le faisceau $\rR^{1}\rho_{X*}(\rho_{X}^{*}(\mathbb{L}))$ de $X_{\fet}$ est le faisceau sur $\Et_{\textnormal{f}/X}$ associé au préfaisceau défini par $Y \mapsto \rH^{1}(Y_{\et},\rho_{X}^{*}(\mathbb{L}))$. La catégorie des $\rho_{X}^{*}(\mathbb{L})$-torseurs sur $Y_{\et}$ est équivalente à la catégorie des $G_{Y}$-fibrés principaux homogènes sur $Y$ pour la topologie étale, où $G_{Y}=G\times_{X}Y$ (\cite{SGAIV} VII 2 b)). Tout $G_{Y}$-fibré principal homogène sur $Y$ est trivialisé par un revêtement étale de $Y$. Donc le faisceau $\rR^{1}\rho_{X*}(\rho_{X}^{*}(\mathbb{L}))$ est trivial.

On a une suite spectrale
\begin{equation}
	\rH^{i}(X_{\fet},\rR^{j}\rho_{X*}(\rho_{X}^{*}(\mathbb{L})))\Rightarrow \rH^{i+j}(X_{\et},\rho_{X}^{*}(\mathbb{L})).
\end{equation}
On en déduit un isomorphisme
\begin{displaymath}
	\rH^{1}(X_{\fet},\rho_{X*}\rho_{X}^{*}(\mathbb{L}))\xrightarrow{\sim} \rH^{1}(X_{\et},\rho_{X}^{*}(\mathbb{L})).
\end{displaymath}
La proposition pour $i=1$ résulte de ce dernier et de l'isomorphisme $\id \xrightarrow{\sim} \rho_{X*}\rho^{*}_{X}$. Comme $\rR^{1}\rho_{X*}(\rho_{X}^{*}(\mathbb{L}))$ est nul, on en déduit par \eqref{suite exacte terme petite} une suite exacte 
\begin{equation}
	0 \to \rH^{2}(X_{\fet},\rho_{X*}(\rho_{X}^{*}(\mathbb{L}))) \to \rH^{2}(X_{\et},\rho_{X}^{*}(\mathbb{L}))\to \rH^{0}(X_{\fet},\rR^{2}\rho_{X*}(\rho_{X}^{*}(\mathbb{L}))).
\end{equation}
La proposition pour $i=2$ résulte alors de ce dernier et de l'isomorphisme $\id \xrightarrow{\sim} \rho_{X*}\rho^{*}_{X}$.


\begin{nothing}\label{categorie a isogenies}
Si $\mathscr{C}$ est une catégorie additive, on désigne par $\mathscr{C}_{\mathbb{Q}}$ et l'on appelle \textit{catégorie des objets de $\mathscr{C}$ à isogénie près} (\cite{AGT} III.6.1), la catégorie ayant mêmes objets que $\mathscr{C}$, et telle que l'ensemble des morphismes soit donné, pour tous $M,N\in \Ob(\mathscr{C})$, par
\begin{equation}
	\Hom_{\mathscr{C}_{\mathbb{Q}}}(M,N)=\Hom_{\mathscr{C}}(M,N)\otimes_{\mathbb{Z}}\mathbb{Q}.
	\label{Hom C Hom CQ}
\end{equation}
On désigne par
\begin{equation}
	\mathscr{C}\to \mathscr{C}_{\mathbb{Q}},\qquad M\mapsto M_{\mathbb{Q}},
	\label{foncteur canonique de isogenie}
\end{equation}
le foncteur canonique. Deux objets de $\mathscr{C}$ sont dits \textit{isogènes} s'ils sont isomorphes dans la catégorie $\mathscr{C}_{\mathbb{Q}}$. Si $\mathscr{C}$ est une catégorie abélienne, la catégorie $\mathscr{C}_{\mathbb{Q}}$ est abélienne et le foncteur canonique \eqref{foncteur canonique de isogenie} est exact (\cite{AGT} III.6.1.4). Tout foncteur additif (resp. exact) entre catégories additives (resp. abéliennes) $\mathscr{C}\to \mathscr{C}'$ induit un foncteur additif (resp. exact) $\mathscr{C}_{\mathbb{Q}}\to \mathscr{C}'_{\mathbb{Q}}$ compatible aux foncteurs canoniques \eqref{foncteur canonique de isogenie} (\cite{AGT} III.6.1.5).

La famille des isogénies de $\mathscr{C}$ permet un calcul de fractions bilatéral (\cite{Il71} I 1.4.2). La catégorie $\mathscr{C}_{\mathbb{Q}}$ s'identifie à la catégorie localisé par rapport aux isogénies et le foncteur \eqref{foncteur canonique de isogenie} est le foncteur de localisation (cf. \cite{AGT} III.6.1).
\end{nothing}

\begin{nothing} \label{not Ext Hi Q}
Soit $(\mathcal{T},A)$ un topos annelé. On désigne par $\Mod_{\mathbb{Q}}(\mathcal{T},A)$ la catégorie des $A$-modules de $\mathcal{T}$ à isogénie près \eqref{categorie a isogenies} dont les objets sont appelés les \textit{$A_{\mathbb{Q}}$-modules}. \'{E}tant donnés deux $A_{\mathbb{Q}}$-modules $M$ et $N$, on note $\Hom_{A_{\mathbb{Q}}}(M,N)$ le groupe des morphismes de $M$ dans $N$ dans $\Mod_{\mathbb{Q}}(\mathcal{T},A)$. Rappelons que la catégorie $\Mod_{\mathbb{Q}}(\mathcal{T},A)$ a suffisamment d'injectifs (\cite{AGT} III.6.1.6) et reprenons les notations de \eqref{definition Ext using Hom derive} pour cette catégorie. Pour tout entier $i\ge 0$ et tous $A_{\mathbb{Q}}$-modules $M$ et $N$, on pose $\Ext_{A_{\mathbb{Q}}}^{i}(M,N)= \Ext_{\Mod_{\mathbb{Q}}(\mathcal{T},A)}^{i}(M,N)$ qui est un $\mathbb{Q}$-espace vectoriel. En particulier, pour tout entier $i\ge 0$ et $A_{\mathbb{Q}}$-module $N$, on pose
\begin{equation}
	\rH^{i}(\mathcal{T},N)=\Ext_{A_{\mathbb{Q}}}^{i}(A_{\mathbb{Q}},N).
	\label{def cohomology classe Q}
\end{equation}

	On note $\Vect_{\mathbb{Q}}$ la catégorie des $\mathbb{Q}$-espaces vectoriels. Le foncteur canonique $\Ab\to \Vect_{\mathbb{Q}}$ induit une équivalence de catégories $\Ab_{\mathbb{Q}}\xrightarrow{\sim} \Vect_{\mathbb{Q}}$. \'{E}tant donné un $A$-module $M$ de $\mathcal{T}$, on a un diagramme commutatif
	\begin{equation}
		\xymatrixcolsep{5pc}\xymatrix{
			\Mod(\mathcal{T},A)\ar[d] \ar[r]^-{\Hom_{A}(M,-)} & \Ab \ar[d]\\
			\Mod_{\mathbb{Q}}(\mathcal{T},A) \ar[r]^-{\Hom_{A_{\mathbb{Q}}}(M_{\mathbb{Q}},-)} & \Vect_{\mathbb{Q}},
		}
	\end{equation}
	où les flèche verticales sont les foncteurs canoniques. Ceux-ci sont exacts et transforment les objets injectifs en des objets injectifs (\cite{AGT} III.6.1.6). Pour tout entier $i\ge 0$ et tous $A$-modules $M$ et $N$, on en déduit par la suite spectrale de Cartan-Leray un isomorphisme
	\begin{equation}
		\Ext_{A}^{i}(M,N)\otimes_{\mathbb{Z}}\mathbb{Q}\xrightarrow{\sim}\Ext_{A_{\mathbb{Q}}}^{i}(M_{\mathbb{Q}},N_{\mathbb{Q}}).
		\label{Ext A otimes Q ExtAQ}
	\end{equation}
En particulier, pour tout entier $i\ge 0$ et tout $A$-module $N$, on a un isomorphisme
\begin{equation}
	\rH^{i}(\mathcal{T},N)\otimes_{\mathbb{Z}}\mathbb{Q}\simeq \rH^{i}(\mathcal{T},N_{\mathbb{Q}}).
	\label{coh N coh class Q}
\end{equation}
\end{nothing}

\begin{nothing}
	Soient $\mathscr{A}$ une catégorie abélienne, $\End(\id_{\mathscr{A}})$ l'anneau des endomorphismes du foncteur identique de $\mathscr{A}$ et $\varphi: \oo\to \End(\id_{\mathscr{A}})$ un homomorphisme. Pour tout objet $M$ de $\mathscr{A}$ et tout $\gamma\in \oo$, on note $\mu_{\gamma}(M)$ l'endomorphisme de $M$ défini par $\varphi(\gamma)$. En particulier, pour tous objets $M$ et $N$ de $\mathscr{A}$, $\Hom_{\mathscr{A}}(M,N)$ est muni d'une structure de $\oo$-module.

	On dit qu'un objet $M$ de $\mathscr{A}$ est \textit{$\alpha$-nul} s'il est annulé par tout élément de $\mm$, i.e. si $\mu_{\gamma}(M)=0$ pour tout $\gamma\in \mm$. Pour toute suite exacte $0\to M''\to M \to M'\to 0$, pour que $M$ soit $\alpha$-nul, il faut et il suffit que $M''$ et $M'$ soient $\alpha$-nuls. On appelle catégorie des $\alpha$-objets de $\mathscr{A}$ et l'on note $\alpha$-$\mathscr{A}$ le quotient de la catégorie $\mathscr{A}$ par la sous-catégorie épaisse formée des objets $\alpha$-nuls (\cite{Ga62} III \S 1; cf. aussi \cite{AG15} 1.4.2). On note
	\begin{equation}
		\alpha:\mathscr{A}\to \alpha\textnormal{-}\mathscr{A}, \qquad M\mapsto \alpha(M)
		\label{foncteur alpha general}
	\end{equation}
	le foncteur canonique; on notera aussi $M^{\alpha}$ au lieu de $\alpha(M)$ lorsqu'il n'y a aucun risque de confusion. La catégorie $\alpha$-$\mathscr{A}$ est abélienne et le foncteur $\alpha$ est exact (\cite{Ga62} III \S 1 prop. 1). On dit qu'un morphisme $f$ de $\mathscr{A}$ est un $\alpha$-isomorphisme (resp. $\alpha$-monomorphisme, resp. $\alpha$-épimorphisme) si $\alpha(f)$ est un isomorphisme (resp. monomorphisme, resp. épimorphisme), autrement dit si son noyau et son conoyau (resp. son noyau, resp. son conoyau) sont $\alpha$-nuls.

	La famille des $\alpha$-isomorphismes de $\mathscr{A}$ permet un calcul de fractions bilatère (\cite{Il71} I 1.4.2). La catégorie $\alpha$-$\mathscr{A}$ s'identifie à la catégorie localisé par rapport aux $\alpha$-isomorphismes et $\alpha$ \eqref{foncteur alpha general} est le foncteur canonique (de localisation).
	\label{Cat Ab almost}
\end{nothing}
\begin{lemma}[\cite{AG15} 1.4.3]
	Les hypothèses étant celles de \ref{Cat Ab almost}, soient, de plus, $f:M\to N$ un morphisme de $\mathscr{A}$, $\gamma\in \oo$ tels que le noyau et le conoyau de $f$ soient annulés par $\gamma$. Alors, il existe un morphisme $g:N\to M$ de $\mathscr{A}$ tel que $g\circ f=\mu_{\gamma^{2}}(M)$ et $f\circ g=\mu_{\gamma^{2}}(N)$.
	\label{lemma 143}
\end{lemma}
%

\begin{nothing}
	Soit $\mathscr{A}$ une catégorie abélienne tensorielle autrement dit, $\mathscr{A}$ est une catégorie abélienne munie d'un foncteur bi-additif $\otimes:\mathscr{A}\times \mathscr{A}\to \mathscr{A}$ et d'un objet unité $A$, vérifiant certaines conditions (\cite{DM82} 1.15), et soit $\varphi: \oo\to \End(A)$ un homomorphisme. On a un homomorphisme canonique $\End(A)\to \End(\id_{\mathscr{A}})$. On peut donc définir la catégorie $\alpha$-$\mathscr{A}$ suivant \ref{Cat Ab almost}. Le produit tensoriel induit un bi-foncteur
	\begin{equation}
		\alpha\textnormal{-}\mathscr{A}\times \alpha\textnormal{-}\mathscr{A}\to \alpha\textnormal{-}\mathscr{A}, \qquad (M,N)\mapsto M\otimes N,
		\label{produit tensoriel alpha}
	\end{equation}
	qui fait de $\alpha\textnormal{-}\mathscr{A}$ une catégorie abélienne tensorielle, dont $A^{\alpha}$ est un objet unité (cf. \cite{AG15} 1.4.4). Le foncteur $\alpha$ induit un homomorphisme $\End(A)\to \End(A^{\alpha})$. Par suite, pour tous objets $M$ et $N$ de $\alpha\textnormal{-}\mathscr{A}$, $\Hom_{\alpha\textnormal{-}\mathscr{A}}(M,N)$ est canoniquement muni d'une structure de $\End(A)$-module.
	\label{Cat Ab Ten almost}
\end{nothing}

\section{Représentations à coefficients dans les modules de type $\alpha$-fini} \label{Alg homologique}
\begin{nothing} \label{alpha o modules}
	Rappelons que $\oo$ désigne le séparé complété $p$-adique de $\mathcal{O}_{\overline{K}}$ et $\mm$ son idéal maximal et que l'on a $\mm^{2}=\mm$ et $v(\mm)=\mathbb{Q}_{>0}$ \eqref{notations 11}. Pour tout nombre réel $r\ge 0$, on note $I_{r}$ l'idéal de $\oo$ formé des éléments $x\in \oo$ tels que $v(x)>r$.
	
	On note $\Mod(\oo)$ la catégorie des $\oo$-modules. Prenant pour $\varphi:\oo\xrightarrow{\sim}\End(\id_{\Mod(\oo)})$ l'isomorphisme canonique, on considère les notions introduites dans \eqref{Cat Ab almost} pour la catégorie $\Mod(\oo)$. Considérons le foncteur
\begin{equation}
	\Mod(\oo) \to \Mod(\oo) \qquad M \mapsto M_{\sharp}=\Hom_{\oo}(\mm,M).
	\label{functor sigma star}
\end{equation}
Pour tout $\oo$-module $M$, on a des $\alpha$-isomorphismes \eqref{Cat Ab almost} canoniques fonctoriels (\cite{AGT} V.2.4)
\begin{eqnarray}
	M&\to& M_{\sharp}, \qquad x\mapsto (m\in \mm\mapsto m x), \label{N to Mstar}\\
	\mm\otimes_{\oo}M&\to& M. \label{mmooM to M}
\end{eqnarray}
\end{nothing}
\begin{prop}[\cite{AGT} V.2.5 et 2.6]
	Soit $f:M\to N$ un morphisme de $\oo$-module.

	\textnormal{(i)} Pour que $f$ soit un $\alpha$-isomorphisme, il faut et il suffit que le morphisme induit $f_{\sharp}:M_{\sharp}\to N_{\sharp}$ soit un isomorphisme.

	\textnormal{(ii)} Pour que $f$ soit un $\alpha$-isomorphisme, il faut et il suffit que le morphisme induit $\mm\otimes_{\oo}M\to \mm\otimes_{\oo}N$ soit un isomorphisme.

	\textnormal{(iii)} Pour tout entier $i\ge 0$, le foncteur $\Ext^{i}_{\oo}(\mm,-)$ envoie les $\alpha$-isomorphismes sur des isomorphismes.

	\textnormal{(iv)} Le foncteur \eqref{functor sigma star} envoie les $\alpha$-monomorphismes en des monomorphismes.
	\label{prop Tsuji alpha iso}
\end{prop}

\begin{rem} \label{m pas projectif}
	Le $\oo$-module $\mm$ n'est pas projectif. En effet, le $\oo$-module $\mm$ est plat et on a $\mm=\varinjlim_{n\in \mathbb{N}}(p^{1/n})$, où $p^{1/n}$ est une racine $n$-ième de $p$ dans $\mathcal{O}_{\overline{K}}$. Mais le système projectif $(\Hom_{\oo}(\oo, (p^{1/n})))_{n\ge 1}$ ne vérifie pas la condition de Mittag-Leffler. L'assertion résulte alors de (\cite{RG71} 3.1.3). 
\end{rem}
\begin{nothing}
	Suivant (\cite{Sch13} 2.2), pour tous $\oo$-modules $M$ et $N$ et tout $\gamma\in \mm$, on dit que $M$ et $N$ sont \textit{$\gamma$-équivalents} et on note $M\approx_{\gamma}N$, s'il existe deux $\oo$-morphismes $f:M\to N$ et $g: N\to M$ tels que $f\circ g=\gamma\id_{N}$ et $g\circ f=\gamma\id_{M}$. On dit que $M$ et $N$ sont \textit{$\alpha$-équivalents} et on note $M\approx N$ si $M$ et $N$ sont $\gamma$-équivalents pour tout $\gamma\in \mm$.
	
	On dit qu'un $\oo$-module $M$ est de \textit{type $\alpha$-fini} si, pour tout $\gamma\in \mm$, il existe un $\oo$-module de type fini $N$ tel que $M\approx_{\gamma}N$. Pour tout entier $n\ge 1$, on dit qu'un $\oo_{n}$-module $M$ est de \textit{type $\alpha$-fini} s'il est de type $\alpha$-fini en tant que $\oo$-module.
	
	On dit qu'une suite de $\oo$-modules est \textit{$\alpha$-exacte} si ses groupes de cohomologie sont $\alpha$-nuls (cf. \cite{AG15} 1.5.2).
	\label{type de alpha fini}
\end{nothing}

\begin{rem}
	(i) D'après \ref{lemma 143}, deux $\oo$-modules $\alpha$-isomorphes sont $\alpha$-équivalents. La réciproque n'est pas vraie. En effet, pour tout nombre réel $r>0$, on a alors $I_{r}\approx \oo$ \eqref{alpha o modules}. Si $r\not\in \mathbb{Q}_{>0}$, $\oo$ et $I_{r}$ ne sont pas $\alpha$-isomorphes (cf. \cite{AG15} 1.8.4).
	
	(ii) La notion de \textit{type $\alpha$-fini} introduite dans \ref{type de alpha fini} est compatible avec la définition générale introduite dans \ref{type alpha fini topos} (cf. \cite{AG15} 1.5.3 et 1.8.5(ii)).
\end{rem}

\begin{theorem}[\cite{Sch13} 2.5] \label{thm classification scholze}	
	Pour tout $\oo$-module de type $\alpha$-fini $M$, il existe une unique suite décroissante de nombres réels positifs $(r_{i})_{i\ge1}$, tendant vers $0$, et un unique entier $n\ge 0$ tels que
	\begin{equation}
		M\approx \oo^{n}\oplus (\oplus_{i\ge 1} \oo/I_{r_{i}}).
		\label{classification M alpha tf}
	\end{equation}
\end{theorem}

\begin{lemma}[\cite{AG15} 1.5.13 et 1.8.7] \label{lemma Gabber exactes tf}
	Soit $0\to M'\to M\to M''\to 0$ une suite $\alpha$-exacte de $\oo$-modules. Pour que $M$ soit de type $\alpha$-fini, il faut et il suffit que $M'$ et $M''$ le soient.
\end{lemma}

\begin{prop} \label{lemma topologie aptf}
	Soit $M$ un $\oo$-module de type $\alpha$-fini.

	\textnormal{(i)} Il existe un sous-$\oo$-module libre de type fini $M^{\circ}$ de $M$ tel que $M/M^{\circ}$ soit annulé par une puissance de $p$.
	
	\textnormal{(ii)} Les sous-modules $(p^{n}M^{\circ})_{n\ge 1}$ forment un système fondamental de voisinages de $0$ pour la topologie $p$-adique de $M$.

	\textnormal{(iii)} La topologie $p$-adique de $M$ est séparée et complète.

	\textnormal{(iv)} Soit $N$ un sous-module de $M$. La topologie $p$-adique de $N$ est induite par la topologie $p$-adique de $M$.
\end{prop}
\textit{Preuve}. (i) D'après \ref{thm classification scholze}, le $\oo$-module $M$ est $\alpha$-équivalent à un $\oo$-module $N=\oo^{n}\oplus (\oplus_{i\ge 1} \oo/I_{r_{i}})$ avec $(r_{i})_{i\ge 1}$ tendant vers $0$. Notons $N^{\circ}$ le sous-$\oo$-module $\oo^{n}$ de $N$. Choisissons un élément $\gamma$ de $\mm$. Il existe des morphismes de $\oo$-modules $f:N\to M$ et $g:M\to N$ tels que $g\circ f=\gamma\id_{N}$ et $f\circ g=\gamma\id_{M}$. Donc, la restriction de $f$ à $N^{\circ}$ est injective. On prend pour $M^{\circ}$ le sous-$\oo$-module libre de type fini $f(N^{\circ})$ de $M$. On a alors $g(M^{\circ})=\gamma N^{\circ}\subset N^{\circ}$. Les morphismes $f$ et $g$ induisent une $\gamma$-équivalence entre $M/M^{\circ}$ et $N/N^{\circ}$:
\begin{equation}
	\bar{f}:N/N^{\circ}\to M/M^{\circ} \qquad \bar{g}: M/M^{\circ}\to N/N^{\circ}.
\end{equation}
Comme $N/N^{\circ}\simeq \oplus_{i\ge 1} \oo/I_{r_{i}}$ est annulé par une puissance de $p$, il en est de même de $M/M^{\circ}$.

(ii) En vertu de (i), il existe un entier $m\ge 1$ tel que $p^{m}(M/M^{\circ})=0$. On en déduit que $p^{m}M\subset M^{\circ}\subset M$. Les systèmes fondamentaux $(p^{n}M^{\circ})_{n\ge 1}$ et $(p^{n}M)_{n\ge 1}$ de voisinages de $0$ sont alors équivalents.

(iii) On a un isomorphisme canonique $M^{\circ}\xrightarrow{\sim} \widehat{M^{\circ}}=\varprojlim M^{\circ}/p^{n}M^{\circ}$. En vertu de (ii), $\widehat{M}=\varprojlim M/p^{n}M^{\circ}$ est le séparé complété $p$-adique de $M$. Pour tout entier $n\ge 1$, on a une suite exacte
\begin{equation}
	0\to M^{\circ}/p^{n}M^{\circ}\to M/p^{n}M^{\circ}\to M/M^{\circ}\to 0.
\end{equation}
Le système projectif $(M^{\circ}/p^{n}M^{\circ})_{n\ge 1}$ vérifie la condition de Mittag-Leffler. La limite projective induit donc une suite exacte
\begin{equation}
	0\to \widehat{M^{\circ}} \to \widehat{M}\to M/M^{\circ}\to 0
\end{equation}
qui s'insère dans un diagramme commutatif
\begin{equation}
	\xymatrix{
		0\ar[r] & M^{\circ} \ar[r] \ar[d]_{\wr} & M \ar[r] \ar[d] & M/M^{\circ} \ar[r] \ar@{=}[d] & 0\\
		0\ar[r] & \widehat{M^{\circ}} \ar[r] & \widehat{M} \ar[r]& M/M^{\circ} \ar[r] & 0
	}
\end{equation}
On en déduit que le morphisme canonique $M\to \widehat{M}$ est un isomorphisme; d'où l'assertion.

(iv) En vertu de \ref{lemma Gabber exactes tf}, $N$ est de type $\alpha$-fini. Il existe un sous-$\oo$-module libre de type fini $N^{\circ}$ de $N$ vérifiant la condition (i). On peut supposer de plus que $N^{\circ}\subset M^{\circ}$. Soit $N^{\circ}_{\sat}=\{x\in M^{\circ}| p^{n}x\in N^{\circ} \textnormal{ pour un entier } n\ge 1\}$ la saturation de $N^{\circ}$ dans $M^{\circ}$. On a alors $(p^{n}M^{\circ})\cap N^{\circ}_{\sat}=p^{n} N^{\circ}_{\sat}$ pour tout $n\ge 1$. Comme $p$ n'est pas un diviseur de zéro dans $M^{\circ}/N_{\sat}^{\circ}$, $M^{\circ}/N_{\sat}^{\circ}$ est $\oo$-plat (\cite{Ab10} 1.9.12). On en déduit que $N^{\circ}_{\sat}$ est de type fini (\cite{Ab10} 1.9.14). Il existe donc un entier $l\ge 1$ tel que $p^{l}N_{\sat}^{\circ}\subset N^{\circ}$. Pour tout entier $n\ge 1$, on a 
\begin{displaymath}
	(p^{n}M^{\circ})\cap N^{\circ}=(p^{n}M^{\circ}\cap N_{\sat}^{\circ})\cap N^{\circ}=(p^{n}N_{\sat}^{\circ})\cap N^{\circ}=p^{n}N_{\sat}^{\circ}.
\end{displaymath}
On en déduit que
\begin{displaymath}
	p^{n}N^{\circ}\subset (p^{n}M^{\circ})\cap N^{\circ}\subset p^{n-l}N^{\circ}.
\end{displaymath}
L'assertion s'ensuit compte tenu de (ii).


\begin{lemma}[\cite{AG15} 1.8.11] \label{lemma 1811}
	Soit $M$ un $\oo$-module tel que $M/pM$ soit de type $\alpha$-fini. Alors, le $\oo$-module $\widehat{M}=\varprojlim M/p^{n}M$ est de type $\alpha$-fini.
\end{lemma}

\begin{lemma}
	\textnormal{(i)} Soit $(M_{n})_{n\ge 1}$ un système projectif de $\oo$-modules $\alpha$-nuls. Alors, les $\oo$-modules $\varprojlim M_{n}$ et $\rR^{1}\varprojlim M_{n}$ sont $\alpha$-nuls.
	
	\textnormal{(ii)} Soient $(M_{n})_{n\ge 1}$ et $(N_{n})_{n\ge 1}$ deux systèmes projectifs de $\oo$-modules et $(f_{n}:M_{n}\to N_{n})_{n\ge 1}$ un système projectif des $\alpha$-isomorphismes. Alors, le foncteur limite projective induit des $\alpha$-isomorphismes $\varprojlim M_{n}\to \varprojlim N_{n}$ et $\rR^{1}\varprojlim M_{n}\to \rR^{1}\varprojlim N_{n}$.
	\label{lemma proj lim alpha nuls}
\end{lemma}

L'assertion (i) est démontrée dans (\cite{GR03} 2.4.2) et l'assertion (ii) résulte de (i).

\begin{definition}[\cite{AGT} V.6.1] \label{alpha plat}
	Soient $R$ une $\oo$-algèbre et $M$ un $R$-module. On dit que $M$ est $\alpha$-plat si il vérifie les conditions équivalentes suivantes:

	(i) Pour tout $R$-module $N$ et tout entier $i\ge 1$, $\Tor_{i}^{R}(M,N)$ est un $\oo$-module $\alpha$-nul.

	(ii) Pour tout $R$-module $N$, $\Tor_{1}^{R}(M,N)$ est un $\oo$-module $\alpha$-nul.

	(iii) Pour tout morphisme $R$-linéaire injectif $f:N_{1}\to N_{2}$, le noyau de $\id_{M}\otimes f:M\otimes_{R}N_{1}\to M\otimes_{R}N_{2}$ est un $\oo$-module $\alpha$-nul.
\end{definition}

\begin{lemma}
	Soient $R$ une $\oo$-algèbre et $M$ un $R$-module $\alpha$-plat. Alors, le $R$-module $\mm\otimes_{\oo}M$ est plat. En particulier, pour tout $\oo$-module $\alpha$-plat $N$, $\mm\otimes_{\oo}N$ est sans-torsion.
	\label{lemma alpha plat plat mm}
\end{lemma}
\textit{Preuve}. En vertu de \ref{prop Tsuji alpha iso}(ii), pour tout $\oo$-module $\alpha$-nul $N$, $\mm\otimes_{\oo}N$ est nul. La $R$-platitude de $\mm\otimes_{\oo}M$ résulte alors de la $\oo$-platitude de $\mm$. On sait qu'un $\oo$-module est plat si et seulement s'il est sans-torsion (\cite{BouAlgcom} chap. VI \S 3.6 lem.1); d'où la deuxième assertion.

\begin{prop}
	Soient $M$ et $N$ deux $\oo_{n}$-modules \eqref{notations 11} (resp. $\oo$-modules) $\alpha$-équivalents. Alors, la $\alpha$-platitude de $M$ est équivalente à celle de $N$.
\end{prop}
\textit{Preuve}. Supposons que $M$ soit $\alpha$-plat. Soit $f:N_{1}\to N_{2}$ un morphisme injectif de $\oo_{n}$-modules. Pour tout $\gamma\in \mm$, il existe des morphismes $\oo$-linéaires $g:M\to N$ et $h:N\to M$ tels que $g\circ h=\gamma\id_{N}$ et $h\circ g=\gamma\id_{M}$. On a alors un diagramme commutatif
\begin{equation}
	\xymatrix{
		N\otimes_{\oo_{n}}N_{1}\ar[r]^{h\otimes \id_{N_{1}}} \ar[d]_{\id_{N}\otimes f} & M\otimes_{\oo_{n}}N_{1}\ar[d]^{\id_{M}\otimes f}\\
		N\otimes_{\oo_{n}}N_{2}\ar[r]^{h\otimes \id_{N_{2}}}& M\otimes_{\oo_{n}}N_{2}
	}
\end{equation}
Comme $\Ker(M\otimes_{\oo_{n}}N_{1}\to M\otimes_{\oo_{n}}N_{2})$ est $\alpha$-nul et $g\circ h=\gamma\id_{N}$, on en déduit que $\Ker(N\otimes_{\oo_{n}}N_{1}\to N\otimes_{\oo_{n}}N_{2})$ est annulé par l'idéal $\gamma\mm$ de $\oo$. La $\alpha$-platitude de $N$ s'ensuit compte tenu du fait que $\mm^{2}=\mm$. La démonstration pour les $\oo$-modules est similaire.

\begin{coro}
	Soient $n$ un entier $\ge 1$ et $M$ un $\oo_{n}$-module (resp. $\oo$-module) de type $\alpha$-fini. Pour que $M$ soit $\alpha$-plat, il faut et il suffit qu'il existe un entier $m\ge 0$ tel que $M\approx \oo_{n}^{m}$ (resp. $M\approx \oo^{m}$).
	\label{coro alpha plat tf}
\end{coro}
\textit{Preuve}. Comme $p^{n}M=0$, d'après \ref{thm classification scholze}, il existe une suite décroissante de nombres réels $(r_{i})_{i\ge 1}$ tendant vers $0$ tels que $r_{i}\le n$ et que $M\approx \oplus_{i\ge 1} \oo/I_{r_{i}}$. Le $\oo$-module $\oo/I_{n}$ est $\alpha$-isomorphe à $\oo_{n}$. Soit $r$ un nombre réel tel que $0<r<n$. Il est clair que $\oo/I_{r}$ n'est pas un $\oo_{n}$-module $\alpha$-plat. L'assertion pour les $\oo_{n}$-modules s'ensuit. La démonstration pour les $\oo$-modules est similaire.

\begin{lemma}
	Soit $(M_{n})_{n\ge 1}$ un système projectif de $\oo$-modules tel que, pour tout entier $n\ge 1$, $M_{n}$ soit un $\oo_{n}$-module $\alpha$-plat de type $\alpha$-fini et que le morphisme canonique $M_{n+1}\otimes_{\oo_{n+1}}\oo_{n}\rightarrow M_{n}$ soit un $\alpha$-isomorphisme. On pose, pour tout entier $n\ge 1$, $N_{n}=\mm\otimes_{\oo}M_{n}$, $\widehat{M}=\varprojlim M_{n}$ et $\widehat{N}=\varprojlim N_{n}$. Alors:

	\textnormal{(i)} On a un $\alpha$-isomorphisme canonique $\widehat{N}\rightarrow \widehat{M}$.

	\textnormal{(ii)} L'homomorphisme canonique $\widehat{N}/p^{n}\widehat{N}\to N_{n}$ est un isomorphisme.

	\textnormal{(iii)} Les $\oo$-modules $\widehat{N}$ et $\widehat{M}$ sont $\alpha$-plats de type $\alpha$-fini.
	\label{lemma aptf P oon}
\end{lemma}
\textit{Preuve}. (i) Pour tout entier $n\ge 1$, le morphisme canonique $N_{n}\to M_{n}$ est un $\alpha$-isomorphisme \eqref{mmooM to M}. En vertu de \ref{lemma proj lim alpha nuls}(ii), la limite projective induit un $\alpha$-isomorphisme canonique $\widehat{N}\to \widehat{M}$.

(ii) D'après \ref{prop Tsuji alpha iso}(ii), le $\alpha$-isomorphisme $M_{n+1}\otimes_{\oo_{n+1}}\oo_{n}\to M_{n}$ induit un isomorphisme $N_{n+1}\otimes_{\oo_{n+1}}\oo_{n}\xrightarrow{\sim} N_{n}$. On notera que, pour tout entier $n\ge 1$, $N_{n}$ est un $\oo_{n}$-module plat \eqref{lemma alpha plat plat mm}. Pour tous entiers $m,n\ge 1$, appliquant le foncteur $N_{m+n}\otimes_{\oo_{m+n}}-$ à la suite exacte $0\to \oo_{m}\xrightarrow{p^{n}} \oo_{m+n}\to \oo_{n}\to 0$, on en déduit une suite exacte:
\begin{equation}
	0\to N_{m}\xrightarrow{p^{n}} N_{m+n}\to N_{n}\to 0.
\end{equation}
Le système projectif $(N_{m})_{m\ge 1}$ vérifie la condition de Mittag-Leffler. On en déduit une suite exacte 
\begin{equation}
	\xymatrix{
		0\ar[r]& \widehat{N}\ar[r]^{p^{n}} & \widehat{N}\ar[r] & N_{n}\ar[r]& 0,
	}
\end{equation}
d'où l'assertion. 

(iii) En vertu de (i) et \ref{coro alpha plat tf}, il suffit de démontrer l'assertion pour $\widehat{N}$. La $\alpha$-finitude de $\widehat{N}$ résulte de (ii) et \ref{lemma 1811}. D'après \ref{thm classification scholze}, $\widehat{N}$ est $\alpha$-équivalent à un $\oo$-module $\oo^{m}\oplus(\oplus_{i\ge 1} \oo/I_{r_{i}})$ avec $(r_{i})_{i\ge 1}$ tendant vers $0$. On en déduit que, pour tout entier $n\ge 1$, on a $\widehat{N}/p^{n}\widehat{N}\approx \oo_{n}^{m}\oplus (\oplus_{i\ge 1} (\oo/I_{r_{i}})/p^{n}(\oo/I_{r_{i}}))$. D'autre part, on notera qu'il existe un entier $l\ge 1$ tel que, pour tout entier $n\ge 1$, $N_{n}\approx \oo_{n}^{l}$ \eqref{coro alpha plat tf}. En vertu de (ii), pour tout entier $i\ge 1$, on a donc $r_{i}=0$; d'où la $\alpha$-platitude de $\widehat{N}$.

\begin{rem}
	Sous les hypothèses de \ref{lemma aptf P oon}. La topologie $p$-adique de $\widehat{N}$ est la limite projective des topologies discrètes sur $N_{n}$.
	\label{topologie p adic limite proj aptf}
\end{rem}

\begin{nothing}
	Soient $V$ un $\mathfrak{C}$-espace vectoriel de dimension finie, $V^{\circ}$ un sous-$\oo$-module libre de type fini de $V$ qui l'engendre sur $\mathfrak{C}$. On appelle \textit{topologie $p$-adique sur $V$} l'unique topologie compatible avec sa structure de groupe additif pour laquelle les sous-groupes $(p^{n}V^{\circ})_{n\ge 1}$ forment un système fondamental de voisinages de $0$. Cette topologie ne dépend pas du choix de $V^{\circ}$ (cf. \cite{AGT} II.2.2). On notera que tout morphisme $\mathfrak{C}$-linéaire de $\mathfrak{C}$-espaces vectoriels de dimension finie est continu pour la topologie $p$-adique.
	\label{topologie R}
\end{nothing}

\begin{nothing}
	Dans la suite de cette section, on se donne un groupe profini $G$. Pour tout entier $n\ge 1$, on munit l'anneau $\oo_{n}=\oo/p^{n}\oo$ (resp. $\oo$, resp. $\mathfrak{C}$) de la topologie discrète (resp. de la topologie $p$-adique, resp. de la topologie $p$-adique \eqref{topologie R}) et de \textit{l'action triviale} de $G$. 
	
	Soit $n$ un entier $\ge 1$. Par \textit{$\oo_{n}$-représentation de $G$}, on sous-entend un $\oo_{n}$-module $M$ muni de la topologie discrète et d'une action continue et $\oo_{n}$-linéaire de $G$. Un morphisme de $\oo_{n}$-représentations de $G$ est une application $\oo_{n}$-linéaire et $G$-équivariante. On désigne par $\Rep_{\oo_{n}}(G)$ la catégorie des $\oo_{n}$-représentations de $G$.

	On note $\breve{\oo}$ le système projectif d'anneaux $(\oo_{n})_{n\ge 1}$. On désigne par $\Rep_{\breve{\oo}}(G)$ la catégorie des systèmes projectifs de $\oo[G]$-modules discrèts $(V_{n})_{n\ge 1}$ tels que, pour tout entier $n\ge 1$, $p^{n}V_{n}=0$. Il est clair que les catégories $\Rep_{\oo_{n}}(G)$ et $\Rep_{\breve{\oo}}(G)$ sont abéliennes.
	
	Par \textit{$\oo$-représentation continue de $G$}, on sous-entend un $\oo$-module \textit{de type $\alpha$-fini} muni de la topologie $p$-adique et d'une action $\oo$-linéaire continue de $G$. Par \textit{$\mathfrak{C}$-représentation continue}, on sous-entend un $\mathfrak{C}$-espace vectoriel de \textit{dimension finie} muni de la topologie $p$-adique \eqref{topologie R} et d'une action $\mathfrak{C}$-linéaire continue de $G$. Un morphisme de $\oo$-représentations (resp. $\mathfrak{C}$-représentations) continues de $G$ est une application $\oo$-linéaire et $G$-équivariante. Une telle application est continue pour la topologie $p$-adique (cf. \ref{topologie R}). On désigne par $\Rep_{\oo}^{\cont}(G)$ (resp. $\Rep_{\mathfrak{C}}^{\cont}(G)$) la catégorie des $\oo$-représentations (resp. $\mathfrak{C}$-représentations) continues de $G$. 
	\label{oo oon reprentations general}
\end{nothing}

\begin{lemma} \label{lemma cat ab}
	Soient $V$ un objet de $\Rep_{\oo}^{\cont}(G)$ et $W$ un sous-$\oo$-module de $V$ invariant par $G$. Alors, $W$ (resp. $V/W$) est un $\oo$-module de type $\alpha$-fini et la restriction de l'action de $G$ sur $W$ (resp. $V/W$) est continue pour la topologie $p$-adique.
\end{lemma}
\textit{Preuve}. La $\alpha$-finitude de $W$ et $V/W$ résulte de \ref{lemma Gabber exactes tf}. La continuité de l'action de $G$ pour la topologie $p$-adique de $W$ s'ensuit compte tenu \ref{lemma topologie aptf}(iv). Soit $\pi:V\to V/W$ la projection. Pour tout entier $n\ge 1$, on a $\pi(p^{n}V)=p^{n}(V/W)$. On en déduit la continuité de l'action de $G$ sur $V/W$.

\begin{prop}\label{cat ab Rep oo cont}
	La catégorie $\Rep_{\oo}^{\cont}(G)$ est abélienne.
\end{prop}
\textit{Preuve}. Cette catégorie est clairement additive. Soit $f:V_{1}\to V_{2}$ un morphisme de $\Rep_{\oo}^{\cont}(G)$. D'après \ref{lemma cat ab}, on définit $\Ker(f)$ (resp. $\Coker(f)$) par le $\oo$-module $\Ker(V_{1}\to V_{2})$ (resp. $\Coker(V_{1}\to V_{2})$) muni de la topologie $p$-adique et de l'action de $G$ induite par $V_{1}$ (resp. $V_{2}$). Il est clair que le morphisme canonique $\Coim(f)\to \Image(f)$ est un isomorphisme.

\begin{nothing} \label{oo oon reprentations}
	Soit $n$ un entier $\ge 1$. On désigne par $\Rep_{\oo_{n}}^{\ltf}(G)$ la sous-catégorie pleine de $\Rep_{\oo_{n}}(G)$ formée des $\oo_{n}$-représentations dont le $\oo_{n}$-module sous-jacent est libre de type fini, par $\Rep_{\oo}^{\ltf}(G)$ la sous-catégorie pleine de $\Rep_{\oo}^{\cont}(G)$ formée des $\oo$-représentations continues dont le $\oo$-module sous-jacent est libre de type fini.
	
	On désigne par $\Rep^{\ltf}_{\breve{\oo}}(G)$ la sous-catégorie pleine de $\Rep_{\breve{\oo}}(G)$ formée des systèmes projectifs $(V_{n})_{n\ge 1}$ tels que, pour tout entier $n\ge 1$, $V_{n}$ soit un objet de $\Rep^{\ltf}_{\oo_{n}}(G)$ et que le morphisme canonique $V_{n+1}\otimes_{\oo_{n+1}}\oo_{n}\xrightarrow{\sim} V_{n}$ soit un isomorphisme de $\oo_{n}$-représentations de $G$.
	
	Le foncteur limite projective induit une équivalence de catégories
\begin{equation}
	\Rep_{\breve{\oo}}^{\ltf}(G)\xrightarrow{\sim} \Rep_{\oo}^{\ltf}(G).
	\label{systeme p adique oo Rep oo}
\end{equation}
En effet, le foncteur $\Rep_{\oo}^{\ltf}(G)\to \Rep_{\breve{\oo}}^{\ltf}(G)$, défini par $V\mapsto (V/p^{n}V)_{n\ge 1}$ est un quasi-inverse de \eqref{systeme p adique oo Rep oo}. Pour tout objet $V$ de $\Rep_{\oo}^{\ltf}(G)$, on note $V[\frac{1}{p}]$ la $\mathfrak{C}$-représentation continue associée. Le foncteur canonique $\Rep_{\oo}^{\ltf}(G)\to \Rep_{\mathfrak{C}}^{\cont}(G)$ induit une équivalence de catégories:
	\begin{equation}
		\Rep_{\oo}^{\ltf}(G)_{\mathbb{Q}}\xrightarrow{\sim} \Rep_{\mathfrak{C}}^{\cont}(G),
		\label{Rep o Q eqi Rep C}
	\end{equation}
	où la source désigne la catégorie des objets de $\Rep_{\oo}^{\ltf}(G)$ à isogénie près (cf. \cite{DW03} prop. 22).
\end{nothing}

\begin{nothing}	\label{alpha iso cat Rep}
	Prenant pour $\varphi$ l'un des homomorphismes canoniques $\oo\to \End(\id_{\Rep_{\oo_{n}}(G)})$, $\oo\to \End(\id_{\Rep_{\breve{\oo}}(G)})$ ou $\oo\to \End(\id_{\Rep_{\oo}^{\cont}(G)})$, on peut considérer les notions introduites dans \ref{Cat Ab almost} pour les catégories abéliennes $\Rep_{\oo_{n}}(G)$, $\Rep_{\breve{\oo}}(G)$ et $\Rep_{\oo}^{\cont}(G)$. Un morphisme de $\Rep_{\oo_{n}}(G)$ (resp. $\Rep_{\oo}^{\cont}(G)$) est un $\alpha$-isomorphisme si et seulement s'il induit un $\alpha$-isomorphisme sur les modules sous-jacents. Un morphisme $(f_{n})_{n\ge 1}$ de $\Rep_{\breve{\oo}}(G)$ est un $\alpha$-isomorphisme si et seulement si, pour tout entier $n\ge 1$, $f_{n}$ est un $\alpha$-isomorphisme de $\Rep_{\oo_{n}}(G)$.
\end{nothing}

\begin{nothing} \label{oo oon rep aptf}
	Soit $n$ un entier $\ge 1$. On désigne par $\Rep^{\aptf}_{\oo_{n}}(G)$ la sous-catégorie pleine de $\Rep_{\oo_{n}}(G)$ formée des $\oo_{n}$-représentations dont le $\oo_{n}$-module sous-jacent est $\alpha$-plat de type $\alpha$-fini \eqref{coro alpha plat tf}, et par $\Rep^{\aptf}_{\oo}(G)$ la sous-catégorie pleine de $\Rep^{\cont}_{\oo}(G)$ formée des $\oo$-représentations continues de $G$ dont le $\oo$-module sous-jacent est $\alpha$-plat de type $\alpha$-fini.
	
	On désigne par $\Rep_{\breve{\oo}}^{\aptf}(G)$ la sous-catégorie pleine de $\Rep_{\breve{\oo}}(G)$ formée des systèmes projectifs de représentations $(V_{n})_{n\ge 1}$ tels que, pour tout entier $n\ge 1$, $V_{n}$ soit un objet de $\Rep_{\oo_{n}}^{\aptf}(G)$ et que le morphisme canonique $V_{n+1}\otimes_{\oo_{n+1}}\oo_{n}\rightarrow V_{n}$ soit un \textit{$\alpha$-isomorphisme} \eqref{alpha iso cat Rep}.

	En munissant $\mm$ de l'action triviale de $G$, d'après \ref{lemma aptf P oon} et \ref{topologie p adic limite proj aptf}, le foncteur limite projective induit un foncteur
\begin{equation}
	F_{1}: \Rep_{\breve{\oo}}^{\aptf}(G)\to \Rep^{\aptf}_{\oo}(G) \qquad (V_{n})_{n\ge 1}\mapsto \varprojlim \mm\otimes_{\oo}V_{n}.
	\label{foncteur aptf oon oo}
\end{equation}
D'autre part, on a un foncteur
\begin{equation}
	F_{2}: \Rep^{\aptf}_{\oo}(G) \to \Rep_{\breve{\oo}}^{\aptf}(G)\qquad V \mapsto (V/p^{n}V)_{n\ge 1}.
	\label{foncteur aptf oo oon}
\end{equation}

En vertu de \ref{lemma topologie aptf}(i)-(ii), on a un foncteur canonique
\begin{equation}
	\Rep^{\aptf}_{\oo}(G)\to \Rep^{\cont}_{\mathfrak{C}}(G) \qquad V\mapsto V[\frac{1}{p}].
	\label{foncteur canonique Rep oo aptf to C}
\end{equation}

	On désigne par $\Rep^{\aptf}_{\oo}(G)_{\mathbb{Q}}$ (resp. $\Rep^{\aptf}_{\breve{\oo}}(G)_{\mathbb{Q}}$) la catégorie des objets de $\Rep^{\aptf}_{\oo}(G)$ (resp. $\Rep^{\aptf}_{\breve{\oo}}(G)$) à isogénie près. 
\end{nothing}

\begin{lemma} 
	\textnormal{(i)} Le foncteur $F_{1}$ envoie les $\alpha$-isomorphismes sur des isomorphismes.

	\textnormal{(ii)} Soit $(V_{n})_{n\ge 1}$ un objet de $\Rep_{\breve{\oo}}^{\aptf}(G)$. On a un isomorphisme canonique fonctoriel
	\begin{displaymath}
		F_{2}(F_{1}( (V_{n})_{n\ge 1}))\xrightarrow{\sim} (\mm\otimes_{\oo}V_{n})_{n\ge 1}
	\end{displaymath}
	et un $\alpha$-isomorphisme canonique $(\mm\otimes_{\oo}V_{n})_{n\ge 1}\to (V_{n})_{n\ge 1}$.
	\label{F1 envoie alpha iso iso}
\end{lemma}
\textit{Preuve}. L'assertion (i) résulte de \ref{prop Tsuji alpha iso}(ii) et l'assertion (ii) résulte de \eqref{mmooM to M} et \ref{lemma aptf P oon}(ii).

\begin{lemma}
	Soit $V$ un objet de $\Rep_{\oo}^{\aptf}(G)$. On munit $\mm\otimes_{\oo}V$ de la topologie $p$-adique et de l'action canonique de $G$. Alors:

	\textnormal{(i)} On a un isomorphisme canonique fonctoriel $\mm\otimes_{\oo}V\xrightarrow{\sim} F_{1}(F_{2}(V))$.	

	\textnormal{(ii)} Le $\alpha$-isomorphisme canonique $\mm\otimes_{\oo}V\to V$ est injectif.
	\label{lemma initial aptf rep}
\end{lemma}
\textit{Preuve}. (i) D'après \ref{lemma topologie aptf}(iii), on a un isomorphisme canonique et fonctoriel de $\oo$-représentations continues $\mm\otimes_{\oo}V\xrightarrow{\sim} \widehat{\mm\otimes_{\oo}V} \simeq \varprojlim \mm\otimes_{\oo}(V/p^{n}V)$; d'où l'assertion.

(ii) Le noyau du morphisme $\mm\otimes_{\oo}V\to V$ est $\alpha$-nul. La proposition résulte du fait que le $\oo$-module $\mm\otimes_{\oo}V$ n'a pas de torsion \eqref{lemma alpha plat plat mm}.

\begin{coro}\label{coro Rep aptf oo aptf oo Q}
	Les foncteurs $F_{1}$ et $F_{2}$ induissent des équivalences de catégories quasi-inverses l'une de l'autre
	\begin{equation}
		\Rep^{\aptf}_{\oo}(G)_{\mathbb{Q}}\xrightarrow{\sim} \Rep^{\aptf}_{\breve{\oo}}(G)_{\mathbb{Q}},\qquad \Rep^{\aptf}_{\breve{\oo}}(G)_{\mathbb{Q}}\xrightarrow{\sim} \Rep^{\aptf}_{\oo}(G)_{\mathbb{Q}}.
	\end{equation}
\end{coro}
\textit{Preuve}. En vertu de \ref{lemma 143}, tout $\alpha$-isomorphisme est une isogénie. L'assertion résulte alors de \ref{F1 envoie alpha iso iso} et \ref{lemma initial aptf rep}.

\begin{prop} \label{prop oo aptf C rep}
	Le foncteur canonique \eqref{foncteur canonique Rep oo aptf to C} induit une équivalence de catégories
	\begin{equation}
		\Rep^{\aptf}_{\oo}(G)_{\mathbb{Q}}\to \Rep_{\mathfrak{C}}^{\cont}(G).
		\label{equivalence oo aptf Crep}
	\end{equation}
\end{prop}
\textit{Preuve}. L'essentielle surjectivité résulte de \eqref{Rep o Q eqi Rep C}. Soient $V_{1}$ et $V_{2}$ deux objets de $\Rep_{\oo}^{\aptf}(G)$. Il suffit de démontrer que le morphisme canonique 
\begin{equation} \label{isomorphisme aptf oo C}
	\Hom_{G}(V_{1},V_{2})\otimes_{\mathbb{Z}}\mathbb{Q}\to \Hom_{G}(V_{1}[\frac{1}{p}],V_{2}[\frac{1}{p}])
\end{equation}
est un isomorphisme. D'après \ref{lemma initial aptf rep}(ii), on a une suite exacte de $\Rep_{\oo}^{\cont}(G)$
\begin{equation}
	0\to \mm\otimes_{\oo}V_{2}\to V_{2}\to k\otimes_{\oo}V_{2}\to 0.
\end{equation}
On en déduit une suite exacte
\begin{equation}
	0\to \Hom_{G}(V_{1},\mm\otimes_{\oo}V_{2})\to \Hom_{G}(V_{1},V_{2})\to \Hom_{G}(V_{1},k\otimes_{\oo}V_{2}).
\end{equation}
Comme $k\otimes_{\oo}V_{2}$ est $\alpha$-nul, on a des isomorphismes $(\mm\otimes_{\oo}V_{2})[\frac{1}{p}]\simeq V_{2}[\frac{1}{p}]$ et $\Hom_{G}(V_{1},\mm\otimes_{\oo}V_{2})\otimes_{\mathbb{Z}}\mathbb{Q}\to \Hom_{G}(V_{1},V_{2})\otimes_{\mathbb{Z}}\mathbb{Q}$. D'après \ref{lemma alpha plat plat mm}, on peut donc se ramener au cas où $V_{2}$ est sans-torsion.

On note $V_{1,\tor}$ le sous-module de torsion de $V_{1}$ et on le munit de la topologie $p$-adique et de l'action continue de $G$ induite par $V_{1}$ \eqref{cat ab Rep oo cont}. On a alors des isomorphismes
\begin{equation}
	\Hom_{G}(V_{1},V_{2})\simeq \Hom_{G}(V_{1}/V_{1,\tor},V_{2})
\end{equation}
et $V_{1}[\frac{1}{p}]\simeq (V_{1}/V_{1,\tor})[\frac{1}{p}]$. On peut donc se ramener au cas où $V_{1}$ est sans-torsion.

Comme $V_{1}$ et $V_{2}$ sont isogènes à des $\oo$-modules libres de type fini \eqref{coro alpha plat tf}, on a un isomorphisme canonique
\begin{equation}
	\Hom_{\oo}(V_{1},V_{2})\otimes_{\mathbb{Z}}\mathbb{Q}\xrightarrow{\sim} \Hom_{\mathfrak{C}}(V_{1}[\frac{1}{p}],V_{2}[\frac{1}{p}]).
	\label{morphisme fully faithful aptf}
\end{equation}
On en déduit l'injectivité de \eqref{isomorphisme aptf oo C}. \'{E}tant donnée une application $\mathfrak{C}$-linéaire et $G$-équivariante $f:V_{1}[\frac{1}{p}]\to V_{2}[\frac{1}{p}]$, il existe une application $\oo$-linéaire $g:V_{1}\to V_{2}$ telle que $g[\frac{1}{p}]=f$. Comme les actions de $G$ sur $\oo$ et $\mathfrak{C}$ sont triviales et les $\oo$-modules sous-jacents à $V_{1}$ et $V_{2}$ n'ont pas de torsion, l'application $g$ est $G$-équivariante; d'où la surjectivité de \eqref{isomorphisme aptf oo C}.

\begin{lemma} \label{lemma extensions sysl}
	Soient $X$ un schéma connexe, $\overline{x}$ un point géométrique de $X$. Pour tout entier $n\ge 1$, on note encore $\oo_{n}$ le faisceau d'anneaux constant de $X_{\fet}$ de valeur $\oo_{n}$ et on note $\breve{\oo}$ l'anneau $(\oo_{n})_{n\ge 1}$ de $X_{\fet}^{\mathbb{N}^{\circ}}$.

	\textnormal{(i)} Pour qu'un $\oo_{n}$-module $\mathbb{L}$ de $X_{\fet}$ soit localement libre de type fini, il faut et il suffit que sa fibre $\mathbb{L}_{\overline{x}}$ en $\overline{x}$ soit un $\oo_{n}$-module libre de type fini.

	\textnormal{(ii)} Pour qu'un $\breve{\oo}$-module $\mathbb{L}=(\mathbb{L}_{n})_{n\ge 1}$ soit localement libre de type fini, il faut et il suffit que $\mathbb{L}$ soit adique et que pour tout entier $n\ge 1$, le $\oo_{n}$-module $\mathbb{L}_{n}$ soit localement libre de type fini.
\end{lemma}
\textit{Preuve}. L'assertion (i) est démontrée dans (\cite{AGT} III.2.11). En calquant la démonstration de (\cite{AGT} III.7.19, 7.20), on vérifie l'assertion (ii). 

\begin{nothing} \label{general notion of rep sysl}
	Conservons les notations de \ref{lemma extensions sysl} et reprenons de plus celles de \ref{oo oon reprentations} pour $G=\pi_{1}(X,\overline{x})$. Soit $n$ un entier $\ge 1$. La restriction du foncteur $\mu_{\overline{x}}:\mathbb{B}_{\pi_{1}(X,\overline{x})}\xrightarrow{\sim} X_{\fet}$ \eqref{foncteur mu} à la sous-catégorie $\Rep_{\oo_{n}}(\pi_1(X,\overline{x}))$ de $\mathbf{B}_{\pi_{1}(X,\overline{x})}$ induit une équivalence de catégories:
	\begin{equation}
		\Rep_{\oo_n}(\pi_1(X,\overline{x})) \xrightarrow{\sim} \Mod(X_{\fet},\oo_n).
		\label{Rep on to Mod on general}
	\end{equation}
	D'après \ref{lemma extensions sysl}(i),	celle-ci induit une équivalence de catégories (\ref{LPtf}, \ref{oo oon reprentations})
	\begin{equation}
		\mu_{\overline{x}}:\Rep_{\oo_n}^{\ltf}(\pi_1(X,\overline{x})) \xrightarrow{\sim} \Sysl(X_{\fet},\oo_n).
		\label{Rep on to Mod on}
	\end{equation}
	D'après \eqref{systeme p adique oo Rep oo} et \ref{lemma extensions sysl}(ii), les foncteurs \eqref{Rep on to Mod on} induisent une équivalence de catégories
\begin{equation}
	\breve{\mu}_{\overline{x}}: \Rep_{\oo}^{\ltf}(\pi_{1}(X,\overline{x})) \xrightarrow{\sim} \Sysl(X^{\mathbb{N}^{\circ}}_{\fet},\breve{\oo}).
	\label{Rep o to Mod o}
\end{equation}
On désigne par $\Sysl_{\mathbb{Q}}(X^{\mathbb{N}^{\circ}}_{\fet},\breve{\oo})$ la catégorie des $\breve{\oo}$-modules localement libres de type fini de $X_{\fet}^{\mathbb{N}^{\circ}}$ à isogénie près. D'après \eqref{Rep o Q eqi Rep C}, le foncteur \eqref{Rep o to Mod o} induit une équivalence de catégories
\begin{equation}
	\breve{\mu}_{\overline{x},\mathbb{Q}}: \Rep_{\mathfrak{C}}^{\cont}(\pi_{1}(X,\overline{x}))\xrightarrow{\sim} \Sysl_{\mathbb{Q}}(X^{\mathbb{N}^{\circ}}_{\fet},\breve{\oo}).
	\label{Rep C to Mod C}
\end{equation}
\end{nothing}

\section{Compléments de géométrie rigide} \label{Geo rigide}
\begin{nothing} \label{XX rig foncteur}	
	Soit $\mathscr{S}=\Spf(\oo)$. Pour tout $\mathscr{S}$-schéma formel de présentation finie $\XX$, on désigne par $\XX^{\rig}$ l'espace rigide cohérent associé (\cite{Ab10} 4.1.6) et par $\XX^{\rig}_{\ad}$ le topos admissible de $\XX^{\rig}$ (\cite{Ab10} 4.4.1). On pose $\Theta=\mathscr{S}^{\rig}$. Pour tout $\mathscr{O}_{\XX}$-module $\mathscr{F}$, on peut lui associer un faisceau $\mathscr{F}^{\rig}$ de $\XX^{\rig}_{\ad}$ (\cite{Ab10} 4.7.4). Le faisceau $(\mathscr{O}_{\XX})^{\rig}$ est un anneau de $\XX_{\ad}^{\rig}$ que l'on note $\mathscr{O}_{\XX^{\rig}}$ (\cite{Ab10} 4.7.5). On a un morphisme canonique de topos annelés (\cite{Ab10} 4.7.5)
	\begin{equation}
		\rho_{\XX}:(\XX_{\ad}^{\rig},\mathscr{O}_{\XX^{\rig}})\to (\XX_{\zar},\mathscr{O}_{\XX}).
	\end{equation}
	Pour tout $\mathscr{O}_{\XX}$-module cohérent $\mathscr{F}$, on a un isomorphisme canonique (\cite{Ab10} 4.7.25)	
	\begin{equation}
		\mathscr{F}^{\rig}\xrightarrow{\sim} \rho_{\XX}^{*}(\mathscr{F}).
		\label{F rig rho XX pullback}
	\end{equation}
\end{nothing}
\begin{nothing}
	Soit $V$ un $\mathfrak{C}$-schéma séparé de type fini. On peut lui associer un $\Theta$-espace rigide quasi-séparé $V^{\an}$ (\cite{Ab10} 7.4.12). On désigne par $V^{\an}_{\ad}$ le topos admissible de $V^{\an}$ (\cite{Ab10} 7.3.4) et $\mathscr{O}_{V^{\an}}$ l'anneau structural de $V^{\an}$ (\cite{Ab10} 7.3.8). On a un morphisme de topos annelés (le morphisme de GAGA \cite{Ab10} (7.6.2.4)):
	\begin{equation}
		\Phi_{V}: (V^{\an}_{\ad},\mathscr{O}_{V^{\an}}) \to (V_{\zar},\mathscr{O}_{V}).
	\end{equation}
	Si $F$ est un $\mathscr{O}_{V}$-module, on pose $F^{\an}=\Phi_{V}^{*}(F)$ (l'image réciproque étant pris au sens des modules). L'anneau $\mathscr{O}_{V^{\an}}$ est cohérent (cf. \cite{Ab10} 7.3.10 et 7.3.11). On désigne par $\Mod^{\coh}(\mathscr{O}_{V})$ (resp. $\Mod^{\coh}(\mathscr{O}_{V^{\an}})$) la catégorie des $\mathscr{O}_{V}$-modules (resp. $\mathscr{O}_{V^{\an}}$-modules) cohérents (\cite{Ab10} 7.3.10).
	\label{V C analytique}
\end{nothing}
\begin{nothing}
	Soient $X$ un $\check{\overline{S}}$-schéma séparé et de présentation finie \eqref{notations 11}, $\widehat{X}$ le schéma formel complété $p$-adique de $X$ et $V=X\otimes_{\oo}\mathfrak{C}$. On a un $\Theta$-morphisme canonique (\cite{Ab10} 7.4.12.2)
	\begin{equation}
		\widehat{X}^{\rig}\to V^{\an}
		\label{morphisme modele rig to an}
	\end{equation}
	Si $X$ est propre sur $\check{\overline{S}}$, celui-ci est un isomorphisme et par suite $V^{\an}$ est cohérent (\cite{Ab10} 7.4.14). Soient $\mathscr{F}$ un $\mathscr{O}_{X}$-module cohérent et $\widehat{\mathscr{F}}$ son complété $p$-adique. On a un isomorphisme fonctoriel (\cite{Ab10} 7.6.7)
	\label{X XX V}
\begin{equation}
	\widehat{\mathscr{F}}^{\rig}\xrightarrow{\sim} (\mathscr{F}|V)^{\an}|\widehat{X}^{\rig}.
	\label{F rig et an}
\end{equation}
\end{nothing}
\begin{prop}
	Soient $V$ un $\mathfrak{C}$-schéma séparé de type fini et $F$ un $\mathscr{O}_{V}$-module cohérent. Si $F^{\an}=0$, $F$ est également nul.
	\label{Fan nul Fnul}
\end{prop}
\textit{Preuve}. Le support de $F$ est un sous-ensemble fermé de $V$. Comme chaque sous-ensemble fermé de $V$ a un point fermé, il suffit de démontrer que, pour tout point fermé $P$ de $V$, $F_{P}$ est nul.

Soient $P$ un $\mathfrak{C}$-point de $V$ et $Q:\Theta\to V^{\an}$ le morphisme d'espaces rigides associé (\cite{Ab10} 7.4.12). D'après (\cite{Ab10} 7.4.12.3), il existe un $\check{\overline{S}}$-schéma $X$ séparé de présentation finie et de fibre générique $V$ tel que $Q$ se factorise en un morphisme $\Theta\to \widehat{X}^{\rig}$ que l'on note encore $Q$, suivi du morphisme canonique $\widehat{X}^{\rig}\to V^{\an}$ \eqref{morphisme modele rig to an}. D'après (\cite{Ab10} 3.3.12 et 4.1.20), il existe un unique point rigide $\mathscr{Q}:\mathscr{S}\to \widehat{X}$ (\cite{Ab10} 3.3.1) tel que $Q=\mathscr{Q}^{\rig}$.

Pour démontrer que $F_{P}$ est nul, on peut se ramener au cas où $X=\Spec(A)$ est affine. On note $\widehat{A}$ le complété $p$-adique de $A$ et on pose $B=A[\frac{1}{p}]$, $B'=\widehat{A}[\frac{1}{p}]$ et $f: B\to B'$ le morphisme canonique. Le point rigide $\mathscr{Q}$ de $\widehat{X}$ correspond à un idéal maximal $\mathfrak{q}$ de $B'$ (cf. \cite{Ab10} 3.3.2). D'après (\cite{Ab10} 1.12.11(i)), $\mathfrak{p}=f^{-1}(\mathfrak{q})$ est l'idéal maximal de $B$ correspondant au point fermé $P$.

Il existe un $\mathscr{O}_{X}$-module cohérent $\mathscr{F}$ tel que $\mathscr{F}\otimes_{\oo}\mathfrak{C}\simeq F$ (\cite{Ab10} 2.6.23). En vertu de \eqref{F rig et an}, on a un isomorphisme $\widehat{\mathscr{F}}^{\rig}\xrightarrow{\sim} F^{\an}|\widehat{X}^{\rig}$. On en déduit que $\widehat{\mathscr{F}}^{\rig}=0$. Il existe un $A$-module cohérent $M$ tel que $\mathscr{F}$ soit isomorphe au $\mathscr{O}_{X}$-module associé à $M$. Notons $\widehat{M}$ le complété $p$-adique de $M$. D'après (\cite{Ab10} 1.12.16), le couple $(A,pA)$ vérifie la condition de Krull (\cite{Ab10} 1.8.25) et on a alors un isomorphisme
\begin{equation}
	M\otimes_{A}\widehat{A}\xrightarrow{\sim}\widehat{M}.
\end{equation}
D'après (\cite{Ab10} 4.8.24), le $B'$-module $\widehat{M}[\frac{1}{p}]$ est nul. On en déduit que
\begin{equation}
	(M[\frac{1}{p}])_{\mathfrak{p}}\otimes_{B_{\mathfrak{p}}}B'_{\mathfrak{q}}\simeq (\widehat{M}[\frac{1}{p}])_{\mathfrak{q}}
\end{equation}
est nul. D'après (\cite{Ab10} 1.12.17), le morphisme canonique $A\to \widehat{A}$ est plat. Il en est même du morphisme local $B_{\mathfrak{p}}\to B'_{\mathfrak{q}}$. En vertu de (\cite{EGAInew} 6.6.1), $B'_{\mathfrak{q}}$ est fidèlement plat sur $B_{\mathfrak{p}}$. Le $B_{\mathfrak{p}}$-module $(M[\frac{1}{p}])_{\mathfrak{p}}\simeq F_{P}$ est alors nul, d'où la proposition.

\begin{coro}
	Sous les hypothèses de \ref{Fan nul Fnul}. Le foncteur
	\begin{equation}
		\Phi_{V}^{*}:\Mod^{\coh}(\mathscr{O}_{V})\to \Mod^{\coh}(\mathscr{O}_{V^{\an}})\qquad F\mapsto F^{\an}
		\label{Phi V star}
	\end{equation}
	est exact et fidèle.
	\label{Phi V fidelement plat coh}
\end{coro}

\textit{Preuve}. Cela résulte de \ref{Fan nul Fnul} et (\cite{Ab10} 7.6.8).

\begin{prop}
	Soient $V$ un schéma séparé de type fini sur $\mathfrak{C}$ et $F$ un $\mathscr{O}_{X}$-module cohérent. Pour que $F$ soit plat, il faut et il suffit que $F^{\an}$ soit un $\mathscr{O}_{V^{\an}}$-module plat.
	\label{Zar Analytique localement libre}
\end{prop}

\textit{Preuve}. Si $F$ est plat, c'est alors un fibré vectoriel sur $V$. On en déduit que $F^{\an}$ est localement libre de type fini et est plat.

Supposons que $F^{\an}$ soit plat. Soient $U$ un ouvert de $V$, $F_{U}$ la restriction de $F$ à $U$ et $F_{U}^{\an}=\Phi_{U}^{*}(F)$ \eqref{Phi V star}. Le $\mathscr{O}_{U^{\an}}$-module $F_{U}^{\an}$ s'identifie à la restriction de $F^{\an}$ sur $U^{\an}$. Comme $U^{\an}_{\ad}$ a suffisamment de points (\cite{Ab10} 7.3.14), $F^{\an}_{U}$ est un $\mathscr{O}_{U^{\an}}$-module plat. Comme le foncteur $\Phi_{U}^{*}$ \eqref{Phi V star} est exact et fidèle, on en déduit que le foncteur
\begin{equation}
	F_{U}\otimes_{\mathscr{O}_{U}}-: \Mod^{\coh}(\mathscr{O}_{U})\to \Mod^{\coh}(\mathscr{O}_{U})
\end{equation}
est exact. La proposition s'ensuit compte tenu de (\cite{Ab10} 1.3.17).

\section{Rappel sur les courbes} \label{courbes rel}
\begin{nothing}
Dans cet article, \textit{une courbe sur un corps $F$} est un schéma séparé de type fini et \textit{géométriquement connexe} sur $F$ de dimension $1$. Une courbe propre sur un corps algébriquement clos est dite \textit{semi-stable} si elle est réduite et a pour seules singularités des points doubles ordinaires. On notera que cette définition est moins restrictive que la définition classique puisqu'on n'impose aucune condition sur l'intersection avec les composantes rationnelles.
\label{courbe sur corps}
\end{nothing}
\begin{nothing}
Soient $T=\Spec(R)$ le spectre d'un anneau de valuation $R$ de hauteur $1$, de corps des fractions $F$, et $\overline{F}$ une clôture algébrique de $F$. On note $\tau$ (resp. $t$) le point générique (resp. fermé) de $T$ et $\overline{\tau}$ le point géométriquement générique correspondant à $\overline{F}$.

Une \textit{$T$-courbe} est un $T$-schéma plat, séparé, de présentation finie et de dimension relative $1$ dont la fibre générique est une courbe \textit{lisse} sur $F$. Une $T$-courbe $f:X\to T$ est \textit{propre (resp. régulière, resp. normale)} si $f$ est propre (resp. $X$ est régulier, resp. $X$ est normal). D'après (\cite{Liu} \S~4 3.8), toute $T$-courbe est irréductible et réduite. D'après la factorisation de Stein (\cite{EGA III} 4.3.1), la fibre spéciale d'une $T$-courbe propre est géométriquement connexe et est donc une $t$-courbe propre (\ref{courbe sur corps}). On dit qu'une $T$-courbe $X_{1}$ \textit{domine} une autre $T$-courbe $X_{2}$ s'il existe un $T$-morphisme dominant de $X_{1}$ dans $X_{2}$.

On dit qu'une $T$-courbe est \textit{semi-stable} si elle est propre et si sa fibre spéciale géométrique est une courbe semi-stable. On dit qu'une $T$-courbe $X$ est \textit{cohomologiquement plate en dimension $0$} si elle est propre et si le faisceau $\mathscr{O}_{X}$ est cohomologiquement plat sur $T$ en dimension $0$, i.e. si le morphisme structural $\lambda: X\to T$ induit un isomorphisme $\mathscr{O}_{T}\xrightarrow{\sim}\lambda_{*}(\mathscr{O}_{X})$ universellement (\cite{EGA III} 7.8.1). On notera qu'une $T$-courbe propre à fibres géométriquement réduites (en particulier une $T$-courbe semi-stable) est cohomologiquement plate en dimension $0$ (cf. \cite{EGA III} 7.8.6).

Une $T$-courbe propre $X$ est appelée \textit{modèle} (resp. $T$-modèle) de sa fibre générique (resp. de sa fibre géométrique générique). Soit $C$ une courbe propre et lisse sur $\overline{F}$. On dit qu'un $T$-modèle $X_{1}$ de $C$ \textit{domine} un autre $T$-modèle $X_{2}$ de $C$, s'il existe un $T$-morphisme de $X_{1}$ dans $X_{2}$ induisant un isomorphisme entre les fibres génériques.
\label{general not of curves}
\end{nothing}

$\hspace*{-1.2em}\bf{\arabic{section}.\stepcounter{theorem}\arabic{theorem}.}$
Pour toute $T$-courbe semi-stable $X$, le foncteur de Picard relatif $\Pic^0_{X/T}$ (\ref{not of Picard}) est représentable par un $T$-schéma semi-abélien (cf. \cite{BLR90} 9.4.1). Si $T=\check{\overline{S}}$ \eqref{notations 11}, on a $\Pic(\check{\overline{S}})=0$ et donc un isomorphisme canonique $\Pic(X)\xrightarrow{\sim} \Pic_{X/\check{\overline{S}}}(\check{\overline{S}})$ (cf. \cite{BLR90} 8.1 Prop.4).

\begin{definition}
	Soient $(T,\tau)$ comme dans \ref{general not of curves} et $X$ une $T$-courbe propre.
	
	(i) On appelle \textit{$\tau$-revêtement de $X$} la donnée d'un $T$-morphisme propre de présentation finie de $T$-courbes $\varphi:Y\to X$ tel que le morphisme $\varphi_{\tau}: Y_{\tau}\to X_{\tau}$ soit un revêtement étale.
	
	(ii) Soit $\varphi:Y\to X$ un $\tau$-revêtement. Si la $T$-courbe propre $Y$ est cohomologiquement plate en dimension $0$ (resp. a des fibres géométriquement réduites, resp. est semi-stable, resp. est régulière), on dit que $\varphi: Y\to X$ est \textit{un $\tau$-revêtement cohomologiquement plat en dimension $0$} (resp. \textit{à fibres géométriquement réduites}, resp. \textit{semi-stable}, resp. \textit{régulier}) et en abrégé un $\tau$-revêtement c.p.d.$0$.
	
	(iii) On dit qu'un $\tau$-revêtement $\varphi:Y\to X$ est \textit{fini} si le morphisme sous-jacent $Y\to X$ est fini.
	
	(iv) Soit $\varphi:Y\to X$ un $\tau$-revêtement. Si $\varphi_{\tau}:Y_{\tau}\to X_{\tau}$ est un revêtement étale et galoisien de groupe des automorphismes $G$ et si l'action de $G$ sur $Y_{\tau}$ s'étend en une $X$-action sur $Y$, on dit que $\varphi:Y\to X$ est un \textit{$\tau$-revêtement galoisien}.
	\label{def of revêtement}
\end{definition}

Soient $\varphi_{1}:Y_{1}\to X$ et $\varphi_{2}: Y_{2}\to X$ deux $\tau$-revêtements de $X$. \textit{Un morphisme de $\varphi_{1}$ dans $\varphi_{2}$} est un $X$-morphisme $\pi:Y_{1}\to Y_{2}$. On notera alors que $Y_{1}$ domine $Y_{2}$ et on dira que $\varphi_{1}$ \textit{domine} $\varphi_{2}$. Supposons de plus, que $\varphi_{1}$ et $\varphi_{2}$ soient galoisiens. Alors, le morphisme $\pi_{\tau}:Y_{1,\tau}\to Y_{2,\tau}$ est un revêtement galoisien. Notons $u: \Aut(Y_{1,\tau}/X_{\tau})\to \Aut(Y_{2,\tau}/X_{\tau})$ l'homomorphisme induit par $\pi_{\tau}$. On dit que $\pi$ est \textit{équivariant} s'il est $u$-équivariant.

\begin{theorem}[Grothendieck; cf. \cite{Ab00} pour $g(C)\ge 2$ et \cite{Liu} \S~10.2 pour $g(C)=1$]
	Pour toute $\overline{K}$-courbe propre et lisse $C$ \eqref{notations 11}, il existe un trait $S'$ fini sur $S$ tel que $C$ admette un $S'$-modèle semi-stable.
	\label{semistable reduction theorem}
\end{theorem}

\begin{nothing}
	Soit $X$ une $S$-courbe propre. Une \textit{désingularisation de $X$} est un modèle régulier de $X_{\eta}$ dominant $X$. On sait qu'il existe une désingularisation de $X$ (cf. \cite{Liu} \S~8 3.44). Si le genre de $X_{\eta}$ est $\ge 1$, grâce au critère de contraction de Castelnuovo (cf. \cite{Liu} \S~9 3.21), on peut construire \textit{une désingularisation minimale $Y$ de $X$} unique à isomorphisme près, c'est-à-dire que toute désingularisation $Y'\to X$ se factorise d'une façon unique à travers $Y$. D'après (\cite{Liu} \S~9 3.13), tout $\eta$-automorphisme de $X_{\eta}$ s'étend en un $S$-automorphisme de la désingularisation minimale de $X$.
	\label{S courbe propriete}
\end{nothing}
\begin{prop}[\cite{Liu06} 2.7 et 2.8] Soit $X$ une $S$-courbe propre. Il existe un trait $S'$ fini sur $S$ tel que $X_{S'}$ soit dominé par un $S'$-modèle semi-stable et régulier $X'$ de $X_{\overline{\eta}}$. Si le genre de $X_{\eta}$ est $\ge 1$, on peut trouver $S'$ de sorte que $X'$ soit la désingularisation minimale de $X_{S'}$.
	\label{semistable rcm}
\end{prop}

\begin{coro}
	Soient $C$ une $\overline{K}$-courbe propre et lisse, $X_{1}$ et $X_{2}$ deux $S$-modèles de $C$. Il existe un trait $S'$ fini sur $S$ et un $S'$-modèle $X_{3}$ semi-stable et régulier de $C$ dominant $X_{1,S'}$ et $X_{2,S'}$.
	\label{semistable rcm coro}
\end{coro}
\textit{Preuve}. Il existe un $S$-modèle $X$ de $C$ dominant $X_{1}$ et $X_{2}$ (\cite{DW03} prop. 20). L'assertion résulte alors de \ref{semistable rcm}.

\begin{nothing}
	On dit qu'un $S$-schéma $X$ a une \textit{réduction semi-stable} si, localement pour la topologie étale, $X$ est isomorphe à un $S$-schéma de la forme
\begin{equation}
	\Spec(\mathcal{O}_{K}[t_1,\dots,t_{m}]/(t_1\dots t_n-\pi)),
	\label{semistable redction form}
\end{equation}
où $\pi$ est une uniformisante de $\mathcal{O}_{K}$, $m$ et $n$ sont des entiers vérifiant $m\ge n\ge 1$. On notera qu'une $S$-courbe propre admet une réduction semi-stable si et seulement si elle est semi-stable et régulière (cf. \cite{Sa04} 1.6).
\label{def reduction semistable}
\end{nothing}
\begin{nothing}
	Soit $X$ une $S$-courbe à fibres géométriquement réduites. D'après (\cite{Liu} \S~4 1.18), pour tout trait $S'$ fini sur $S$, le schéma $X_{S'}$ est normal. On en déduit par passage à la limite projective que le schéma $\overline{X}=X\times_{S}\overline{S}$ \eqref{basic notation} est normal.
	\label{X bar normal}
\end{nothing}
\begin{prop} \label{domine par bon revêtement}
	Soit $X$ une $S$-courbe propre.

	\textnormal{(i)} \'{E}tant donnés un nombre fini de $\eta$-revêtements (resp. $\eta$-revêtements galoisiens) $\varphi_i:Y_i\to X$, il existe un trait $(S',\eta')$ fini sur $(S,\eta)$, un $\eta'$-revêtement semi-stable, régulier et galoisien $\varphi:Y\to X_{S'}$ et pour chaque $\varphi_i$, un morphisme (resp. un morphisme équivariant) de $\varphi$ dans $\varphi_{i}\times_{S}S'$.
	
	\textnormal{(ii)} Pour tout $\eta$-revêtement \textrm{fini} (resp. \textrm{fini et galoisien}) $\varphi:Y\to X$, il existe un trait $(S',\eta')$ fini sur $(S,\eta)$ tel que $\varphi\times_{S}S'$ soit dominé par un $\eta'$-revêtement \textrm{fini}, à fibres géométriquement réduites (resp. \textrm{fini, galoisien}, à fibres géométriquement réduites) de $X_{S'}$.
	
	\textnormal{(iii)} Pour tout $\eta$-revêtement $\varphi:Y\to X$, il existe un trait $(S',\eta')$ fini sur $(S,\eta)$, un $S'$-modèle semi-stable $X'$ de $X_{\overline{\eta}}$ dominant $X_{S'}$ et un $\eta'$-revêtement \textrm{semi-stable, régulier, fini et galoisien} $\varphi':Y'\to X'$ qui s'insère dans un diagramme commutatif
	\begin{equation}
		\xymatrix{
			Y'\ar[r]^{\varphi'} \ar[d] & X'\ar[d]\\
			Y_{S'}\ar[r]^{\varphi_{S'}} & X_{S'}.
		}
		\label{revêtement galoisien domine revêtement}
	\end{equation}
\end{prop}

\textit{Preuve}. (i) Si le genre de $X_{\eta}$ est nul, $X_{\overline{\eta}}$ n'admet aucun revêtement étale non-trivial. Dans ce cas, l'énoncé (i) résulte de \ref{semistable rcm} appliqué aux produits des $Y_{i}$ au-dessus de $X$. Pour $g(X_{\eta})\ge 1$, cela résulte de \ref{semistable rcm} et des propriétés de la désingularisation minimale, appliquées au produit des $Y_{i}$ au-dessus de $X$ (cf. \cite{DW03} Thm. 12 et \cite{DW05} Thm. 4).

(ii) D'après (\cite{Epp73} Thm. 2.0), il existe un trait $(S',\eta')$ fini sur $(S,\eta)$ tel que la normalisation $Y'$ de $Y_{S'}$ ait une fibre spéciale géométriquement réduite. Comme $Y$ est excellent, $Y'$ est fini sur $Y$. Supposons de plus que $\varphi:Y\to X$ soit galoisien. Comme $Y'$ est la clôture intégrale de $Y$ dans $Y\times_{S}\eta'$, l'action de $\Aut(Y_{\eta}/X_{\eta})$ sur $Y$ s'étend en une action sur $Y'$, i.e. $Y'\to X_{S'}$ est aussi un $\eta'$-revêtement galoisien.

(iii) D'après (i), quitte à remplacer $K$ par une extension finie, $\varphi$ est dominé par un $\eta$-revêtement semi-stable, régulier et galoisien $Y'\to X$. Notons $G=\Aut(Y'_{\eta}/X_{\eta})$. L'action de $G$ sur $Y_{\eta}'$ s'étend en une $X$-action sur $Y'$. Comme $Y'$ est projective (\cite{Lich68} Thm. 2.8), le quotient $Y'/G$ existe (\cite{SGAI} V.1); notons le $X'$. Calquant la démonstration de (\cite{Liu06} 3.8) (cf. aussi \cite{Ray90} Prop. 5), on vérifie que $X'$ est un modèle semi-stable de $X_{\eta}$ dominant $X$. Donc $Y'\to X'$ est un $\eta$-revêtement semi-stable, régulier, fini et galoisien. Comme le diagramme
\begin{equation}
		\xymatrix{
			Y'_{\eta}\ar[r] \ar@{=}[d] & X'_{\eta}\ar@{=}[d]\\
			Y_{\eta}\ar[r] & X_{\eta},
		}
\end{equation}
est commutatif, il en est de même de \eqref{revêtement galoisien domine revêtement} en vertu de (\cite{EGA I} 7.2.2.1).

\begin{coro} Soit $X$ une $\overline{S}$-courbe propre.
	
	\textnormal{(i)} \'{E}tant donnés un nombre fini de $\overline{\eta}$-revêtements de $X$, il existe un $\overline{\eta}$-revêtement semi-stable et galoisien de $X$ qui les domine tous.

	\textnormal{(ii)} Tout $\overline{\eta}$-revêtement $\varphi:Y\to X$ fini et galoisien est dominé par un $\overline{\eta}$-revêtement fini, galoisien et à fibres géométriquement réduites.
	\label{domine par bon rev bar}
\end{coro}

\textit{Preuve}. D'après la théorie de la descente (\cite{EGA IV} 8.8.3 et 8.10.5), quitte à remplacer $S$ par une extension finie, tout $\overline{\eta}$-revêtement (resp. $\overline{\eta}$-revêtement fini, resp. $\overline{\eta}$-revêtement fini et galoisien) de $X$ se descend en un $\eta$-revêtement (resp. $\eta$-revêtement fini, resp. $\eta$-revêtement fini et galoisien). Le corollaire résulte alors de \ref{domine par bon revêtement}.

%
%


\begin{nothing} \label{model curve vect}
	Soient $X$ une $\check{\overline{S}}$-courbe projective et de présentation finie et $C=X_{\check{\overline{\eta}}}$ \eqref{notations 11}. L'anneau $\oo$ est universellement cohérent (\cite{Ab10} 1.12.15). Par suite, l'anneau $\mathscr{O}_{X}$ est cohérent. On note $\Mod^{\coh}(\mathscr{O}_{X})$ (resp. $\Mod^{\coh}(\mathscr{O}_{C})$) la catégorie des $\mathscr{O}_{X}$-modules (resp. $\mathscr{O}_{C}$-modules) cohérents et $\Mod^{\coh}_{\mathbb{Q}}(\mathscr{O}_{X})$ la catégorie de $\mathscr{O}_{X}$-modules cohérents à isogénie près. D'après (\cite{Ab10} 2.6.23), tout $\mathscr{O}_{C}$-module cohérent se prolonge en un $\mathscr{O}_{X}$-module cohérent. \'{E}tant donnés deux objets $\mathcal{F}$ et $\mathcal{F}'$ de $\Mod^{\coh}(\mathscr{O}_{X})$, on a un isomorphisme
\begin{equation}
	\Hom_{\mathscr{O}_{X}}(\mathcal{F},\mathcal{F}')\otimes_{\mathbb{Z}}\mathbb{Q}\xrightarrow{\sim}\Hom_{\mathscr{O}_{C}}(\mathcal{F}_{\check{\overline{\eta}}},\mathcal{F}_{\check{\overline{\eta}}}').
	\label{Hom OX Hom OC pf}
\end{equation}
Celui-ci implique que $\mathcal{F}$ et $\mathcal{F}'$ sont isogènes dans $\Mod^{\coh}(\mathscr{O}_{X})$ si et seulement si leurs fibres génériques sont isomorphes sur $C$. Donc, le foncteur canonique $\Mod^{\coh}(\mathscr{O}_{X})\to \Mod^{\coh}(\mathscr{O}_{C})$ induit une équivalence de catégories
\begin{equation}
	\Mod^{\coh}_{\mathbb{Q}}(\mathscr{O}_{X})\xrightarrow{\sim} \Mod^{\coh}(\mathscr{O}_{C}).
	\label{module coh X Q to C}
\end{equation}
\end{nothing}
\begin{nothing} \label{jmathXX X coh}
	Conservons les notations de \ref{model curve vect}. On pose $\mathscr{S}=\Spf(\oo)$ et on note $\XX$ le schéma formel complété $p$-adique de $X$ qui est propre et de présentation finie sur $\mathscr{S}$. En particulier, $\mathscr{O}_{\XX}$ est cohérent (\cite{Ab10} 2.8.1). Pour tout $\mathscr{O}_{X}$-module $\mathcal{F}$, on note $\widehat{\mathcal{F}}$ son complété $p$-adique. On note $\Mod^{\coh}(\mathscr{O}_{\XX})$ la catégorie des $\mathscr{O}_{\XX}$-modules cohérents. Comme $X$ est projectif et de présentation finie sur $\check{\overline{S}}$, le foncteur
\begin{equation}
	\Mod^{\coh}(\mathscr{O}_{X}) \to \Mod^{\coh}(\mathscr{O}_{\XX}),\qquad \mathcal{F} \mapsto \widehat{\mathcal{F}}
	\label{Grothendieck existence formel}
\end{equation}
induit une équivalence de catégories (\cite{Ab10} 2.13.8).

Pour tout $\mathscr{O}_{\XX}$-module $\mathscr{F}$, on pose
\begin{equation}
	\mathscr{F}[\frac{1}{p}]=\mathscr{F}\otimes_{\mathbb{Z}_{p}}\mathbb{Q}_{p}.
	\label{F inverse p}
\end{equation}
On note $\Mod^{\coh}(\mathscr{O}_{\XX}[\frac{1}{p}])$ la catégorie des $\mathscr{O}_{\XX}[\frac{1}{p}]$-modules cohérents et $\Mod^{\coh}_{\mathbb{Q}}(\mathscr{O}_{\XX})$ la catégorie des $\mathscr{O}_{\XX}$-modules cohérents à isogénie près. Alors le foncteur canonique
\begin{equation}
	\Mod^{\coh}(\mathscr{O}_{\XX})\to \Mod^{\coh}(\mathscr{O}_{\XX}[\frac{1}{p}]), \qquad \mathscr{F}\mapsto \mathscr{F}[\frac{1}{p}],
\end{equation}
induit une équivalence de catégories (\cite{AGT} III.6.16)
\begin{equation}
	\Mod^{\coh}_{\mathbb{Q}}(\mathscr{O}_{\XX})\xrightarrow{\sim} \Mod^{\coh}(\mathscr{O}_{\XX}[\frac{1}{p}]).
	\label{Mod OXX Q to Mod OXXp}
\end{equation}
Choisissons des quasi-inverses de \eqref{Grothendieck existence formel} et \eqref{Mod OXX Q to Mod OXXp}. On désigne par $\jmath_{\XX}$ le foncteur composé de l'équivalence de catégories \eqref{module coh X Q to C} et des quasi-inverses de \eqref{Grothendieck existence formel}, \eqref{Mod OXX Q to Mod OXXp}
\begin{equation}
	\jmath_{\XX}: \Mod^{\coh}(\mathscr{O}_{\XX}[\frac{1}{p}]) \xrightarrow{\sim} \Mod^{\coh}(\mathscr{O}_{C}).
	\label{jmathXX coh}
\end{equation}
\end{nothing}

\begin{lemma} \label{lemma LPtf VB}
	Conservons les notations de \ref{jmathXX X coh}. Pour tout $\mathscr{O}_{\XX}[\frac{1}{p}]$-module localement projectif de type fini $\mathscr{F}$ \eqref{LPtf}, $\jmath_{\XX}(\mathscr{F})$ \eqref{jmathXX coh} est un fibré vectoriel sur $C$.
\end{lemma}
\textit{Preuve}. Posons $F=\jmath_{\XX}(\mathscr{F})$. On désigne par $\XX^{\rig}$ l'espace rigide cohérent associé au schéma formel $\XX$ \eqref{XX rig foncteur} et $C^{\an}$ l'espace rigide quasi-séparé associé au $\mathfrak{C}$-schéma propre $C$ \eqref{V C analytique}. On a un isomorphisme $\XX^{\rig}\xrightarrow{\sim} C^{\an}$ \eqref{X XX V} et un morphisme de topos annelés (\cite{Ab10} 4.7.5)
\begin{equation}
	\rho_{\XX}:(\XX_{\ad}^{\an},\mathscr{O}_{\XX^{\rig}})\to (\XX_{\zar},\mathscr{O}_{\XX}[\frac{1}{p}]).
\end{equation}
En vertu de \eqref{F rig rho XX pullback} et \eqref{F rig et an}, on a un isomorphisme
\begin{equation}
	\rho_{\XX}^{*}(\mathscr{F})\xrightarrow{\sim} F^{\an}.
\end{equation}
Le $\mathscr{O}_{C^{\an}}$-module $F^{\an}$ est alors localement projectif de type fini. Par suite, c'est un $\mathscr{O}_{C^{\an}}$-module plat \eqref{LPtf}. D'après \ref{Zar Analytique localement libre}, $F$ est plat. Le lemme s'ensuit.

\begin{lemma} \label{isomorphisme section inverse p}
	Conservons les notations de \ref{jmathXX X coh}. 
	
	\textnormal{(i)} Pour tout $\mathscr{O}_{\XX}$-module $\mathscr{H}$, on a un isomorphisme canonique
	\begin{equation}
		\Gamma(\XX,\mathscr{H}[\frac{1}{p}])\simeq \Gamma(\XX,\mathscr{H})\otimes_{\mathbb{Z}}\mathbb{Q}.
		\label{isomorphisme section inverse p eq}
	\end{equation}

	\textnormal{(ii)} Soient $\mathscr{F}$ un $\mathscr{O}_{\XX}[\frac{1}{p}]$-module cohérent et $F=\jmath_{\XX}(\mathscr{F})$. Pour tout entier $i\ge 0$, on a un isomorphisme canonique
	\begin{equation}
		\rH^{i}(\XX,\mathscr{F})\simeq \rH^{i}(C,F).
	\end{equation}
\end{lemma}
\textit{Preuve}. (i) Par platitude de $\mathscr{O}_{\XX}$ sur $\oo$, il suffit de démontrer l'isomorphisme \eqref{isomorphisme section inverse p eq} pour les ouverts affines. Dans ce cas, cela résulte de (\cite{Ab10} 2.10.5).

(ii) Il existe un $\mathscr{O}_{X}$-module cohérent $\mathcal{F}$ tel que $\mathscr{F}\simeq \widehat{\mathcal{F}}[\frac{1}{p}]$ et que $F\simeq \mathcal{F}_{\check{\overline{\eta}}}$. D'après (\cite{Ab10} 2.13.2), on a un isomorphisme canonique
\begin{equation}
	\rH^{i}(X,\mathcal{F})\simeq \rH^{i}(\XX,\widehat{\mathcal{F}}).
\end{equation}
On déduit par (i) un isomorphisme $\rH^{i}(\XX,\widehat{\mathcal{F}})\otimes_{\mathbb{Z}}\mathbb{Q}\simeq \rH^{i}(\XX,\mathscr{F})$. Par platitude, on a un isomorphisme canonique
\begin{equation}
	\rH^{i}(X,\mathcal{F})\otimes_{\mathbb{Z}}\mathbb{Q}\simeq \rH^{i}(C,F).
\end{equation}
L'assertion s'ensuit.

\section{Correspondance de Deninger-Werner} \label{Cor DW}
On rappelle que le corps résiduel $k$ de $\mathcal{O}_{K}$ est une clôture algébrique de $\mathbb{F}_{p}$ \eqref{notations 11}.
\begin{definition}
	(i) Soient $C$ une courbe propre et lisse sur un corps algébriquement clos de caractéristique $p$ et $\Fr_{C}:C\to C$ le Frobenius absolu de $C$. On dit qu'un fibré vectoriel $F$ sur $C$ est \textit{fortement semi-stable} si $(\Fr_C^{n})^{*}(F)$ est un fibré vectoriel semi-stable pour tout entier $n\ge 0$.

	(ii) Soient $C$ une courbe propre sur un corps algébriquement clos de caractéristique $p$ et $\pi:\widetilde{C}\to C$ la normalisation du sous-schéma réduit sous-jacent à $C$. On dit qu'un fibré vectoriel $F$ sur $C$ est \textit{fortement semi-stable de degré 0} si $\pi^*(F)$ est fortement semi-stable de degré $0$ sur chaque composante irréductible de $\widetilde{C}$.
	\label{fortement semistable}
\end{definition}
On notera qu'un fibré vectoriel semi-stable sur une courbe propre et lisse de genre $g\ge 2$ sur un corps de caractéristique $p$ n'est pas nécessairement fortement semi-stable (cf. \cite{Gi73}).

\begin{theorem}[\cite{DW05} Thm. 16, Thm. 17 et \cite{To07} 2.4]
	Soient $X$ une $\overline{S}$-courbe propre, $\check{X}=X\times_{\overline{S}}\check{\overline{S}}$ \eqref{notations 11} et $\mathcal{F}$ un fibré vectoriel sur $\check{X}$. Les conditions suivantes sont équivalentes:

	\textnormal{(i)} La fibre spéciale $\mathcal{F}_{s}$ de $\mathcal{F}$ est fortement semi-stable de degré $0$ sur $X_s$.

	\textnormal{(ii)} Pour tout $n\ge 1$, il existe un $\overline{\eta}$-revêtement $\varphi:X'\to X$ \textnormal{\eqref{def of revêtement}} tel que $\varphi_n^{*}(\mathcal{F}_n)$ soit libre de type fini, où $\mathcal{F}_{n}$ (resp. $\varphi_{n}$) est la réduction modulo $p^{n}$ de $\mathcal{F}$ (resp. $\varphi$) \eqref{basic notation}.

	\textnormal{(iii)} Il existe un entier $n\ge 1$ et un $\overline{\eta}$-revêtement $\varphi:X'\to X$ tels que $\varphi_n^{*}(\mathcal{F}_n)$ soit libre de type fini.
	\label{DW fortement semistable}
\end{theorem}

\begin{definition}
	Conservons les notations de \ref{DW fortement semistable}. On appelle \textit{fibré vectoriel de Deninger-Werner sur $\check{X}$} tout fibré vectoriel sur $\check{X}$ vérifiant les conditions équivalentes de \ref{DW fortement semistable}. On désigne par $\BB^{\DW}_{\check{X}}$ la sous-catégorie pleine de $\Vect_{\check{X}}$ \eqref{notations vect bundle} formée des fibrés vectoriels de Deninger-Werner.
	\label{def of categorie DW}
\end{definition}
\begin{prop}
	Soit $X$ une $\overline{S}$-courbe semi-stable. Un fibré en droites sur $\check{X}$ est de Deninger-Werner si et seulement si sa classe appartient à $\Pic^{0}_{\check{X}/\check{\overline{S}}}(\check{\overline{S}})$ \eqref{not of Picard}.
	\label{DW deg 1 semistable}
\end{prop}
\textit{Preuve}. D'après (\cite{DW05} Thm. 12), la catégorie $\BB^{\DW}_{\check{X}}$ contient tous les fibrés en droites dont les classes sont dans $\Pic^0_{\check{X}/\check{\overline{S}}}(\check{\overline{S}})$. Inversement, pour tout fibré en droites $L\in \BB^{\DW}_{\check{X}}$, l'image inverse de la fibre spéciale $L_{s}$ sur la normalisation de chaque composante irréductible de $X_{s}$ est de degré $0$. D'après (\cite{DW05} Thm.13), la fibre générique $L_{\check{\overline{\eta}}}$ de $L$ est donc aussi de degré $0$. Ceci implique que $L\in \Pic^{0}_{\check{X}/\check{\overline{S}}}(\check{\overline{S}})$ en vertu de (\cite{BLR90} 9.1.13).

\begin{prop}
	Soient $X$ une $\overline{S}$-courbe propre et lisse, $n$ un entier $\ge 1$ et $L$ un fibré en droites dont la classe appartient à $\Pic^0_{\check{X}/\check{\overline{S}}}(\check{\overline{S}})$. Alors, il existe un $\overline{\eta}$-revêtement fini, galoisien et à fibres géométriquement réduites de $X$ trivialisant $L_{n}$.
	\label{DW deg 1}
\end{prop}
\textit{Preuve}. Quitte à remplacer $K$ par une extension finie, il existe une $S$-courbe propre et lisse $Y$ telle que $X\simeq \overline{Y}$ et un point $x\in Y(S)$. On sait que $\Pic_{Y/S}^{0}$ est un schéma abélien sur $S$ (\cite{BLR90} 9.4.4). La classe de $L$ est contenue dans $\Pic_{Y/S}^{0}(\check{\overline{S}})$. On a un morphisme de $S$-schémas
\begin{equation}
	j_{x}:Y\to \Pic_{Y/S}^{0}
	\label{jx}
\end{equation}
défini, pour tout $S$-schéma $T$, par
\begin{equation}
	Y(T)\to \Pic^{0}_{Y/S}(T) \qquad y\mapsto [\mathscr{O}_{Y_{T}}(y-x)].
	\label{jx morphisme}
\end{equation}
Notons $A=(\Pic^0_{Y/S})^{\vee}$ le schéma abélien dual de $\Pic_{Y/S}^{0}$ et considérons le morphisme d'Albanese défini par $x\in Y(S)$,
\begin{equation}
	\varphi_x: Y \xrightarrow{j_{x}} \Pic^{0}_{Y/S} \xrightarrow{\Phi} A,
	\label{morphisme d'Albanese}
\end{equation}
où $\Phi$ est la polarisation canonique (\cite{MB85} 2.6.4). D'après (\cite{MB85} 2.7.9), le morphisme
\begin{equation}
	A^{\vee}\simeq \Pic^0_{A/S} \xrightarrow{\varphi_x^{*}} \Pic^0_{Y/S}\simeq A^{\vee},
	\label{-id map}
\end{equation}
induit par $\varphi_x$, s'identifie à $-\id_{A^{\vee}}$. Il existe donc un fibré en droites $L'$ sur $\check{\overline{A}}=A\times_{S}\check{\overline{S}}$ dont la classe appartient à $\Pic^{0}_{A/S}(\check{\overline{S}})$ tel que $\varphi_{x}^{*}(L')\simeq L$. Par descente (\cite{EGA IV} 8.8.3), le groupe $\Pic_{A/S}(k)\simeq \Pic_{A_{s}/k}(k)$ est de torsion. Il existe un entier $N$ tel que
\begin{equation}
	N^{*}[L'_{s}]=0 \qquad \textrm{dans} ~ \Pic^0_{A/S}(k),
	\label{line bundle trivial mod p}
\end{equation}
où on a encore noté $N:\Pic^0_{A/S}\to \Pic^0_{A/S}$ le morphisme de la multiplication par $N$. Le groupe
\begin{displaymath}
	\Ker (\Pic^0_{A/S}(\oo_{n})\to \Pic^0_{A/S}(k))
\end{displaymath}
est de $p^n$-torsion en vertu de (\cite{Ta67} 2.4). On a donc
\begin{equation}
	(p^nN)^{*}[L_n']=0 \qquad \textrm{dans} ~ \Pic^0_{A/S}(\oo_{n}).
	\label{line bundle trivial mod pn}
\end{equation}

On désigne par $Y^{(n)}$ le $T$-schéma défini par le diagramme cartésien
\begin{equation}
	\xymatrix{
		Y^{(n)} \ar[d] \ar[r] &A\ar[d]^{p^nN}\\
		Y \ar[r]^{\varphi_x} &A}
	\label{Jacobien method}
\end{equation}
En vertu de \eqref{line bundle trivial mod pn}, l'image réciproque de $L_n$ par $Y^{(n)}_n\to Y_n$ est triviale. D'après la théorie du corps de classes géométrique (\cite{Se75} VI 11), $Y^{(n)}_{\eta}$ est géométriquement connexe et est un revêtement étale et galoisien de $Y_{\eta}$. D'après \eqref{Jacobien method}, l'action du groupe $\Aut(Y^{(n)}_{\eta}/Y_{\eta})$ s'étend en une $Y$-action sur $Y^{(n)}$. Donc $Y^{(n)} \to Y$ est un $\eta$-revêtement fini et galoisien. La proposition résulte alors de \ref{domine par bon rev bar}(ii).

\begin{nothing} \label{homomorphisme psi conjugaison}
	Soient $X$ une $\overline{S}$-courbe propre, $\overline{x}$ un $\overline{\eta}$-point de $X_{\overline{\eta}}$ et $n$ un entier $\ge 1$. On se propose de construire un foncteur de $\BB^{\DW}_{\check{X}}$ dans $\Rep_{\oo_n}^{\ltf}(\pi_1(X_{\overline{\eta}},\overline{x}))$ \eqref{oo oon reprentations}. Soit $\mathcal{F}$ un fibré vectoriel de Deninger-Werner sur $\check{X}$. Comme $X$ est propre, le $\overline{\eta}$-point $\overline{x}$ se prolonge en une section $\overline{S}\to X$ que l'on note encore (abusivement) $\overline{x}$. On note $\overline{x}_n: \overline{S}_n\to X_{n}$ la réduction modulo $p^{n}$ de $\overline{x}$. On définit le $\oo_n$-module sous-jacent à la représentation associée à $\mathcal{F}$ par
\begin{equation}
	\mathbb{V}_{n}(\mathcal{F})=\Gamma(\overline{S}_{n},\overline{x}_n^*(\mathcal{F}_{n})),
	\label{def of Fxn}
\end{equation}
qui est un $\oo_n$-module libre de type fini.

	D'après \ref{domine par bon rev bar}(i), il existe un $\overline{\eta}$-revêtement c.p.d.$0$ et galoisien $\varphi:X'\to X$ de groupe des automorphismes $G$ tel que $\varphi_n^{*}(\mathcal{F}_n)$ soit libre de type fini. Soit $(C_{i},\overline{y}_{i})_{i\in I}$ un revêtement universel galoisien de $(X_{\overline{\eta}},\overline{x})$ \eqref{rep and sysl}. Choisissons un point géométrique $\overline{y}\in X'_{\overline{\eta}}(\overline{\eta})$ au-dessus de $\overline{x}$. Il existe $i\in I$ et un épimorphisme
\begin{displaymath}
	\xi_{i}:C_{i}\to X'_{\overline{\eta}}\qquad \overline{y}_{i} \mapsto \overline{y}.
\end{displaymath}
Celui-ci induit un épimorphisme
\begin{displaymath}
	\xi_{i}: \Aut(C_{i}/X_{\overline{\eta}}) \to G
\end{displaymath}
défini pour tout $\sigma\in \Aut(C_{i}/X_{\overline{\eta}})$, par $\xi_{i}(\sigma)(\overline{y})=\xi_{i}\circ \sigma(\overline{y}_{i})$. Le composé des homomorphismes
\begin{equation}
	\xi_{\overline{y}}: \pi_{1}(X_{\overline{\eta}},\overline{x}) \twoheadrightarrow \Aut(C_{i}/X_{\overline{\eta}})^{\op} \xrightarrow{\xi_{i}} G^{\op}
	\label{psi y homo}
\end{equation}
est indépendant du choix de $i\in I$ mais dépend du choix de $\overline{y}$ au-dessus de $\overline{x}$. En effet, soient $\overline{y}$ et $\overline{y}'$ deux $\overline{\eta}$-points de $X_{\overline{\eta}}'$ au-dessus de $\overline{x}$. Il existe un élément $h\in G$ tel que $h(\overline{y})=\overline{y}'$. D'après (\cite{DW03} Thm.13 (19)), les homomorphismes $\xi_{\overline{y}}, \xi_{\overline{y}'}: \pi_{1}(C,\overline{x}) \to G^{\op}$ sont reliés par:
\begin{equation}
	\xi_{\overline{y}'}(\gamma)=h\xi_{\overline{y}}(\gamma)h^{-1}, \qquad \forall \gamma\in \pi_{1}(C,\overline{x}).
	\label{conjugaison de psi}
\end{equation}
\end{nothing}
\begin{nothing}
	Conservons les hypothèses et notations de \ref{homomorphisme psi conjugaison}. On construit une action $\rho_{n,\mathcal{F}}$ de $\pi_1(X_{\overline{\eta}},\overline{x})$ sur $\mathbb{V}_{n}(\mathcal{F})$. L'action de $G$ sur $X'_{\overline{\eta}}$ s'étend en une $X$-action sur $X'$. Tout $g\in G$ induit donc un automorphisme
\begin{equation}
	g_{n}^{*}: \Gamma(X'_n,\varphi_n^{*}(\mathcal{F}_{n})) \to \Gamma(X'_n,\varphi_n^{*}(\mathcal{F}_{n})).
	\label{l'action de G sur section}
\end{equation}
Le $\overline{\eta}$-point $\overline{y}$ se prolonge en un $\overline{S}$-point de $X'$ au-dessus de $\overline{x}\in X(\overline{S})$ que l'on note encore (abusivement) $\overline{y}$. Comme le morphisme structural $\lambda:X'\to \overline{S}$ vérifie $\lambda_{*}(\mathscr{O}_{X'})=\mathscr{O}_{\overline{S}}$ universellement, le morphisme $\overline{y}_n:\overline{S}_n\to X'_n$ déduit de $\overline{y}$ induit un isomorphisme:
\begin{equation}
	\overline{y}_n^{*}:\Gamma(X'_n,\varphi_n^{*}(\mathcal{F}_{n}))\xrightarrow{\sim}\Gamma(\overline{S}_n, \overline{y}_n^*(\varphi_n^*(\mathcal{F}_{n})))=\mathbb{V}_{n}(\mathcal{F}).
	\label{pullback to a point rho}
\end{equation}
Pour tout $\gamma\in \pi_{1}(X_{\overline{\eta}},\overline{x})$, posant $g=\xi_{\overline{y}}(\gamma)$, on définit un automorphisme $\rho_{n,\mathcal{F}}(\gamma)$ de $\mathbb{V}_{n}(\mathcal{F})$ par
\begin{equation}
	\rho_{n,\mathcal{F}}(\gamma): \mathbb{V}_{n}(\mathcal{F})\xrightarrow{(\overline{y}_{n}^{*})^{-1}} \Gamma(X'_{n},\varphi_{n}^{*}(\mathcal{F}_{n})) \xrightarrow{g_{n}^{*}} \Gamma(X'_n,\varphi_n^{*}(\mathcal{F}_{n})) \xrightarrow{\overline{y}_{n}^{*}} \mathbb{V}_{n}(\mathcal{F}).
	\label{def of rho n F gamma}
\end{equation}
On notera que $\rho_{n,\mathcal{F}}$ est l'homomorphisme composé
\begin{equation}
	\rho_{n,\mathcal{F}}: \pi_{1}(X_{\overline{\eta}},\overline{x}) \xrightarrow{\xi_{\overline{y}}} G^{\op} \to \Aut_{\oo_{n}}\Gamma(X'_{n},\varphi_{n}^{*}(\mathcal{F}_{n})) \xrightarrow[\sim]{\overline{y}_{n}^{*}} \Aut_{\oo_{n}} \mathbb{V}_{n}(\mathcal{F}).
	\label{def of rho n F}
\end{equation}
On définit la $\oo_{n}$-représentation de $\pi_{1}(X_{\overline{\eta}},\overline{x})$ associée à $\mathcal{F}$ par $(\mathbb{V}_{n}(\mathcal{F}), \rho_{n,\mathcal{F}})$ que l'on note encore $\mathbb{V}_{n}(\mathcal{F})$.
\label{l'action du groupe fondamentale}	
\end{nothing}

\begin{prop}[\cite{DW03} Thm. 13] \label{def of DW n independant}
	Sous les hypothèses de \ref{homomorphisme psi conjugaison}, la $\oo_{n}$-représentation $\mathbb{V}_{n}(\mathcal{F})$ de $\pi_{1}(X_{\overline{\eta}},\overline{x})$ est indépendante des choix du $\overline{\eta}$-revêtement c.p.d.$0$ \textnormal{(\ref{def of revêtement}(ii))} et galoisien $\varphi:X'\to X$ et du point $\overline{y}$ de $X'_{\overline{\eta}}$ au-dessus de $\overline{x}$, à isomorphisme près.
\end{prop}


\begin{nothing} \label{definition of functor rho}
	Conservons les hypothèses et notations de \ref{homomorphisme psi conjugaison}. \'{E}tant donnés un morphisme $f:\mathcal{F}\to \mathcal{F}'$ de $\BB^{\DW}_{\check{X}}$ et un entier $n\ge 1$, il existe un $\overline{\eta}$-revêtement c.p.d.$0$ et galoisien $\varphi:X'\to X$ trivialisant $\mathcal{F}_n$ et $\mathcal{F}_n'$ (cf. \ref{domine par bon rev bar}(i)). L'isomorphisme \eqref{pullback to a point rho} est clairement fonctoriel en $\mathcal{F}$. D'après \ref{def of DW n independant}, le morphisme $f$ induit un morphisme de $\oo_n$-représentations de $\pi_1(X_{\overline{\eta}},\overline{x})$
\begin{equation}
	\mathbb{V}_n(f): \mathbb{V}_n(\mathcal{F})\to \mathbb{V}_n(\mathcal{F}').
	\label{rhon f}
\end{equation}
La correspondance $\mathcal{F}\mapsto \mathbb{V}_{n}(\mathcal{F})$ définit ainsi un foncteur qu'on note \eqref{oo oon reprentations}
\begin{equation}
	\mathbb{V}_n: \BB^{\DW}_{\check{X}} \to \Rep_{\oo_n}^{\ltf}(\pi_1(X_{\overline{\eta}},\overline{x})).
	\label{rho VVn}
\end{equation}
\end{nothing}
\begin{prop}[\cite{DW03} Thm. 14 et 16]	\label{propriete DW foncteur}
	\textnormal{(i)} Les foncteurs $\{\mathbb{V}_n\}_{n\ge 1}$ \eqref{rho VVn} sont compatibles et ils définissent un foncteur $\oo$-linéaire et exact \eqref{oo oon reprentations}
	\begin{equation}
		\mathbb{V}: \BB^{\DW}_{\check{X}}\to \Rep_{\oo}^{\ltf}(\pi_1(X_{\overline{\eta}},\overline{x})),
		\label{rho entier}
	\end{equation}
	qui commute aux produits tensoriels, aux passages aux duaux et aux homomorphismes internes. On le note aussi $\mathbb{V}_{\check{X}}$ pour signifier que la construction dépend de $X$.
	
	\textnormal{(ii)} Pour tout morphisme dominant $f:X'\to X$ de $\overline{S}$-modèles de $X_{\overline{\eta}}$, l'image inverse d'un fibré vectoriel de Deninger-Werner sur $X$ par $f$ est de Deninger-Werner sur $X'$. De plus, on a un isomorphisme de foncteurs
	\begin{equation}
		\mathbb{V}_{\check{X}}\xrightarrow{\sim} \mathbb{V}_{\check{X}'}\circ f^{*}
		\label{XX XX}
	\end{equation}
\end{prop}

$\hspace*{-1.2em}\bf{\arabic{section}.\stepcounter{theorem}\arabic{theorem}.}$
Dans la suite de cette section, on fixe une courbe propre et lisse $C$ sur $\overline{K}$, $\overline{x}$ un $\overline{\eta}$-point de $C$ et on pose $\check{C}=C\otimes_{\overline{K}}\mathfrak{C}$.

\begin{definition}
	(i) On dit qu'un fibré vectoriel $F$ sur $\check{C}$ est \textit{de Deninger-Werner} s'il existe un $\overline{S}$-modèle $X$ de $C$ et un fibré vectoriel de Deninger-Werner $\mathcal{F}$ sur $\check{X}$ muni d'un isomorphisme $F\simeq \mathcal{F}_{\check{\overline{\eta}}}$. On désigne par $\BB^{\DW}_{\check{C}}$ la sous-catégorie pleine de $\Vect_{\check{C}}$ formée des fibrés vectoriels de Deninger-Werner.

	(ii) On dit qu'un fibré vectoriel $F$ sur $\check{C}$ est \textit{potentiellement de Deninger-Werner} s'il existe un revêtement étale et connexe $\pi:C'\to C$ tel que, posant $\check{\pi}=\pi\otimes_{\overline{K}}\mathfrak{C}$, $\check{\pi}^{*}(F)$ soit de Deninger-Werner sur $\check{C}'=C'\otimes_{\overline{K}}\mathfrak{C}$. On désigne par $\BB^{\DW}_{\check{C}\sharp}$ la sous-catégorie pleine de $\Vect_{\check{C}}$ formée des fibrés vectoriels potentiellement de Deninger-Werner.

	\label{def of categorie DW rat}
\end{definition}

\begin{theorem}[\cite{DW05} Thm. 11, 12, 13] 
	\textnormal{(i)} Les sous-catégories $\BB^{\DW}_{\check{C}}$ et $\BB^{\DW}_{\check{C}\sharp}$ de $\Vect_{\check{C}}$ sont stables par extensions.

	\textnormal{(ii)} La catégorie $\BB^{\DW}_{\check{C}\sharp}$ contient tous les fibrés en droites de degré $0$ sur $\check{C}$.

	\textnormal{(iii)} Tout fibré vectoriel de $\BB_{\check{C}\sharp}^{\DW}$ est semi-stable de degré $0$.
	\label{DW cat property}
\end{theorem}
\begin{nothing}
	Soit $X$ un $\overline{S}$-modèle de $C$. D'après \ref{model curve vect}, on a un foncteur pleinement fidèle $\jmath_{\check{X}}:\BB_{\check{X},\mathbb{Q}}^{\DW}\to \BB_{\check{C}}^{\DW}$. D'après \eqref{Rep o Q eqi Rep C}, \eqref{XX XX} et (\cite{DW03} Cor. 21), on obtient un foncteur $\mathfrak{C}$-linéaire, exact \eqref{oo oon reprentations general}:
\begin{equation}
	\mathbb{V}_{\check{C}}: \BB^{\DW}_{\check{C}}\to \Rep_{\mathfrak{C}}^{\cont}(\pi_1(C,\overline{x}))
	\label{rho C check}
\end{equation}
qui commute aux produits tensoriels, aux duaux et aux homomorphismes internes (cf. \cite{DW05} Thm. 28).
\label{rho CC DW}
\end{nothing}
\section{Topos annelé de Faltings} \label{Topos de Faltings}
\begin{nothing}
	Dans cette section, on se donne un $S$-schéma plat, séparé et de type fini $X$ tel que $\overline{X}=X\times_{S}\overline{S}$ \eqref{basic notation} soit normal. On notera que $\overline{X}$ est localement irréductible dans le sens de (\cite{AGT} III.3.1). En effet, comme $\overline{X}$ est $\overline{S}$-plat, ses points génériques sont les points génériques du schéma $X_{\overline{\eta}}$, qui est noethérien; l'ensemble des points génériques de $\overline{X}$ est donc fini (cf. \cite{AGT} III.3.2(ii)). On note $\hbar: \overline{X}\to X$ et $h:X_{\overline{\eta}}\to X$ les morphismes canoniques. On désigne par $E$ la catégorie des morphismes de schémas $V\to U$ au-dessus du morphisme canonique $X_{\overline{\eta}}\to X$, i.e. des diagrammes commutatifs
\begin{equation}
	\xymatrix{
		V\ar[r] \ar[d] & U \ar[d]\\
		X_{\overline{\eta}}\ar[r] & X
	}
	\label{V U audessus X}
\end{equation}
tels que le morphisme $U\to X$ soit étale et que le morphisme $V\to U_{\overline{\eta}}$ soit fini et étale (\cite{AGT} VI.10.1). Elle est fibrée au-dessus de la catégorie $\Et_{/X}$ des $X$-schémas étales, par le foncteur
\begin{equation}
	\pi: E\to \Et_{/X}, \qquad (V\to U)\mapsto U.
	\label{def of pi}
\end{equation}
La fibre de $\pi$ au-dessus d'un $X$-schéma étale $U$ est la catégorie $\Et_{\textnormal{f}/U_{\overline{\eta}}}$ des schémas finis étales au-dessus de $U_{\overline{\eta}}$, que l'on équipe de la topologie étale \eqref{fet et}. La catégorie fibrée $\pi$ est alors un site fibré. (\cite{SGAIV} VI 7.2.1)

On munit $E$ de la topologie co-évanescente (\cite{AGT} VI.5.3), c'est-à-dire, la topologie engendrée par les recouvrements $\{(V_i\to U_i)\to (V\to U)\}_{i\in I}$ des deux types suivants:

(v) $U_i=U$ pour tout $i \in I$, et $(V_i\to V)_{i\in I}$ est un recouvrement;	

(c) $(U_i\to U)_{i\in I}$ est un recouvrement et $V_i\simeq U_i\times_{U}V$ pour tout $i\in I$.

On appelle \textit{topos de Faltings de $X$} et l'on note $\widetilde{E}$, le topos des faisceaux d'ensembles sur $E$. Si $F$ est un préfaisceau sur $E$, on note $F^{a}$ le faisceau associé. Signalons la description commode et simple suivante de $\widetilde{E}$.
\label{basic topos de Faltings}
\end{nothing}

\begin{prop}[\cite{AGT} VI.5.10]
	La donnée d'un faisceau $F$ sur $E$ est équivalente à la donnée pour tout objet $U$ de $\Et_{/X}$ d'un faisceau $F_U$ de $U_{\overline{\eta},\fet}$ \eqref{fet et} et pour tout morphisme $f : U'\to U$ de $\Et_{/X}$ d'un morphisme $F_U\to (f_{\overline{\eta}})_{\fet*}(F_{U'})$, ces morphismes étant soumis à des relations de compatibilité, tels que, pour toute famille couvrante $(f_n:U_n\to U)_{n\in \Sigma}$ de $\Et_{/X}$, si pour tout $(m,n)\in \Sigma^2$, on pose $U_{mn} = U_m \times_U U_n$ et on note $f_{mn} : U_{mn} \to U$ le morphisme canonique, la suite de morphismes de faisceaux de $U_{\overline{\eta},\fet}$
	\begin{equation}
		F_U\to \prod_{n\in \Sigma}(f_{n,\overline{\eta}})_{\fet*}(F_{U_n})\rightrightarrows \prod_{n,m\in \Sigma}(f_{mn,\overline{\eta}})_{\fet*}(F_{U_{mn}})
		\label{suite exacte faisceau E}
	\end{equation}
soit exacte.
	\label{def of sheaf on E}
\end{prop}

En vertu de la proposition précédente, on peut identifier tout faisceau $F$ sur $E$ au foncteur $\{U\to F_U\}$ associé, où $U\in\Ob(\Et_{/X})$ et $F_U\in \Ob(U_{\overline{\eta},\fet})$ est la restriction de $F$ à la fibre de $\pi$ au-dessus de $U$.\\

$\hspace*{-1.2em}\bf{\arabic{section}.\stepcounter{theorem}\arabic{theorem}.}$
Pour tout $U\in \Ob(\Et_{/X})$, on a un foncteur canonique:
\begin{equation}
	\alpha_{U !}:\Et_{\textnormal{f}/U_{\overline{\eta}}}\to E, \qquad V \mapsto (V \to U).
	\label{def of alpha!}
\end{equation}
Le foncteur $\alpha_{X !}:\Et_{\textnormal{f}/X_{\overline{\eta}}}\to E$ est continu et exact à gauche (\cite{AGT} VI.5.32). Il définit donc un morphisme de topos (\cite{AGT} VI.(10.6.3))
\begin{equation}
	\beta: \widetilde{E}\to X_{\overline{\eta},\fet}.
	\label{foncteur topos beta}
\end{equation}
Le foncteur
\begin{equation}
	\sigma^+:\Et_{/X} \to E, \qquad U\mapsto (U_{\overline{\eta}} \to U)
	\label{foncteur site sigma}
\end{equation}
est continu et exact à gauche (\cite{AGT} VI.5.32). Il définit donc un morphisme de topos (\cite{AGT} VI.(10.6.4)):
\begin{equation}
	\sigma: \widetilde{E}\to X_{\et}.
	\label{foncteur topos sigma}
\end{equation}
Par ailleurs, le foncteur
\begin{equation}
	\Psi^{+}: E\to \Et_{/X_{\overline{\eta}}}, \qquad (V\to U)\mapsto V
	\label{foncteur site Phi}
\end{equation}
est continu et exact à gauche (\cite{AGT} VI.10.7). Il définit donc un morphisme de topos:
\begin{equation}
	\Psi:X_{\overline{\eta},\et}\to \widetilde{E}.
	\label{foncteur topos Phi}
\end{equation}
On vérifie aussitôt qu'on a un isomorphisme canonique
\begin{equation}
	\beta\circ \Psi\xrightarrow{\sim}\rho_{X_{\overline{\eta}}},
	\label{beta Psi commute}
\end{equation}
où $\rho_{X_{\overline{\eta}}}: X_{\overline{\eta},\et}\to X_{\overline{\eta},\fet}$ est le morphisme canonique \eqref{topos varsigma X}. On en déduit par adjonction un morphisme
\begin{equation}
	\beta^*\rightarrow\Psi_{*}\rho_{X_{\overline{\eta}}}^{*},
	\label{beta Psi}
\end{equation}
qui est un isomorphisme en vertu de  (\cite{AGT} VI.10.9(iii)).\\

$\hspace*{-1.2em}\bf{\arabic{section}.\stepcounter{theorem}\arabic{theorem}.}$
On désigne par $D$ la catégorie des morphismes de schémas $V\to U$ au-dessus du morphisme canonique $h: X_{\overline{\eta}}\to X$, i.e. des diagrammes commutatifs
\begin{equation}
	\xymatrix{
		V\ar[r] \ar[d] & U \ar[d]\\
		X_{\overline{\eta}} \ar[r] & X,
	}
	\label{V to U coevanescente}
\end{equation}
tels que les flèches verticales soient étales (cf. \cite{AGT} VI.4.1). Elle est fibrée au-dessus de la catégorie $\Et_{/X}$ des $X$-schémas étales, par le foncteur
\begin{equation}
	\varpi:D\to \Et_{/X} \qquad (V\to U)\mapsto U.	
	\label{functor fibre D}
\end{equation}
La fibre de $\varpi$ au-dessus d'un schéma étale $U$ est la catégorie $\Et_{/U_{\overline{\eta}}}$ des schémas étales au-dessus de $U_{\overline{\eta}}$, que l'on équipe de la topologie étale. La catégorie fibrée $\varpi$ est alors un site fibré (cf. \cite{SGAIV} VI 7.2.1).

On munit $D$ de la topologie co-évanescente (\cite{AGT} VI.5.3), c'est-à-dire la topologie engendrée par les recouvrements $\{(V_i\to U_i)\to (V\to U)\}_{i\in I}$ des deux types suivants:

(v) $U_i=U$ pour tout $i \in I$, et $(V_i\to V)_{i\in I}$ est un recouvrement;	

(c) $(U_i\to U)_{i\in I}$ est un recouvrement et $V_i\simeq U_i\times_{U}V$ pour tout $i\in I$.

Le topos des faisceaux d'ensembles sur $D$ est appelé \textit{le topos co-évanescent} et est noté $X_{\et}\overleftarrow{\times}_{X_{\et}}X_{\overline{\eta},\et}$. Celui-ci vérifie une propriété universelle qui explique cette notation (cf. \cite{AGT} VI.3.7 et VI.3.12). Les foncteurs
\begin{eqnarray}
	&\rmp_1^{+}&:\Et_{/X} \to D,\qquad U\mapsto (U_{\overline{\eta}}\to U),\\
	&\rmp_2^{+}&:\Et_{/X_{\overline{\eta}}} \to D,\qquad V\mapsto (V\to X),
	\label{p1p2 site}
\end{eqnarray}
sont exacts à gauche et continus. Ils définissent donc deux morphismes de topos (\cite{SGAIV} IV 4.9.2)
\begin{eqnarray}
	&\rmp_1&: X_{\et}\overleftarrow{\times}_{X_{\et}}X_{\overline{\eta},\et} \to X_{\et},\\
	&\rmp_2&: X_{\et}\overleftarrow{\times}_{X_{\et}}X_{\overline{\eta},\et} \to X_{\overline{\eta},\et}.
	\label{p1p2 topos}
\end{eqnarray}
\begin{nothing} \label{site topos coevanesente}
Tout objet de $E$ est naturellement un objet de $D$. On définit ainsi un ($\Et_{/X}$)-foncteur cartésien pleinement fidèle
\begin{equation}
	\rho^{+}:E\to D
	\label{coevanescent Faltings site}
\end{equation}
dont la fibre au-dessus d'un schéma étale $U$ sur $X$ est le foncteur d'injection canonique $\Et_{\textnormal{f}/U_{\overline{\eta}}}\to \Et_{/U_{\overline{\eta}}}$. Le foncteur \eqref{coevanescent Faltings site} est continu et exact à gauche (\cite{AGT} VI.10.15). Il définit donc un morphisme de topos
\begin{equation}
	\rho: X_{\et}\overleftarrow{\times}_{X_{\et}}X_{\overline{\eta},\et} \to \widetilde{E}.
	\label{coevanescent Faltings topos}
\end{equation}
Il résulte aussitôt des définitions que les carrés du diagramme
\begin{equation}
	\xymatrix{
		X_{\et} \ar@{=}[d] & X_{\et}\overleftarrow{\times}_{X_{\et}}X_{\overline{\eta},\et} \ar[l]_-{\rmp_1} \ar[r]^-{\rmp_2} \ar[d]^{\rho} & X_{\overline{\eta},\et} \ar[d]^{\rho_{X_{\overline{\eta}}}}\\
		X_{\et} & \widetilde{E} \ar[l]_{\sigma} \ar[r]^{\beta} & X_{\overline{\eta},\fet},
	}
	\label{diagramme commutatif coevanescente Faltings}
\end{equation}
où $\rho_{X_{\overline{\eta}}}$ est le morphisme \eqref{topos varsigma X}, sont commutatifs à isomorphismes canoniques près.
\end{nothing}

\begin{nothing}
D'après (\cite{AGT} VI.3.11), la donnée d'un point de $X_{\et}\overleftarrow{\times}_{X_{\et}}X_{\overline{\eta},\et}$ est équivalente à la donnée d'une paire de points géométriques $\overline{x}$ de $X$ et $\overline{y}$ de $X_{\overline{\eta}}$ et d'une flèche de spécialisation $u: \overline{y} \to X_{(\overline{x})}$, où $X_{(\overline{x})}$ est le localisé strict de $X$ en $\overline{x}$. Un tel point sera noté $(\overline{y}\rightsquigarrow \overline{x})$. On désigne par $\rho(\overline{y}\rightsquigarrow \overline{x})$ son image par $\rho$, qui est donc un point de $\widetilde{E}$. D'après (\cite{AGT} VI.5.30 et VI.10.18), lorsque $(\overline{y}\rightsquigarrow \overline{x})$ décrit la famille des points de $X_{\et}\overleftarrow{\times}_{X_{\et}}X_{\overline{\eta},\et}$, la famille des foncteurs fibres de $X_{\et}\overleftarrow{\times}_{X_{\et}}X_{\overline{\eta},\et}$ (resp. $\widetilde{E}$) associés aux points $(\overline{y}\rightsquigarrow \overline{x})$ (resp. $\rho(\overline{y}\rightsquigarrow \overline{x})$) est conservative.
\label{points topos Faltings}
\end{nothing}

\begin{nothing}
	Reprenons les notations de \ref{fet et}. On désigne par
	\begin{equation}
		\varpi_{\scoh}: D_{\scoh}\to \Et_{\scoh/X}, \qquad \pi_{\scoh}: E_{\scoh}\to \Et_{\scoh/X}
		\label{varpi pi scoh}
	\end{equation}
	les sites fibrés déduits de $\varpi$ \eqref{functor fibre D} et $\pi$ \eqref{def of pi} par changement de base par le foncteur d'injection canonique
	\begin{equation}
		\Et_{\scoh/X}\to \Et_{/X},
		\label{Et scoh to Et}
	\end{equation}
	et par $\iota_{1}: D_{\scoh}\to D$ et $\iota_{2}: E_{\scoh}\to E$ les projections canoniques. D'après (\cite{AGT} VI.5.21), si l'on munit $D_{\scoh}$ (resp. $E_{\scoh}$) de la topologie co-évanescente définie par $\varpi_{\scoh}$ (resp. $\pi_{\scoh}$) et si l'on note $\widetilde{D}_{\scoh}$ (resp. $\widetilde{E}_{\scoh}$) le topos des faisceaux d'ensembles sur $D_{\scoh}$ (resp. $E_{\scoh}$), le foncteur $\iota_{1}$ (resp. $\iota_{2}$) induit par restriction une équivalence de catégories $X_{\et}\overleftarrow{\times}_{X_{\et}}X_{\overline{\eta},\et}\xrightarrow{\sim}\widetilde{D}_{\scoh}$ (resp. $\widetilde{E}\xrightarrow{\sim}\widetilde{E}_{\scoh}$).
	\label{Faltings coevanescent scoh}
\end{nothing}

\begin{nothing}
Pour tout $(V\to U) \in \Ob(E)$, on note $\overline{U}^{V}$ la clôture intégrale de $\overline{U}$ dans $V$. On désigne par $\overline{\mathscr{B}}$ le préfaisceau d'anneaux sur $E$ défini pour tout $(V\to U) \in \Ob(E)$, par
\begin{equation}
	\overline{\mathscr{B}}(V\to U)=\Gamma(\overline{U}^{V},\mathscr{O}_{\overline{U}^{V}}).
	\label{def of bar B}
\end{equation}
Pour tout $U\in \Ob(\Et_{/X})$, on pose \eqref{def of alpha!}
\begin{equation}
	\overline{\mathscr{B}}_{U}=\overline{\mathscr{B}}\circ \alpha_{U!}.
	\label{def of bar B UU}
\end{equation}
D'après (\cite{AGT} III.8.16), $\overline{\mathscr{B}}$ est un faisceau pour la topologie co-évanescente de $E$. On dira dans la suite que $\overline{\mathscr{B}}$ est l'anneau de $\widetilde{E}$ associé à $\overline{X}$. Pour tout entier $n\ge 1$ et tout $U\in \Ob(\Et_{/X})$, on pose
\begin{eqnarray}
	&&\overline{\mathscr{B}}_n=\overline{\mathscr{B}}/p^n\overline{\mathscr{B}}, \label{def of Bn} \\	
	&&\overline{\mathscr{B}}_{U,n}=\overline{\mathscr{B}}_{U}/p^{n}\overline{\mathscr{B}}_{U}.
	\label{bar B U n}
\end{eqnarray}
La correspondance $\{U\mapsto \overline{\mathscr{B}}_{U,n}\}$ forme naturellement un préfaisceau sur $E$ dont le faisceau associé est canoniquement isomorphe au à $\overline{\mathscr{B}}_{n}$ (cf. \cite{AGT} VI.8.2 et VI.8.9).
\label{ring B bar}
\end{nothing}

$\hspace*{-1.2em}\bf{\arabic{section}.\stepcounter{theorem}\arabic{theorem}.}$
Considérons l'homomorphisme
\begin{equation}
	\hbar_{*}(\mathscr{O}_{\overline{X}})\to \sigma_{*}(\overline{\mathscr{B}}),
	\label{ring B OX}
\end{equation}
où $\hbar: \overline{X}\to X$ est la projection canonique, défini pour tout $U\in \Ob(\Et_{/X})$ par l'homomorphisme canonique
\begin{equation}
	\Gamma(\overline{U},\mathscr{O}_{\overline{U}})\to \Gamma(\overline{U}^{U_{\overline{\eta}}},\mathscr{O}_{\overline{U}^{U_{\overline{\eta}}}}).
\end{equation}
Sauf mention explicite du contraire, on considère $\sigma:\widetilde{E}\to X_{\et}$ \eqref{foncteur topos sigma} comme un morphisme de topos annelés, respectivement par $\overline{\mathscr{B}}$ et $\hbar_{*}(\mathscr{O}_{\overline{X}})$. Nous utilisons pour les modules la notation $\sigma^{-1}$ pour désigner l'image inverse au sens des faisceaux abéliens et nous réservons la notation $\sigma^{*}$ pour l'image inverse au sens des modules.

\begin{nothing}
	Comme $X_{\eta}$ est un ouvert de $X_{\et}$, i.e. un sous-objet de l'objet final $X$ de $\Et_{/X}$ (\cite{SGAIV} IV 8.3), $\sigma^*(X_{\eta})$ est un ouvert de $\widetilde{E}$. On note $\gamma:\widetilde{E}_{/\sigma^{*}(X_{\eta})}\to \widetilde{E}$ le morphisme de localisation de $\widetilde{E}$ en $\sigma^*(X_{\eta})$. On désigne par $\widetilde{E}_s$ le sous-topos fermé de $\widetilde{E}$ complémentaire de l'ouvert $\sigma^*(X_{\eta})$, i.e. la sous-catégorie pleine de $\widetilde{E}$ formée des faisceaux $F$ tels que $\gamma^*(F)$ soit un objet final de $\widetilde{E}_{/\sigma^*(X_{\eta})}$ (\cite{SGAIV} IV 9.3.5). On note
\begin{equation}
	\delta:\widetilde{E}_s\to \widetilde{E}
	\label{delta}
\end{equation}
le plongement canonique, c'est-à-dire le morphisme de topos tel que $\delta_{*}: \widetilde{E}_s\to \widetilde{E}$ soit le foncteur d'injection canonique. D'après  (\cite{AGT} III.9.7), pour tout entier $n\ge 1$, l'anneau $\overline{\mathscr{B}}_n$ \eqref{def of Bn} est un objet de $\widetilde{E}_s$.
\label{sous-topos ferme}
\end{nothing}
\begin{nothing} \label{Es Bn topos Faltings}
Soit $n$ un entier $\ge 1$. Notons $d: X_s \to X$ et $d_{n}:X_{s}\to X_{n}$ les morphismes canoniques \eqref{basic notation}. Il existe un morphisme de topos (\cite{AGT} III.9.8):
\begin{equation}
	\sigma_s: \widetilde{E}_s\to X_{s,\et}
	\label{morphisme de topos sigma s}
\end{equation}
unique à isomorphisme près tel que le diagramme
\begin{equation}
	\xymatrix{
		\widetilde{E}_s \ar[r]^{\sigma_s} \ar[d]_{\delta}& X_{s,\et} \ar[d]^{d}\\
		\widetilde{E} \ar[r]^{\sigma} & X_{\et},
	}
	\label{morphisme de topos commutative special}
\end{equation}
soit commutatif à isomorphisme près.

Pour tout entier $n\ge 1$, comme le corps résiduel de $\mathcal{O}_{\overline{K}}$ est algébriquement clos, on a une immersion fermée $\overline{d}_{n} :X_{s} \to \overline{X}_{n}$ qui relèvent $d_{n}$. Comme $\overline{d}_{n}$ est un homéomorphisme universel, on peut considérer le faisceau $\mathscr{O}_{\overline{X}_n}$ de $\overline{X}_{n,\et}$ (ou de $\overline{X}_{n,\zar}$) comme un faisceau de $X_{s,\et}$ (ou de $X_{s,\zar}$). Alors, l'homomorphisme canonique $\sigma^{-1}(\hbar_*(\mathscr{O}_{\overline{X}}))\to \overline{\mathscr{B}}$ \eqref{ring B OX} induit un homomorphisme d'anneaux sur $\widetilde{E}_s$ (\cite{AGT} III.9.9):
\begin{equation}
	\sigma_s^{-1}(\mathscr{O}_{\overline{X}_n})\to \overline{\mathscr{B}}_n.
	\label{morphisme anneau sigma n}
\end{equation}
Par suite, le morphisme de topos $\sigma_s$ \eqref{morphisme de topos sigma s} est sous-jacent à un morphisme de topos annelés que l'on note
\begin{equation}
	\sigma_n: (\widetilde{E}_s, \overline{\mathscr{B}}_n)\to (X_{s,\et},\mathscr{O}_{\overline{X}_n}).
	\label{morphisme de topos sigma n}
\end{equation}
\end{nothing}
\begin{nothing}
	On désigne par $\breve{\overline{\mathscr{B}}}$ l'anneau $(\overline{\mathscr{B}}_n)_{n\ge 1}$ de $\widetilde{E}^{\mathbb{N}^{\circ}}_{s}$ \eqref{limite projective de topos} et par $\mathscr{O}_{\breve{\overline{X}}}$ l'anneau $(\mathscr{O}_{\overline{X}_n})_{n\ge 1}$ de $X_{s,\et}^{\mathbb{N}^{\circ}}$ (ou de $X_{s,\zar}^{\mathbb{N}^{\circ}}$). Les morphismes $(\sigma_n)_{n\ge 1}$ \eqref{morphisme de topos sigma n} induisent un morphisme de topos annelés:
\begin{equation}
	\breve{\sigma}: (\widetilde{E}^{\mathbb{N}^{\circ}}_{s},\breve{\overline{\mathscr{B}}} ) \to (X_{s,\et}^{\mathbb{N}^{\circ}}, \mathscr{O}_{\breve{\overline{X}}}).
	\label{morphisme de topos sigma limit}
\end{equation}
On note $\XX$ le schéma formel complété $p$-adique de $\overline{X}$. On désigne par:
\begin{equation}
	\breve{u}: (X_{s,\et}^{\mathbb{N}^{\circ}}, \mathscr{O}_{\breve{\overline{X}}}) \to (X_{s,\zar}^{\mathbb{N}^{\circ}}, \mathscr{O}_{\breve{\overline{X}}})
	\label{morphisme de topos et to zar}
\end{equation}
le morphisme canonique, par
\begin{equation}
	\lambda: (X_{s,\zar}^{\mathbb{N}^{\circ}}, \mathscr{O}_{\breve{\overline{X}}}) \to (X_{s,\zar}, \mathscr{O}_{\XX})
	\label{morphisme de topos zar limit}
\end{equation}
le morphisme de topos annelés dont le foncteur image directe est la limite projective \eqref{limite projective de topos}, et par
\begin{equation}
	\rT: (\widetilde{E}^{\mathbb{N}^{\circ}}_{s},\breve{\overline{\mathscr{B}}} ) \to (X_{s,\zar}, \mathscr{O}_{\XX})
	\label{morphisme de topos T}
\end{equation}
le morphisme composé de topos annelés $\lambda\circ \breve{u}\circ \breve{\sigma}$.
\label{morphisme sigma T}
\end{nothing}

\begin{nothing}
Si $A$ est un anneau, on note encore $A$ le faisceau constant de valeur $A$ de $X_{\overline{\eta},\fet}$. Comme l'anneau $\overline{\mathscr{B}}(X_{\overline{\eta}}\to X)=\Gamma(\overline{X},\mathscr{O}_{\overline{X}})$ est une $\mathcal{O}_{\overline{K}}$-algèbre, il existe un homomorphisme canonique d'anneaux $\mathcal{O}_{\overline{K}}\to \beta_{*}(\overline{\mathscr{B}})$ de $X_{\overline{\eta},\fet}$. Celui-ci induit, pour tout entier $n\ge 1$, un homomorphisme canonique $\oo_n\to \beta_{*}(\overline{\mathscr{B}}_n)$ de $X_{\overline{\eta},\fet}$. On définit le morphisme de topos annelés
\begin{equation}
	\beta_n: (\widetilde{E}_s,\overline{\mathscr{B}}_n)\to (X_{\overline{\eta},\fet},\oo_n)
	\label{topos beta n 2}
\end{equation}
par le morphisme de topos composé $\beta\circ\delta$ (\eqref{foncteur topos beta} et \eqref{delta}) et l'homomorphisme canonique $\oo_n\to \beta_{*}(\overline{\mathscr{B}}_n)$. On note $\breve{\oo}$ l'anneau $(\oo_n)_{n\ge 1}$ de $X_{\overline{\eta},\fet}^{\mathbb{N}^{\circ}}$ et
\begin{equation}
	\breve{\beta}: (\widetilde{E}_s^{\mathbb{N}^{\circ}},\breve{\overline{\mathscr{B}}})\to ( X_{\overline{\eta},\fet}^{\mathbb{N}^{\circ}}, \breve{\oo})
	\label{topos annele beta}
\end{equation}
le morphisme de topos annelés induit par les morphismes $(\beta_n)_{n\ge 1}$ \eqref{topos beta n 2}.
\label{topos annele morphisme beta}
\end{nothing}
$\hspace*{-1.2em}\bf{\arabic{section}.\stepcounter{theorem}\arabic{theorem}.}$
Pour qu'un faisceau $F$ de $\widetilde{E}$ soit un objet de $\widetilde{E}_{s}$, il faut et il suffit que, pour tout point $(\overline{y}\rightsquigarrow \overline{x})$ de $X_{\et}\overleftarrow{\times}_{X_{\et}}X_{\overline{\eta},\et}$ (\ref{site topos coevanesente}) tel que $\overline{x}$ soit au-dessus de $\eta$, la fibre $F_{\rho(\overline{y}\rightsquigarrow \overline{x})}$ de $F$ en $\rho(\overline{y}\rightsquigarrow \overline{x})$ soit un singleton (cf. \cite{AGT} III.9.6). Soient $n$ un entier $\ge 1$ et $\mathbb{L}$ un $\oo_n$-module de $X_{\overline{\eta},\fet}$. Comme les foncteurs fibres commutent au produit tensoriel (\cite{SGAIV} IV 13.4) et que $\overline{\mathscr{B}}_{n}$ est un objet de $\widetilde{E}_{s}$ (\cite{AGT} III.9.7), on en déduit que le faisceau $\beta^{*}(\mathbb{L})\otimes_{\oo_{n}}\overline{\mathscr{B}}_{n}$ de $\widetilde{E}$ est un objet de $\widetilde{E}_{s}$. Par adjonction, on a donc un isomorphisme de $\widetilde{E}$
\begin{equation}
	\beta^{*}(\mathbb{L})\otimes_{\oo_{n}}\overline{\mathscr{B}}_{n}\xrightarrow{\sim} \delta_{*}(\beta_n^{*}(\mathbb{L})).
	\label{beta pullback E Es meme}
\end{equation}
\begin{nothing}
On note $\Mod(\overline{\mathscr{B}}_n)$ la catégorie des $\overline{\mathscr{B}}_n$-modules de $\widetilde{E}_s$ et $\Mod^{\tf}(\overline{\mathscr{B}}_n)$ la sous-catégorie pleine formée des $\overline{\mathscr{B}}_n$-modules de type fini. On note $\Mod(\breve{\overline{\mathscr{B}}})$ la catégorie des $\breve{\overline{\mathscr{B}}}$-modules de $\widetilde{E}_s^{\mathbb{N}^{\circ}}$, $\Mod^{\atf}(\breve{\overline{\mathscr{B}}})$ la sous-catégorie pleine formée des $\breve{\overline{\mathscr{B}}}$-modules adiques de type fini \eqref{limite projective de topos} et $\Mod_{\mathbb{Q}}(\breve{\overline{\mathscr{B}}})$ (resp. $\Mod_{\mathbb{Q}}^{\atf}(\breve{\overline{\mathscr{B}}})$) la catégorie des objets de $\Mod(\breve{\overline{\mathscr{B}}})$ (resp. $\Mod^{\atf}(\breve{\overline{\mathscr{B}}})$) à isogénie près \eqref{categorie a isogenies}.
\label{notations of Bmodules E}
\end{nothing}

\begin{nothing}
Supposons, dans la suite de cette section, que $X$ soit un $S$-schéma propre à réduction semi-stable de fibre géométrique générique connexe \eqref{def reduction semistable}. On notera que $X$ est adéquat dans le sens de (\cite{AGT} III.4.7) et par suite que $\overline{X}$ est normal (\cite{AGT} III.4.2 (iii)). D'après (\cite{Ach14} 6.1), pour tout point géométrique $\overline{x}$ de $X$, il existe un voisinage étale $U$ de $\overline{x}$ dans $X$ dont la fibre générique $U_{\eta}$ est un schéma $K(\pi,1)$ (cf. \cite{AG15} 1.2.2). On en déduit, pour tout faisceau abélien de torsion localement constant constructible $\mathscr{F}$ de $X_{\overline{\eta},\et}$ et tout entier $i\ge 1$, que $\rR^i\Psi_*(\mathscr{F})=0$ \eqref{foncteur topos Phi} (cf. \cite{Ach14} 9.6). On a donc un isomorphisme
\begin{equation}
	\rH^i(\widetilde{E}, \Psi_*(\mathscr{F}))\xrightarrow{\sim} \rH^i(X_{\overline{\eta},\et},\mathscr{F})
	\label{K pi 1 iso}
\end{equation}
déduit de la suite spectrale de Cartan-Leray. 
\end{nothing}
\begin{nothing}
	Soit $\mathbb{L}$ un faisceau abélien de torsion de $X_{\overline{\eta},\fet}$. Le faisceau $\rho_{X_{\overline{\eta}}}^{*}(\mathbb{L})$ \eqref{topos varsigma X} est isomorphe à une limite inductive de faisceaux abéliens de torsion localement constants et constructibles de $X_{\overline{\eta},\et}$ (cf. \cite{AGT} VI.9.20). D'après (\cite{AGT} VI.10.5 et VI.10.10), les topos $X_{\overline{\eta},\et}$ et $\widetilde{E}$ et le morphisme de topos $\Psi$ sont cohérents. Par \eqref{K pi 1 iso} et passage à la limite (\cite{SGAIV} VI 5.1), on a un isomorphisme
\begin{equation}
	\rH^i(\widetilde{E}, \Psi_*(\rho_{X_{\overline{\eta}}}^{*}(\mathbb{L})))\xrightarrow{\sim} \rH^i(X_{\overline{\eta},\et},\rho_{X_{\overline{\eta}}}^{*}(\mathbb{L}))
	\label{K pi 1 iso L}
\end{equation}
déduit de la suite spectrale de Cartan-Leray. Pour tout entier $i\ge 0$, considérons le morphisme composé
\begin{equation}
	\rH^i(X_{\overline{\eta},\fet},\mathbb{L})\xrightarrow{\beta^{*}}\rH^i(\widetilde{E},\beta^{*}(\mathbb{L}))\xrightarrow{\sim}\rH^i(\widetilde{E}, \Psi_*(\rho_{X_{\overline{\eta}}}^*(\mathbb{L})))\xrightarrow{\sim}\rH^i(X_{\overline{\eta},\et},\rho^*_{X_{\overline{\eta}}}(\mathbb{L})),
	\label{suite de cohomologie}
\end{equation}
où la deuxième flèche est induite par l'isomorphisme $\beta^{*}\xrightarrow{\sim}\Psi_*\rho^*_{X_{\overline{\eta}}}$ \eqref{beta Psi} et la troisième flèche est l'isomorphisme \eqref{K pi 1 iso L}. D'après \ref{general morphism of topos}, le morphisme composé \eqref{suite de cohomologie} s'identifie au morphisme induit par $\rho^*_{X_{\overline{\eta}}}$:
\begin{equation}
	\rho_{X_{\overline{\eta}}}^{*}: \rH^i(X_{\overline{\eta},\fet},\mathbb{L})\to \rH^i(X_{\overline{\eta},\et},\rho^*_{X_{\overline{\eta}}}(\mathbb{L})).
	\label{varsigma pullback F}
\end{equation}
qui est un isomorphisme pour $i=0,1$, et est un monomorphisme pour $i=2$ (cf. \ref{H1 et fet}). Comme les deux dernière flèches de \eqref{suite de cohomologie} sont des isomorphismes, on en déduit le résulte suivant.
\end{nothing}
\begin{lemma} \label{coro beta pullback}
	Soit $\mathbb{L}$ un faisceau abélien de torsion de $X_{\overline{\eta},\fet}$. L'homomorphisme
	\begin{equation}
		\beta^*:\rH^i(X_{\overline{\eta},\fet},\mathbb{L})\to \rH^i(\widetilde{E},\beta^{*}(\mathbb{L})).
		\label{beta pullback}
	\end{equation}
est un isomorphisme si $i=0,1$, et est un monomorphisme si $i=2$.
\end{lemma}

\begin{theorem}[Faltings; \cite{AG15} 2.4.16]
	Soient $n$ un entier $\ge 1$ et $\mathbb{L}$ un $\oo_{n}$-module localement libre de type fini de $X_{\overline{\eta},\fet}$. Le morphisme $\oo_{n}$-linéaire canonique
\begin{equation}
	\rH^i(\widetilde{E},\beta^{*}(\mathbb{L}))\to \rH^i(\widetilde{E}, \beta^*(\mathbb{L})\otimes_{\oo_{n}}\overline{\mathscr{B}}_{n})
	\label{iso de Faltings}
\end{equation}
est un $\alpha$-isomorphisme \eqref{alpha o modules}.
\label{thm iso de Faltings}
\end{theorem}

D'après (\cite{AGT} VI.9.20), il existe une extension finie $F$ de $\mathbb{Q}_{p}$ contenue dans $\overline{K}$ et, notant $\mathcal{O}_{F}$ l'anneau de valuation de $F$, un $(\mathcal{O}_{F}/p^{n}\mathcal{O}_{F})$-module localement libre de type fini $\mathbb{L}'$ de $X_{\overline{\eta},\fet}$ tels que $\mathbb{L}\simeq\mathbb{L}'\otimes_{\mathcal{O}_{F}/p^{n}\mathcal{O}_{F}}\oo_{n}$. Comme $\oo_{n}$ est plat sur $\mathcal{O}_{F}/p^{n}\mathcal{O}_{F}$, il suffit de démontrer que le morphisme $\oo_{n}$-linéaire canonique
\begin{equation}
	\rH^i(\widetilde{E},\beta^{*}(\mathbb{L}'))\otimes_{\mathcal{O}_{F}/p^{n}\mathcal{O}_{F}}\oo_{n}\to \rH^i(\widetilde{E}, \beta^*(\mathbb{L}')\otimes_{\mathcal{O}_{F}/p^{n}\mathcal{O}_{F}}\overline{\mathscr{B}}_{n})	
\end{equation}
est un $\alpha$-isomorphisme. La même démonstration de (\cite{AG15} 2.4.16) s'applique alors aux $(\mathcal{O}_{F}/p^{n}\mathcal{O}_{F})$-modules localement libres de type fini de $X_{\overline{\eta},\fet}$, en remplaçant le morphisme de Frobenius par une puissance de ce dernier.


\begin{prop}
	Soient $n$ un entier $\ge 1$ et $\mathbb{L}$ un $\oo_n$-module localement libre de type fini de $X_{\overline{\eta},\fet}$. Le morphisme $\oo_{n}$-linéaire induit par $\beta_n$ \eqref{topos beta n 2}
	\begin{equation}
		\beta_n^{*}:\rH^i(X_{\overline{\eta},\fet},\mathbb{L})\to \rH^i(\widetilde{E}_{s},\beta_n^{*}(\mathbb{L}))
		\label{presque iso beta n}
	\end{equation}
	est un $\alpha$-isomorphisme si $i=0,1$, et est un $\alpha$-monomorphisme si $i=2$ \eqref{alpha o modules}.
\label{presque iso faltings modifie n}
\end{prop}
\textit{Preuve}. Comme le foncteur $\delta_{*}:\widetilde{E}_{s}\to \widetilde{E}$ \eqref{delta} est exact, on a un isomorphisme
\begin{displaymath}
	\delta^{*}: \rH^{i}(\widetilde{E},\beta^{*}(\mathbb{L}_{n})\otimes_{\oo_{n}}\overline{\mathscr{B}}_{n})\xrightarrow{\sim}\rH^{i}(\widetilde{E}_{s},\beta_{n}^{*}(\mathbb{L}_{n}))
\end{displaymath}
en vertu de \eqref{beta pullback E Es meme}. Le morphisme \eqref{presque iso beta n} s'identifie au morphisme composé
\begin{displaymath}
	\rH^i(X_{\overline{\eta},\fet},\mathbb{L}) \xrightarrow{\beta^{*}} \rH^i(\widetilde{E},\beta^{*}(\mathbb{L}))\to\rH^{i}(\widetilde{E},\beta^{*}(\mathbb{L})\otimes_{\oo_{n}}\overline{\mathscr{B}}_{n})\xrightarrow{\delta^{*}} \rH^{i}(\widetilde{E}_{s},\beta_{n}^{*}(\mathbb{L}_{n})),
\end{displaymath}
où la première flèche est le morphisme \eqref{beta pullback} et la deuxième flèche est le $\alpha$-isomorphisme \eqref{iso de Faltings}. La proposition s'ensuit compte tenu de \ref{coro beta pullback} et \ref{thm iso de Faltings}.

\begin{coro}
Soit $\mathbb{L}=(\mathbb{L}_n)_{n\ge 1}$ un $\breve{\oo}$-module localement libre de type fini de $X_{\overline{\eta},\fet}^{\mathbb{N}^{\circ}}$. Pour $i=0,1$, le morphisme $\oo$-linéaire induit par $\breve{\beta}$ \eqref{topos annele beta}
	\begin{equation}
		\breve{\beta}^{*}:\rH^i(X_{\overline{\eta},\fet}^{\mathbb{N}^{\circ}},\mathbb{L})\to\rH^i(\widetilde{E}^{\mathbb{N}^{\circ}}_{s},\breve{\beta}^{*}(\mathbb{L}))
		\label{presque iso beta}
	\end{equation}
	est un $\alpha$-isomorphisme.
	\label{presque iso faltings modifie}
\end{coro}

\textit{Preuve}. Pour tout entier $i\ge 0$, on a un diagramme commutatif (cf. \cite{AGT} VI.7.10)
\begin{equation}
	\xymatrix{
		0\ar[r] & \rR^1\varprojlim_{n\ge 1} \rH^{i-1}(X_{\overline{\eta},\fet},\mathbb{L}_{n})\ar[r] \ar[d]_{(\beta_{n}^{*})}& \rH^i(X_{\overline{\eta},\fet}^{\mathbb{N}^{\circ}},\mathbb{L})\ar[r] \ar[d]_{\breve{\beta}^{*}} & \varprojlim_{n\ge 1}\rH^{i}(X_{\overline{\eta},\fet},\mathbb{L}_n)\ar[r] \ar[d]^{(\beta_{n}^{*})} & 0\\
		0\ar[r] & \rR^1\varprojlim_{n\ge 1} \rH^{i-1}(\widetilde{E}_{s},\beta_{n}^{*}(\mathbb{L}_{n}))\ar[r] & \rH^i(\widetilde{E}_{s}^{\mathbb{N}^{\circ}},\breve{\beta}^{*}(\mathbb{L}))\ar[r]  & \varprojlim_{n\ge 1}\rH^{i}(\widetilde{E}_{s},\beta_{n}^{*}(\mathbb{L}_{n}))\ar[r]& 0
}
	\label{suite exact de Jannsen}
\end{equation}
où l'on a posé $\rH^{-1}(X_{\overline{\eta},\fet},\mathbb{L}_n)=0$ et $\rH^{-1}(\widetilde{E}_{s},\beta_{n}^{*}(\mathbb{L}_{n}))=0$ pour tout $n\ge 1$. L'assertion résulte alors de \ref{lemma proj lim alpha nuls}(ii) et \ref{presque iso faltings modifie n}. \\
%

$\hspace*{-1.2em}\bf{\arabic{section}.\stepcounter{theorem}\arabic{theorem}.}$
On désigne par $\widetilde{E}_{pt}$ le topos de Faltings associé au morphisme canonique $h_{S}: \overline{\eta}\to S$. Pour tout $U\in \Ob(\Et_{/S})$, le site $\Et_{\textnormal{f}/U_{\overline{\eta}}}$ est équivalent au site $\Et_{\coh/U_{\overline{\eta}}}$ \eqref{fet et}. D'après (\cite{AGT} VI.5.11), le foncteur $\rho^{+}$ induit une équivalence de topos \eqref{coevanescent Faltings topos}:
\begin{equation}
	\rho : S_{\et}\overleftarrow{\times}_{S_{\et}}\overline{\eta}_{\et} \xrightarrow{\sim} \widetilde{E}_{pt},
	\label{coevanescent Faltings iso trait}
\end{equation}
qui s'insère dans le diagramme commutatif \eqref{diagramme commutatif coevanescente Faltings}
\begin{equation}
	\xymatrix{
		S_{\et} \ar@{=}[d] & S_{\et}\overleftarrow{\times}_{S_{\et}}\overline{\eta}_{\et} \ar[l]_-{\rmp_1} \ar[r]^-{\rmp_2} \ar[d]^{\rho} & \overline{\eta}_{\et} \ar[d]^{\wr}\\
		S_{\et} & \widetilde{E}_{pt} \ar[l]_{\sigma} \ar[r]^{\beta} & \overline{\eta}_{\fet}.
	}
	\label{diagramme commutatif coevanescente Faltings 2}
\end{equation}

$\hspace*{-1.2em}\bf{\arabic{section}.\stepcounter{theorem}\arabic{theorem}.}$
Soient $a: \overline{S}\to S$, $h_{S}:\overline{\eta}\to S$ les morphismes canoniques et $a^{+}: \Et_{\scoh/S}\to \Et_{\scoh/\overline{S}}$, $h_{S}^{+}:\Et_{\scoh/S} \to \Et_{\coh/\overline{\eta}}$ les foncteurs de changement de base associés \eqref{fet et}. On note $D_{\scoh}$ et $E_{\scoh}$ les sites introduits dans \ref{Faltings coevanescent scoh} relativement au morphisme $h_{S}$. On observera que $\overline{\eta}_{\et}$ étant un topos ponctuel, il existe un unique morphisme de topos
\begin{equation}
	b: \overline{S}_{\et} \to \overline{\eta}_{\et},
	\label{morphisme de topos b}
\end{equation}
tel que $b^{*}(\overline{\eta})=\overline{S}$ (\cite{SGAIV} IV 4.3). Il est induit par l'unique foncteur exact à gauche
\begin{displaymath}
	b^{+}:\Et_{\coh/\overline{\eta}}\to \Et_{\scoh/\overline{S}}.
\end{displaymath}

\begin{lemma}
	\textnormal{(i)} Pour tout $U\in \Ob(\Et_{\scoh/S})$, il existe un $\overline{S}$-morphisme fonctoriel en $U$
\begin{equation}
	\overline{U}\to b^{+}(U_{\overline{\eta}}),
	\label{2 morphisme site}
\end{equation}
induisant un isomorphisme au-dessus de $\overline{\eta}$.

\textnormal{(ii)} Pour tout $(V\to U)\in \Ob(E_{\scoh})$, il existe un isomorphisme de $\overline{S}$-schémas fonctoriel en $(V\to U)$
\begin{equation}
	\overline{U}\times_{b^{+}(U_{\overline{\eta}})}b^{+}(V)\xrightarrow{\sim} \overline{U}^{V},
	\label{u pullback iso UV normalisation}
\end{equation}
où $\overline{U}^{V}$ désigne la clôture intégrale de $\overline{U}$ dans $V$ \eqref{ring B bar}.
	\label{morphisme u+}
\end{lemma}

\textit{Preuve}. D'après (\cite{EGA IV} 18.10.8), tout objet $U$ de $\Et_{\scoh/S}$ est la somme de deux sous-schémas ouverts $U_{1}$ et $U_{2}$, où $U_{1}$ est fini sur $S$ et $U_{2}$ est fini sur $\eta$. On peut traiter indépendamment le cas où $U$ est fini sur $S$ et celui où $U$ est fini sur $\eta$.

Si $U$ est fini sur $S$, $U$ est isomorphe à une somme de copies de $S$. Par suite, on a un isomorphisme canonique $\overline{U}\xrightarrow{\sim} b^{+}(\overline{U}_{\overline{\eta}})$, d'où le lemme dans ce cas.

Si $U$ est fini sur $\eta$, on a $\overline{U}\simeq U_{\overline{\eta}}$. La construction de \eqref{2 morphisme site} est immédiate dans ce cas puisque le morphisme composé
\begin{displaymath}
	\Et_{\coh/\overline{\eta}}\xrightarrow{b^{+}}\Et_{\scoh/\overline{S}} \to \Et_{\coh/\overline{\eta}},
\end{displaymath}
où la seconde flèche est déduite du changement de base, est l'identité. De plus, on a $\overline{U}\times_{b^+(U_{\overline{\eta}})}b^{+}(V)\simeq U_{\overline{\eta}}\times_{U_{\overline{\eta}}}V=V$ et $\overline{U}^{V}\simeq V$; d'où l'isomorphisme \eqref{u pullback iso UV normalisation} dans ce cas.
\begin{nothing}
On déduit de \eqref{2 morphisme site} un morphisme de foncteurs
\begin{equation}
	t^{+}:a^{+} \to b^{+} h^{+}_{S},
	\label{2 morphisme site t plus}
\end{equation}
et donc un $2$-morphisme
\begin{displaymath}
	t:h_{S}b\to a.
\end{displaymath}
D'après la propriété universelle des topos co-évanescents (\cite{AGT} VI.3.7), le triplet $(a,b,t)$ induit un morphisme de topos
\begin{equation}
	u:\overline{S}_{\et}\to S_{\et}\overleftarrow{\times}_{S_{\et}}\overline{\eta}_{\et},
	\label{morphisme topos u}
\end{equation}
qui s'insère dans le diagramme commutatif
\begin{equation}
	\xymatrix{
		& \overline{S}_{\et} \ar[ld]_{a} \ar[d]^{u} \ar[rd]^{b}& \\
		S_{\et} \ar@{=}[rd] & S_{\et}\overleftarrow{\times}_{S_{\et}}\overline{\eta}_{\et} \ar[l]_{\rmp_1} \ar[r]^{\rmp_2} & \overline{\eta}_{\et} \ar[ld]^{h_{S}}\\
		& S_{\et} &
	}
	\label{diagramme commutatif bar S coevanescente}
\end{equation}
On vérifie aussitôt que
\begin{equation}
	u(\overline{\eta})=(\id:\overline{\eta}\rightsquigarrow \overline{\eta})\qquad \textnormal{et} \qquad u(s)=(\overline{\eta}\rightsquigarrow s),
	\label{u maps points to points}
\end{equation}
où $\overline{\eta}\rightsquigarrow s$ est l'unique flèche de spécialisation de $\overline{\eta}$ vers $s$. En vertu de (\cite{AGT} VI.3.7.2), $u$ est induit par le foncteur
\begin{equation}
	u^{+}: D_{\scoh}\to \Et_{\scoh/\overline{S}} \qquad (V\to U) \mapsto \overline{U}\times_{b^{+}(U_{\overline{\eta}})}b^{+}(V),
	\label{morphism of site u}
\end{equation}
où la flèche $\overline{U}\to b^{+}(U_{\overline{\eta}})$ est défini dans \eqref{2 morphisme site}. On en déduit que $u^{*}( (\overline{\eta}\to \eta)^{a})=\overline{\eta}$. Le topos $\eta_{\et}\overleftarrow{\times}_{\eta_{\et}}\overline{\eta}_{\et}$ est canoniquement isomorphe au topos localisé de $S_{\et}\overleftarrow{\times}_{S_{\et}}\overline{\eta}_{\et}$ en $(\overline{\eta}\to \eta)^{a}$ (\cite{AGT} VI.3.14). Notons $\widetilde{D}_{s}$ le sous-topos fermé de $S_{\et}\overleftarrow{\times}_{S_{\et}}\overline{\eta}_{\et}$ complémentaire de l'ouvert $(\overline{\eta}\to \eta)^{a}$. Le morphisme $u$ \eqref{morphisme topos u} induit deux morphismes de topos (\cite{SGAIV} IV 5.11 et 9.4.3)
\begin{equation}
	v: \overline{\eta}_{\et} \to \eta_{\et}\overleftarrow{\times}_{\eta_{\et}}\overline{\eta}_{\et}\qquad \textnormal{et} \qquad u_{s}: s_{\et} \to \widetilde{D}_{s}.
	\label{morphisme us v}
\end{equation}
On identifie le morphisme $v$ (resp. $u_{s}$) au point $(\id:\overline{\eta}\rightsquigarrow \overline{\eta})$ (resp. $(\overline{\eta}\rightsquigarrow s)$) de $\widetilde{D}$.
\end{nothing}
\begin{prop}
	Le morphisme $u$ \eqref{morphisme topos u} est une équivalence de topos.
	\label{S bar et et coevanescent}
\end{prop}

D'après (\cite{Il14} 2.3.2), le topos $S_{\et}\overleftarrow{\times}_{S_{\et}}\overline{\eta}_{\et}$ est local de centre $(\overline{\eta}\rightsquigarrow s)$ (\cite{SGAIV} VI 8.4.6). En vertu de \eqref{u maps points to points}, le morphisme de topos $u$ est local (\cite{Il14} 2.1). La proposition résulte alors de \ref{v eta bar}, \ref{us s et} et \ref{lemma de topos ouvert ferme} ci-dessous.

\begin{lemma}
	Le morphisme $v$ \eqref{morphisme us v} est une équivalence de topos.
	\label{v eta bar}
\end{lemma}
\textit{Preuve}. On désigne par $D'_{\scoh}$ le site défini par la sous-catégorie pleine $(D_{\scoh})_{/(\overline{\eta}\to \eta)}$ de $D_{\scoh}$ munie de la topologie induite via le foncteur d'inclusion. En vertu de (\cite{AGT} VI.10.14), le topos $\eta_{\et}\overleftarrow{\times}_{\eta_{\et}}\overline{\eta}_{\et}$ est canoniquement isomorphe au topos des faisceaux d'ensembles sur $D'_{\scoh}$. D'après \eqref{morphism of site u}, le morphisme $v$ est induit par le foncteur
\begin{equation}
	v^{+}: D'_{\scoh} \to \Et_{\coh/\overline{\eta}} \qquad (V\to U)\mapsto V.
	\label{morphisme site v}
\end{equation}
On désigne par $\mathscr{C}_{\overline{\eta}}$ la catégorie des voisinages du point $v$ de $\eta_{\et}\overleftarrow{\times}_{\eta_{\et}}\overline{\eta}_{\et}$ dans le site $D'_{\scoh}$. Les objets de $\mathscr{C}_{\overline{\eta}}$ sont les couples formés d'un objet $(V\to U)$ de $D'_{\scoh}$ et d'un $\overline{\eta}$-morphisme $m$ de $\overline{\eta}$ dans $V$. On désigne par $\mathscr{D}_{\overline{\eta}}$ la sous-catégorie pleine de $\mathscr{C}_{\overline{\eta}}$ formée des objets de la forme $(\overline{\eta}\to U,\id_{\overline{\eta}})$. Il est clair que $\mathscr{D}_{\overline{\eta}}^{\circ}$ forme une sous-catégorie cofinale de $\mathscr{C}_{\overline{\eta}}^{\circ}$.

Pour tout faisceau $F$ de $\eta_{\et}\overleftarrow{\times}_{\eta_{\et}}\overline{\eta}_{\et}$, considérons le morphisme d'adjonction
\begin{equation}
	v^{*}(F)\to v^{*}v_{*}v^{*}(F).
	\label{v vvv adjonction F}
\end{equation}
Celui-ci s'identifie à la limite inductive sur la catégorie $\mathscr{D}_{\overline{\eta}}^{\circ}$ des morphismes
\begin{equation}
	F(\overline{\eta}\to U) \to v_{*}(v^{*}(F))(\overline{\eta}\to U).
	\label{section to fibre canonique}
\end{equation}
D'après \eqref{morphisme site v}, on a un isomorphisme $v_{*}(v^{*}(F))(\overline{\eta}\to U)\simeq v^{*}(F)(\overline{\eta})=v^{*}(F)$. Le morphisme \eqref{section to fibre canonique} s'identifie alors au morphisme canonique $F(\overline{\eta}\to U)\to v^{*}(F)$. En prenant la limite inductive, le morphisme \eqref{v vvv adjonction F} est donc un isomorphisme. Comme le point $v$ est conservatif pour le topos $\eta_{\et}\overleftarrow{\times}_{\eta_{\et}}\overline{\eta}_{\et}$ (\cite{AGT} VI.5.28), on en déduit que le morphisme d'adjonction
\begin{equation}
	\id \to v_{*}v^{*}
	\label{adjonction v}
\end{equation}
est un isomorphisme. Celui-ci implique que $v^{*}$ est pleinement fidèle (\cite{SGAIV} I 6.4). Il est clair que le foncteur fibre $v^{*}$ est essentiellement surjectif, d'où le lemme.

\begin{lemma}
	Soit $\mathcal{T}$ un topos local de point central $s:\Pt\to \mathcal{T}$. On suppose que le point $s$ est conservatif pour $\mathcal{T}$. Alors, $s$ est une équivalence de topos.
	\label{lemma topos ponctuel}
\end{lemma}
\textit{Preuve}. Notons $\varepsilon:\mathcal{T}\to \Pt$ la projection canonique (\cite{SGAIV} IV 4.3). Rappelons que, pour tout faisceau $F$ de $\mathcal{T}$, le morphisme naturel $\Gamma(\mathcal{T},F)\to s^{*}(F)$ est bijectif. Celui-ci induit un isomorphisme canonique $\varepsilon_{*}\xrightarrow{\sim} s^{*}$. On en déduit un isomorphisme
\begin{equation}
	\varepsilon_{*}s_{*}\xrightarrow{\sim} s^{*}s_{*}.
	\label{adjonction and epsilon s}
\end{equation}
On notera que le morphisme composé $\varepsilon\circ s: \Pt \to \mathcal{T} \to \Pt$ s'identifie au morphisme identique. On en déduit que le composé de \eqref{adjonction and epsilon s} et du morphisme d'adjonction
\begin{equation}
	\varepsilon_{*}s_{*}\xrightarrow{\sim} s^{*}s_{*}\to \id
\end{equation}
s'identifie au morphisme identique. Par suite, le morphisme d'adjonction
\begin{equation}
	s^{*}s_{*}\xrightarrow{\sim} \id.
	\label{adjonction isomorphisme s}
\end{equation}
est un isomorphisme. D'autre part, le morphisme composé déduit des morphismes d'adjonction $s^{*}\to s^{*}s_{*}s^{*}\to s^{*}$ s'identifie au morphisme identique (\cite{GR14} 1.1.12). On en déduit un isomorphisme canonique
\begin{equation}
	s^{*}\xrightarrow{\sim}s^{*}s_{*}s^{*}.
	\label{coadjonction compose isomorphisme s}
\end{equation}
Comme le foncteur fibre $s^{*}$ est conservatif, le morphisme d'adjonction
\begin{equation}
	\id \xrightarrow{\sim} s_{*}s^{*}
	\label{coadjonction isomorphisme s}
\end{equation}
est un isomorphisme. Donc $s$ induit une équivalence de topos.

\begin{lemma}
	Le morphisme $u_{s}$ \eqref{morphisme us v} est une équivalence de topos.
	\label{us s et}
\end{lemma}
\textit{Preuve}. Les points $v$ et $u_{s}$ sont conservatifs pour le topos $S_{\et}\overleftarrow{\times}_{S_{\et}}\overline{\eta}_{\et}$ (cf. \ref{points topos Faltings}). Par \ref{v eta bar} et (\cite{SGAIV} IV 9.7.3), on en déduit que le point $u_{s}$ est conservatif pour $\widetilde{D}_{s}$. D'après \ref{lemma topos ponctuel}, il suffit de démontrer que le topos $\widetilde{D}_{s}$ est local de centre $u_{s}$. Considérons le morphisme composé de topos
\begin{equation}
	\xymatrix{
		s_{\et} \ar[r]^-{u_{s}} & \widetilde{D}_{s} \ar[r]^-{i} & S_{\et}\overleftarrow{\times}_{S_{\et}}\overline{\eta}_{\et} \ar[r]^-{\varepsilon}& \Pt,
	}
	\label{morphisme compose set Ds Pt}
\end{equation}
où $i$ est le morphisme d'injection canonique et $\varepsilon$ est la projection canonique. Rappelons (\cite{SGAIV} IV 9.5.8) que le morphisme d'adjonction
\begin{equation}
	i^{*}i_{*}\xrightarrow{\sim} \id
	\label{sous topos ferme i adjonction}
\end{equation}
est un isomorphisme. Comme $S_{\et}\overleftarrow{\times}_{S_{\et}}\overline{\eta}_{\et}$ est local de centre $i\circ u_{s}$ \eqref{u maps points to points}, on a un isomorphisme canonique
\begin{equation}
	\varepsilon_{*}\xrightarrow{\sim} u_{s}^{*}i^{*}.
\end{equation}
On en déduit que le composé
\begin{equation}
	\varepsilon_{*}i_{*}\to u^{*}_{s}i^{*}i_{*}\to u^{*}_{s}
\end{equation}
est un isomorphisme. Le foncteur ``section globale'' $\Gamma(\widetilde{D}_{s},-)$ est donc isomorphe au foncteur fibre $u_{s}^{*}$, i.e. le topos $\widetilde{D}_{s}$ est local de centre $u_{s}$.

\begin{lemma}
	Soient $f:\mathcal{T}'\to \mathcal{T}$ un morphisme de topos, $U$ un ouvert de $\mathcal{T}$ et $U'=f^{*}(U)$. Notons $j:\mathcal{T}_{/U}\to \mathcal{T}$ (resp. $j':\mathcal{T}'_{/U'}\to \mathcal{T}'$) le morphisme de localisation de $\mathcal{T}$ (resp. $\mathcal{T}'$) en $U$ (resp. $U'$) et $i:\mathcal{T}_{s}\to \mathcal{T}$ (resp. $i':\mathcal{T}_{s}'\to \mathcal{T}'$) le sous-topos fermé de $\mathcal{T}$ (resp. $\mathcal{T}'$) complémentaire de l'ouvert $U$ (resp. $U'$) \textnormal{(\cite{SGAIV} IV 9.3.5)}. Le morphisme $f$ induit deux morphismes de topos \textnormal{(\cite{SGAIV} IV 5.11 et 9.4.3)}
	\begin{equation}
		f_{/U}:\mathcal{T}'_{/U'}\to \mathcal{T}_{/U} \qquad\textnormal{et} \qquad f_{s}:\mathcal{T}_{s}'\to \mathcal{T}_{s}.
		\label{fo fs general}
	\end{equation}
	Supposons que $\mathcal{T}_{s}$ (resp. $\mathcal{T}_{s}'$) soit ponctuel, $\mathcal{T}$ (resp. $\mathcal{T}'$) soit local de centre $i$ (resp. $i'$) et que $f_{/U}$ soit une équivalence de topos. Alors, $f$ est une équivalence de topos.
	\label{lemma de topos ouvert ferme}
\end{lemma}
\textit{Preuve}. Posons $\rho'=i'^{*}\circ j'_{*}$ et $\rho=i^{*}\circ j_{*}$. Le topos $\mathcal{T}$ (resp. $\mathcal{T}'$) est équivalent au topos obtenu par recollement de $(\mathcal{T}_{/U},\mathcal{T}_{s},\rho)$ (resp. $(\mathcal{T}_{/U'}',\mathcal{T}_{s}',\rho')$) (cf. \cite{SGAIV} IV 9.5.4). Considérons le diagramme
\begin{equation}
	\xymatrix{
		\mathcal{T}_{/U'}'\ar[r]^{j_{*}'} \ar[d]_{f_{/U*}} & \mathcal{T}' \ar[r]^{i'^{*}} \ar[d]_{f_{*}} & \mathcal{T}_{s}' \ar[d]^{f_{s*}} \\
		\mathcal{T}_{/U}\ar[r]^{j_{*}} & \mathcal{T} \ar[r]^{i^{*}} & \mathcal{T}_{s}
	}
	\label{commutatif rho et rho'}
\end{equation}
dont le carré de gauche est clairement commutatif. Le foncteur $i^{*}$ (resp. $i'^{*}$) est isomorphe au foncteur ``section globale'' $\Gamma(\mathcal{T},-)$ (resp. $\Gamma(\mathcal{T}',-)$) et le foncteur $f_{s*}$ est isomorphe au foncteur identique. On en déduit que le carré de droit de \eqref{commutatif rho et rho'} est aussi commutatif. Les foncteurs $\rho$ et $\rho'$ sont donc compatibles. Par suite, $f_{*}$ induit une équivalence de catégories en vertu de (\cite{SGAIV} IV 9.5.2).

\begin{nothing}
	Soit $\overline{\mathscr{B}}_{pt}$ l'anneau de $\widetilde{E}_{pt}$ associé à $\overline{S}$ \eqref{ring B bar}. Par \eqref{def of bar B}, \eqref{u pullback iso UV normalisation} et \eqref{morphism of site u}, on a un isomorphisme d'anneaux de $\widetilde{E}_{pt}$.
	\begin{equation}
		\rho_{*}u_{*}(\mathscr{O}_{\overline{S}})\xrightarrow{\sim} \overline{\mathscr{B}}_{pt}.
		\label{OSbar ring B}
	\end{equation}
Par les équivalences de topos \eqref{coevanescent Faltings iso trait} et \eqref{morphisme topos u}:
	\begin{equation}
		\overline{S}_{\et}\xrightarrow{u} S_{\et}\overleftarrow{\times}_{S_{\et}}\overline{\eta}_{\et} \xrightarrow{\rho} \widetilde{E}_{pt},
		\label{equivalence of topos}
	\end{equation}
	l'anneau $\mathscr{O}_{\overline{S}}$ de $\overline{S}_{\et}$ s'identifie alors à l'anneau $\overline{\mathscr{B}}_{pt}$ de $\widetilde{E}_{pt}$ et pour tout entier $n\ge 1$, l'anneau $\mathscr{O}_{\overline{S}_{n}}$ de $\overline{S}_{\et}$ s'identifie alors à l'anneau $\overline{\mathscr{B}}_{pt,n}$ de $\widetilde{E}_{pt}$ \eqref{def of Bn}. D'après \eqref{diagramme commutatif coevanescente Faltings 2} et \ref{us s et}, le sous-topos fermé $s_{\et}$ de $\overline{S}_{\et}$ s'identifie au sous-topos fermé $\widetilde{E}_{pt,s}$ \eqref{delta} de $\widetilde{E}_{pt}$.

	D'après \eqref{diagramme commutatif coevanescente Faltings 2} et \eqref{diagramme commutatif bar S coevanescente}, le morphisme de topos $a:\overline{S}_{\et}\to S_{\et}$ induit une équivalence de topos annelés
	\begin{equation}
		a_{n}: (s_{\et},\mathscr{O}_{\overline{S}_{n}}) \xrightarrow{\sim} (s_{\et},\mathscr{O}_{\overline{S}_{n}}),
		\label{equivalence de topos anneles an}
	\end{equation}
	qui s'identifie au morphisme de topos annelés $\sigma_{n}$ \eqref{morphisme de topos sigma n} de $\widetilde{E}_{pt}$. De plus, le morphisme de topos $b:\overline{S}_{\et}\to \overline{\eta}_{\fet}$ induit une équivalence de topos annelés
	\begin{equation}
		b_{n}:(s_{\et},\mathscr{O}_{\overline{S}_{n}}) \xrightarrow{\sim} (\overline{\eta}_{\fet},\oo_{n})
		\label{topos bn equi}
	\end{equation}
	qui s'identifie au morphisme de topos annelés $\beta_{n}$ \eqref{topos beta n 2} de $\widetilde{E}_{pt}$.
	\label{not topos Faltings point}
\end{nothing}

\section{Correspondance de Deninger-Werner via le topos annelé de Faltings} \label{DW via Faltings}
Dans cette section, on se donne une $S$-courbe semi-stable $X$ et un $\overline{\eta}$-point $\overline{x}$ de $X_{\overline{\eta}}$. On pose $C=X_{\overline{\eta}}$. Rappelons que $\overline{X}=X\times_{S}\overline{S}$ est un schéma normal (cf. \ref{S courbe propriete}). On reprend les notations de \S~\ref{Topos de Faltings} pour $X$.

\begin{nothing}
	Pour tout $\overline{\eta}$-point $\overline{y}$ de $C$. Il existe un trait $S'$ fini sur $S$, un morphisme $h_{S'}:\overline{\eta}\to S'$ et un $S$-morphisme $y:S'\to X$ qui s'insèrent dans un diagramme commutatif de schémas
	\begin{equation}
		\xymatrix{
			\overline{\eta}\ar[r]^{h_{S'}} \ar[d]_{\overline{y}} & S' \ar[d]^{y}\\
			C \ar[r]& X
		}
		\label{diagramme y}
	\end{equation}
	On identifie le topos annelé $(\overline{S}_{\et},\mathscr{O}_{\overline{S}})$ au topos annelé de Faltings associé au morphisme $h_{S'}$ et au schéma $\overline{S}$ (cf. \ref{not topos Faltings point}). Pour tout entier $n\ge 1$, le diagramme \eqref{diagramme y} induit par fonctorialité (\cite{AGT} III.8.20) un morphisme de topos annelés qu'on note
\begin{equation}
	[\overline{y}]_{n}: (s_{\et},\mathscr{O}_{\overline{S}_n})\to (\widetilde{E}_{s},\overline{\mathscr{B}}_{n}).
	\label{fibre en point faltings topos}
\end{equation}
Par fonctorialité du topos de Faltings (\cite{AGT} VI.10.12), celui-ci est indépendant des choix du trait $S'$ fini sur $S$ et du diagramme \eqref{diagramme y}.
\label{pullback by a point Faltings}
\end{nothing}
\begin{prop} \label{DW for trivial Faltings modules}
	Supposons $X$ régulier. Soient $n$ un entier $\ge 1$, $M$ un $\overline{\mathscr{B}}_{n}$-module libre de type fini et $\overline{y}$ un point géométrique de $C$. Le morphisme $[\overline{y}]_{n}$ \eqref{fibre en point faltings topos} et le foncteur $(-)_{\sharp}$ \eqref{functor sigma star} induisent un isomorphisme fonctoriel
\begin{equation}
	([\overline{y}]_{n}^{*})_{\sharp}: \Gamma(\widetilde{E}_{s},M)_{\sharp}\xrightarrow{\sim} \Gamma(s_{\et},[\overline{y}]_{n}^{*}(M))_{\sharp}.
	\label{pullback yn iso}
\end{equation}
\end{prop}
\textit{Preuve}. Comme $M$ est libre de type fini, il existe un $\oo_{n}$-module libre de type fini $\mathbb{L}$ de $C_{\fet}$ tel que $M\simeq\beta_{n}^{*}(\mathbb{L})$. D'après (\cite{AGT} VI.(10.12.6)), on a un diagramme commutatif de topos annelés
\begin{equation}
	\xymatrix{
		(s_{\et},\mathscr{O}_{\overline{S}_{n}}) \ar[r]^{ [\overline{y}]_{n}}\ar[d]_{b_n}& (\widetilde{E}_s,\overline{\mathscr{B}}_{n})\ar[d]^{\beta_n}\\
		(\overline{\eta}_{\fet},\oo_n)\ar[r]^{\overline{y}}& (C_{\fet}, \oo_n)
	}
	\label{commutatif topos Psi fibre}
\end{equation}
où $b_{n}$ est l'équivalence de topos annelés \eqref{topos bn equi}, induisant un diagramme commutatif:
\begin{equation}
	\xymatrix{
		\Gamma(C_{\fet},\mathbb{L}) \ar[r]^{\overline{y}^*} \ar[d]_{\beta_n^*} &\Gamma(\overline{\eta}_{\fet},\overline{y}^*(\mathbb{L}))\ar[d]^{b_{n}^{*}}\\
		\Gamma(\widetilde{E}_{s}, M)\ar[r]^-{[\overline{y}]_{n}^{*}} & \Gamma(s_{\et}, [\overline{y}]_{n}^{*} (M))
	}
	\label{le morphisme yn topos Faltings Mn}
\end{equation}
où $\beta_{n}^{*}$ induit un $\alpha$-isomorphisme en vertu de \ref{presque iso faltings modifie n} et $\overline{y}^{*}$ induit un isomorphisme puisque $\mathbb{L}$ est un faisceau constant et $C$ est connexe. On en déduit un $\alpha$-isomorphisme
	\begin{equation}
		[\overline{y}]_{n}^{*}: \Gamma(\widetilde{E}_{s},M)\to \Gamma(s_{\et},[\overline{y}]_{n}^{*}(M)).
	\label{pullback yn presque iso}
	\end{equation}
L'isomorphisme \eqref{pullback yn iso} résulte de \ref{prop Tsuji alpha iso}(i).

\begin{nothing}	\label{fonctoriel Faltings topos}
	Soient $(S',\eta')$ un trait fini sur $(S,\eta)$, $\varphi:X'\to X_{S'}$ un $\eta'$-revêtement semi-stable. On note encore $s$ le point fermé de $S'$ et on fixe un $\eta$-morphisme $\overline{\eta}\to \eta'$ et par suite un $S$-morphisme $\overline{S}\to S'$. On pose $X'_{s}=X'\times_{S'}s$, $\overrightharpoon{X}'=X'\times_{S'}\overline{S}$ et $C'=X'\times_{S'}\overline{\eta}$. On désigne par
	\begin{equation}
		\phi:X'\xrightarrow{\varphi} X_{S'}\to X,~ \overline{\varphi}: \overrightharpoon{X}'\to \overline{X},~ \hbar':\overrightharpoon{X}'\to X' \textnormal{ et } h':C'\to X'
		\label{notation relative to varphi}
	\end{equation}
les morphismes canoniques. On a un diagramme commutatif induit par $\varphi$
	\begin{equation}
		\xymatrix{
			C' \ar[r] \ar[d]_{\varphi_{\overline{\eta}}} & \overrightharpoon{X}' \ar[r]^{\hbar'} \ar[d]^{\overline{\varphi}} & X'\ar[d]^{\phi} \\
			C \ar[r] & \overline{X} \ar[r]^{\hbar} & X.
		}
		\label{C' C X' X}
	\end{equation}
	Comme $X'$ est une $S'$-courbe semi-stable, le schéma $\overrightharpoon{X}'$ est normal \eqref{S courbe propriete}. On désigne par $(\widetilde{E}',\overline{\mathscr{B}}')$ le topos annelé de Faltings associé au schéma $X'$ au-dessus de $S'$ (cf. \ref{basic topos de Faltings} et \ref{ring B bar}) et par
	\begin{equation}
		\Phi: (\widetilde{E}',\overline{\mathscr{B}}') \to (\widetilde{E}, \overline{\mathscr{B}})
		\label{morphisme E' to E}
	\end{equation}
	le morphisme induit par fonctorialité par \eqref{C' C X' X} (\cite{AGT} III.8.20). D'après (\cite{AGT} VI.(10.12.6)), les carrés du diagramme
	\begin{equation}
		\xymatrix{
			X'_{\et} \ar[d]_{\phi_{\et}} &\widetilde{E}' \ar[l]_{\sigma'} \ar[r]^{\beta'} \ar[d]^{\Phi} & C'_{\fet} \ar[d]^{\varphi_{\overline{\eta}}} \\
			X_{\et} & \widetilde{E} \ar[l]_{\sigma} \ar[r]^{\beta} & C_{\fet},
		}
		\label{fonctorialite Phi}
	\end{equation}
où $\beta'$ et $\sigma'$ sont les morphismes canoniques \eqref{foncteur topos beta} et \eqref{foncteur topos sigma} relatifs au $S'$-schéma $X'$, sont commutatifs à isomorphismes canoniques près.

On note $\widetilde{E}'_{s}$ le sous-topos fermé de $\widetilde{E}'$ complémentaire de l'ouvert $\sigma'^{*}(X_{\eta'}')$ \eqref{sous-topos ferme},
	\begin{equation}
		\delta':\widetilde{E}'_{s} \to \widetilde{E}'
		\label{delta'}
	\end{equation}
	le plongement canonique \eqref{delta} et
	\begin{equation}
		\sigma_{s}': \widetilde{E}_{s}' \to X_{s,\et}'
		\label{sigma s}
	\end{equation}
	le morphisme canonique de topos \eqref{morphisme de topos sigma s}. D'après \eqref{fonctorialite Phi}, on a un isomorphisme $\Phi^*(\sigma^*(X_{\eta}))\simeq \sigma'^{*}(X_{\eta'}')$. En vertu de (\cite{SGAIV} IV 9.4.3), il existe donc un morphisme de topos
	\begin{equation}
		\Phi_{s} :\widetilde{E}'_{s}\to \widetilde{E}_{s}
		\label{Phi s}
	\end{equation}
	unique à isomorphisme près tel que le diagramme
	\begin{equation}
		\xymatrix{
			\widetilde{E}_{s}' \ar[r]^{\Phi_{s}} \ar[d]_{\delta'} &\widetilde{E}_{s} \ar[d]^{\delta} \\
			\widetilde{E}' \ar[r]^{\Phi} & \widetilde{E}
		}
		\label{Phi s Phi commute}
	\end{equation}
	soit commutatif. Il résulte de \eqref{fonctorialite Phi} et de (\cite{SGAIV} IV 9.4.3) que le diagramme de morphismes de topos
	\begin{equation}
	\xymatrix{
		\widetilde{E}'_{s} \ar[r]^{\Phi_{s}} \ar[d]_{\sigma_{s}'} & \widetilde{E}_{s} \ar[d]^{\sigma_{s}}\\
		X_{s,\et}' \ar[r]^{\varphi_{s}} & X_{s,\et}
		}
		\label{Phi s varphi s}
	\end{equation}
	est commutatif à isomorphisme près.

	Pour tout entier $n\ge 1$, on pose $\overline{\mathscr{B}}'_{n}=\overline{\mathscr{B}}'/p^{n}\overline{\mathscr{B}}'$, et on désigne par
	\begin{equation}
		\sigma'_{n}: (\widetilde{E}'_{s},\overline{\mathscr{B}}_{n}') \to (X_{s,\et}',\mathscr{O}_{\overrightharpoon{X}_{n}'})
		\label{sigma'n}
	\end{equation}
	le morphisme de topos annelés induit par $\sigma'$ \eqref{morphisme de topos sigma n}, et par
	\begin{equation}
		\beta_{n}': (\widetilde{E}'_{s},\overline{\mathscr{B}}_{n}') \to (C_{\fet}',\oo_n)
		\label{beta'n}
	\end{equation}
	le morphisme de topos annelés induit par $\beta'$ \eqref{topos beta n 2}. L'homomorphisme canonique $\Phi^{-1}(\overline{\mathscr{B}})\to \overline{\mathscr{B}}'$ induit un homomorphisme $\Phi_{s}^{*}(\overline{\mathscr{B}}_{n}) \to \overline{\mathscr{B}}_n'$. Le morphisme $\Phi_{s}$ est donc sous-jacent à un morphisme de topos annelés, que l'on note
	\begin{equation}
		\Phi_{n}: (\widetilde{E}_{s}',\overline{\mathscr{B}}'_{n}) \to (\widetilde{E}_s,\overline{\mathscr{B}}_s).
		\label{Phi n}
	\end{equation}
	Il résulte de \eqref{Phi s varphi s} et de la définition de \eqref{morphisme anneau sigma n} que le diagramme de morphismes de topos annelés
	\begin{equation}
		\xymatrix{
			(\widetilde{E}_{s}', \overline{\mathscr{B}}_n') \ar[r]^{\Phi_n} \ar[d]_{\sigma_n'} & (\widetilde{E}_s, \overline{\mathscr{B}}_n) \ar[d]^{\sigma_n} \\
			(X_{s,\et}',\mathscr{O}_{\overrightharpoon{X}_n'}) \ar[r]^{\overline{\varphi}_{n}}& (X_{s,\et},\mathscr{O}_{\overline{X}_n}),
		}
		\label{fonctorialite sigma n}
	\end{equation}
	où $\overline{\varphi}_{n}$ est la réduction modulo $p^{n}$ de $\overline{\varphi}$, est commutatif à isomorphisme près. D'après \eqref{fonctorialite Phi}, \eqref{Phi s Phi commute} et la définition de \eqref{topos beta n 2}, le diagramme de morphismes de topos annelés
	\begin{equation}
		\xymatrix{
			(\widetilde{E}_s', \overline{\mathscr{B}}_n') \ar[r]^{\Phi_n}  \ar[d]_{\beta_n'} & (\widetilde{E}_s, \overline{\mathscr{B}}_n) \ar[d]^{\beta_n} \\
			(C_{\fet}',\oo_n) \ar[r]^{\varphi_{\overline{\eta}}} &(C_{\fet},\oo_n)
		}
		\label{fonctorialite beta n}
	\end{equation}
	est commutatif à isomorphisme près.

Les morphismes $(\Phi_{n})_{n\ge 1}$ définissent un morphisme de topos annelés
\begin{equation}
	\breve{\Phi}:(\widetilde{E}_s'^{\mathbb{N}^{\circ}},\breve{\overline{\mathscr{B}}}')\to (\widetilde{E}_s^{\mathbb{N}^{\circ}},\breve{\overline{\mathscr{B}}}).
	\label{Phi tout niveau}
\end{equation}

\end{nothing}
\begin{nothing}	\label{G invariant galoisien Faltings topos}
	Conservons les notations de \ref{fonctoriel Faltings topos} et supposons de plus que le $\eta'$-revêtement semi-stable $\varphi:X'\to X_{S'}$ soit galoisien. Posons $G=\Aut(X'_{\eta'}/X_{\eta'})$, qui est isomorphe à $\Aut(C'/C)$. Tout $g\in G$ s'étend en un $X$-automorphisme de $X'$. Par suite, celui-ci induit pour tout entier $n\ge 1$, un morphisme de topos annelés $g_{n}:(\widetilde{E}'_{s},\overline{\mathscr{B}}_{n}')\to (\widetilde{E}'_{s},\overline{\mathscr{B}}_{n}')$ tel que $\Phi_{n}\simeq \Phi_{n}\circ g_{n}$ et que pour tous $g,h\in G$, on ait un isomorphisme
	\begin{displaymath}
		(gh)_{n}^{*}\simeq h_{n}^{*}\circ g_{n}^{*}.
	\end{displaymath}
\end{nothing}

\begin{definition}
(i) On dit qu'un $\overline{\mathscr{B}}_n$-module $M_n$ est \textit{potentiellement libre de type fini} s'il est de type fini, et s'il existe un trait $(S',\eta')$ fini sur $(S,\eta)$ et un $\eta'$-revêtement semi-stable $\varphi:X'\to X_{S'}$, tel que, reprenant les notations de \ref{fonctoriel Faltings topos}, $\Phi_{n}^{*}(M_n)$ soit un $\overline{\mathscr{B}}'_n$-module libre de type fini.

(ii) On dit qu'un $\breve{\overline{\mathscr{B}}}$-module $M=(M_n)_{n\ge 1}$ est \textit{potentiellement libre de type fini} si $M$ est adique de type fini \eqref{limite projective de topos}, et si, pour tout entier $n\ge 1$, le $\overline{\mathscr{B}}_n$-module $M_n$ est potentiellement libre de type fini.
\label{def pltf}
\end{definition}

On désigne par $\Mod^{\lt}(\overline{\mathscr{B}}_n)$ (resp. $\Mod^{\lt}(\breve{\overline{\mathscr{B}}})$) la sous-catégorie pleine de $\Mod(\overline{\mathscr{B}}_n)$ (resp. $\Mod(\breve{\overline{\mathscr{B}}})$), formée des $\overline{\mathscr{B}}_n$-modules potentiellement libres de type fini (resp. $\breve{\overline{\mathscr{B}}}$-modules potentiellement libres de type fini). On désigne par $\Mod^{\lt}_{\mathbb{Q}}(\breve{\overline{\mathscr{B}}})$ la catégorie des objets de $\Mod^{\lt}(\breve{\overline{\mathscr{B}}})$ à isogénie près, qui est une sous-catégorie pleine de $\Mod^{\atf}_{\mathbb{Q}}(\breve{\overline{\mathscr{B}}})$ \eqref{notations of Bmodules E}. Les objets de $\Mod^{\lt}_{\mathbb{Q}}(\breve{\overline{\mathscr{B}}})$ sont appelés \textit{les $\breve{\overline{\mathscr{B}}}_{\mathbb{Q}}$-modules potentiellement libres de type fini}.

\begin{nothing}
	Dans cette section, on reprend les notations de \ref{oo oon reprentations} pour le groupe pro-fini $\pi_{1}(C,\overline{x})$. Soit $n\ge1$ un entier. On note abusivement $\beta_{n}^{*}$, $\breve{\beta}^{*}$ et $\breve{\beta}^{*}_{\mathbb{Q}}$ les foncteurs composés suivants:
\begin{eqnarray}
	&&\beta_{n}^{*}: \Rep_{\oo_{n}}^{\ltf}(\pi_{1}(C,\overline{x}))\xrightarrow{\mu_{\overline{x}}}\Sysl(C_{\fet},\oo_{n})\xrightarrow{\beta_{n}^{*}} \Mod(\overline{\mathscr{B}}_{n}),\label{betan star rep}\\
	&&\breve{\beta}^{*}: \Rep_{\oo}^{\ltf}(\pi_{1}(C,\overline{x})) \xrightarrow{\breve{\mu}_{\overline{x}}} \Sysl(C_{\fet}^{\mathbb{N}^{\circ}},\breve{\oo}) \xrightarrow{\breve{\beta}^{*}} \Mod(\breve{\overline{\mathscr{B}}}),\label{beta star rep}\\
	&&\breve{\beta}_{\mathbb{Q}}^{*}: \Rep_{\mathfrak{C}}^{\cont}(\pi_{1}(C,\overline{x})) \xrightarrow{\breve{\mu}_{\overline{x},\mathbb{Q}}} \Sysl_{\mathbb{Q}}(C_{\fet}^{\mathbb{N}^{\circ}},\breve{\oo})\xrightarrow{\breve{\beta}_{\mathbb{Q}}^{*}} \Mod_{\mathbb{Q}}(\breve{\overline{\mathscr{B}}}), \label{beta Q rep}
\end{eqnarray}
où $\mu_{\overline{x}}$, $\breve{\mu}_{\overline{x}}$ et $\breve{\mu}_{\overline{x},\mathbb{Q}}$ sont les équivalences de catégories définies dans \ref{general notion of rep sysl}.
\label{notation beta n rep}
\end{nothing}
\begin{prop}
	Soit $n$ un entier $\ge 1$. Le foncteur $\beta_{n}^{*}$ (resp. $\breve{\beta}^{*}$, resp. $\breve{\beta}^{*}_{\mathbb{Q}}$) est exact et il se factorise à travers la sous-catégorie $\Mod^{\lt}(\overline{\mathscr{B}}_n)$ (resp. $\Mod^{\lt}(\breve{\overline{\mathscr{B}}})$, resp. $\Mod^{\lt}_{\mathbb{Q}}(\breve{\overline{\mathscr{B}}})$).
	\label{betan lt}
\end{prop}
\textit{Preuve}. L'anneau $\overline{\mathscr{B}}$ est plat sur $\mathcal{O}_{\overline{K}}$ (\cite{AGT} III.9.2) et par suite $\overline{\mathscr{B}}_{n}$ est une $\oo_{n}$-algèbre plate, d'où l'exactitude du foncteur $\beta_{n}^{*}$. On en déduit que les foncteurs $\breve{\beta}^{*}$ et $\breve{\beta}_{\mathbb{Q}}^{*}$ sont exacts.

\'{E}tant donné un objet $\mathbb{L}$ de $\Sysl(C_{\fet},\oo_n)$, il existe un revêtement étale et galoisien $\phi: C'\to C$ tel que $\phi^{*}(\mathbb{L})$ soit un $\oo_{n}$-module libre de type fini de $C_{\fet}'$. D'après \ref{domine par bon revêtement}(i), quitte à remplacer $S$ par une extension finie, il existe un $\eta$-revêtement semi-stable $\varphi:X' \to X$ de fibre géométrique générique $\phi:C'\to C$. Reprenons les notations de \ref{fonctoriel Faltings topos} dans ce cas. D'après \eqref{fonctorialite beta n}, on a un isomorphisme
\begin{equation}
	\Phi_n^*(\beta_n^*(\mathbb{L}))\simeq \beta_n'^*(\varphi_{\overline{\eta}}^*(\mathbb{L})).
	\label{pullback trivial de Bn rhon}
\end{equation}
Par suite, $\beta_{n}^{*}(\mathbb{L})$ est potentiellement libre de type fini. L'assertion pour le foncteur $\breve{\beta}^{*}$ s'ensuit compte tenu de \ref{lemma extensions sysl}(ii).

\begin{prop}
	Supposons $X$ régulier. Le foncteur \eqref{beta Q rep} est pleinement fidèle.	
	\label{beta pullback pleinement fidele}
\end{prop}
\textit{Preuve}. \'{E}tant donnés deux objets $\mathbb{L}_{1}, \mathbb{L}_{2}$ de $\Sysl(C_{\fet}^{\mathbb{N}^{\circ}},\breve{\oo})$, on pose
\begin{displaymath}
	\mathbb{L}=\FHom_{\breve{\oo}}(\mathbb{L}_{1},\mathbb{L}_{2}),
\end{displaymath}
qui est encore un $\breve{\oo}$-module localement libre de type fini de $C_{\fet}^{\mathbb{N}^{\circ}}$. Comme le $\breve{\oo}$-module $\mathbb{L}_{1}$ est localement libre de type fini, on a un isomorphisme canonique
\begin{equation}
	\breve{\beta}^{*}(\mathbb{L})\xrightarrow{\sim}\FHom_{\breve{\overline{\mathscr{B}}}}(\breve{\beta}^{*}(\mathbb{L}_{1}),\breve{\beta}^{*}(\mathbb{L}_{2})).
	\label{Vn presentation finie}
\end{equation}
D'après \ref{presque iso faltings modifie} et \eqref{Vn presentation finie}, le foncteur $\breve{\beta}$ induit un isomorphisme:
\begin{equation}
	\Hom_{\breve{\oo}}(\mathbb{L}_{1},\mathbb{L}_{2})\otimes_{\mathbb{Z}}\mathbb{Q} \xrightarrow{\sim} \Hom_{\breve{\overline{\mathscr{B}}}}(\breve{\beta}^{*}(\mathbb{L}_{1}),\breve{\beta}^{*}(\mathbb{L}_{2}))\otimes_{\mathbb{Z}}\mathbb{Q};
	\label{Gamma Hom beta}
\end{equation}
d'où la pleine fidélité de \eqref{beta Q rep}.

\begin{lemma}\label{general notations of VnM}
	Soit $M$ un $\overline{\mathscr{B}}_n$-module potentiellement libre de type fini. Le $\oo_{n}$-module 
	\begin{equation}
		\mathscr{W}_{n}(M)=\Gamma(s_{\et}, [\overline{x}]_n^{*}(M))_{\sharp},
		\label{def rho Mn x}
	\end{equation}
	où $[\overline{x}]_{n}: (s_{\et},\mathscr{O}_{\overline{X}_{n}})\to (\widetilde{E}_{s},\overline{\mathscr{B}}_{n})$ est le morphisme de topos annelés induit par $\overline{x}$ \eqref{fibre en point faltings topos} et $(-)_{\sharp}$ est le foncteur \eqref{functor sigma star}, est $\alpha$-plat de type $\alpha$-fini \textnormal{(\ref{type de alpha fini}, \ref{alpha plat})}. 
\end{lemma}
\textit{Preuve}. D'après \ref{domine par bon revêtement}(i), il existe un trait $(S',\eta')$ fini sur $(S,\eta)$ et un $\eta'$-revêtement semi-stable, régulier et galoisien $\varphi:X'\to X_{S'}$ tel que, fixant un $\eta$-morphisme $\overline{\eta}\to \eta'$ et reprenant les notations de \ref{fonctoriel Faltings topos}, $\Phi_{n}^{*}(M)$ soit libre de type fini et que chaque $\overline{\eta}$-point de $C'$ au-dessus de $\overline{x}$ se descend en un $\eta'$-point de $X_{\eta'}'$.

Choisissons un $\overline{\eta}$-point $\overline{y}$ de $C'$ au-dessus de $\overline{x}$, et un $\eta'$-point de $X_{\eta'}'$ qui le descend. On note $y$ le $S'$-point de $X'$ correspondant et $x:S'\to X$ le $S$-morphisme induit par $\varphi$. Le diagramme commutatif de couples de schémas
\begin{displaymath}
	\xymatrix{
		(\overline{\eta}\to S') \ar[dr]_-{(\overline{x},x)} \ar[r]^-{(\overline{y},y)} & (C'\to X') \ar[d]^{(\varphi_{\overline{\eta}},\phi)} \\
	& (C\to X)}
\end{displaymath}
induit un diagramme commutatif de topos annelés:
\begin{displaymath}
	\xymatrix{
		(s_{\et},\mathscr{O}_{\overline{S}_n}) \ar[dr]_{[\overline{x}]_n} \ar[r]^{[\overline{y}]_{n}} & (\widetilde{E}_{s}',\overline{\mathscr{B}}_{n}') \ar[d]^{\Phi_n}\\
		& (\widetilde{E}_{s},\overline{\mathscr{B}}_{n})}
\end{displaymath}
On en déduit un isomorphisme de $\oo_{n}$-modules $\mathscr{W}_{n}(M)\simeq \Gamma(s_{\et},[\overline{y}]_{n}^{*}(\Phi_n^*(M)))_{\sharp}$. Comme $\Phi_n^*(M)$ est un $\overline{\mathscr{B}}_n'$-module libre de type fini, l'assertion s'ensuit compte tenu de \ref{coro alpha plat tf}.

\begin{nothing}
	On fixe un entier $\ge 1$. Dans la suite de cette section, on se propose de construire un foncteur de $\Mod^{\lt}(\overline{\mathscr{B}}_n)$ dans $\Rep_{\oo_{n}}^{\aptf}(\pi_{1}(C,\overline{x}))$ \eqref{oo oon rep aptf}. \'{E}tant donné un $\overline{\mathscr{B}}_n$-module potentiellement libre de type fini $M$, on définit le $\oo_n$-module $\mathscr{W}_{n}(M)$ sous-jacent à la $\oo_{n}$-représentation associée à $M$ par \eqref{def rho Mn x}. Reprenant les notations de la preuve de \ref{general notations of VnM}, le morphisme $[\overline{y}]_{n}$ induit un isomorphisme $\oo_{n}$-linéaire (cf. \ref{DW for trivial Faltings modules})
\begin{equation}
	([\overline{y}]_{n}^{*})_{\sharp}:\Gamma(\widetilde{E}'_{s}, \Phi^*_n (M))_{\sharp} \xrightarrow{\sim} \mathscr{W}_{n}(M).
	\label{fibre yn iso}
\end{equation}
	
	On construit une action $\varrho_{M}$ de $\pi_1(C,\overline{x})$ sur $\mathscr{W}_{n}(M)$ comme suit. Posons $G=\Aut(X'_{\eta'}/X_{\eta'})$ qui est isomorphe à $\Aut(C'/C)$. D'après \ref{G invariant galoisien Faltings topos}, tout $g\in G$ induit un morphisme de topos annelés $g_{n}:(\widetilde{E}'_{s},\overline{\mathscr{B}}_{n}')\to (\widetilde{E}'_{s},\overline{\mathscr{B}}_{n}')$ et par suite un automorphisme
	\begin{equation}
		(g_{n}^{*})_{\sharp}: \Gamma(\widetilde{E}'_{s}, \Phi^*_n (M))_{\sharp} \xrightarrow{\sim} \Gamma(\widetilde{E}'_{s}, \Phi^*_n (M))_{\sharp}.
		\label{Gamma n autormorphisme}
	\end{equation}
	On a un homomorphisme surjectif \eqref{psi y homo}
	\begin{equation}
		\xi_{\overline{y}}: \pi_{1}(C,\overline{x}) \to G^{\op}.
		\label{psi bar y 2}
	\end{equation}
	Pour tout $\gamma\in \pi_{1}(C,\overline{x})$, posant $g=\xi_{\overline{y}}(\gamma)$, on définit un automorphisme $\varrho_{M}(\gamma)$ de $\mathscr{W}_{n}(M)$ par
	\begin{equation}
		\varrho_{M}(\gamma): \mathscr{W}_{n}(M)\xrightarrow{(([\overline{y}]_{n}^{*})_{\sharp})^{-1}} \Gamma(\widetilde{E}'_{s}, \Phi^*_n(M))_{\sharp} \xrightarrow{(g_{n}^{*})_{\sharp}} \Gamma(\widetilde{E}'_{s}, \Phi^*_n(M))_{\sharp} \xrightarrow{([\overline{y}]_{n}^{*})_{\sharp}} \mathscr{W}_{n}(M).
		\label{def of varrho M}
	\end{equation}
	On notera que $\varrho_{M}$ est l'homomorphisme composé
\begin{equation}
	\varrho_{M}: \pi_{1}(C,\overline{x}) \xrightarrow{\xi_{\overline{y}}} G^{\op} \to \Aut_{\oo_{n}}\Gamma(\widetilde{E}'_{s}, \Phi^*_n(M))_{\sharp} \xrightarrow[\sim]{([\overline{y}]_{n}^{*})_{\sharp}} \Aut_{\oo_{n}} \mathscr{W}_{n}(M).
	\label{def of varrho M compose}
\end{equation}
On définit la $\oo_{n}$-représentation de $\pi_{1}(C,\overline{x})$ associée à $M$ par $(\mathscr{W}_{n}(M),\varrho_{M})$ qu'on note simplement $\mathscr{W}_{n}(M)$.
	\label{construction of varrho}
\end{nothing}

\begin{prop}
	\'{E}tant donné un $\overline{\mathscr{B}}_n$-module potentiellement libre de type fini $M$, la $\oo_{n}$-représentation $\mathscr{W}_{n}(M)$ de $\pi_{1}(C,\overline{x})$ est indépendante des choix du trait $(S',\eta')$ fini sur $(S,\eta)$, du $\eta'$-revêtement semi-stable, régulier et galoisien $\varphi: X'\to X_{S'}$ et du point $\overline{y}$ de $C'(\overline{\eta})$ au-dessus de $\overline{x}$ à isomorphisme près.
	\label{indepandence varrhoM}
\end{prop}
\textit{Preuve}. Le $\oo_{n}$-module \eqref{def rho Mn x} sous-jacent à $\mathscr{W}_{n}(M)$ est clairement indépendant des données ci-dessus. Il suffit de démontrer que l'isomorphisme \eqref{def of varrho M} est indépendant des mêmes données.

Reprenons les notations de \ref{general notations of VnM} et soient $\overline{y}$ et $\overline{y}'$ deux $\overline{\eta}$-points de $C'$ au-dessus de $\overline{x}$. Il existe un élément $h\in G\simeq\Aut(C'/C)$ tel que $h(\overline{y})=\overline{y}'$. Avec les notations de \ref{homomorphisme psi conjugaison}, on a deux homomorphismes
\begin{equation}
	\xi_{\overline{y}}: \pi_{1}(C,\overline{x}) \to G^{\op} \qquad \xi_{\overline{y}'}: \pi_{1}(C,\overline{x}) \to G^{\op}.
	\label{2 psi morphismes}
\end{equation}
reliés par la relation \eqref{conjugaison de psi}:
\begin{equation}
	\xi_{\overline{y}'}(\gamma)=h\xi_{\overline{y}}(\gamma)h^{-1} \qquad \forall \gamma\in \pi_{1}(C,\overline{x}).
	\label{conjugaison de psi 2}
\end{equation}
On déduit par fonctorialité un diagramme commutatif
\begin{equation}
	\xymatrix{
		\Gamma(\widetilde{E}'_{s}, \Phi_{n}^{*} (M))_{\sharp} \ar[d]_{([\overline{y}']^{*}_{n})_{\sharp}} \ar[r]^{(h_{n}^{*})_{\sharp}} & \Gamma(\widetilde{E}'_{s}, \Phi_{n}^{*} (M))_{\sharp} \ar[ld]^{([\overline{y}]_{n}^{*})_{\sharp}}\\
		\mathscr{W}_{n}(M). &
	}
\end{equation}
Pour tout $\gamma\in \pi_{1}(C,\overline{x})$, posant $g=\xi_{\overline{y}}(\gamma)$ et $g'=\xi_{\overline{y}'}(\gamma)$, on en déduit que $g'=hgh^{-1}$ dans $G$ et que
\begin{equation}
	([\overline{y}']^{*}_{n})_{\sharp}\circ (g_n'^{*})_{\sharp} \circ ([\overline{y}']^{*}_{n})_{\sharp}^{-1} \simeq ([\overline{y}]_{n}^{*})_{\sharp}\circ (h_{n}^{*})_{\sharp}\circ(g_{n}'^{*})_{\sharp}\circ (h_{n}^{*})_{\sharp}^{-1} \circ (([\overline{y}]_{n}^{*})_{\sharp})^{-1} =([\overline{y}]_{n}^{*})_{\sharp}\circ (g_{n}^{*})_{\sharp} \circ (([\overline{y}]_{n}^{*})_{\sharp})^{-1}.
	\label{independance of y}
\end{equation}
Ceci implique que l'automorphisme $\varrho_{M}(\gamma)$ \eqref{def of varrho M} est indépendant du choix du point $\overline{y}$ au-dessus de $\overline{x}$.

Soient $(S'',\eta'')$ un trait fini sur $(S,\eta)$ et $\phi:X''\to X_{S''}$ un $\eta''$-revêtements semi-stable, régulier et galoisien trivialisant $M$. D'après \ref{domine par bon revêtement}(i), on peut supposer que $S''$ domine $S'$ et qu'il existe un morphisme équivariant de $\phi$ dans $\varphi\times_{S'}S''$ (cf. \ref{def of revêtement}). Notons-le $\psi: X''\to X'\times_{S'}S''$ et soit $u: \Aut(X''_{\eta''}/X_{\eta''})\to \Aut(X'_{\eta''}/X_{\eta''})$ l'homomorphisme induit par $\psi_{\eta''}$. Alors, le morphisme équivariant $\psi$ induit un morphisme $u$-équivariant de topos annelés de Faltings:
\begin{equation}
	\Psi_{n}:(\widetilde{E}''_{s}, \overline{\mathscr{B}}''_{n})\to (\widetilde{E}'_{s},\overline{\mathscr{B}}_{n}'),
\end{equation}
où la source désigne le topos annelé de Faltings associé au $S''$-schéma $X''$, qui est muni d'une action de $\Aut(X_{\eta''}''/X_{\eta''})$ \eqref{G invariant galoisien Faltings topos}. On en déduit, pour tout $g\in \Aut(X_{\eta''}''/X_{\eta''})$, un diagramme commutatif
\begin{equation}
	\xymatrix{
		\Gamma(\widetilde{E}'_{s}, \Phi^*_n(M))_{\sharp} \ar[r]^{(u(g)_{n}^{*})_{\sharp}} \ar[d]_{(\Psi_{n}^{*})_{\sharp}} & \Gamma(\widetilde{E}'_{s}, \Phi^*_n(M))_{\sharp} \ar[d]^{(\Psi_{n}^{*})_{\sharp}} \\
		\Gamma(\widetilde{E}''_{s}, \Psi_{n}^{*}(\Phi^*_n(M)))_{\sharp} \ar[r]^{(g_{n}^{*})_{\sharp}} & \Gamma(\widetilde{E}''_{s}, \Psi_{n}^{*}(\Phi^*_n(M)))_{\sharp}
	}
\end{equation}
La proposition s'ensuit compte tenu de \eqref{def of varrho M}.

\begin{nothing} \label{fonctorialite de varrho}
	Soit $f:M\to M'$ un morphisme de $\overline{\mathscr{B}}_n$-modules potentiellement libres de type fini. D'après \ref{domine par bon revêtement}(i), il existe un trait $(S',\eta')$ fini sur $(S,\eta)$ et un $\eta'$-revêtement semi-stable, régulier et galoisien $\varphi: X'\to X_{S'}$ tels que, reprenant les notations de \ref{fonctoriel Faltings topos}, $\Phi_{n}^{*}(M)$ et $\Phi_{n}^{*}(M')$ soient libres de type fini. L'automorphisme \eqref{def of varrho M} est fonctoriel en $M$. Par \ref{indepandence varrhoM}, on peut associer à $f:M\to M'$ un morphisme de $\oo_n$-représentations de $\pi_1(C,\overline{x})$
\begin{equation}
	\mathscr{W}_n(f):\mathscr{W}_{n}(M)\to \mathscr{W}_{n}(M').
	\label{varrho n f}
\end{equation}
La correspondance $M\mapsto \mathscr{W}_{n}(M)$ définit ainsi un foncteur $\oo_{n}$-linéaire
\begin{equation}
	\mathscr{W}_{n}:\Mod^{\lt}(\overline{\mathscr{B}}_{n})\to \Rep_{\oo_{n}}^{\aptf}(\pi_{1}(C,\overline{x})).
	\label{foncteur varrho n}
\end{equation}
\end{nothing}
\begin{lemma} \label{rho M entier}
	Soit $M=(M_{n})_{n\ge 1}$ un $\breve{\overline{\mathscr{B}}}$-module potentiellement libre de type fini. Pour tout entier $n\ge 1$, le morphisme canonique $M_{n+1}\to M_{n}$ induit un $\alpha$-isomorphisme fonctoriel de $\oo_{n}$-représentations \eqref{alpha iso cat Rep} 
\begin{equation}
	\mathscr{W}_{n+1}(M_{n+1})\otimes_{\oo_{n+1}}\oo_{n}\to \mathscr{W}_{n}(M_{n}).
	\label{Wn+1 mod pn Wn iso}
\end{equation}
\end{lemma}
\textit{Preuve}. Soit $n$ un entier $\ge 1$. D'après \ref{domine par bon revêtement}(i), il existe un trait $(S',\eta')$ fini sur $(S,\eta)$ et un $\eta'$-revêtement semi-stable, régulier et galoisien $\varphi: X'\to X_{S'}$ trivialisant $M_{n+1}$. On a un isomorphisme de $\oo_{n}$-modules 
\begin{equation}
	\Gamma(s_{\et},[\overline{x}]_{n+1}^{*}(M_{n+1}))\otimes_{\oo_{n+1}}\oo_{n}\xrightarrow{\sim} \Gamma(s_{\et},[\overline{x}]_{n}^{*}(M_{n})). 
\end{equation}
En vertu de \eqref{N to Mstar}, celui-ci induit un $\alpha$-isomorphisme canonique de $\oo_{n}$-modules
\begin{equation}
	\mathscr{W}_{n+1}(M_{n+1})\otimes_{\oo_{n+1}}\oo_{n}\to \mathscr{W}_{n}(M_{n}).
\end{equation}
D'après \ref{indepandence varrhoM} et la fonctorialité de \eqref{def of varrho M}, ce dernier est $\pi_{1}(C,\overline{x})$-équivariant. On vérifie aussitôt que \eqref{Wn+1 mod pn Wn iso} est fonctoriel en $M$ (cf. \ref{fonctorialite de varrho}).\\

$\hspace*{-1.2em}\bf{\arabic{section}.\stepcounter{theorem}\arabic{theorem}.}$
En vertu de \eqref{foncteur aptf oon oo} et \ref{rho M entier}, les foncteurs $\mathscr{W}_{n}$ induisent un foncteur $\oo$-linéaire
	\begin{equation}
		\mathscr{W}: \Mod^{\lt}(\breve{\overline{\mathscr{B}}})\to \Rep_{\oo}^{\aptf}(\pi_{1}(C,\overline{x}))\qquad M=(M_{n})_{n\ge 1}\mapsto \varprojlim \mm\otimes_{\oo}\mathscr{W}_{n}(M_{n}).
		\label{varrho entier}
	\end{equation}
	On le note aussi $\mathscr{W}_{X}$ pour signifier que la construction dépend de $X$. D'après \ref{prop oo aptf C rep}, le foncteur $\mathscr{W}$ s'étend en un foncteur $\mathfrak{C}$-linéaire:
\begin{equation}
	\mathscr{W}_{\mathbb{Q}}: \Mod^{\lt}_{\mathbb{Q}}(\breve{\overline{\mathscr{B}}})\to \Rep_{\mathfrak{C}}^{\cont}(\pi_{1}(C,\overline{x})).
	\label{varrho Q}
\end{equation}

\begin{prop} \label{WX WX' compatible}
	Soient $(S',\eta')$ un trait fini sur $(S,\eta)$, $X'$ un $S'$-modèle semi-stable de $C$ dominant $X_{S'}$. Reprenant les notations de \ref{fonctoriel Faltings topos} pour $\varphi:X'\to X_{S'}$, l'image inverse par $\breve{\Phi}$ \eqref{Phi tout niveau} d'un $\breve{\overline{\mathscr{B}}}$-module potentiellement libre de type fini est un $\breve{\overline{\mathscr{B}}}'$-module potentiellement libre de type fini. De plus, on a un isomorphisme de foncteurs
	\begin{equation}
		\mathscr{W}_{X}\xrightarrow{\sim}\mathscr{W}_{X'}\circ \breve{\Phi}^{*}.
		\label{varrho modeles compatible}
	\end{equation}
\end{prop}
\textit{Preuve}. Soient $n$ un entier $\ge 1$ et $M_{n}$ un $\overline{\mathscr{B}}_{n}$-module potentiellement libre de type fini. D'après \ref{domine par bon revêtement}(i), quitte à remplacer $S'$ par une extension finie, il existe un $\eta'$-revêtement semi-stable, régulier et galoisien $\psi: X''\to X'$ tel que $\Psi_{n}^{*}(\Phi_{n}^{*}(M_{n}))$ soit un $\overline{\mathscr{B}}_{n}''$-module libre de type fini, où $(\widetilde{E}_{s}'',\overline{\mathscr{B}}_{n}'')$ désigne le topos annelé de Faltings associé à la $S'$-courbe semi-stable $X''$, $\Psi_{n}$ le morphisme de topos annelés induit par $\psi$ \eqref{Phi n}. Donc $\Phi_{n}^{*}(M_{n})$ est un $\overline{\mathscr{B}}_{n}'$-module potentiellement libre de type fini. On en déduit la première proposition.

Notons $\mathscr{W}_{X,n}(M_{n})$ et $\mathscr{W}_{X',n}(\Phi_{n}^{*}(M_{n}))$ les $\oo_{n}$-représentations associées à $M_{n}$ et $\Phi_{n}^{*}(M_{n})$, respectivement. On notera que $\varphi\circ\psi$ est aussi un $\eta'$-revêtement semi-stable, régulier et galoisien de $X_{S'}$. D'après \eqref{fibre yn iso}, on a un isomorphisme de $\oo_{n}$-modules 
\begin{equation}
	\mathscr{W}_{X,n}(M_{n})\xrightarrow{\sim}\mathscr{W}_{X',n}(\Phi_{n}^{*}(M_{n})).
	\label{iso WXn WX'n}
\end{equation}
En vertu de \eqref{def of varrho M} et \ref{indepandence varrhoM}, celui-ci est $\pi_{1}(C,\overline{x})$-équivariant. On vérifie aussitôt que l'isomorphisme \eqref{iso WXn WX'n} est fonctoriel en $M_{n}$ (cf. \ref{fonctorialite de varrho}). L'isomorphisme \eqref{varrho modeles compatible} s'ensuit.

\begin{prop}
	Soient $n$ un entier $\ge 1$ et $V$ un objet de $\Rep_{\oo_{n}}^{\ltf}(\pi_{1}(C,\overline{x}))$. Alors, on a un $\alpha$-isomorphisme canonique et fonctoriel de $\oo_{n}$-représentations
	\begin{equation}
		V\rightarrow \mathscr{W}_{n}(\beta_{n}^{*}(V))
		\label{varrho n B n}
	\end{equation}
	où $\beta_{n}^{*}$ est le foncteur \eqref{betan star rep}.
	\label{rep generalise rep vrai n}
\end{prop}

\textit{Preuve}. Soit $\phi:C'\to C$ un revêtement étale et galoisien tel que l'action de $\pi_1(C,\overline{x})$ sur $V$ se factorise à travers $\Aut(C'/C)^{\op}$. D'après \ref{domine par bon revêtement}(i), il existe un trait $(S',\eta')$ fini sur $(S,\eta)$ et un $\eta'$-revêtement semi-stable, régulier et galoisien $\varphi:X'\to X_{S'}$ de fibre géométrique générique $\varphi_{\overline{\eta}}=\phi: C'\to C$. Reprenant les notations de \ref{fonctoriel Faltings topos} et \ref{G invariant galoisien Faltings topos} pour $\varphi:X'\to X_{S'}$, tout $g\in \Aut(C'/C)$ induit un diagramme commutatif de topos annelés (cf. \eqref{fonctorialite beta n})
\begin{equation}
	\xymatrix{
		(\widetilde{E}_s',\overline{\mathscr{B}}'_n) \ar[r]^{\beta_n'} \ar[d]_{g_{n}} & (C'_{\fet},\oo_n) \ar[d]^{g}\\
		(\widetilde{E}_s',\overline{\mathscr{B}}'_n) \ar[r]^{\beta_n'}& (C'_{\fet},\oo_n)
	}
	\label{fonctorialite de G beta}
\end{equation}
On note $\mathbb{L}=\mu_{\overline{x}}(V)$ le $\oo_n$-module localement libre de type fini de $C_{\fet}$ associé à $V$ \eqref{Rep on to Mod on}. Choisissons un $\overline{\eta}$-point $\overline{y}$ de $C'$ au-dessus $\overline{x}$. On désigne par 
\begin{equation}
	\tau(V): V\xrightarrow{\sim} \mathscr{W}_{n}(\beta_{n}^{*}(\mathbb{L}))
	\label{FV definition}
\end{equation}
le $\alpha$-isomorphisme de $\oo_{n}$-modules défini par le composé
\begin{equation}
	V\xrightarrow{\sim} \Gamma(C_{\fet}',\phi^{*}(\mathbb{L})) \rightarrow \Gamma(\widetilde{E}_s',\beta_n'^*(\phi^{*}(\mathbb{L}))) \rightarrow \Gamma(\widetilde{E}_s',\beta_n'^*(\phi^{*}(\mathbb{L})))_{\sharp}\xrightarrow{\sim} \Gamma(\widetilde{E}_s',\Phi_n^*(\beta_{n}^{*}(\mathbb{L})))_{\sharp}\xrightarrow{([\overline{y}]_{n}^{*})_{\sharp}} \mathscr{W}_{n}(\beta_{n}^{*}(\mathbb{L}))
	\label{V vers section globale}
\end{equation}
où la première flèche est un isomorphisme puisque $\phi^{*}(\mathbb{L})$ est constant et $C'$ est connexe, la deuxième flèche est le $\alpha$-isomorphisme induit par $\beta'^{*}_{n}$ \eqref{presque iso beta n}, la troisième flèche est un $\alpha$-isomorphisme en vertu de \eqref{N to Mstar}, la quatrième flèche provient de \eqref{fonctorialite beta n} et la dernière flèche est l'isomorphisme \eqref{pullback yn presque iso}. D'après \eqref{fonctorialite de G beta}, pour tout $g\in \Aut(C'/C)$, le morphisme \eqref{V vers section globale} s'insère dans un diagramme commutatif:
\begin{equation}
	\xymatrix{
		V\ar[r]^-{\sim} \ar[d]_{g} & \Gamma(C_{\fet}',\phi^{*}(\mathbb{L}))\ar[d]^{g^{*}} \ar[r] & \Gamma(\widetilde{E}_s',\beta_n'^*(\phi^{*}(\mathbb{L})))_{\sharp} \ar[d]^{(g_{n}^*)_{\sharp}}  \ar[r]^-{\sim} &  \Gamma(\widetilde{E}_s',\Phi_n^*(\beta_{n}^{*}(\mathbb{L})))_{\sharp} \ar[r]^-{([\overline{y}]_{n}^*)_{\sharp}} \ar[d]^{(g_{n}^*)_{\sharp}}& \mathscr{W}_{n}(\beta_{n}^{*}(\mathbb{L})) \\
		V\ar[r]^-{\sim} & \Gamma(C_{\fet}',\phi^{*}(\mathbb{L})) \ar[r] & \Gamma(\widetilde{E}_s',\beta_n'^*(\phi^{*}(\mathbb{L})))_{\sharp} \ar[r]^-{\sim} &  \Gamma(\widetilde{E}_s',\Phi_n^*(\beta_{n}^{*}(\mathbb{L}))_{\sharp} \ar[r]^-{([\overline{y}]_{n}^*)_{\sharp}}& \mathscr{W}_{n}(\beta_{n}^{*}(\mathbb{L}))
	}
	\label{comparaison de descente}
\end{equation}
Celui-ci implique que le $\alpha$-isomorphisme $\tau(V)$ \eqref{FV definition} est compatible aux actions de $\pi_1(C,\overline{x})$ en vertu de \eqref{def of varrho M}. Par suite, on obtient un $\alpha$-isomorphisme de $\Rep_{\oo_n}^{\aptf}(\pi_1(C,\overline{x}))$
\begin{equation}
	\tau(V): V\to \mathscr{W}_n(\beta_{n}^{*}(V)).
	\label{Transformation naturelle F(V)}
\end{equation}

Soit $f:V'\to V$ un morphisme de $\Rep_{\oo_n}^{\ltf}(\pi_1(C,\overline{x}))$. Choisissons un revêtement étale et galoisien $\phi:C'\to C$ trivialisant les représentations $V$ et $V'$, un trait $(S',\eta')$ fini sur $(S,\eta)$ et un $\eta'$-revêtement semi-stable, régulier et galoisien $\varphi:X'\to X_{S'}$ de fibre générique géométrique $\phi:C'\to C$. Le diagramme commutatif \eqref{comparaison de descente} étant fonctoriel en $V$, on en déduit un diagramme commutatif
\begin{equation}
	\xymatrixcolsep{4pc}\xymatrix{
		V'\ar[r]^{f} \ar[d]_{\tau(V')} & V \ar[d]^{\tau(V)} \\
		\mathscr{W}_n(\beta_{n}^{*}(V')) \ar[r]^{(\mathscr{W}_n \circ \beta_{n}^{*})(f)} & \mathscr{W}_n(\beta_{n}^{*}(V))
	}
\end{equation}
d'où la fonctorialité du $\alpha$-isomorphisme \eqref{varrho n B n}.

\begin{coro}
	\textnormal{(i)} Soit $V$ un objet de $\Rep^{\ltf}_{\oo}(\pi_{1}(C,\overline{x}))$. On a un isomorphisme canonique fonctoriel
	\begin{equation}
		\mm\otimes_{\oo}V\xrightarrow{\sim} \mathscr{W}(\breve{\beta}^{*}(V)).
	\end{equation}

	\textnormal{(ii)} Le foncteur composé $\mathscr{W}_{\mathbb{Q}}\circ \breve{\beta}^{*}_{\mathbb{Q}}$ est isomorphe au foncteur identique.
	\label{rep generalise rep vrai}
\end{coro}
\textit{Preuve}. Pour tout entier $n\ge 1$, on pose $V_{n}=V/p^{n}V$. D'après \ref{rep generalise rep vrai n}, on a un $\alpha$-isomorphisme canonique et fonctoriel $(V_{n})_{n\ge 1}\to (\mathscr{W}_{n}(\beta_{n}^{*}(V_{n})))_{n\ge 1}$ de $\Rep_{\breve{\oo}}^{\aptf}(\pi_{1}(C,\overline{x}))$ \eqref{oo oon rep aptf}. L'assertion (i) résulte alors de \ref{F1 envoie alpha iso iso}(i) et \ref{lemma initial aptf rep}(ii). L'assertion (ii) résulte de (i) et de l'isomorphisme canonique $(\mm\otimes_{\oo}V)[\frac{1}{p}]\to V[\frac{1}{p}]$. \\

$\hspace*{-1.2em}\bf{\arabic{section}.\stepcounter{theorem}\arabic{theorem}.}$
Posons $\check{\overline{X}}=X\times_{S}\check{\overline{S}}$ \eqref{notations 11}. Pour tout entier $n\ge 1$ et tout $\mathscr{O}_{\check{\overline{X}}}$-module $\mathcal{F}$, on considère $\mathcal{F}_{n}=\mathcal{F}\otimes_{\mathscr{O}_{S}}\mathscr{O}_{S_{n}}$ comme un faisceau de $X_{s,\et}$. On pose $\breve{\mathcal{F}}=(\mathcal{F}_{n})_{n\ge 1}$ que l'on considère aussi comme un faisceau de $X_{s,\et}^{\mathbb{N}^{\circ}}$. D'après (\cite{AGT} III.7.18), on a un isomorphisme $\breve{\sigma}^{*}(\breve{\mathcal{F}})\simeq (\sigma_{n}^{*}(\mathcal{F}_{n}))_{n\ge 1}$, où $\sigma_{n}$ et $\breve{\sigma}$ sont des morphismes de topos annelés \eqref{morphisme de topos sigma n} et \eqref{morphisme de topos sigma limit}. Soit $n$ un entier $\ge 1$. Reprenant les notations de \ref{def of categorie DW}, on note abusivement $\sigma_{n}^{*}$ et $\breve{\sigma}^{*}$ les foncteurs
\begin{eqnarray}
	\sigma_{n}^{*}: \BB_{\check{\overline{X}}}^{\DW} \to \Mod(\overline{\mathscr{B}}_n), \qquad \mathcal{F}\mapsto \sigma_{n}^{*}(\mathcal{F}_{n}) \label{sigma n star DW}\\
	\breve{\sigma}^{*}: \BB_{\check{\overline{X}}}^{\DW} \to \Mod^{\atf}(\breve{\overline{\mathscr{B}}}), \qquad \mathcal{F}\mapsto \breve{\sigma}^{*}(\breve{\mathcal{F}}) \label{breve sigma star DW}.
\end{eqnarray}
%
\begin{prop}
	Soit $n$ un entier $\ge 1$. Alors, le foncteur $\sigma_{n}^{*}$ \eqref{sigma n star DW} (resp. $\breve{\sigma}^{*}$ \eqref{breve sigma star DW}) se factorise à travers la sous-catégorie $\Mod^{\lt}(\overline{\mathscr{B}}_{n})$ (resp. $\Mod^{\lt}(\breve{\overline{\mathscr{B}}})$).
\end{prop}
\textit{Preuve}. Soient $n$ un entier $\ge 1$ et $\mathcal{F}$ un fibré vectoriel de Deninger-Werner sur $\check{\overline{X}}$. D'après \ref{domine par bon revêtement}(i), quitte à remplacer $S$ par une extension finie, il existe un $\eta$-revêtement semi-stable $\varphi: X'\to X$ tel que $\overline{\varphi}_{n}$ trivialise $\mathcal{F}_n$. Reprenant les notations de \ref{fonctoriel Faltings topos} pour $\varphi$. D'après \eqref{fonctorialite sigma n}, on a un isomorphisme
\begin{equation}
	\Phi_{n}^{*}(\sigma_{n}^{*}(\mathcal{F}_{n}))\simeq \sigma'^{*}_{n}(\overline{\varphi}_{n}^{*}(\mathcal{F}_{n})).
\end{equation}
Le $\overline{\mathscr{B}}_{n}$-module $\sigma_n^*(\mathcal{F}_n)$ est donc potentiellement libre de type fini. Comme $\breve{\mathcal{F}}$ est adique, on en déduit que le $\breve{\overline{\mathscr{B}}}$-module $\breve{\sigma}^{*}(\breve{\mathcal{F}})$ est potentiellement libre de type fini.

\begin{prop}
	Soient $n$ un entier $\ge 1$ et $\mathcal{F}$ un fibré vectoriel de Deninger-Werner sur $\check{\overline{X}}$. On a un $\alpha$-isomorphisme canonique et fonctoriel
	\begin{equation}
		\mathbb{V}_{n}(\mathcal{F})\rightarrow \mathscr{W}_{n}(\sigma_{n}^{*}(\mathcal{F}_{n}))
		\label{commutatif rho varrho n}
	\end{equation}
	où $\mathbb{V}_{n}$ est le foncteur de Deninger-Werner \eqref{rho VVn}.
	\label{2 def DW sont meme}
\end{prop}
\textit{Preuve}. D'après \ref{domine par bon revêtement}(i), quitte à remplacer $S$ par une extension finie, il existe un $\eta$-revêtement semi-stable, régulier et galoisien $\varphi:X' \to X$ tel que $\overline{\varphi}_{n}$ trivialise $\mathcal{F}_n$. Reprenant les notations de \ref{pullback by a point Faltings} et \ref{fonctoriel Faltings topos}, on a un diagramme commutatif de topos annelés:
\begin{equation}
	\xymatrix{
		(s_{\et}, \mathscr{O}_{\overline{S}_{n}}) \ar@/_2pc/[dd]_{a_{n}} \ar[d]^{\wr} \ar[r]_{[\overline{y}]_{n}} \ar@/^2pc/[rr]^{[\overline{x}_n]} & (\widetilde{E}'_s,\overline{\mathscr{B}}_{n}')\ar[d]^{\sigma'_n} \ar[r]_{\Phi_n} & (\widetilde{E}_s,\overline{\mathscr{B}}_{n}) \ar[d]^{\sigma_n}\\
		(s_{\et}, \mathscr{O}_{\overline{S}_n})\ar[r]^{y_n}\ar[d]^{\wr} & (X'_{s,\et},\mathscr{O}_{\overline{X}_n'}) \ar[r]^{\overline{\varphi}_n} \ar[d]^{u_n'}& (X_{s,\et},\mathscr{O}_{\overline{X}_n}) \ar[d]^{u_{n}}\\
		(s_{\zar}, \mathscr{O}_{\overline{S}_n})\ar[r]^{y_n} \ar@/_2pc/[rr]_{\overline{x}_n} & (X'_{s,\zar},\mathscr{O}_{\overline{X}_n'}) \ar[r]^{\overline{\varphi}_n}& (X_{s,\zar},\mathscr{O}_{\overline{X}_n})
	}
	\label{diagramme commutatif topos Faltings et}
\end{equation}
où $u_{n}$ et $u_{n}'$ sont les morphismes de topos annelés canoniques. Rappelons \eqref{def of Fxn} et \eqref{def rho Mn x} que les $\oo_{n}$-modules sous-jacents aux représentations $\mathbb{V}_{n}(\mathcal{F})$ et $\mathscr{W}_n(\sigma_n^*(\mathcal{F}_n))$ sont définis par
\begin{equation}
	\mathbb{V}_{n}(\mathcal{F})=\Gamma(s_{\zar},\overline{x}_{n}^{*}(\mathcal{F}_{n})) \quad \textnormal{et}\quad \mathscr{W}_{n}(\sigma_{n}^{*}(\mathcal{F}_{n}))=\Gamma(s_{\et}, [x_n]^{*}(\sigma_n^*(\mathcal{F}_n)))_{\sharp}.
	\label{def rho sigma Fn}
\end{equation}
Comme on a $[\overline{x}_n]^{*}(\sigma_n^*(\mathcal{F}_n))\simeq a_{n}^{*}(\overline{x}_{n}^{*}(\mathcal{F}_n))$ \eqref{diagramme commutatif topos Faltings et}, l'équivalence de topos annelés $a_{n}$ et le foncteur $(-)_{\sharp}$ \eqref{functor sigma star} induisent un $\alpha$-isomorphisme de $\oo_{n}$-modules
\begin{equation}
	a_{n}^{*}: \mathbb{V}_{n}(\mathcal{F})\rightarrow \mathscr{W}_{n}(\sigma_{n}^{*}(\mathcal{F}_{n})).
	\label{isomorphisme entre fibre xn}
\end{equation}

Posons $C'=X_{\overline{\eta}}'$ et $G=\Aut(C'/C)$. Reprenant les notations de \ref{G invariant galoisien Faltings topos}, les carrés commutatifs de droite de \eqref{diagramme commutatif topos Faltings et} sont compatibles aux actions de $G$ sur $(\widetilde{E}_{s}',\overline{\mathscr{B}}_{n}')$, $(X_{s,\et}',\mathscr{O}_{\check{\overline{X}}'_{n}})$ et $(X_{s,\zar}',\mathscr{O}_{\check{\overline{X}}'_{n}})$. On en déduit par le foncteur $(-)_{\sharp}$, pour tout $g\in G$, un diagramme commutatif:
\begin{equation}
	\xymatrix{
		\mathbb{V}_{n}(\mathcal{F}) \ar[r]^-{(\overline{y}_{n}^*)^{-1}} \ar[d]_{a_{n}^*} & \Gamma(X_{s,\zar}',\overline{\varphi}_n^{*}(\mathcal{F}_n)) \ar[d]^-{(u_{n}'\sigma_n')^*} \ar[r]^{g_{n}^{*}}& \Gamma(X_{s,\zar}',\overline{\varphi}_n^{*}(\mathcal{F}_n)) \ar[d]^-{(u_{n}'\sigma_n')^*} \ar[r]^-{\overline{y}_{n}^*} & \mathbb{V}_{n}(\mathcal{F}) \ar[d]^{a_{n}^{*}}\\
		\mathscr{W}_{n}(\sigma_{n}^{*}(\mathcal{F}_{n})) \ar[r]^-{(([\overline{y}]_{n}^{*})_{\sharp})^{-1}}  & \Gamma(\widetilde{E}'_s, (\sigma_n')^*(\varphi_n^*(\mathcal{F}_n)))_{\sharp} \ar[r]^{(g_{n}^{*})_{\sharp}} &\Gamma(\widetilde{E}'_s, (\sigma_n')^*(\varphi_n^*(\mathcal{F}_n)))_{\sharp} \ar[r]^-{([\overline{y}]_{n}^{*})_{\sharp}} &  \mathscr{W}_{n}(\sigma_{n}^{*}(\mathcal{F}_{n}))}
	\label{compatible of rho Fn rho sigma Fn}
\end{equation}
Celui-ci implique que le $\alpha$-isomorphisme \eqref{isomorphisme entre fibre xn} est compatible aux actions de $\pi_{1}(C,\overline{x})$ en vertu de \eqref{def of rho n F gamma} et \eqref{def of varrho M}. On en déduit un $\alpha$-isomorphisme de $\oo_{n}$-représentations de $\pi_{1}(C,\overline{x})$
\begin{equation}
	\tau(\mathcal{F}): \mathbb{V}_{n}(\mathcal{F})\rightarrow\mathscr{W}_{n}(\sigma_{n}^{*}(\mathcal{F})).
	\label{lambda F definition}
\end{equation}

Soit $f:\mathcal{F}'\to \mathcal{F}$ un morphisme de $\BB^{\DW}_{\check{\overline{X}}}$. D'après \ref{domine par bon revêtement}(i), quitte à remplacer $S$ par une extension finie, il existe un $\eta$-revêtement semi-stable et régulier $\varphi:X' \to X$ tels que $\overline{\varphi}_{n}$ trivialise $\mathcal{F}_n$ et $\mathcal{F}_{n}'$. Le diagramme commutatif \eqref{compatible of rho Fn rho sigma Fn} est fonctoriel en $\mathcal{F}$. On en déduit par \ref{def of DW n independant} et \ref{indepandence varrhoM} un diagramme commutatif
\begin{equation}
	\xymatrixcolsep{4pc}\xymatrix{
		\mathbb{V}_{n}(\mathcal{F}') \ar[r]^{\mathbb{V}_{n}(f)} \ar[d]_{\tau(\mathcal{F}')} & 	\mathbb{V}_{n}(\mathcal{F}) \ar[d]^{\tau(\mathcal{F})}\\
		\mathscr{W}_{n}(\sigma_{n}^{*}(\mathcal{F}'))\ar[r]^{\mathscr{W}_{n}\circ\sigma_{n}^{*}(f)} & \mathscr{W}_{n}(\sigma_{n}^{*}(\mathcal{F}'))
	}
\end{equation}
d'où la fonctorialité du $\alpha$-isomorphisme \eqref{commutatif rho varrho n}.

\begin{coro}
	Soit $\mathcal{F}$ un fibré vectoriel de Deninger-Werner sur $\check{\overline{X}}$. On a un isomorphisme canonique et fonctoriel de $\oo$-représentations continues de $\pi_{1}(C,\overline{x})$
	\begin{equation}
		\mm\otimes_{\oo}\mathbb{V}(\mathcal{F})\xrightarrow{\sim} \mathscr{W}(\breve{\sigma}^{*}(\breve{\mathcal{F}}))
		\label{commutatif rho varrho entier}
	\end{equation}
	où $\mathbb{V}$ est le foncteur de Deninger-Werner \eqref{rho entier}.
	\label{compatible of 2 rho}
\end{coro}
\textit{Preuve}. D'après \ref{2 def DW sont meme}, on a un $\alpha$-isomorphisme canonique et fonctoriel $(\mathbb{V}_{n}(\mathcal{F}))_{n\ge 1}\to (\mathscr{W}_{n}(\sigma_{n}^{*}(\mathcal{F}_{n})))_{n\ge 1}$ de $\Rep_{\breve{\oo}}^{\aptf}(\pi_{1}(C,\overline{x}))$ \eqref{oo oon rep aptf}. L'assertion résulte alors de \ref{F1 envoie alpha iso iso}(i) et \ref{lemma initial aptf rep}(ii).

\section{Un énoncé de descente pour les fibrés de Deninger-Werner} \label{Descente Gal}
\begin{nothing}
	Soient $\mathcal{T}$ un topos et $G$ un groupe fini. \textit{Une action de $G$ sur $\mathcal{T}$} est la donnée pour tout $g\in G$ d'un morphisme de topos $g:\mathcal{T}\to\mathcal{T}$ tels que $(\id)^{*}=\id_{\mathcal{T}}$ et que pour tous $g,h\in G$, on ait un isomorphisme
	\begin{equation}
		c_{g,h}:g^{*}h^{*}\xrightarrow{\sim}(hg)^{*},
	\end{equation}
	vérifiant des relations de cocycle pour la composition (\cite{SGAI} VI 7.4). Supposons que $\mathcal{T}$ soit muni d'une telle action de $G$. \textit{Un objet $G$-équivariant} de $\mathcal{T}$ est la donnée d'un faisceau $\mathscr{F}$ de $\mathcal{T}$ et pour tout $g\in G$ d'un isomorphisme
	\begin{equation}
		\tau_{g}^{\mathscr{F}}: \mathscr{F}\xrightarrow{\sim}g^{*}(\mathscr{F}),
		\label{systeme de G actions}
	\end{equation}
	tels que $\tau_{\id}^{\mathscr{F}}=\id_{\mathscr{F}}$ et que pour tous $g,h\in G$, on ait
	\begin{equation}
		\tau_{gh}^{\mathscr{F}}=c_{h,g}\circ h^{*}(\tau_{g}^{\mathscr{F}})\circ \tau_{h}^{\mathscr{F}}.
	\end{equation}
	Une telle donnée est appelée \textit{une action de $G$} (à droite) sur $\mathscr{F}$. \'{E}tant donnés deux objets $G$-équivariants $\mathscr{F}_{1}$ et $\mathscr{F}_{2}$, un morphisme $f: \mathscr{F}_{1}\to \mathscr{F}_{2}$ de $\mathcal{T}$ est dit \textit{$G$-équivariant}, si pour tout $g\in G$, le diagramme
	\begin{equation}
		\xymatrix{
			\mathscr{F}_{1} \ar[r]^{f} \ar[d]_{\tau_{g}^{\mathscr{F}_{1}}} & \mathscr{F}_{2}\ar[d]^{\tau_{g}^{\mathscr{F}_{2}}}\\
			g^{*}(\mathscr{F}_{1}) \ar[r]^{g^{*}(f)} & g^{*}(\mathscr{F}_{2})
		}
	\end{equation}
	est commutatif. On désigne par $\mathcal{T}_{G}$ la catégorie des objets $G$-équivariants de $\mathcal{T}$.
	
	Soit $A$ un anneau de $\mathcal{T}$ muni d'une action de $G$ compatible avec sa structure d'anneau, i.e. tel que les isomorphismes $\tau^{A}_{g}$ \eqref{systeme de G actions} soient des isomorphismes d'anneaux. \textit{Un $A$-module $G$-équivariant} est la donnée d'un $A$-module $M$ de $\mathcal{T}$ et d'une action de $G$ sur $M$ compatible avec sa structure de $A$-module.
	\label{notation G action}
\end{nothing}

\begin{nothing} \label{topos de Faltings star}
	Dans la suite de cette section, on se donne une $S$-courbe semi-stable $X$ et un revêtement étale et galoisien $C'$ de $C=X_{\overline{\eta}}$. On pose $G=\Aut(C'/C)$ et on note $\psi:C'\to C$ le morphisme canonique. D'après \ref{X bar normal}, le schéma $\overline{X}$ est normal. On reprend les notations de \S~\ref{Topos de Faltings} pour $X$. On note $E^{\sharp}$ le site défini par la catégorie $E_{/(C'\to X)}$ \eqref{basic topos de Faltings} munie de la topologie induite par le foncteur canonique
\begin{equation}
	\gamma: E_{/(C'\to X)}\to E.
	\label{gamma C'}
\end{equation}
Le site $E^{\sharp}$ est alors canoniquement équivalent au site de Faltings associé au morphisme canonique $C'\to X$ (cf. \cite{AGT} VI.10.1 et VI.10.14). On désigne par $\widetilde{E}^{\sharp}$ le topos localisé de $\widetilde{E}$ en $\beta^{*}(C')=(C'\to X)^{a}$ qui est canoniquement équivalent au topos des faisceaux d'ensembles sur le site $E^{\sharp}$ (\cite{AGT} VI.10.14) et par
	\begin{equation}
		\Phi: \widetilde{E}^{\sharp} \to \widetilde{E}
		\label{Phi star}
	\end{equation}
	le morphisme de localisation. En vertu de (\cite{SGAIV} III 2.5), celui-ci s'identifie au morphisme de topos induit par fonctorialité du topos de Faltings (\cite{AGT} VI.10.12). On désigne par $\widetilde{E}^{\sharp}_{s}$ le sous-topos fermé de $\widetilde{E}^{\sharp}$ complémentaire de l'ouvert $\Phi^{*}(\sigma^{*}(X_{\eta}))$ de $\widetilde{E}^{\sharp}$ et par
	\begin{equation}
		\delta^{\sharp}: \widetilde{E}_{s}^{\sharp}\to \widetilde{E}^{\sharp}
		\label{delta sharp}
	\end{equation}
le prolongement canonique. En vertu de (\cite{SGAIV} IV 9.4.3), il existe un morphisme de topos canonique
	\begin{equation}
		\Phi_{s}: \widetilde{E}^{\sharp}_{s}\to \widetilde{E}_{s}.
		\label{Phis sharp}
	\end{equation}
	On pose $\overline{\mathscr{B}}^{\sharp}=\Phi^{*}(\overline{\mathscr{B}})$ \eqref{ring B bar} et, pour tout entier $n\ge 1$, $\overline{\mathscr{B}}^{\sharp}_{n}=\overline{\mathscr{B}}^{\sharp}/p^{n}\overline{\mathscr{B}}^{\sharp}$ qui est isomorphe à l'anneau $\Phi^{*}(\overline{\mathscr{B}}_{n})$. Comme $\overline{\mathscr{B}}_{n}$ est un objet de $\widetilde{E}_{s}$ \eqref{sous-topos ferme}, on en déduit que l'anneau $\overline{\mathscr{B}}_{n}^{\sharp}$ appartient à $\widetilde{E}_{s}^{\sharp}$. Le morphisme de topos \eqref{Phis sharp} est donc sous-jacent à un morphisme de topos annelés qu'on note
\begin{equation}
	\Phi_{n}: (\widetilde{E}^{\sharp}_{s}, \overline{\mathscr{B}}^{\sharp}_{n}) \to (\widetilde{E}_{s}, \overline{\mathscr{B}}_{n}).
	\label{Phi star n}
\end{equation}
\end{nothing}

\begin{nothing}
	Par fonctorialité du topos de Faltings (\cite{AGT} VI.10.12), l'action de $G$ sur $C'$ induit une action de $G$ sur le topos $\widetilde{E}^{\sharp}$ \eqref{notation G action}. En particulier, tout $g\in G$ induit un morphisme $g:\widetilde{E}^{\sharp}\to \widetilde{E}^{\sharp}$ tel que
	\begin{equation}
		\Phi\xrightarrow{\sim} \Phi\circ g.
		\label{G action sur Esharp}
	\end{equation}
	On en déduit que $g^{*}(\Phi^{*}(\sigma^{*}(X_{\eta})))\simeq \Phi^{*}(\sigma^{*}(X_{\eta}))$. L'action de $G$ sur $\widetilde{E}^{\sharp}$ induit alors une action de $G$ sur $\widetilde{E}_{s}^{\sharp}$. Pour tout $g\in G$, on a donc un morphisme $g_{s}:\widetilde{E}_{s}^{\sharp}\to \widetilde{E}_{s}^{\sharp}$ tel que
	\begin{equation}
		\Phi_{s}\xrightarrow{\sim} \Phi_{s}\circ g_{s}.
		\label{G action sur Esharps}
	\end{equation}
	On désigne par $\widetilde{E}_{G}^{\sharp}$ la catégorie des objets $G$-équivariants de $\widetilde{E}^{\sharp}$ \eqref{notation G action}. D'après \eqref{G action sur Esharp}, l'image inverse de $\Phi$ induit un foncteur
	\begin{equation}
		\Phi^{*}: \widetilde{E} \to \widetilde{E}^{\sharp}_{G}.
		\label{Phi E Esharp G}
	\end{equation}
	En particulier, l'anneau $\overline{\mathscr{B}}_{n}^{\sharp}$ est muni d'une action de $G$ compatible avec sa structure d'anneau. Pour tout entier $n\ge 1$ et tout $g\in G$, le morphisme de topos $g_{s}$ est sous-jacent à un morphisme de topos annelés
	\begin{displaymath}
		g_{n}:(\widetilde{E}_{s}^{\sharp},\overline{\mathscr{B}}_{n}^{\sharp})\to (\widetilde{E}_{s}^{\sharp},\overline{\mathscr{B}}_{n}^{\sharp})
	\end{displaymath}
	tel qu'on ait un isomorphisme canonique
\begin{equation}
	\Phi_{n}\xrightarrow{\sim} \Phi_{n}\circ g_{n}.
\end{equation}
L'image inverse de $\Phi_{n}$ induit alors un foncteur
	\begin{equation}
		\Phi^{*}_{n}: \Mod(\overline{\mathscr{B}}_{n})\to \Mod_{G}(\overline{\mathscr{B}}^{\sharp}_{n}),
		\label{Phin E Esharp G Mod}
	\end{equation}
	où $\Mod_{G}(\overline{\mathscr{B}}_{n}^{\sharp})$ désigne la sous-catégorie pleine de $\widetilde{E}_{G}^{\sharp}$ formée des $\overline{\mathscr{B}}_{n}^{\sharp}$-modules $G$-équivariants de $\widetilde{E}_{s}^{\sharp}$ \eqref{notation G action}.
\end{nothing}

\begin{prop}
	Le foncteur $\Phi^{*}:\widetilde{E}\to \widetilde{E}^{\sharp}_{G}$ \eqref{Phi E Esharp G} est une équivalence de catégories.
	\label{lemme de descente}
\end{prop}

\textit{Preuve}. Considérons
\begin{equation}
	\mathbf{F}\to E
	\label{FAIS E}
\end{equation}
le topos fibré associé au site $E$ (cf. \cite{Gi71} II 3.4.1): la catégorie fibre de $\mathbf{F}$ au-dessus de $(V\to U)\in \Ob(E)$ est le topos localisé $\widetilde{E}_{/(V\to U)^{a}}$ de $\widetilde{E}$ en $(V\to U)^{a}$, et pour tout morphisme $h:(V'\to U')\to (V\to U)$ de $E$ le foncteur image inverse $h^{*}: \widetilde{E}_{/(V\to U)^{a}}\to \widetilde{E}_{/(V'\to U')^{a}}$ est le foncteur image inverse par le morphisme de localisation par rapport à $h$. On rappelle que $\mathbf{F}$ est un champ au-dessus de $E$ (cf. \cite{Gi71} II 3.4.4).

Le morphisme $\varpi=(\psi,\id):(C'\to X) \to (C\to X)$ forme un recouvrement de $(C\to X)$ dans $E$ (cf. \ref{basic topos de Faltings}(v)) que l'on note (abusivement) encore $\varpi$. La catégorie $\widetilde{E}=\mathbf{F}(C\to X)$ est donc équivalente à la catégorie $\mathbf{F}(\varpi)$ des données de descente relativement au recouvrement $\varpi$.

Comme le morphisme $\psi: C'\to C$ est un torseur sous $G$ de $\Et_{\textnormal{f}/C}$, on a un isomorphisme dans $E$
\begin{equation}
	(C'\to X)\times G \simeq (C'\to X)\times_{(C\to X)} (C'\to X).
	\label{G torsor}
\end{equation}
On note encore $G$ le faisceau constant de $\widetilde{E}$ de valeur $G$. Le morphisme $(C'\to X)^{a}\to (C\to X)^{a}$ est alors un torseur sous $G$ de $\widetilde{E}$ (\cite{Gi71} III 1.4.1). Par la descente galoisienne (cf. \cite{BLR90} 6.2.B), la catégorie des données de descente $\mathbf{F}(\varpi)$ relativement au recouvrement $\varpi$ est équivalente à la catégorie $\widetilde{E}^{\sharp}_{G}$, d'où la proposition.

\begin{coro}
	Le foncteur $\Phi_{n}^{*}:\Mod(\overline{\mathscr{B}}_{n})\to \Mod_{G}(\overline{\mathscr{B}}_{n}^{\sharp})$ \eqref{Phin E Esharp G Mod} est une équivalence de catégories.
	\label{lemme de descente Bn}
\end{coro}

\begin{nothing}
	Pour tout objet $U$ de $\Et_{/X}$, on pose $U'_{\overline{\eta}}=U\times_{X}C'$ et on note $\varphi_{U}: U'_{\overline{\eta}}\to U_{\overline{\eta}}$ le morphisme canonique. Considérons le site fibré
	\begin{equation}
		\pi^{\sharp}=\pi\circ \gamma: E^{\sharp}\to \Et_{/X},
		\label{site fibre pi sharp}
	\end{equation}
	dont la fibre au-dessus de $U$ est le site $\Et_{\textnormal{f}/U_{\overline{\eta}}'}$ \eqref{morphisme can sharp}. Pour tout faisceau $\mathscr{G}$ de $\widetilde{E}^{\sharp}$, on note $\mathscr{G}_{U}$ sa restriction à $\Et_{\textnormal{f}/U_{\overline{\eta}}'}$. Le topos $\widetilde{E}^{\sharp}$ est canoniquement équivalent au topos de Faltings associé au morphisme canonique $(C'\to X)$. La donnée d'un faisceau $\mathscr{G}$ de $\widetilde{E}^{\sharp}$ est alors équivalente à la donnée pour tout objet $U$ de $\Et_{/X}$ d'un faisceau $\mathscr{G}_{U}$ de $U_{\overline{\eta},\fet}'$ vérifiant des conditions de comptabilité et de recollement (cf. \cite{AGT} (VI.5.11.1)). Comme $\varphi_{U}$ est un torseur sous $G$ pour la topologie étale de $U_{\overline{\eta}}$, le topos $U_{\overline{\eta},\fet}'$ est muni d'une action de $G$ \eqref{notation G action} induite par fonctorialité. D'après (\cite{AGT} VI.5.11), la donnée d'une action de $G$ sur $\mathscr{G}$ est équivalente à la donnée pour tout objet $U$ de $\Et_{/X}$ d'une action de $G$ sur $\mathscr{G}_{U}$ compatible aux morphismes de restriction.
\label{F to FU G action}
\end{nothing}

\begin{nothing}
	Soit $\mathscr{F}=\{U\mapsto \mathscr{F}_{U}\}_{U\in \Ob(\Et_{/X})}$ un préfaisceau sur $E$ tel que pour tout objet $U$ de $\Et_{/X}$, $\mathscr{F}_{U}$ soit un faisceau de $U_{\overline{\eta},\fet}$. Pour tout objet $U$ de $\Et_{/X}$, la flèche $\varphi_{U}:U'_{\overline{\eta}}\to U_{\overline{\eta}}$ est un objet de $\Et_{\textnormal{f}/U_{\overline{\eta}}}$. Donc, $\varphi_{U,\fet}^{*}$ est un foncteur de restriction. On en déduit un isomorphisme canonique de préfaisceaux sur $E^{\sharp}$
	\begin{equation}
		\mathscr{F}\circ \gamma \xrightarrow{\sim} \{U\mapsto \varphi_{U,\fet}^{*}(\mathscr{F}_{U})\},
		\label{resctriction prefaisceau}
	\end{equation}
	où $\gamma$ est le foncteur \eqref{gamma C'}. Le faisceau associé à $\mathscr{F}\circ \gamma$ est isomorphe à $\Phi^{*}(\mathscr{F}^{a})$. Le faisceau $\varphi_{U,\fet}^{*}(\mathscr{F}_{U})$ de $U_{\overline{\eta},\fet}'$ est muni d'une action de $G$ induite par $\varphi_{U,\fet}^{*}$ fonctorielle en $U$. D'après \ref{F to FU G action}, ils induisent une action de $G$ sur $\Phi^{*}(\mathscr{F}^{a})$ qui coïncide avec celle de $\Phi^{*}(\mathscr{F}^{a})$ induite par $\Phi^{*}$ \eqref{Phi E Esharp G}.
	\label{les actions de G sur Esharp prefaisceau}
\end{nothing}

\begin{nothing} \label{categorie P affine}
	On désigne par $\mathbf{P}$ la sous-catégorie pleine de $\Et_{/X}$ formée des schémas affines. On munit $\mathbf{P}$ de la topologie induite par celle de $\Et_{/X}$. Alors, $\mathbf{P}$ est une famille topologiquement génératrice du site $\Et_{/X}$. On désigne par
\begin{equation}
	\pi_{\mathbf{P}}^{\sharp}: E_{\mathbf{P}}^{\sharp} \to \mathbf{P}
	\label{pi P star}
\end{equation}
le site fibré déduit de $\pi^{\sharp}$ \eqref{site fibre pi sharp} par changement de base par $\mathbf{P}\to \Et_{/X}$. On munit $E_{\mathbf{P}}^{\sharp}$ de la topologie co-évanescente définie par $\pi_{\mathbf{P}}^{\sharp}$ et on note $\widetilde{E}_{\mathbf{P}}^{\sharp}$ le topos des faisceaux d'ensembles sur $E_{\mathbf{P}}$. D'après (\cite{AGT} VI.5.21 et VI.5.22), la topologie de $E_{\mathbf{P}}^{\sharp}$ est induite par celle de $E^{\sharp}$ au moyen du foncteur de projection canonique $E_{\mathbf{P}}^{\sharp}\to E^{\sharp}$, et celui-ci induit par restriction une équivalence de catégories
\begin{equation}
	\widetilde{E}^{\sharp}\xrightarrow{\sim} \widetilde{E}^{\sharp}_{\mathbf{P}}.
	\label{equivalence P Faltings}
\end{equation}
On désigne par $\widehat{E}_{\mathbf{P}}^{\sharp}$ la catégorie des préfaisceaux d'ensembles sur $E_{\mathbf{P}}$. Comme $E_{\mathbf{P}}^{\sharp}$ est une sous-catégorie topologiquement génératrice de $E^{\sharp}$, le foncteur ``faisceau associé'' sur $E_{\mathbf{P}}^{\sharp}$ induit un foncteur que l'on note aussi
\begin{equation}
	\widehat{E}_{\mathbf{P}}^{\sharp} \to \widetilde{E}^{\sharp},\qquad \mathscr{F}\to \mathscr{F}^{a}.
\end{equation}
Soient $\mathscr{G}=\{W\mapsto \mathscr{G}_{W}\}$ ($W\in \Ob(\Et_{/X})$) un préfaisceau sur $E^{\sharp}$, $\mathscr{G}_{\mathbf{P}}=\{U\mapsto \mathscr{G}_{U}\}$ ($U\in \Ob(\mathbf{P})$) l'objet de $\widehat{E}_{\mathbf{P}}^{\sharp}$ obtenu en restreignant $\mathscr{G}$ à $E_{\mathbf{P}}^{\sharp}$. On en déduit un isomorphisme canonique de $\widetilde{E}^{\sharp}$
\begin{equation}
	\mathscr{G}_{\mathbf{P}}^{a} \xrightarrow{\sim} \mathscr{G}^{a}.
\end{equation}
\end{nothing}

\begin{nothing}	\label{ring B sharp}
	Dans la suite de cette section, on suppose, de plus, qu'il existe un trait $(S',\eta')$ fini sur $(S,\eta)$ et un $\eta'$-revêtement fini, galoisien et à fibres géométriquement réduites $\varphi:X'\to X_{S'}$ \eqref{def of revêtement} tel que, en fixant un $\eta$-morphisme $\overline{\eta}\to \eta'$, $X'_{\overline{\eta}}\simeq C'$ et que $\varphi_{\overline{\eta}}:X_{\overline{\eta}}'\to X_{\overline{\eta}}$ s'identifie au morphisme $\psi$. On pose $\overrightharpoon{X}'=X'\times_{S'}\overline{S}$ et on note
\begin{equation}
	\hbar: \overline{X}\to X\quad \textnormal{et} \quad \overline{\varphi}:\overrightharpoon{X}'\to \overline{X}
		\label{morphisme can sharp}
\end{equation}
les morphismes canoniques.

	Pour tout $(V\to U)\in \Ob(E^{\sharp})$, on pose $U'=U\times_{X}X'$, $\overrightharpoon{U}'=U'\times_{X'}\overrightharpoon{X}'$ et $\overline{U}^{V}$ (resp. $\overrightharpoon{U}'^{V}$) la clôture intégrale de $\overline{U}$ (resp. $\overrightharpoon{U}'$) dans $V$. Comme $\overrightharpoon{U}'$ est fini sur $\overline{U}$, on a un isomorphisme canonique $\overrightharpoon{U}'^{V}\xrightarrow{\sim} \overline{U}^{V}$. Comme $\Phi^{*}$ est le foncteur de restriction, on a $\Phi^{*}(\overline{\mathscr{B}})(V\to U)\simeq \overline{\mathscr{B}}(V\to U)$. On en déduit par \eqref{def of bar B} un isomorphisme fonctoriel en $(V\to U)$
	\begin{equation}
		\overline{\mathscr{B}}^{\sharp}(V\to U)\xrightarrow{\sim} \Gamma(\overrightharpoon{U}'^{V},\mathscr{O}_{\overrightharpoon{U}'^{V}}).
		\label{B star def}
	\end{equation}
	Pour tout $U\in \Ob(\Et_{/X})$, on pose $\overline{\mathscr{B}}^{\sharp}_{U}=\varphi_{U,\fet}^{*}(\overline{\mathscr{B}}_{U})$ \eqref{def of bar B UU} et $\overline{\mathscr{B}}_{U,n}^{\sharp}=\overline{\mathscr{B}}_{U}^{\sharp}/p^{n}\overline{\mathscr{B}}_{U}^{\sharp}$ qui est isomorphe à $\varphi_{U,\fet}^{*}(\overline{\mathscr{B}}_{U,n})$ \eqref{bar B U n}. En vertu de \eqref{B star def}, on a un homomorphisme d'anneaux de $U'_{\overline{\eta},\fet}$
	\begin{equation}
		\mathscr{O}_{\overrightharpoon{X}'}(\overrightharpoon{U}') \to \overline{\mathscr{B}}^{\sharp}_{U},
		\label{OX prime to Bstar}
	\end{equation}
	où la source désigne le faisceau constant de valeur $\mathscr{O}_{\overrightharpoon{X}'}(\overrightharpoon{U}')$. Comme l'action de $G$ sur $C'$ s'étend en une action de $G$ sur $\overrightharpoon{X}'$, l'anneau $\mathscr{O}_{\overrightharpoon{X}'}(\overrightharpoon{U}')$ est muni d'une action de $G$ telle que l'homomorphisme \eqref{OX prime to Bstar} soit $G$-équivariant. D'après \ref{ring B bar} et \eqref{resctriction prefaisceau}, l'anneau $\overline{\mathscr{B}}_{n}^{\sharp}$ est isomorphe au faisceau sur $E^{\sharp}$ associé au préfaisceau
\begin{equation}
	\{U\to \overline{\mathscr{B}}_{U,n}^{\sharp}\} \qquad U\in \Ob(\Et_{/X}).
	\label{B star n description}
\end{equation}
\end{nothing}
$\hspace*{-1.2em}\bf{\arabic{section}.\stepcounter{theorem}\arabic{theorem}.}$
Comme $\widetilde{E}^{\sharp}$ est isomorphe au topos de Faltings associé au morphisme canonique $C'\to X$, on a un morphisme canonique (\cite{AGT} (VI.10.6.3))
\begin{equation}
	\beta^{\sharp}: \widetilde{E}^{\sharp}\to C_{\fet}'
	\label{beta star}
\end{equation}
qui rend commutatif le diagramme (cf. \cite{AGT} (VI.10.12.6))
\begin{displaymath}
	\xymatrix{
		\widetilde{E}^{\sharp} \ar[r]^{\Phi} \ar[d]_{\beta^{\sharp}} & \widetilde{E}\ar[d]^{\beta}\\
		C_{\fet}' \ar[r]^{\varphi_{\overline{\eta}}}& C_{\fet}.
	}
\end{displaymath}
Soit $n$ un entier $\ge 1$. On désigne par $\beta_{n}^{\sharp}: (\widetilde{E}_{s}^{\sharp}, \overline{\mathscr{B}}_{n}^{\sharp})\to (C'_{\fet},\oo_{n})$ le morphisme de topos annelés défini par le composé $\beta^{\sharp}\circ \delta^{\sharp}$ \eqref{delta sharp} et le morphisme canonique $\oo_{n}\to \beta_{*}^{\sharp}(\overline{\mathscr{B}}_{n}^{\sharp})$. Par des arguments similaires à \ref{fonctoriel Faltings topos}, on en déduit un diagramme commutatif
\begin{equation}
	\xymatrix{
		(\widetilde{E}_{s}^{\sharp}, \overline{\mathscr{B}}_n^{\sharp}) \ar[r]^{\Phi_n} \ar[d]_{\beta_n^{\sharp}} & (\widetilde{E}_s, \overline{\mathscr{B}}_n)\ar[d]^{\beta_n}\\
		(C_{\fet}',\oo_n) \ar[r]^{\varphi_{\overline{\eta}}}& (C_{\fet},\oo_n),
	}
	\label{diagramme com 2 Faltings et fet}
\end{equation}
à isomorphisme près, qui commute aux actions de $G$ sur $(\widetilde{E}_{s}^{\sharp}, \overline{\mathscr{B}}_n^{\sharp})$ et $(C_{\fet}',\oo_n)$.

\begin{theorem}
	Soient $X$ une $S$-courbe semi-stable, $\overline{x}$ un point géométrique de $X_{\overline{\eta}}$ et $\mathcal{F}$ un fibré vectoriel sur $\check{\overline{X}}=X\times_{S}\check{\overline{S}}$ \eqref{notations 11}. Supposons qu'il existe un entier $n\ge 1$, un trait $(S',\eta')$ fini sur $(S,\eta)$ et un $\eta'$-revêtement fini, galoisien et à fibres géométriquement réduites $\varphi:X'\to X_{S'}$ tels que, fixant un $\eta$-morphisme $\overline{\eta}\to \eta'$ et reprenant les notations de \eqref{morphisme can sharp}, $\overline{\varphi}_{n}^{*}(\mathcal{F}_{n})$ soit libre de type fini. Alors, on a un isomorphisme $\overline{\mathscr{B}}_{n}$-linéaire
	\begin{equation}
		\gamma_{n}:\sigma_{n}^{*}(\mathcal{F}_{n})\xrightarrow{\sim} \beta_{n}^{*} (\mathbb{V}_{n}(\mathcal{F})),
		\label{Weil Tate admissible mod pn}
	\end{equation}
	où $\mathbb{V}_{n}(\mathcal{F})$ est la $\oo_{n}$-représentation de $\pi_{1}(X_{\overline{\eta}},\overline{x})$ associée à $\mathcal{F}$ \eqref{l'action du groupe fondamentale}, $\sigma_{n}$ est le morphisme de topos annelés \eqref{morphisme de topos sigma n} et $\beta_{n}^{*}$ est le foncteur \eqref{betan star rep}.
	\label{Mod Weil-Tate mod pn}
\end{theorem}

\textit{Preuve}. Pour tout objet $U$ de $\Et_{/X}$, on pose $U'=U\times_{X}X'$, $U_{s}'=U'\times_{S'}s$ et $U_{\overline{\eta}}'=U'\times_{S'}\overline{\eta}$. D'après \eqref{resctriction prefaisceau} et (\cite{AGT} VI.5.34(ii) et VI.8.9), le $\overline{\mathscr{B}}_n^{\sharp}$-module $\Phi_{n}^*(\sigma_n^*(\mathcal{F}_n))$ est isomorphe au faisceau associé au préfaisceau:
\begin{equation}
	\{U\to \mathcal{F}_n(U_{s})\otimes_{\mathscr{O}_{\overline{X}_n}(U_{s})}\overline{\mathscr{B}}_{U,n}^{\sharp}\} \qquad U\in \Ob(\Et_{/X}).
	\label{def of sigma Fn}
\end{equation}
où $\mathcal{F}_n(U_{s})$ et $\mathscr{O}_{\overline{X}_n}(U_{s})$ sont considérés comme des faisceaux constants de $U_{\overline{\eta},\fet}'$. D'après \eqref{OX prime to Bstar}, l'anneau $\overline{\mathscr{B}}_{U,n}^{\sharp}$ de $U_{\overline{\eta},\fet}'$ est une $\mathscr{O}_{\overrightharpoon{X}'_{n}}(U_{s}')$-algèbre. Pour tout objet $U$ de $\mathbf{P}$ \eqref{categorie P affine}, on en déduit un isomorphisme $G$-équivariant de faisceaux de $U'_{\overline{\eta},\fet}$
\begin{equation}
	\mathcal{F}_n(U_{s})\otimes_{\mathscr{O}_{\overline{X}_n}(U_{s})}\overline{\mathscr{B}}_{U,n}^{\sharp}\xrightarrow{\sim} \overline{\varphi}^*_n(\mathcal{F}_n)(U'_{s})\otimes_{\mathscr{O}_{\overrightharpoon{X}'_n}(U'_{s})}\overline{\mathscr{B}}_{U,n}^{\sharp},
	\label{mod pullback by sigma}
\end{equation}
où les actions de $G$ sur $\mathcal{F}_n(U_{s})$ et $\mathscr{O}_{\overline{X}_n}(U_{s})$ sont triviales et $\overline{\varphi}_{n}^{*}(\mathcal{F}_n)(U'_{s})$ et $\mathscr{O}_{\overrightharpoon{X}'_n}(U'_{s})$ sont considérés comme des faisceaux constants de $U_{\overline{\eta},\fet}'$ munis des actions de $G$ induites par $\overline{\varphi}_{n}^{*}$. On note $M_{U}$ le but de \eqref{mod pullback by sigma}.

On pose $\mathbb{L}_{n}(\mathcal{F})=\mu_{\overline{x}}(\mathbb{V}_{n}(\mathcal{F}))$ le $\oo_{n}$-module localement libre de type fini de $C_{\fet}$ associé à $\mathbb{V}_{n}(\mathcal{F})$ \eqref{Rep on to Mod on}. D'après (\cite{AGT} VI.8.9 et IV.10.9), le $\overline{\mathscr{B}}_n^{\sharp}$-module $\Phi_{n}^*(\beta_{n}^{*}(\mathbb{L}_{n}(\mathcal{F})))\simeq (\beta_n^{\sharp})^{*}(\varphi_{\overline{\eta}}^{*}(\mathbb{L}_{n}(\mathcal{F})))$ est isomorphe au faisceau associé au préfaisceau
\begin{equation}
	\{U\to h_{U,\fet}'^{*}(\varphi_{\overline{\eta},\fet}^{*}(\mathbb{L}_{n}(\mathcal{F})))\otimes_{\oo_n}\overline{\mathscr{B}}^{\sharp}_{U,n}\} \qquad U\in \Ob(\Et_{/X}),
	\label{mod pullback by beta}
\end{equation}
où $h_{U}':U_{\overline{\eta}}'\to C'$ est le morphisme canonique. On pose $N_{U}=h_{U,\fet}'^{*}(\varphi_{\overline{\eta},\fet}^{*}(\mathbb{L}_{n}(\mathcal{F})))\otimes_{\oo_n}\overline{\mathscr{B}}^{\sharp}_{U,n}$ qui est muni d'une action de $G$ fonctorielle en $U$ (cf. \ref{les actions de G sur Esharp prefaisceau}).

Comme $\overline{\varphi}_{n}^{*}(\mathcal{F}_{n})$ est un $\mathscr{O}_{\overrightharpoon{X}'_{n}}$-module libre de type fini, d'après \eqref{def of rho n F gamma} et \ref{def of DW n independant}, on a un isomorphisme $G$-équivariant de $\oo_{n}$-modules de $C_{\fet}'$
\begin{equation}
	\overline{\varphi}_{n}^{*}(\mathcal{F}_{n})(X'_{s})\simeq \varphi_{\overline{\eta},\fet}^{*}(\mathbb{L}_{n}(\mathcal{F})),
	\label{DW condition Faltings}
\end{equation}
où $\overline{\varphi}_{n}^{*}(\mathcal{F}_{n})(X'_{s})$ est considéré comme le faisceau constant de $C_{\fet}'$ muni de l'action de $G$ induite par $\overline{\varphi}_{n}^{*}$. On a un isomorphisme $G$-équivariant et fonctoriel en $U$ de $\mathscr{O}_{\overrightharpoon{X}'_n}(U'_{s})$-modules
\begin{equation}
	\overline{\varphi}^*_n(\mathcal{F}_n)(U'_{s}) \xrightarrow{\sim} \overline{\varphi}_n^*(\mathcal{F}_n)(X'_{s})\otimes_{\oo_{n}}\mathscr{O}_{\overrightharpoon{X}'_n}(U'_{s}).
\end{equation}
On en déduit par \eqref{DW condition Faltings} un isomorphisme $G$-équivariant et fonctoriel en $U$ de $\mathscr{O}_{\overrightharpoon{X}'_n}(U'_{s})$-modules de $U_{\overline{\eta},\fet}'$
\begin{equation}
	\overline{\varphi}^*_n(\mathcal{F}_n)(U'_{s})\simeq h_{U,\fet}'^{*}(\varphi_{\overline{\eta},\fet}^{*}(\mathbb{L}_{n}(\mathcal{F})))\otimes_{\oo_n}\mathscr{O}_{\overrightharpoon{X}'_n}(U'_{s}).
	\label{DW via 2 pullback}
\end{equation}
En vertu de \eqref{mod pullback by sigma} et \eqref{mod pullback by beta}, on en déduit, pour tout objet $U$ de $\mathbf{P}$, un isomorphisme $G$-équivariant et fonctoriel en $U$ de $\overline{\mathscr{B}}_{U,n}^{\sharp}$-modules de $U_{\overline{\eta},\fet}'$
\begin{equation}
	M_{U}\xrightarrow{\sim}N_{U}.
	\label{MU iso NU}
\end{equation}
D'après \ref{les actions de G sur Esharp prefaisceau} et \ref{categorie P affine}, on en déduit un isomorphisme $G$-équivariant de $\overline{\mathscr{B}}_n^{\sharp}$-modules
\begin{equation}
	\Phi_n^{*} (\sigma_n^{*}(\mathcal{F}_n))\xrightarrow{\sim} \Phi_n^*(\beta_{n}^{*}(\mathbb{L}_{n}(\mathcal{F}))).
	\label{donne de descente Galois pour Weil-Tate}
\end{equation}
Le théorème résulte alors de \ref{lemme de descente Bn}.

\begin{rem}
	Pour tout entier $1\le m \le n$, on a un isomorphisme $\overline{\mathscr{B}}_{m}$-linéaire $\gamma_{m}: \sigma_{m}^{*}(\mathcal{F}_{m})\xrightarrow{\sim} \beta_{m}^{*} (\mathbb{V}_{m}(\mathcal{F}))$. En vertu de la preuve, les isomorphismes $\gamma_{n}$ et $\gamma_{m}$ sont compatibles.
	\label{rem 912}
\end{rem}
\section{Fibrés vectoriels et représentations de Weil-Tate} \label{FV Rep WT}
\begin{nothing} \label{notations general section WT}
	Soient $X$ une $S$-courbe semi-stable, $\overline{x}$ un point géométrique de $X_{\overline{\eta}}$. On reprend les notations de \S~\ref{Topos de Faltings} pour $X$. Rappelons qu'on a des morphismes de topos annelés (cf. \ref{morphisme sigma T} et \ref{topos annele morphisme beta})
\begin{eqnarray}
	(X_{s,\et},\mathscr{O}_{\overline{X}_{n}})\xleftarrow{\sigma_{n}} &(\widetilde{E}_{s},\overline{\mathscr{B}}_{n})&\xrightarrow{\beta_{n}}(X_{\overline{\eta},\fet},\oo_{n}), \qquad \forall n \ge 1\\
	(X_{s,\et}^{\mathbb{N}^{\circ}},\mathscr{O}_{\breve{\overline{X}}})\xleftarrow{\breve{\sigma}} &(\widetilde{E}_{s}^{\mathbb{N}^{\circ}},\breve{\overline{\mathscr{B}}})&\xrightarrow{\breve{\beta}}(X_{\overline{\eta},\fet}^{\mathbb{N}^{\circ}},\breve{\oo}).
\end{eqnarray}
On reprend les notations de \ref{notation beta n rep} pour les foncteurs $\beta_{n}^{*}$, $\breve{\beta}^{*}$ et $\breve{\beta}^{*}_{\mathbb{Q}}$.

On désigne par $\XX$ le schéma formel complété $p$-adique de $\overline{X}$. Rappelons que $\XX$ est de présentation finie sur $\mathscr{S}=\Spf(\oo)$ et qu'on a un morphisme de topos annelés \eqref{morphisme de topos T}
\begin{equation}
	\rT: (\widetilde{E}_{s}^{\mathbb{N}^{\circ}},\breve{\overline{\mathscr{B}}})\to (X_{s,\zar},\mathscr{O}_{\XX}).
	\label{moprhisme de topos anneles T}
\end{equation}
Pour tout $\mathscr{O}_{\XX}$-module cohérent $\mathscr{F}$ de $X_{s,\zar}$ et tout entier $n\ge 1$, on note $\mathscr{F}_{n}=\mathscr{F}/p^{n}\mathscr{F}$ que l'on considère comme faisceau de $X_{s,\zar}$ ou de $X_{s,\et}$. D'après (\cite{AGT} (III.11.1.12)), on a un isomorphisme canonique
\begin{equation}
	\breve{\sigma}^{*}( (\mathscr{F}_{n})_{n\ge 1})\simeq \rT^{*}(\mathscr{F}).
	\label{compare sigma T}
\end{equation}
Le foncteur $\rT_{*}$ induit un foncteur additif et exact à gauche que l'on note encore
\begin{equation}
	\rT_{*}:\Mod_{\mathbb{Q}}(\breve{\overline{\mathscr{B}}})\to \Mod(\mathscr{O}_{\XX}[\frac{1}{p}]).
\end{equation}
D'après \eqref{Mod OXX Q to Mod OXXp}, le foncteur $\rT^{*}$ induit un foncteur additif que l'on note encore
\begin{equation}
	\rT^{*}: \Mod^{\coh}(\mathscr{O}_{\XX}[\frac{1}{p}])\to \Mod^{\atf}_{\mathbb{Q}}(\breve{\overline{\mathscr{B}}}).
	\label{T star rationel}
\end{equation}
\end{nothing}

\begin{definition} \label{B associe}
	\textnormal{(i)} Soient $n$ un entier $\ge 1$, $\mathcal{F}_{n}$ un fibré vectoriel sur $\overline{X}_{n}$ et $\mathcal{V}_{n}$ une $\oo_{n}$-représentation de $\pi_{1}(X_{\overline{\eta}},\overline{x})$ \eqref{oo oon reprentations general}. On dit que $\mathcal{F}_{n}$ et $\mathcal{V}_{n}$ sont $\overline{\mathscr{B}}_{n}$\textit{-associés} s'il existe un isomorphisme $\overline{\mathscr{B}}_{n}$-linéaire
	\begin{equation}
		\sigma_n^*(\mathcal{F}_n)\xrightarrow{\sim} \beta_{n}^{*}(\mathcal{V}_{n}).
		\label{module pn admissible}
	\end{equation}

	\textnormal{(ii)} Soient $\mathscr{F}$ un $\mathscr{O}_{\XX}[\frac{1}{p}]$-module localement projectif de type fini \eqref{LPtf} et $V$ une $\mathfrak{C}$-représentation continue de $\pi_{1}(X_{\overline{\eta}},\overline{x})$ \eqref{oo oon reprentations}. On dit que $\mathscr{F}$ et $V$ sont \textit{$\breve{\overline{\mathscr{B}}}_{\mathbb{Q}}$-associés} s'il existe un isomorphisme $\breve{\overline{\mathscr{B}}}_{\mathbb{Q}}$-linéaire
	\begin{equation}
		\rT^{*}(\mathscr{F})\xrightarrow{\sim} \breve{\beta}_{\mathbb{Q}}^{*}(V).
		\label{B admissible}
	\end{equation}
\end{definition}

\begin{definition}
	\textnormal{(i)} On dit qu'un $\mathscr{O}_{\XX}[\frac{1}{p}]$-module localement projectif de type fini $\mathscr{F}$ est \textit{de Weil-Tate} s'il est $\breve{\overline{\mathscr{B}}}_{\mathbb{Q}}$-associé à une $\mathfrak{C}$-représentation continue de $\pi_{1}(X_{\overline{\eta}},\overline{x})$.

	\textnormal{(ii)} On dit qu'une $\mathfrak{C}$-représentation continue $V$ de $\pi_{1}(X_{\overline{\eta}},\overline{x})$ est \textit{de Weil-Tate relativement à $X$} si elle est $\breve{\overline{\mathscr{B}}}_{\mathbb{Q}}$-associée à un $\mathscr{O}_{\XX}[\frac{1}{p}]$-module localement projectif de type fini.
	\label{Weil-Tate X}
\end{definition}

On désigne par $\Mod^{\WT}(\mathscr{O}_{\XX}[\frac{1}{p}])$ la sous-catégorie pleine de $\Mod^{\coh}(\mathscr{O}_{\XX}[\frac{1}{p}])$ formée des $\mathscr{O}_{\XX}[\frac{1}{p}]$-modules de Weil-Tate et par $\Rep_{\mathfrak{C}}^{\WT/X}(\pi_{1}(X_{\overline{\eta}},\overline{x}))$ la sous-catégorie pleine de $\Rep_{\mathfrak{C}}^{\cont}(\pi_{1}(X_{\overline{\eta}},\overline{x}))$ formée des $\mathfrak{C}$-représentations de Weil-Tate relativement à $X$.

\begin{rem} \label{FV de WT sur X}
	On dit qu'un fibré vectoriel $\mathcal{F}$ sur $\check{\overline{X}}=X\times_{S}\check{\overline{S}}$ \eqref{notations 11} est \textit{de Weil-Tate} si le $\mathscr{O}_{\XX}[\frac{1}{p}]$-module $\widehat{\mathcal{F}}[\frac{1}{p}]$ est de Weil-Tate. D'après \eqref{compare sigma T}, il revient au même de dire qu'il existe une $\mathfrak{C}$-représentation continue $V$ de $\pi_{1}(X_{\overline{\eta}},\overline{x})$ et un isomorphisme $\breve{\overline{\mathscr{B}}}_{\mathbb{Q}}$-linéaire
	\begin{equation}
		(\breve{\sigma}^{*}( (\mathcal{F}_{n})_{n\ge 1}))_{\mathbb{Q}}\xrightarrow{\sim} \breve{\beta}_{\mathbb{Q}}^{*}(V).
		\label{sigma beta B associe}
	\end{equation}
\end{rem}

\begin{prop}
	Soient $X$ une $S$-courbe semi-stable, $\overline{x}$ un point géométrique de $X_{\overline{\eta}}$, $n$ un entier $\ge 1$ et $\mathcal{F}$ un fibré vectoriel de Deninger-Werner sur $\check{\overline{X}}$ \eqref{def of categorie DW}. Alors, il existe un trait $S'$ fini sur $S$, un $S'$-modèle semi-stable et régulier $X'$ de $X_{\overline{\eta}}$ \eqref{Notations C X LPft Vect} et un $\eta'$-revêtement $\varphi:X'\to X_{S'}$ tels que $\overline{\varphi}^{*}_{n}(\mathcal{F}_{n})$ et $\mathbb{V}_{\check{\overline{X}},n}(\mathcal{F})$ soient $\overline{\mathscr{B}}_{n}'$-associés, où $(\widetilde{E}_{s}',\overline{\mathscr{B}}_{n}')$ désigne le topos annelé de Faltings associé au schéma $X'$ au-dessus de $S'$ et $\mathbb{V}_{\check{\overline{X}},n}(\mathcal{F})$ est la $\oo_{n}$-représentation de $\pi_{1}(X_{\overline{\eta}},\overline{x})$ associée à $\mathcal{F}$ par \eqref{rho VVn}.
	\label{DW module pn admissible}
\end{prop}
\textit{Preuve}. Quitte à remplacer $S$ par une extension finie, d'après \ref{domine par bon revêtement}(iii), il existe un $S$-modèle semi-stable $X'$ de $X_{\overline{\eta}}$ et un $\eta$-revêtement $\varphi:X'\to X$ tels que $\overline{\varphi}_{n}^{*}(\mathcal{F}_{n})$ soit trivialisé par un $\eta$-revêtement semi-stable, fini et galoisien $Y'\to X'$. D'après \ref{Mod Weil-Tate mod pn}, $\overline{\varphi}^{*}_{n}(\mathcal{F}_{n})$ et $\mathbb{V}_{\check{\overline{X}}',n}(\overline{\varphi}^{*}(\mathcal{F}))$ sont $\overline{\mathscr{B}}_{n}'$-associés. D'après \eqref{semistable rcm coro}, \eqref{fonctorialite sigma n} et \eqref{fonctorialite beta n}, on peut supposer que $X$ est régulier. Alors, la proposition s'ensuit compte tenu de \ref{propriete DW foncteur}(ii).

\begin{prop}
	Soient $X$ une $S$-courbe propre et lisse, $\overline{x}$ un point géométrique de $X_{\overline{\eta}}$, $\mathcal{L}$ un fibré en droites de classe appartenant à $\Pic^{0}_{\check{\overline{X}}/\check{\overline{S}}}(\check{\overline{S}})$. Pour tout entier $n\ge 1$, on a associé à $\mathcal{L}$ une $\oo_{n}$-représentation $\mathbb{V}_{n}(\mathcal{L})$ de $\pi_{1}(X_{\overline{\eta}},\overline{x})$ \eqref{l'action du groupe fondamentale} et une $\oo$-représentation continue $\mathbb{V}(\mathcal{L})$ de $\pi_{1}(X_{\overline{\eta}},\overline{x})$ \eqref{rho entier}. Alors, $\mathcal{L}_{n}$ et $\mathbb{V}_{n}(\mathcal{L})$ sont $\overline{\mathscr{B}}_{n}$-associés; et $\widehat{\mathcal{L}}[\frac{1}{p}]$ et $\mathbb{V}(\mathcal{L})[\frac{1}{p}]$ \eqref{oo oon reprentations} sont $\breve{\overline{\mathscr{B}}}_{\mathbb{Q}}$-associés.
	\label{bonne reduction fibre en droite}	
\end{prop}
\textit{Preuve}. Soit $n$ un entier $\ge 1$. D'après \ref{DW deg 1} et \ref{Mod Weil-Tate mod pn}, on a un isomorphisme
\begin{equation}
	\gamma_{n}:\sigma_{n}^{*}(\mathcal{L}_{n})\xrightarrow{\sim} \beta_{n}^{*}(\mathbb{V}_{n}(\mathcal{L})),
	\label{Bn associe line bundle}
\end{equation}
d'où le premier énoncé. En vertu de \ref{rem 912}, les isomorphismes $(\gamma_{n})_{n\ge 1}$ \eqref{Bn associe line bundle} sont compatibles. La proposition s'ensuit compte tenu \eqref{compare sigma T}.

\begin{prop}
	Soient $X$ une $S$-courbe semi-stable et régulière et $\overline{x}$ un point géométrique de $X_{\overline{\eta}}$.
	
	\textnormal{(i)} La restriction du foncteur $\rT^{*}$ \eqref{T star rationel} à $\Mod^{\WT}(\mathscr{O}_{\XX}[\frac{1}{p}])$ se factorise à travers la sous-catégorie pleine $\Mod^{\lt}_{\mathbb{Q}}(\breve{\overline{\mathscr{B}}})$ \textnormal{\eqref{def pltf}} de $\Mod^{\atf}_{\mathbb{Q}}(\breve{\overline{\mathscr{B}}})$.
	
	\textnormal{(ii)} Pour tout objet $\mathscr{F}$ de $\Mod^{\WT}(\mathscr{O}_{\XX}[\frac{1}{p}])$, on a un isomorphisme $\breve{\overline{\mathscr{B}}}_{\mathbb{Q}}$-linéaire, canonique et fonctoriel
	\begin{equation}
		\rT^{*}(\mathscr{F}) \xrightarrow{\sim} \breve{\beta}_{\mathbb{Q}}^{*}(\mathscr{W}_{\mathbb{Q}}(\rT^{*}(\mathscr{F}))),
		\label{isomorphisme VX WT}
	\end{equation}
	où $\mathscr{W}_{\mathbb{Q}}: \Mod^{\lt}_{\mathbb{Q}}(\breve{\overline{\mathscr{B}}})\to \Rep^{\cont}_{\mathfrak{C}}(\pi_{1}(X_{\overline{\eta}},\overline{x}))$ est le foncteur \eqref{varrho Q}.
	\label{isomorphisme VX WT lemma}
\end{prop}
\textit{Preuve}. Soit $\mathscr{F}$ un objet de $\Mod^{\WT}(\mathscr{O}_{\XX}[\frac{1}{p}])$. Choisissons un objet $V$ de $\Rep_{\mathfrak{C}}^{\cont}(\pi_{1}(X_{\overline{\eta}},\overline{x}))$ et un isomorphisme $\breve{\overline{\mathscr{B}}}_{\mathbb{Q}}$-linéaire
\begin{equation}
	\tau:\rT^{*}(\mathscr{F}) \simeq \breve{\beta}^{*}_{\mathbb{Q}}(V).
	\label{isomorphisme BQ associe fibre}
\end{equation}
L'assertion (i) résulte alors de \ref{betan lt}. D'après \ref{rep generalise rep vrai}(ii), on a un isomorphisme $\breve{\overline{\mathscr{B}}}_{\mathbb{Q}}$-linéaire canonique et fonctoriel en $V$
\begin{equation}
	\breve{\beta}^{*}_{\mathbb{Q}}(V) \xrightarrow{\sim} \breve{\beta}_{\mathbb{Q}}^{*}(\mathscr{W}_{\mathbb{Q}}(\breve{\beta}^{*}_{\mathbb{Q}}(V))).
	\label{adjonction V beta}
\end{equation}
On en déduit compte tenu de \eqref{isomorphisme BQ associe fibre} et \eqref{adjonction V beta} un isomorphisme $\breve{\overline{\mathscr{B}}}_{\mathbb{Q}}$-linéaire	
\begin{displaymath}
		\rT^{*}(\mathscr{F}) \xrightarrow{\sim} \breve{\beta}_{\mathbb{Q}}^{*}(\mathscr{W}_{\mathbb{Q}}(\rT^{*}(\mathscr{F}))).
\end{displaymath}
Montrons que cet isomorphisme ne dépend pas des choix de $V$ et de l'isomorphisme $\tau$ et qu'il est fonctoriel en $\mathscr{F}$. Il suffit de démontrer que, pour tous objets $V$ et $V'$ de $\Rep_{\mathfrak{C}}^{\cont}(\pi_{1}(X_{\overline{\eta}},\overline{x}))$ et tout morphisme $\breve{\overline{\mathscr{B}}}_{\mathbb{Q}}$-linéaire $\theta: \breve{\beta}_{\mathbb{Q}}^{*}(V)\to \breve{\beta}_{\mathbb{Q}}^{*}(V')$, le diagramme
\begin{equation}
	\xymatrix{
		\breve{\beta}^{*}_{\mathbb{Q}}(V) \ar[r]^-{\sim} \ar[d]_{\theta} & \breve{\beta}_{\mathbb{Q}}^{*}(\mathscr{W}_{\mathbb{Q}}(\breve{\beta}^{*}_{\mathbb{Q}}(V))) \ar[d]^{\breve{\beta}_{\mathbb{Q}}^{*}(\mathscr{W}_{\mathbb{Q}}(\theta))}\\
		\breve{\beta}^{*}_{\mathbb{Q}}(V') \ar[r]^-{\sim} & \breve{\beta}_{\mathbb{Q}}^{*}(\mathscr{W}_{\mathbb{Q}}(\breve{\beta}^{*}_{\mathbb{Q}}(V')))
	}
	\label{diagramme induit par theta V}
\end{equation}
est commutatif. D'après la pleine fidélité du foncteur $\breve{\beta}_{\mathbb{Q}}^{*}$ \eqref{beta pullback pleinement fidele}, $\theta$ est l'image d'un morphisme $\vartheta:V\to V'$ de $\Rep^{\cont}_{\mathfrak{C}}(\pi_{1}(X_{\overline{\eta}},\overline{x}))$ par $\breve{\beta}_{\mathbb{Q}}^{*}$. La commutativité de \eqref{diagramme induit par theta V} résulte alors de la fonctorialité de \eqref{adjonction V beta}.

\begin{coro}
	Sous les hypothèses de \ref{isomorphisme VX WT lemma}, pour tout objet $\mathscr{F}$ de $\Mod^{\WT}(\mathscr{O}_{\XX}[\frac{1}{p}])$, la $\mathfrak{C}$-représentation continue $\mathscr{W}_{\mathbb{Q}}(\rT^{*}(\mathscr{F}))$ est de Weil-Tate relativement à $X$. On désigne par $\mathscr{V}_{\XX}$ le foncteur
	\begin{equation}
		\mathscr{V}_{\XX}: \Mod^{\WT}(\mathscr{O}_{\XX}[\frac{1}{p}]) \to \Rep_{\mathfrak{C}}^{\WT/X}(\pi_{1}(X_{\overline{\eta}},\overline{x})),\qquad \mathscr{F}\mapsto \mathscr{W}_{\mathbb{Q}}(\rT^{*}(\mathscr{F})).
		\label{definition foncteur VXX}
	\end{equation}
	Alors, $\mathscr{V}_{\XX}(\mathscr{F})$ et $\mathscr{F}$ sont $\breve{\overline{\mathscr{B}}}_{\mathbb{Q}}$-associés.
\end{coro}

\begin{prop}
	Soient $X$ une $S$-courbe semi-stable et régulière. On note $\LPtf(\mathscr{O}_{\XX}[\frac{1}{p}])$ la sous-catégorie pleine de $\Mod^{\coh}(\mathscr{O}_{\XX}[\frac{1}{p}])$ formée des $\mathscr{O}_{\XX}[\frac{1}{p}]$-modules localement projectifs de type fini.
	
	\textnormal{(i)} Pour tout objet $\mathscr{F}$ de $\LPtf(\mathscr{O}_{\XX}[\frac{1}{p}])$, le morphisme d'adjonction \textnormal{(cf. \cite{AGT} III.12.1)}
	\begin{equation}
		\mathscr{F}\to \rT_{*}(\rT^{*}(\mathscr{F}))
		\label{formule projectif T}
	\end{equation}
	est un isomorphisme.

	\textnormal{(ii)} La restriction du foncteur $\rT^{*}$ \eqref{T star rationel} à $\LPtf(\mathscr{O}_{\XX}[\frac{1}{p}])$ est pleinement fidèle.
\label{general property of T pushback}
\end{prop}

\textit{Preuve}. (i) D'après (\cite{AGT} III.12.4(ii)), on a un isomorphisme déduit du morphisme d'adjonction
\begin{equation}
	\mathscr{F}\otimes_{\mathscr{O}_{\XX}[\frac{1}{p}]}\rT_{*}(\breve{\overline{\mathscr{B}}}_{\mathbb{Q}})\xrightarrow{\sim} \rT_{*}(\rT^{*}(\mathscr{F})).
\end{equation}
L'isomorphisme \eqref{formule projectif T} résulte alors de l'isomorphisme (\cite{AGT} III.11.8)
\begin{equation}
	\mathscr{O}_{\XX}[\frac{1}{p}]\xrightarrow{\sim} \rT_{*}(\breve{\overline{\mathscr{B}}}_{\mathbb{Q}}).
\end{equation}

(ii) Soient $\mathscr{G}_{1}$ et $\mathscr{G}_{2}$ deux $\mathscr{O}_{\XX}$-modules cohérents de $X_{s,\zar}$ tels que $\mathscr{G}_{1}[\frac{1}{p}]$ et $\mathscr{G}_{2}[\frac{1}{p}]$ soient localement projectifs de type fini. Notons $\mathscr{G}=\FHom_{\mathscr{O}_{\XX}}(\mathscr{G}_1,\mathscr{G}_2)$. Alors, $\mathscr{G}[\frac{1}{p}]$ est aussi localement projectif de type fini \eqref{LPtf}. Comme $\mathscr{G}_{1}$ est de présentation finie, on a un isomorphisme canonique
\begin{equation}
	\rT^{*}(\mathscr{G})\xrightarrow{\sim}\FHom_{\breve{\overline{\mathscr{B}}}}(\rT^{*}(\mathscr{G}_1),\rT^{*}(\mathscr{G}_2)).
	\label{iso of hom T}
\end{equation}
Considérons le morphisme composé
\begin{equation}
	f: \Gamma(X_{s,\zar},\mathscr{G})\to  \Gamma(X_{s,\zar}, \rT_{*}(\rT^{*}(\mathscr{G}))) \xrightarrow{\sim} \Gamma(\widetilde{E}^{\mathbb{N}^{\circ}}_s,\rT^{*}(\mathscr{G})).
	\label{iso global section pullback T}
\end{equation}
D'après (i) et \ref{isomorphisme section inverse p}(i), le morphisme $f\otimes_{\mathbb{Z}}\mathbb{Q}$ est un isomorphisme. On en déduit la pleine fidélité de $\rT^{*}$.

\begin{prop} \label{formule projectif T lemma}
	Soient $X$ une $S$-courbe semi-stable et régulière et $\overline{x}$ un point géométrique de $X_{\overline{\eta}}$. Pour toute $\mathfrak{C}$-représentation $V$ de $\pi_{1}(X_{\overline{\eta}},\overline{x})$ de Weil-Tate relativement à $X$, le $\mathscr{O}_{\XX}[\frac{1}{p}]$-module $\rT_{*}(\breve{\beta}_{\mathbb{Q}}^{*}(V))$ est de Weil-Tate. De plus, on a un isomorphisme $\breve{\overline{\mathscr{B}}}_{\mathbb{Q}}$-linéaire, canonique et fonctoriel
	\begin{equation}
		\breve{\beta}_{\mathbb{Q}}^{*}(V) \xrightarrow{\sim}\rT^{*}(\rT_{*}(\breve{\beta}_{\mathbb{Q}}^{*}(V))).
		\label{adjonction T WT V}
	\end{equation}
\end{prop}
\textit{Preuve}. Soit $V$ une $\mathfrak{C}$-représentation de $\pi_{1}(X_{\overline{\eta}},\overline{x})$ de Weil-Tate relativement à $X$. Choisissons un $\mathscr{O}_{\XX}[\frac{1}{p}]$-module de Weil-Tate $\mathscr{F}$ et un isomorphisme $\breve{\overline{\mathscr{B}}}_{\mathbb{Q}}$-linéaire
\begin{equation}
	\tau: \breve{\beta}_{\mathbb{Q}}^{*}(V)\simeq \rT^{*}(\mathscr{F}).
	\label{condition BQ associe v}
\end{equation}
D'après \eqref{formule projectif T}, $\rT_{*}(\breve{\beta}_{\mathbb{Q}}^{*}(V))$ est de Weil-Tate.

En vertu de \eqref{formule projectif T}, on a un isomorphisme $\breve{\overline{\mathscr{B}}}_{\mathbb{Q}}$-linéaire, canonique et fonctoriel en $\mathscr{F}$
\begin{equation}
	\rT^{*}(\mathscr{F})\xrightarrow{\sim} \rT^{*}(\rT_{*}(\rT^{*}(\mathscr{F}))).
	\label{isomorphisme adjonction T}
\end{equation}
On en déduit compte tenu de \eqref{condition BQ associe v} et \eqref{isomorphisme adjonction T} un isomorphisme $\breve{\overline{\mathscr{B}}}_{\mathbb{Q}}$-linéaire
\begin{displaymath}
		\breve{\beta}_{\mathbb{Q}}^{*}(V) \xrightarrow{\sim}\rT^{*}(\rT_{*}(\breve{\beta}_{\mathbb{Q}}^{*}(V))).
\end{displaymath}
Montrons que cet isomorphisme ne dépend pas des choix de $\mathscr{F}$ et de l'isomorphisme $\tau$ et qu'il est fonctoriel en $V$. Il suffit de montrer que, pour tous $\mathscr{O}_{\XX}[\frac{1}{p}]$-modules localement projectifs de type fini $\mathscr{F}$, $\mathscr{F}'$ et tout morphisme $\breve{\overline{\mathscr{B}}}_{\mathbb{Q}}$-linéaire $\theta: \rT^{*}(\mathscr{F})\to \rT^{*}(\mathscr{F}')$, le diagramme
\begin{equation}
	\xymatrix{
		\rT^{*}(\mathscr{F}) \ar[r]^-{\sim} \ar[d]_{\theta} & \rT^{*}(\rT_{*}(\rT^{*}(\mathscr{F}))) \ar[d]^{\rT^{*}(\rT_{*}(\theta))} \\
		\rT^{*}(\mathscr{F}')\ar[r]^-{\sim} & \rT^{*}(\rT_{*}(\rT^{*}(\mathscr{F}')))
	}
	\label{commutatif Bassocie pour F}
\end{equation}
est commutatif. D'après la pleine fidélité du foncteur $\rT^{*}$ \ref{general property of T pushback}(ii), $\theta$ est l'image d'un morphisme $\vartheta:\mathscr{F}\to \mathscr{F}'$ de $\LPtf(\mathscr{O}_{\XX}[\frac{1}{p}])$ par $\rT^{*}$. La commutativité de \eqref{commutatif Bassocie pour F} résulte de la fonctorialité de \eqref{isomorphisme adjonction T} en $\mathscr{F}$.

\begin{coro}
	Sous les hypothèses de \ref{formule projectif T lemma}. Pour toute $\mathfrak{C}$-représentation $V$ de Weil-Tate relativement à $X$, le $\mathscr{O}_{\XX}[\frac{1}{p}]$-module $\rT_{*}(\breve{\beta}_{\mathbb{Q}}^{*}(V))$ est de Weil-Tate. On désigne par $\mathscr{T}_{\XX}$ le foncteur
\begin{equation}
	\mathscr{T}_{\XX}:\Rep_{\mathfrak{C}}^{\WT/X}(\pi_{1}(X_{\overline{\eta}},\overline{x}))\to \Mod^{\WT}(\mathscr{O}_{\XX}[\frac{1}{p}]),\qquad V\mapsto \rT_{*}(\breve{\beta}_{\mathbb{Q}}^{*}(V)).
	\label{definition foncteur TXX}
\end{equation}
Alors, $V$ et $\mathscr{T}_{\XX}(V)$ sont $\breve{\overline{\mathscr{B}}}_{\mathbb{Q}}$-associés.
\end{coro}

\begin{prop}
	Soient $X$ une $S$-courbe semi-stable et régulière et $\overline{x}$ un point géométrique de $X_{\overline{\eta}}$. Les foncteurs $\mathscr{V}_{\XX}$ \eqref{definition foncteur VXX} et $\mathscr{T}_{\XX}$ \eqref{definition foncteur TXX} sont des équivalence de catégories quasi-inverses l'une de l'autre.
	\label{equivalence WT XX}
\end{prop}
\textit{Preuve}. D'après \ref{isomorphisme VX WT lemma}(ii), pour tout objet $\mathscr{F}$ de $\Mod^{\WT}(\mathscr{O}_{\XX}[\frac{1}{p}])$, on a un isomorphisme $\breve{\overline{\mathscr{B}}}_{\mathbb{Q}}$-linéaire, canonique et fonctoriel
\begin{equation}
	\rT^{*}(\mathscr{F})\xrightarrow{\sim} \breve{\beta}_{\mathbb{Q}}^{*}(\mathscr{V}_{\XX}(\mathscr{F})).
	\label{B associe F VF}
\end{equation}
Composant l'isomorphisme canonique \eqref{formule projectif T} et l'image de \eqref{B associe F VF} par $\rT_{*}$, on déduit un isomorphisme fonctoriel
\begin{equation}
	\mathscr{F}\xrightarrow{\sim} \mathscr{T}_{\XX}(\mathscr{V}_{\XX}(\mathscr{F})).
	\label{F TXX VXX}
\end{equation}

D'autre part, d'après \eqref{adjonction T WT V}, pour tout objet $V$ de $\Rep^{\WT/X}_{\mathfrak{C}}(\pi_{1}(X_{\overline{\eta}},\overline{x}))$, on a un isomorphisme $\breve{\overline{\mathscr{B}}}_{\mathbb{Q}}$-linéaire, canonique et fonctoriel
\begin{equation}
	\breve{\beta}^{*}_{\mathbb{Q}}(V)\xrightarrow{\sim} \rT^{*}(\mathscr{T}_{\XX}(V)).
	\label{B associe V TV}
\end{equation}
Composant l'isomorphisme canonique $V\xrightarrow{\sim} \mathscr{W}_{\mathbb{Q}}(\breve{\beta}_{\mathbb{Q}}^{*}(V))$ (\ref{rep generalise rep vrai}(ii)) et l'image de \eqref{B associe V TV} par $\mathscr{W}_{\mathbb{Q}}$, on déduit un isomorphisme fonctoriel
\begin{equation}
	V\xrightarrow{\sim} \mathscr{V}_{\XX}(\mathscr{T}_{\XX}(V)).
	\label{V VXX TXX}
\end{equation}

\begin{nothing} \label{notations morphisme 2 models}
	Soient $X$ une $S$-courbe semi-stable et régulière, $(S',\eta')$ un trait fini sur $(S,\eta)$ et $X'$ un modèle semi-stable et régulier de $X_{\eta'}$ dominant $X_{S'}$. Reprenons les notations de \ref{notations general section WT} pour $X$. On fixe un $\eta$-morphisme $\overline{\eta}\to \eta'$ et par suite un $S$-morphisme $\overline{S}\to S'$. On désigne par $\XX'$ le schéma formel complété $p$-adique de $X\times_{S'}\overline{S}$ et par $\mathfrak{g}:\XX'\to \XX$ le morphisme de schémas formels induit par la projection canonique $X'\to X_{S'}$. Reprenant les notations de \ref{fonctoriel Faltings topos}, on désigne par
\begin{equation}
	\rT': (\widetilde{E}'^{\mathbb{N}^{\circ}}_{s},\breve{\overline{\mathscr{B}}}')\to (X'_{s,\zar},\mathscr{O}_{\XX'}).
\end{equation}
le morphisme \eqref{moprhisme de topos anneles T} relatif à $X'$. On a alors un diagramme commutatif à isomorphismes près (cf. \eqref{fonctorialite beta n} et \cite{AGT} III.(9.11.15)).
\begin{equation}
	\xymatrix{
		(X'_{s,\zar}, \mathscr{O}_{\XX'} ) \ar[d]_{\mathfrak{g}} & (\widetilde{E}'^{\mathbb{N}^{\circ}}_{s},\breve{\overline{\mathscr{B}}}') \ar[l]_-{\rT'} \ar[r]^{\breve{\beta}'} \ar[d]^{\breve{\Phi}} & ( X^{\mathbb{N}^{\circ}}_{\overline{\eta},\fet},\breve{\oo}) \ar@{=}[d] \\
		(X_{s,\zar},\mathscr{O}_{\XX}) & (\widetilde{E}^{\mathbb{N}^{\circ}}_{s},\breve{\overline{\mathscr{B}}}) \ar[l]_-{\rT} \ar[r]^{\breve{\beta}} & (X_{\overline{\eta},\fet}^{\mathbb{N}^{\circ}},\breve{\oo})
		}
	\label{fonctorialite breve Phi T}
\end{equation}
\end{nothing}

\begin{prop}
	Conservons les notations de \ref{notations morphisme 2 models}. Alors:
	
	\textnormal{(i)} L'image inverse d'un $\mathscr{O}_{\XX}[\frac{1}{p}]$-module de Weil-Tate par $\mathfrak{g}$ est un $\mathscr{O}_{\XX'}[\frac{1}{p}]$-module de Weil-Tate.

	\textnormal{(ii)} Toute $\mathfrak{C}$-représentation continue de Weil-Tate relativement à $X$ est de Weil-Tate relativement à $X'$.

	\textnormal{(iii)} On a des diagrammes commutatifs à isomorphismes près
	\begin{equation}
	\xymatrix{
		\Mod^{\WT}(\mathscr{O}_{\XX}[\frac{1}{p}]) \ar[d]_{\mathfrak{g}^{*}} \ar[r]^-{\mathscr{V}_{\XX}} & \Rep_{\mathfrak{C}}^{\WT/X}(\pi_{1}(X_{\overline{\eta}},\overline{x})) \ar@{^{(}->}[d]  \\
		\Mod^{\WT}(\mathscr{O}_{\XX'}[\frac{1}{p}])   \ar[r]^-{\mathscr{V}_{\XX'}} & \Rep_{\mathfrak{C}}^{\WT/X'}(\pi_{1}(X_{\overline{\eta}},\overline{x}))
	}
	\label{compatible model foncteur V}
	\end{equation}
	et
	\begin{equation}
	\xymatrix{
		\Rep_{\mathfrak{C}}^{\WT/X}(\pi_{1}(X_{\overline{\eta}},\overline{x})) \ar@{^{(}->}[d] \ar[r]^-{\mathscr{T}_{\XX}} & \Mod^{\WT}(\mathscr{O}_{\XX}[\frac{1}{p}]) \ar[d]_{\mathfrak{g}^{*}} \\
		\Rep_{\mathfrak{C}}^{\WT/X'}(\pi_{1}(X_{\overline{\eta}},\overline{x})) \ar[r]^-{\mathscr{T}_{\XX'}} & \Mod^{\WT}(\mathscr{O}_{\XX'}[\frac{1}{p}])
	}
	\label{compatible model foncteur T}
	\end{equation}
	où $\mathscr{T}_{\XX'}$ et $\mathscr{V}_{\XX'}$ sont les foncteurs \eqref{definition foncteur TXX} et \eqref{definition foncteur VXX} relatifs à $X'$.
	\label{lemma compatible TV formal models}
\end{prop}
\textit{Preuve}. Soient $\mathscr{F}$ un $\mathscr{O}_{\XX}[\frac{1}{p}]$-module de Weil-Tate et $V$ une $\mathfrak{C}$-représentation de Weil-Tate associée. Il existe donc un isomorphisme $\breve{\overline{\mathscr{B}}}_{\mathbb{Q}}$-linéaire
\begin{equation}
	\rT^{*}(\mathscr{F})\xrightarrow{\sim} \breve{\beta}_{\mathbb{Q}}^{*}(V).
\end{equation}
En vertu de \eqref{fonctorialite breve Phi T}, on en déduit un isomorphisme $\breve{\overline{\mathscr{B}}}_{\mathbb{Q}}'$-linéaire
\begin{equation}
	\rT'^{*}(\mathfrak{g}^{*}(\mathscr{F})) \xrightarrow{\sim} \breve{\beta}'^{*}_{\mathbb{Q}}(V),
\end{equation}
d'où les assertions (i) et (ii).

Comme $\mathscr{V}_{\XX}=\rT^{*}\circ \mathscr{W}_{\mathbb{Q}}$, la commutativité du diagramme \eqref{compatible model foncteur V} résulte de \ref{WX WX' compatible} et \eqref{fonctorialite breve Phi T}. La commutativité du diagramme \eqref{compatible model foncteur T} s'ensuit par \ref{equivalence WT XX}, d'où l'assertion (iii).

\begin{nothing} \label{Notations C X LPft Vect}
	Dans la suite de cette section, on fixe une $\overline{K}$-courbe propre et lisse $C$ et un point géométrique $\overline{x}$ de $C$. On pose $\check{C}=C\otimes_{\overline{K}}\mathfrak{C}$. Soit $S'$ un trait fini sur $S$. Un \textit{$S'$-modèle de $C$} est la donnée d'un $S$-morphisme $\overline{S}\to S'$ et d'une $S'$-courbe propre $X$ munie d'un isomorphisme $C\simeq X\times_{S'}\overline{\eta}$. Soit $X$ un $S'$-modèle semi-stable et régulier de $C$. Reprenant les notations de \ref{jmathXX X coh}, on a un foncteur pleinement fidèle en vertu de \eqref{lemma LPtf VB}
\begin{equation}
	\jmath_{\XX}: \LPtf(\mathscr{O}_{\XX}[\frac{1}{p}])\to \Vect_{\check{C}}.
	\label{jmath LPtf Vect}
\end{equation}
\end{nothing}
\begin{definition}
	(i) On dit qu'un fibré vectoriel $F$ sur $\check{C}$ est \textit{de Weil-Tate} s'il existe un trait $S'$ fini sur $S$, un $S'$-modèle semi-stable et régulier $X$ de $C$ et un $\mathscr{O}_{\XX}[\frac{1}{p}]$-module de Weil-Tate $\mathscr{F}$ tel que $\jmath_{\XX}(\mathscr{F})\simeq F$ \eqref{jmath LPtf Vect}.

	(ii) On dit qu'une $\mathfrak{C}$-représentation continue $V$ de $\pi_{1}(C,\overline{x})$ est \textit{de Weil-Tate} s'il existe un trait $S'$ fini sur $S$, un $S'$-modèle semi-stable $X$ de $C$ tels que $V$ soit de Weil-Tate relativement à $X$.
	\label{Weil-Tate}
\end{definition}
On désigne par $\BB_{\check{C}}^{\WT}$ la sous-catégorie pleine de $\Vect_{\check{C}}$ formée des fibrés vectoriels de Weil-Tate et par $\SyWT$ la sous-catégorie pleine de $\Rep_{\mathfrak{C}}^{\cont}(\pi_{1}(C,\overline{x}))$ formée des $\mathfrak{C}$-représentations continues de Weil-Tate. Soit $X$ un $S$-modèle semi-stable et régulier de $C$. La restriction du foncteur $\jmath_{\XX}$ \eqref{jmath LPtf Vect} à $\Mod^{\WT}(\mathscr{O}_{\XX}[\frac{1}{p}])$ induit un foncteur pleinement fidèle que l'on note encore
\begin{equation}
	\jmath_{\XX}: \Mod^{\WT}(\mathscr{O}_{\XX}[\frac{1}{p}])\to \BB_{\check{C}}^{\WT}.
	\label{jmath XX check C WT}
\end{equation}
Soient $S'$ un trait fini sur $S$ et $X'$ un $S'$-modèle semi-stable et régulier de $C$ dominant $X_{S'}$ et reprenons les notations de \ref{notations morphisme 2 models}. En vertu de \ref{lemma compatible TV formal models}(i) et de la définition de $\jmath_{\XX}$ \eqref{jmathXX coh}, on vérifie que le diagramme
	\begin{equation}
		\xymatrix{
			\Mod^{\WT}(\mathscr{O}_{\XX}[\frac{1}{p}]) \ar[rd]^{\jmath_{\XX}} \ar[d]_{\mathfrak{g}^{*}} & \\
			\Mod^{\WT}(\mathscr{O}_{\XX'}[\frac{1}{p}]) \ar[r]^-{\jmath_{\XX'}} & \BB_{\check{C}}^{\WT},
		}
	\end{equation}
	où $\jmath_{\XX'}$ désigne le foncteur \eqref{jmath XX check C WT} relatif à $X'$, est commutatif.

\begin{lemma}
	\textnormal{(i)} Pour tout morphisme $f$ de $\BB_{\check{C}}^{\WT}$, il existe un trait $S'$ fini sur $S$ et un $S'$-modèle semi-stable et régulier $X$ de $C$ tels que $f$ soit contenu dans l'image du foncteur $\jmath_{\XX}$ \eqref{jmath XX check C WT}.
	
	\textnormal{(ii)} Pour tout morphisme $f$ de $\SyWT$, il existe un trait $S'$ fini sur $S$ et un $S'$-modèle semi-stable $X$ de $C$ tels que $f$ provienne de la sous-catégorie pleine $\Rep^{\WT/X}_{\mathfrak{C}}(\pi_{1}(C,\overline{x}))$ de $\Rep^{\WT}_{\mathfrak{C}}(\pi_{1}(C,\overline{x}))$.

	\label{lemma models}
\end{lemma}
\textit{Preuve}. (i) Soient $F$ et $F'$ deux fibrés vectoriels de Weil-Tate sur $\check{C}$. D'après \ref{semistable rcm coro}, il existe un trait $S'$ fini sur $S$ et un $S'$-modèle $X$ semi-stable et régulier de $C$ tels que $F$ et $F'$ soient les images des $\mathscr{O}_{\XX}[\frac{1}{p}]$-modules localement projectifs de type fini par $\jmath_{\XX}$ \eqref{jmath XX check C WT}. L'assertion résulte alors de la pleine fidélité de \eqref{jmath XX check C WT}.

En calquant la démonstration de (i), on vérifie l'assertion (ii).


\begin{nothing}	\label{foncteur DW pour WT fibres}
	En vertu de \ref{lemma compatible TV formal models}, \ref{Weil-Tate} et \ref{lemma models}, les foncteurs $\mathscr{V}_{\XX}$ fournissent un foncteur
	\begin{equation}
		\mathscr{V}_{\check{C}}: \BB_{\check{C}}^{\WT}\to \SyWT.
		\label{varrho check C WT}
	\end{equation}
De même, les foncteurs $\jmath_{\XX}\circ \mathscr{T}_{\XX}$ fournissent un foncteur:
	\begin{equation}
		\mathscr{T}_{\check{C}}: \SyWT \to \BB_{\check{C}}^{\WT}.
		\label{foncteur T WT}
	\end{equation}
\end{nothing}

\begin{theorem}
	Les foncteurs $\mathscr{V}_{\check{C}}$ \eqref{varrho check C WT} et $\mathscr{T}_{\check{C}}$ \eqref{foncteur T WT} sont des équivalences de catégories quasi-inverses l'une de l'autre.
	\label{theorem WT}
\end{theorem}
\textit{Preuve}. Pour tout objet $V$ de $\Rep_{\mathfrak{C}}^{\WT}(\pi_{1}(C,\overline{x}))$, quitte à remplacer $S$ par une extension finie, il existe un $S$-modèle semi-stable et régulier $X$ de $C$ tel que $V$ soit de Weil-Tate relativement à $X$. D'après \eqref{V VXX TXX}, on a un isomorphisme
\begin{equation}
	V\xrightarrow{\sim} \mathscr{V}_{\check{C}}(\mathscr{T}_{\check{C}}(V))
\end{equation}
Il résulte de la fonctorialité de \eqref{V VXX TXX} et de \ref{lemma models}(ii) que celui-ci est fonctoriel en $V$.

De même, pour tout objet $F$ de $\BB_{\check{C}}$, quitte à remplacer $S$ par une extension finie, il existe un $S$-modèle semi-stable et régulier $X$ de $C$ et un objet $\mathscr{F}$ de $\Mod^{\WT}(\mathscr{O}_{\XX}[\frac{1}{p}])$ tel que $F\simeq \jmath_{\XX}(\mathscr{F})$. D'après \eqref{F TXX VXX}, on a un isomorphisme
\begin{equation}
	F\xrightarrow{\sim} \mathscr{T}_{\check{C}}(\mathscr{V}_{\check{C}}(F)).
\end{equation}
Il résulte de la fonctorialité de \eqref{F TXX VXX} et de \ref{lemma models}(i) que celui-ci est fontoriel en $F$.

\begin{theorem}
	Tout fibré vectoriel de Deninger-Werner sur $\check{C}$ est de Weil-Tate. De plus, la restriction du foncteur $\mathscr{V}_{\check{C}}:\BB_{\check{C}}^{\WT}\to \Rep_{\mathfrak{C}}^{\cont}(\pi_{1}(C,\overline{x}))$ \eqref{varrho check C WT} à la sous-catégorie $\BB_{\check{C}}^{\DW}$ s'identifie au foncteur de Deninger-Werner $\mathbb{V}_{\check{C}}$ \eqref{rho C check}.
	\label{DW implique WT}
\end{theorem}

On présentera la démonstration dans \ref{DW implique WT coro} et \ref{foncteur DW WT pf}.

\begin{coro}
	Le foncteur de Deninger-Werner $\mathbb{V}_{\check{C}}$ \eqref{rho C check} est pleinement fidèle.
\end{coro}
\textit{Preuve}. Cela résulte de \ref{theorem WT} et \ref{DW implique WT}.

\section{Correspondance de Simpson $p$-adique pour les représentations de Weil-Tate} \label{Simpson p adique}
\begin{nothing}
On munit $S$ de la structure logarithmique $\mathscr{M}_{S}$ définie par son point fermé (\cite{AGT} II.6.1), autrement dit, $\mathscr{M}_{S}=j_{*}(\mathscr{O}_{\eta}^{\times})\cap \mathscr{O}_{S}$, où $j:\eta \to S$ est l'injection canonique. On munit $\overline{S}$ et $\check{\overline{S}}$ \eqref{notations 11} des structures logarithmiques $\mathscr{M}_{\overline{S}}$ et $\mathscr{M}_{\check{\overline{S}}}$ images inverses de $\mathscr{M}_{S}$ (cf. \cite{AGT} II.5.10).

On associe fonctoriellement à toute $\mathbb{Z}_{(p)}$-algèbre $A$ l'anneau
\begin{equation}
	A^{\flat}=\varprojlim_{x\mapsto x^p}A/pA
	\label{perfectoid}
\end{equation}
et un homomorphisme $\theta:W(A^{\flat})\to \widehat{A}$ de l'anneau $W(A^{\flat})$ des vecteurs de Witt de $A^{\flat}$ dans le séparé complété $p$-adique $\widehat{A}$ de $A$ (cf. \cite{AGT} II.9.3).

On désigne par $\underline{p}$ l'élément de $(\mathcal{O}_{\overline{K}})^{\flat}$ \eqref{perfectoid} induit par la suite $(p^{\frac{1}{pn}})_{n\ge0}$ \eqref{notations 11} et on pose
\begin{equation}
	\xi=[\underline{p}]-p \in W((\mathcal{O}_{\overline{K}})^{\flat}),
	\label{def of xi}
\end{equation}
où $[~]$ est le représentant multiplicatif. On pose
\begin{equation}
	\mathscr{A}_2(\mathcal{O}_{\overline{K}})=W((\mathcal{O}_{\overline{K}})^{\flat})/\Ker(\theta)^2.
	\label{def of A2OK}
\end{equation}
On a une suite exacte (\cite{AGT} II.9.5(iv))
\begin{equation}
	\xymatrix{
		0\ar[r] & W( (\mathcal{O}_{\overline{K}})^{\flat})\ar[r]^{\cdot \xi}& W( (\mathcal{O}_{\overline{K}})^{\flat})\ar[r]^{\theta} &\oo\to 0.
	}
\end{equation}
Elle induit une suite exacte
\begin{equation}
	0\to \oo \xrightarrow{\cdot \xi}\mathscr{A}_2(\mathcal{O}_{\overline{K}})\xrightarrow{\theta}\oo \to 0,
	\label{suite exacte de A2S}
\end{equation}
où $\cdot \xi$ est induit par la multiplication par $\xi$ dans $\mathscr{A}_2(\mathcal{O}_{\overline{K}})$. L'idéal $\Ker(\theta)$ est un $\oo$-module libre de base $\xi$. On le note $\xi\oo$ et on note $\xi^{-1}\oo$ le $\oo$-module dual de $\xi\oo$. Pour tout $\oo$-module $M$, on désigne les $\oo$-modules $M\otimes_{\oo}\xi\oo$ et $M\otimes_{\oo}\xi^{-1}\oo$ simplement par $\xi M$ et $\xi^{-1}M$ respectivement.

On pose
\begin{equation}
	\mathscr{A}_2(\overline{S})=\Spec(\mathscr{A}_2(\mathcal{O}_{\overline{K}}))
	\label{def of A2S}
\end{equation}
que l'on munit de la structure logarithmique $\mathscr{M}_{\mathscr{A}_2(\overline{S})}$ définie dans (\cite{AGT} II.9.8). L'homomorphisme $\theta$ induit alors une immersion fermée exacte (\cite{AGT} II.5.22)
\begin{equation}
	i_{\overline{S}}: (\check{\overline{S}},\mathscr{M}_{\check{\overline{S}}})\to (\mathscr{A}_2(\overline{S}),\mathscr{M}_{\mathscr{A}_2(\overline{S})}).
	\label{immersion fermee exacte}
\end{equation}
\end{nothing}
\begin{nothing}
Dans la suite de cette section, on se donne une $S$-courbe semi-stable et régulière $X$ que l'on munit de la structure logarithmique $\mathscr{M}_{X}$ définie par sa fibre spécial $X_s$. Le morphisme de schémas logarithmiques $f:(X,\mathscr{M}_{X})\to (S,\mathscr{M}_S)$ est alors adéquat au sens de (\cite{AGT} III.4.7). On notera que $X_{\eta}$ est le sous-schéma ouvert maximal de $X$ où la structure logarithmique $\mathscr{M}_{X}$ est triviale. Pour alléger les notations, on pose
\begin{equation}
	\widetilde{\Omega}_{X/S}^{1}=\Omega_{(X,\mathscr{M}_X)/(S,\mathscr{M}_S)}^1.
	\label{omega XS log}
\end{equation}

On munit $\check{\overline{X}}=X\times_{S}\check{\overline{S}}$ \eqref{basic notation} de la structure logarithmique $\mathscr{M}_{\check{\overline{X}}}$ image inverse de $\mathscr{M}_X$. Posons $\mathcal{F}=\FHom_{\mathscr{O}_{\check{\overline{X}}}}(\widetilde{\Omega}_{X/S}^{1}\otimes_{\mathscr{O}_{S}}\mathscr{O}_{\check{\overline{S}}},\mathscr{O}_{\check{\overline{X}}})$. Comme $X$ est une $S$-courbe, on a alors $\rH^{2}(X_{\check{\overline{\eta}}},\mathscr{F}_{\check{\overline{\eta}}})=0$ et $\rH^{2}(X_{s},\mathscr{F}_{s})=0$. Compte tenu de (\cite{Mu70} \S~5 Cor.2), on en déduit $\rH^{2}(\check{\overline{X}},\mathscr{F})=0$. D'après (\cite{Ka89} Prop. 3.14), il existe une $(\mathscr{A}_2(\overline{S}),\mathscr{M}_{\mathscr{A}_2(\overline{S})})$-déformation lisse $(\widetilde{X},\mathscr{M}_{\widetilde{X}})$ de $(\check{\overline{X}},\mathscr{M}_{\check{\overline{X}}})$, c'est-à-dire un diagramme commutatif
\begin{equation}
	\xymatrix{
		(\check{\overline{X}},\mathscr{M}_{\check{\overline{X}}})\ar[r] \ar[d] & (\widetilde{X},\mathscr{M}_{\widetilde{X}}) \ar[d]\\\
		(\check{\overline{S}},\mathscr{M}_{\check{\overline{S}}})\ar[r]^-{i_{\check{\overline{S}}}} & (\mathscr{A}_2(\overline{S}),\mathscr{M}_{\mathscr{A}_2(\overline{S})}).
	}
	\label{A2 deformation}
\end{equation}
On fixe une telle déformation $(\widetilde{X},\mathscr{M}_{\widetilde{X}})$ \eqref{A2 deformation} dans la suite de cette section.
\end{nothing}
\begin{nothing}
	On note $\XX$ le schéma formel complété $p$-adique de $\overline{X}$. On désigne par $\xi^{-1}\widetilde{\Omega}^1_{\XX/\mathscr{S}}$ le complété $p$-adique du $\mathscr{O}_{\overline{X}}$-module $\xi^{-1} \widetilde{\Omega}_{\overline{X}/\overline{S}}^{1}=\xi^{-1} \widetilde{\Omega}_{X/S}^{1}\otimes_{\mathscr{O}_{X}}\mathscr{O}_{\overline{X}}$ et on pose, pour tout entier $j\ge 0$, $\xi^{-j}\widetilde{\Omega}^{j}_{\XX/\mathscr{S}}=\wedge^{j}_{\mathscr{O}_{\XX}}(\xi^{-1}\widetilde{\Omega}^1_{\XX/\mathscr{S}})$. On sous-entend par $\mathscr{O}_{\XX}[\frac{1}{p}]$-module de Higgs à coefficients dans $\xi^{-1}\widetilde{\Omega}^1_{\XX/\mathscr{S}}$, un $\mathscr{O}_{\XX}[\frac{1}{p}]$-module de Higgs à coefficients dans $\xi^{-1}\widetilde{\Omega}^{1}_{\XX/\mathscr{S}}[\frac{1}{p}]$ \eqref{topos module higgs}. On désigne par $\mathbf{MH}(\mathscr{O}_{\XX}[\frac{1}{p}],\xi^{-1}\widetilde{\Omega}_{\XX/\mathscr{S}}^1)$ la catégorie de tels modules et par $\mathbf{MH}^{\coh}(\mathscr{O}_{\XX}[\frac{1}{p}],\xi^{-1}\widetilde{\Omega}_{\XX/\mathscr{S}}^1)$ la sous-catégorie pleine formée des modules de Higgs dont le $\mathscr{O}_{\XX}[\frac{1}{p}]$-module sous-jacent est cohérent. On appelle \textit{$\mathscr{O}_{\XX}[\frac{1}{p}]$-fibré de Higgs} à coefficients dans $\xi^{-1}\widetilde{\Omega}^1_{\XX/\mathscr{S}}$ tout $\mathscr{O}_{\XX}[\frac{1}{p}]$-module de Higgs à coefficients dans $\xi^{-1}\widetilde{\Omega}^1_{\XX/\mathscr{S}}$ dont le $\mathscr{O}_{\XX}[\frac{1}{p}]$-module sous-jacent est localement projectif de type fini \eqref{LPtf}.
\end{nothing}
\begin{nothing}	\label{AG associe}
	On reprend les notations de \S~\ref{Topos de Faltings} pour le $S$-schéma $X$. Soient $\mathscr{M}$ un $\breve{\overline{\mathscr{B}}}_{\mathbb{Q}}$-module adique et de type fini et $\mathscr{N}$ un $\mathscr{O}_{\XX}[\frac{1}{p}]$-fibré de Higgs à coefficients dans $\xi^{-1}\widetilde{\Omega}^1_{\XX/\mathscr{S}}$. Dans (\cite{AGT} III.12.10(ii)), Abbes et Gros définissent la propriété pour $\mathscr{M}$ et $\mathscr{N}$ d'\textit{être associés}. On dit que $\mathscr{M}$ est \textit{de Dolbeault} s'il admet un fibré de Higgs associé et que $\mathscr{N}$ est \textit{soluble} s'il admet un $\breve{\overline{\mathscr{B}}}_{\mathbb{Q}}$-module associé (cf. \cite{AGT} III.12.11).

On désigne par $\Mod^{\Dolb}_{\mathbb{Q}}(\breve{\overline{\mathscr{B}}})$ la sous-catégorie pleine de $\Mod^{\atf}_{\mathbb{Q}}(\breve{\overline{\mathscr{B}}})$ formée des $\breve{\overline{\mathscr{B}}}_{\mathbb{Q}}$-modules de Dolbeault, et par $\mathbf{MH}^{\sol}(\mathscr{O}_{\XX}[\frac{1}{p}],\xi^{-1}\widetilde{\Omega}_{\XX/\mathscr{S}}^1)$ la sous-catégorie pleine de $\mathbf{MH}(\mathscr{O}_{\XX}[\frac{1}{p}],\xi^{-1}\widetilde{\Omega}_{\XX/\mathscr{S}}^1)$ formée des $\mathscr{O}_{\XX}[\frac{1}{p}]$-fibrés de Higgs solubles à coefficients dans $\xi^{-1}\widetilde{\Omega}_{\XX/\mathscr{S}}^1$. Il existe un foncteur (cf. \cite{AGT} III.12.18)
	\begin{equation}
		\mathscr{H}: \Mod_{\mathbb{Q}}^{\Dolb}(\breve{\overline{\mathscr{B}}})\to \mathbf{MH}^{\sol}(\mathscr{O}_{\XX}[\frac{1}{p}],\xi^{-1}\widetilde{\Omega}_{\XX/\mathscr{S}}^1),
		\label{foncteur H dolb}
	\end{equation}
et un foncteur dans l'autre sens (cf. \cite{AGT} III.12.23)
	\begin{equation}
		\mathscr{W}:\mathbf{MH}^{\sol}(\mathscr{O}_{\XX}[\frac{1}{p}],\xi^{-1}\widetilde{\Omega}_{\XX/\mathscr{S}}^1)\to \Mod_{\mathbb{Q}}^{\Dolb}(\breve{\overline{\mathscr{B}}}).
		\label{foncteur W}
	\end{equation}
\end{nothing}

\begin{theorem}[\cite{AGT} III.12.26]
	Les foncteurs \eqref{foncteur H dolb} et \eqref{foncteur W} sont des équivalences de catégories quasi-inverses l'une de l'autre.
	\label{equivalence of categories Simpson}
\end{theorem}

%
\begin{lemma} \label{lemma fibre vect Dolb}
	Soit $\mathscr{F}$ un $\mathscr{O}_{\XX}[\frac{1}{p}]$-module localement projectif de type fini. Alors, le $\breve{\overline{\mathscr{B}}}_{\mathbb{Q}}$-module $\rT^{*}(\mathscr{F})$ \eqref{T star rationel} et le fibré de Higgs $(\mathscr{F},0)$ sont associés \eqref{AG associe}, et on a un isomorphisme canonique et fonctoriel
	\begin{equation}
		(\mathscr{F},0)\xrightarrow{\sim} \mathscr{H}(\rT^{*}(\mathscr{F})).
		\label{F0 iso to HTF}
	\end{equation}
\end{lemma}
\textit{Preuve}. Démontrons que $\rT^{*}(\mathscr{F})$ et $(\mathscr{F},0)$ sont \textit{$r$-associés} dans le sens de (\cite{AGT} III.12.10(i)) pour tout nombre rationnel $r>0$. On désigne par $\breve{\mathscr{C}}^{(r)}$ la \textit{$\breve{\overline{\mathscr{B}}}$-algèbre de Higgs-Tate d'épaisseur $r$ associée à $(\widetilde{X},\mathscr{M}_{\widetilde{X}})$} (cf. \cite{AGT} III.10.31), par $\breve{d}^{(r)}$ la $\breve{\overline{\mathscr{B}}}$-dérivation de $\breve{\mathscr{C}}^{(r)}$ défini dans (\cite{AGT} (III.12.7.1)), par $\Xi^{r}$ la catégorie des $p^{r}$-isoconnexions intégrables relativement à l'extension $\breve{\mathscr{C}}^{(r)}/\breve{\overline{\mathscr{B}}}$ (cf. \cite{AGT} III.6.10) et par $\Xi^{r}_{\mathbb{Q}}$ la catégorie des objets de $\Xi^{r}$ à isogénie près \eqref{categorie a isogenies}. Dans la suite, on va construire un isomorphisme canonique fonctoriel
\begin{equation}
	\rT^{r+}(\mathscr{F},0)\xrightarrow{\sim} \mathfrak{S}^{r}(\rT^{*}(\mathscr{F})),
	\label{iso r associe}
\end{equation}
où $\mathfrak{S}^{r}: \Mod^{\atf}_{\mathbb{Q}}(\breve{\overline{\mathscr{B}}}) \to \Xi_{\mathbb{Q}}^{r}$ et $\rT^{r+}:\mathbf{MH}^{\coh}(\mathscr{O}_{\XX}[\frac{1}{p}],\xi^{-1}\widetilde{\Omega}_{\XX/\mathscr{S}}^1) \to \Xi_{\mathbb{Q}}^{r}$ sont des foncteurs (\cite{AGT} (III.12.7.3)) et (\cite{AGT} (III.12.7.8)). En vertu de \eqref{Mod OXX Q to Mod OXXp}, il existe un $\mathscr{O}_{\XX}$-module cohérent $\mathscr{G}$ tel que $\mathscr{F}\simeq \mathscr{G}[\frac{1}{p}]$. Par définition de $\mathfrak{S}^{r}$ (\cite{AGT} (III.12.7.2)), on a
\begin{equation}
	\mathfrak{S}^{r}(\rT^{*}(\mathscr{F}))=(\breve{\mathscr{C}}^{(r)}\otimes_{\breve{\overline{\mathscr{B}}}}\rT^{*}(\mathscr{G}),\breve{\mathscr{C}}^{(r)}\otimes_{\breve{\overline{\mathscr{B}}}}\rT^{*}(\mathscr{G}), \id, p^r \breve{d}^{(r)}\otimes \id)_{\mathbb{Q}}.
	\label{Sr on TF}
\end{equation}
Compte tenu de la définition de $\rT^{r+}$ (\cite{AGT} (III.12.7.6)), on a
\begin{equation}
	\rT^{r+}(\mathscr{F},0)=(\breve{\mathscr{C}}^{(r)}\otimes_{\breve{\overline{\mathscr{B}}}}\rT^{*}(\mathscr{G}),\breve{\mathscr{C}}^{(r)}\otimes_{\breve{\overline{\mathscr{B}}}}\rT^{*}(\mathscr{G}), \id\otimes_{\breve{\overline{\mathscr{B}}}}\rT^{*}(\id), p^r \breve{d}^{(r)}\otimes \rT^{*}(\id))_{\mathbb{Q}}.
	\label{Tr on F0}
\end{equation}
On prend pour \eqref{iso r associe} l'isomorphisme induit par \eqref{Sr on TF} et \eqref{Tr on F0}. D'après (\cite{AGT} III.12.17(i)), on en déduit un isomorphisme
\begin{displaymath}
	(\mathscr{F},0)\xrightarrow{\sim} \mathscr{H}(\rT^{*}(\mathscr{F})).
\end{displaymath}
La fonctorialité de \eqref{F0 iso to HTF} se déduit de celle de \eqref{iso r associe}. Le lemme s'ensuit.

\begin{prop} \label{Higgs nul et Weil-Tate}
	Soient $\overline{x}$ un point géométrique de $X_{\overline{\eta}}$, $V$ une $\mathfrak{C}$-représentation continue de $\pi_{1}(X_{\overline{\eta}},\overline{x})$.
	
	\textnormal{(i)} Les conditions suivantes sont équivalentes:

	\qquad \textnormal{(a)} Le $\breve{\overline{\mathscr{B}}}_{\mathbb{Q}}$-module $\breve{\beta}^{*}_{\mathbb{Q}}(V)$ est de Dolbeault \eqref{AG associe} et le fibré de Higgs $\mathscr{H}(\breve{\beta}^{*}_{\mathbb{Q}}(V))$ est de champ nul.

	\qquad \textnormal{(b)} $V$ est de Weil-Tate relativement à $X$ \eqref{Weil-Tate X}.

	\textnormal{(ii)} Si les conditions de \textnormal{(i)} sont remplies, on a un isomorphisme canonique fonctoriel de fibrés de Higgs
	\begin{equation}
		\mathscr{H}(\breve{\beta}^{*}_{\mathbb{Q}}(V))\xrightarrow{\sim}(\mathscr{T}_{\XX}(V),0),
		\label{WT Dolb champs 0}
	\end{equation}
	où $\mathscr{T}_{\XX}$ est le foncteur \eqref{definition foncteur TXX}.
\end{prop}
\textit{Preuve}. Supposons d'abord que la condition \textnormal{(a)} soit remplie. On note $\mathscr{F}$ le $\mathscr{O}_{\XX}[\frac{1}{p}]$-module localement projectif de type fini sous-jacent à $\mathscr{H}(\breve{\beta}^{*}_{\mathbb{Q}}(V))$. D'après \ref{lemma fibre vect Dolb}, on a un isomorphisme canonique
\begin{equation}
	(\mathscr{F},0)\xrightarrow{\sim} \mathscr{H}(\rT^{*}(\mathscr{F})).
\end{equation}
En vertu de \ref{equivalence of categories Simpson}, on en déduit un isomorphisme de $\breve{\overline{\mathscr{B}}}_{\mathbb{Q}}$-modules
\begin{equation}
	\breve{\beta}_{\mathbb{Q}}^{*}(V)\xrightarrow{\sim}\rT^{*}(\mathscr{F}),
\end{equation}
d'où la condition (b).

Supposons que la condition (b) soit vraie. D'après \ref{formule projectif T lemma}, le $\mathscr{O}_{\XX}[\frac{1}{p}]$-module $\mathscr{T}_{\XX}(V)$ est localement projectif de type fini. De plus, on a un isomorphisme fonctoriel en $V$ \eqref{adjonction T WT V}
\begin{equation}
	\breve{\beta}_{\mathbb{Q}}^{*}(V)\xrightarrow{\sim} \rT^{*}(\mathscr{T}_{\XX}(V)).
\end{equation}
La condition (a) et l'assertion (ii) résultent alors de \ref{lemma fibre vect Dolb} appliqué à $\mathscr{F}=\mathscr{T}_{\XX}(V)$.

\begin{nothing}	\label{preparation pour HT decomposition}
	Soient $C$ une $\overline{K}$-courbe propre et lisse, $\overline{x}$ un point géométrique de $C$ et $V$ une $\mathfrak{C}$-représentation de Weil-Tate de $\pi_{1}(C,\overline{x})$. Choisissons un $\breve{\oo}$-module localement libre de type fini $\mathbb{L}$ de $C_{\fet}^{\mathbb{N}^{\circ}}$ dont la classe d'isogénie correspond à $V$ via \eqref{Rep C to Mod C}. On pose, pour tout entier $i\ge 0$, $\rH^{i}(C_{\fet}^{\mathbb{N}^{\circ}},V)=\rH^{i}(C_{\fet}^{\mathbb{N}^{\circ}},\mathbb{L})\otimes_{\mathbb{Z}}\mathbb{Q}$. On notera que cette définition ne dépend pas du choix de $\mathbb{L}$ \eqref{coh N coh class Q}.
\end{nothing}

\begin{prop} \label{HT ss prop}
	Conservons les notations de \ref{preparation pour HT decomposition} et soit $\mathscr{T}(V)$ le fibré vectoriel sur $\check{C}$ associé à $V$. On a alors une suite exacte
	\begin{equation}
		0\to \rH^{1}(\check{C}_{\zar},\mathscr{T}(V))\to \rH^{i}(C_{\fet}^{\mathbb{N}^{\circ}},V) \to \rH^{0}(\check{C}_{\zar},\mathscr{T}(V)\otimes_{\mathscr{O}_{\check{C}}}(\xi^{-1}\Omega_{\check{C}/\mathfrak{C}}^{1}))\to 0.
		\label{HT ss}
	\end{equation}
\end{prop}
\textit{Preuve}. Quitte à remplacer $K$ par une extension finie, il existe un $S$-modèle semi-stable et régulier $X$ de $C$ tel que $V$ soit de Weil-Tate relativement à $X$ \eqref{Weil-Tate X}. Reprenant les notations précédentes pour $X$, on a un isomorphisme canonique de $\breve{\overline{\mathscr{B}}}_{\mathbb{Q}}$-modules
\begin{equation}
	\breve{\beta}_{\mathbb{Q}}^{*}(V)\simeq \rT^{*}(\mathscr{T}_{\XX}(V)).
	\label{BQ associe HT decom}
\end{equation}

Reprenant les notations de \ref{not Ext Hi Q} pour le topos annelé $(\widetilde{E}_{s}^{\mathbb{N}^{\circ}},\breve{\overline{\mathscr{B}}})$, d'après \eqref{presque iso faltings modifie}, le foncteur $\breve{\beta}_{\mathbb{Q}}^{*}$ \eqref{beta Q rep} induit un isomorphisme canonique
\begin{equation}
	\rH^{1}(C_{\fet}^{\mathbb{N}^{\circ}},V)\simeq \rH^{1}(\widetilde{E}_{s}^{\mathbb{N}^{\circ}},\breve{\beta}_{\mathbb{Q}}^{*}(V)).
	\label{Faltings Q iso HT}
\end{equation}

Par ailleurs, choisissant un $\mathscr{O}_{\XX}$-module cohérent $\mathscr{F}$ tel que $\mathscr{F}[\frac{1}{p}]\simeq \mathscr{T}_{\XX}(V)$ \eqref{Mod OXX Q to Mod OXXp}, on a une suite exacte
\begin{equation}
	\rH^{i}(\XX,\rR^{j}\rT_{*}(\rT^{*}(\mathscr{F})))\Rightarrow \rH^{i+j}(\widetilde{E}_{s}^{\mathbb{N}^{\circ}}, \rT^{*}(\mathscr{F})).
\end{equation}
Compte tenu de \ref{isomorphisme section inverse p}(i), on en déduit une suite spectrale
\begin{equation} \label{ss HT XX}
	\rH^{i}(\XX,\rR^{j}\rT_{*}(\rT^{*}(\mathscr{T}_{\XX}(V))))\Rightarrow \rH^{i+j}(\widetilde{E}_{s}^{\mathbb{N}^{\circ}}, \rT^{*}(\mathscr{T}_{\XX}(V))).
\end{equation}
D'après \eqref{WT Dolb champs 0} et (\cite{AGT} III.12.34), on a un isomorphisme canonique de $\mathscr{O}_{\XX}[\frac{1}{p}]$-module
\begin{equation}
	\rR^{j}\rT_{*}(\rT^{*}(\mathscr{T}_{\XX}(V)))\simeq \mathscr{T}_{\XX}(V)\otimes_{\mathscr{O}_{\XX}}\xi^{-j}\widetilde{\Omega}_{\XX/\mathscr{S}}^{j}.
\end{equation}
Compte tenu de \ref{isomorphisme section inverse p}(ii), on a un isomorphisme $\mathfrak{C}$-linéaire
\begin{equation}
	\rH^{i}(\XX, \mathscr{T}_{\XX}(V)\otimes_{\mathscr{O}_{\XX}}(\xi^{-j}\widetilde{\Omega}_{\XX/\mathscr{S}}^{j}))\simeq \rH^{i}(\check{C}_{\zar}, \mathscr{T}(V)\otimes_{\mathscr{O}_{\check{C}}}(\xi^{-j}\Omega_{\check{C}/\mathfrak{C}}^{j})).
\end{equation}
Comme $C$ est une courbe sur $\overline{K}$, ce dernier est nul si $i\ge 2$. On en déduit par \eqref{ss HT XX} une suite exacte de $\mathfrak{C}$-espaces vectoriels
\begin{equation}
	0\to \rH^{1}(\check{C}_{\zar},\mathscr{T}(V))\to \rH^{1}(\widetilde{E}_{s}^{\mathbb{N}^{\circ}}, \rT^{*}(\mathscr{T}_{\XX}(V))) \to \rH^{0}(\check{C}_{\zar},\mathscr{T}(V)\otimes_{\mathscr{O}_{\check{C}}}(\xi^{-1}\Omega_{\check{C}/\mathfrak{C}}^{1}))\to 0.
	\label{suite exacte en Es N}
\end{equation}
La proposition s'ensuit compte tenu de \eqref{BQ associe HT decom}, \eqref{Faltings Q iso HT} et \eqref{suite exacte en Es N}.

\section{Faisceaux de $\alpha$-modules} \label{Faisceaux alpha}
\begin{nothing} \label{notations alpha A modules}
	Pour toute $\oo$-algèbre $A$, on désigne par $\Mod(A)$ la catégorie abélienne tensorielle des $A$-modules. Suivant \ref{Cat Ab almost}, prenant pour $\phi: \oo\to A=\End(\id_{\Mod(A)})$ l'homomorphisme structural, on appelle \textit{catégorie des $\alpha$-$A$-modules} et l'on note $\alpha$-$\Mod(A)$ le quotient de la catégorie $\Mod(A)$ par la sous-catégorie épaisse des $A$-modules $\alpha$-nuls. On désigne par
	\begin{equation}
		\alpha:\Mod(A)\to \alpha\textnormal{-}\Mod(A)\qquad M\mapsto M^{\alpha},
		\label{foncteur alpha ModA}
	\end{equation}
	le foncteur canonique \eqref{foncteur alpha general}. Pour tous $A$-modules $M$ et $N$, on a un isomorphisme canonique fonctoriel (\cite{AG15} (1.4.7.1))
	\begin{equation}
		\Hom_{\alpha\textnormal{-}\Mod(A)}(M^{\alpha},N^{\alpha})\xrightarrow{\sim} \Hom_{\Mod(A)}(\mm\otimes_{\oo}M,N).
		\label{Hom isomorphisme alpha modules}
	\end{equation}
	On désigne par $\sigma_{*}$ le foncteur
	\begin{equation}
		\sigma_{*}:\alpha\textnormal{-}\Mod(A)\to \Mod(A),\qquad P\mapsto \Hom_{\alpha\textnormal{-}\Mod(A)}(A^{\alpha},P),
		\label{foncteur varsigma star}
	\end{equation}
et par $\sigma_{!}$ le foncteur
\begin{equation}
	\sigma_{!}: \alpha\textnormal{-}\Mod(A) \to \Mod(A), \qquad P\mapsto \mm\otimes_{\oo}\sigma_{*}(P).
	\label{foncteur varsigma !}
\end{equation}
En vertu de \eqref{Hom isomorphisme alpha modules}, pour tout $A$-module $M$, on a un isomorphisme canonique fonctoriel
\begin{equation}
	\sigma_{*}(M^{\alpha})\xrightarrow{\sim} \Hom_{\oo}(\mm,M).
	\label{Hom oo mm M}
\end{equation}
\end{nothing}
\begin{prop}[\cite{AG15} 1.4.8]
	\textnormal{(i)} Le foncteur $\sigma_{*}$ est un adjoint à droite du foncteur $\alpha$ \eqref{foncteur alpha ModA}.

	\textnormal{(ii)} Le morphisme d'adjonction $\alpha\circ \sigma_{*}\to \id$ est un isomorphisme.	

	\textnormal{(iii)} Le foncteur $\sigma_{!}$ est un adjoint à gauche du foncteur $\alpha$.

	\textnormal{(iv)} Le morphisme d'adjonction $\id\to \alpha\circ \sigma_{!}$ est un isomorphisme.
	\label{propriete adjoint alpha}
\end{prop}

\begin{coro}[\cite{AG15} 1.4.9]
	Les foncteurs $\alpha$ et $\sigma_{!}$ sont exacts et le foncteur $\sigma_{*}$ est exact à gauche.
	\label{exactitude alpha AMod}
\end{coro}

\begin{nothing}	
On appelle \textit{$\alpha$-$\oo$-algèbre} (ou \textit{$\oo^{\alpha}$-algèbre}) un monoïde unitaire commutatif de $\alpha$-$\Mod(\oo)$ (\cite{AG15} 1.4.11). On désigne par $\Alg(\oo)$ la catégorie des $\oo$-algèbres et par $\alpha$-$\Alg(\oo)$ la catégorie des $\alpha$-$\oo$-algèbres. Le foncteur $\alpha$ étant monoïdal, il induit un foncteur que l'on note encore
	\begin{equation}
		\alpha:\Alg(\oo)\to \alpha\textnormal{-}\Alg(\oo).
		\label{alpha o algebra}
	\end{equation}
	Compte tenu de l'isomorphisme canonique $A^{\alpha}\xrightarrow{\sim}A^{\alpha}\otimes_{A^{\alpha}}A^{\alpha}$, le foncteur $\sigma_{*}$ \eqref{foncteur varsigma star} induit un foncteur que l'on note encore
	\begin{equation}
		\sigma_{*}:\alpha\textnormal{-}\Alg(\oo)\to \Alg(\oo), \qquad P\mapsto \Hom_{\alpha\textnormal{-}\Mod(A)}(A^{\alpha},P).
		\label{sigma o algebra}
	\end{equation}
	En vertu de (\cite{AG15} 1.4.12), celui-ci est un adjoint à droite du foncteur $\alpha$ \eqref{alpha o algebra} et le morphisme d'adjonction $\alpha\circ\sigma_{*}\to \id $ est un isomorphisme.
	\label{Alg alpha sigma}
\end{nothing}
\begin{nothing}
	Soit $A$ une $\oo$-algèbre. On pose $A^{\alpha}=\alpha(A)$ \eqref{alpha o algebra} et on désigne par $\Mod(A^{\alpha})$ la catégorie des $A^{\alpha}$-modules unitaires de $\alpha$-$\Mod(\oo)$. Le foncteur $\alpha:\Mod(\oo)\to \alpha\textnormal{-}\Mod(\oo)$ \'{e}tant monoïdal, il induit un foncteur
	\begin{equation}
		\alpha_{A}:\Mod(A)\to \Mod(A^{\alpha}).
	\end{equation}
	Celui-ci transforme les $\alpha$-isomorphismes en des isomorphismes. Il induit une équivalence de catégories (\cite{AG15} 1.4.13)
	\begin{equation}
		\alpha\textnormal{-}\Mod(A)\xrightarrow{\sim} \Mod(A^{\alpha}).
	\end{equation}

	On pose $A'=\sigma_{*}(A^{\alpha})$ \eqref{sigma o algebra}. D'après \ref{Alg alpha sigma}, on a un morphisme canonique de $\oo$-algèbres $\lambda:A\to A'$ qui induit un isomorphisme $A^{\alpha}\xrightarrow{\sim} A'^{\alpha}$. Pour tous $\oo^{\alpha}$-modules $P$ et $Q$, on a un morphisme $\oo$-linéaire canonique
	\begin{equation}
		\Hom_{\alpha\textnormal{-}\Mod(\oo)}(\oo^{\alpha},P)\otimes_{\oo}\Hom_{\alpha\textnormal{-}\Mod(\oo)}(\oo^{\alpha},Q)\to \Hom_{\alpha\textnormal{-}\Mod(\oo)}(\oo^{\alpha},P\otimes_{\oo^{\alpha}}Q)
	\end{equation}
	défini par fonctorialité et composition. Par suite, le foncteur $\sigma_{*}$ \eqref{foncteur varsigma star} induit un foncteur
	\begin{equation}
		\tau'_{*}:\Mod(A^{\alpha})\to \Mod(A').
	\end{equation}
Composant avec le foncteur induit par $\lambda: A\to A'$, on obtient un foncteur
\begin{equation}
		\tau_{*}:\Mod(A^{\alpha})\to \Mod(A).
	\label{le foncteur tau A module}
\end{equation}
	On désigne par $\tau_{!}$ le foncteur
	\begin{equation}
		\tau_{!}:\Mod(A^{\alpha})\to \Mod(A) \qquad P\mapsto \mm\otimes_{\oo}\tau_{*}(P).
	\end{equation}
	Le foncteur $\tau_{*}$ (resp. $\tau_{!}$) est un adjoint à droite (resp. gauche) de $\alpha_{A}$; les morphismes d'adjonction $\alpha_{A}\circ\tau_{*}\to \id$ et $\id \to \alpha_{A}\circ \tau_{!}$ sont des isomorphismes (cf. \cite{AG15} 1.4.13).
	\label{BMod unitaire}
\end{nothing}

\begin{nothing}
	Dans la suite de cette section, $\mathscr{C}$ désigne un site, $\widehat{\mathscr{C}}$ la catégorie des préfaisceaux d'ensembles sur $\mathscr{C}$ et $\widetilde{\mathscr{C}}$ le topos des faisceaux d'ensembles sur $\mathscr{C}$. On note $\oo_{\widehat{\mathscr{C}}}$ (resp. $\mm_{\widehat{\mathscr{C}}}$) le préfaisceau constant sur $\mathscr{C}$ de valeur $\oo$ (resp. $\mm$) et $\oo_{\widetilde{\mathscr{C}}}$ (resp. $\mm_{\widetilde{\mathscr{C}}}$) le faisceau associé. On désigne par $\Mod(\oo_{\widehat{\mathscr{C}}})$ (resp. $\Mod(\oo_{\widetilde{\mathscr{C}}})$) la catégorie abélienne tensorielle des ($\oo_{\widehat{\mathscr{C}}}$)-modules de $\widehat{\mathscr{C}}$ (resp. ($\oo_{\widetilde{\mathscr{C}}}$)-modules de $\widetilde{\mathscr{C}}$). Prenant pour $\oo\to \End(\oo_{\widehat{\mathscr{C}}})$ (resp. $\oo \to \End(\oo_{\widetilde{\mathscr{C}}})$) l'homomorphisme canonique, on appelle \textit{catégorie des $\alpha$-$\oo_{\widehat{\mathscr{C}}}$-modules} (resp. \textit{$\alpha$-$\oo_{\widetilde{\mathscr{C}}}$-modules}) et l'on note $\alpha$-$\Mod(\oo_{\widehat{\mathscr{C}}})$ (resp. $\alpha$-$\Mod(\oo_{\widetilde{\mathscr{C}}})$) le quotient de $\Mod(\oo_{\widehat{\mathscr{C}}})$ (resp. $\Mod(\oo_{\widetilde{\mathscr{C}}})$) par la sous-catégorie épaisse des modules $\alpha$-nuls.
\end{nothing}
\begin{nothing}	\label{faisceaux de alpha oo modules}
	On appelle \textit{catégorie des préfaisceaux de $\alpha$-$\oo$-modules sur $\mathscr{C}$} et l'on note $\alpha$-$\widehat{\mathscr{C}}$ la catégorie des préfaisceaux sur $\mathscr{C}$ à valeurs dans la catégorie $\alpha$-$\Mod(\oo)$, i.e., la catégorie des foncteurs de $\mathscr{C}^{\circ}$ à valeurs dans $\alpha$-$\Mod(\oo)$ (\cite{AG15} 1.4.16).

	On dit qu'un préfaisceau de $\alpha$-$\oo$-modules $F$ sur $\mathscr{C}$ est \textit{séparé} (resp. un \textit{faisceau}) si pour tout objet $X$ de $\mathscr{C}$ et tout crible couvrant $\mathscr{R}$ de $X$, le morphisme canonique
	\begin{equation}
		F(X)\to \varprojlim_{(Y,u)\in (\mathscr{C}_{/\mathscr{R}})^{\circ}} F(Y)
	\end{equation}
	est un monomorphisme (resp. isomorphisme) (cf. \cite{AG15} 1.4.21 et 1.4.24). On note $\alpha$-$\widetilde{\mathscr{C}}$ la sous-catégorie pleine de $\alpha$-$\widehat{\mathscr{C}}$ des faisceaux de $\alpha$-$\oo$-modules sur $\mathscr{C}$.

	Le foncteur $\alpha:\Mod(\oo)\to \alpha$-$\Mod(\oo)$ \eqref{foncteur alpha ModA} définit un foncteur exact et monoïdal
	\begin{equation}
		\widehat{\alpha}: \Mod(\oo_{\widehat{\mathscr{C}}}) \to \alpha\textnormal{-}\widehat{\mathscr{C}}.
		\label{alpha oo modules presheaf}
	\end{equation}
	D'après (\cite{AG15} 1.4.9(i)), celui-ci transforme les faisceaux de $\oo$-modules en des faisceaux de $\alpha$-$\oo$-modules. Il induit donc un foncteur
\begin{equation}
	\widetilde{\alpha}:\Mod(\oo_{\widetilde{\mathscr{C}}})\to \alpha\textnormal{-}\widetilde{\mathscr{C}}.
	\label{alpha oo modules sheaf}
\end{equation}
Ce dernier transforme les $\alpha$-isomorphismes en des isomorphismes (\cite{AG15} 1.4.23(iii)). Il induit une équivalence de catégories (\cite{AG15} 1.4.35)
\begin{equation}
	\alpha\textnormal{-}\Mod(\oo_{\widetilde{\mathscr{C}}})\xrightarrow{\sim}\alpha\textnormal{-}\widetilde{\mathscr{C}}.
	\label{equivalence alpha o mod et alpha faisceau}
\end{equation}
On munit $\alpha$-$\widetilde{\mathscr{C}}$ de la structure de catégorie abélienne tensorielle déduite de celle de $\alpha$-$\Mod(\oo_{\widetilde{\mathscr{C}}})$ \eqref{Cat Ab Ten almost} via l'équivalence \eqref{equivalence alpha o mod et alpha faisceau}.
	
	En vertu de (\cite{AG15} 1.4.31), le foncteur canonique d'inclusion $\iota:\alpha\textnormal{-}\widetilde{\mathscr{C}}\to \alpha\textnormal{-}\widehat{\mathscr{C}}$ admet un adjoint à gauche	
	\begin{equation}
		\overline{\textsf{a}}:\alpha\textnormal{-}\widehat{\mathscr{C}}\to \alpha\textnormal{-}\widetilde{\mathscr{C}},
		\label{faisceau associe almost}
	\end{equation}
	dit foncteur ``faisceau de $\alpha$-$\oo$-modules associé''. On désigne par $i:\Mod(\oo_{\widetilde{\mathscr{C}}})\to \Mod(\oo_{\widehat{\mathscr{C}}})$ le foncteur canonique et par $\textsf{a}:\Mod(\oo_{\widehat{\mathscr{C}}})\to \Mod(\oo_{\widetilde{\mathscr{C}}})$ le foncteur ``faisceau de $\oo$-modules associé''. On a alors des diagrammes commutatifs à isomorphismes près (\cite{AG15} 1.4.35, 1.4.36)
	\begin{equation}
		\xymatrix{
			\Mod(\oo_{\widehat{\mathscr{C}}})\ar[r]^{\textsf{a}} \ar[d]_{\widehat{\alpha}} & \Mod(\oo_{\widetilde{\mathscr{C}}})\ar[d]^{\widetilde{\alpha}} \\
			\alpha\textnormal{-}\widehat{\mathscr{C}}\ar[r]^{\overline{\textsf{a}}}& \alpha\textnormal{-}\widetilde{\mathscr{C}}.
		} \qquad
		\xymatrix{
			\Mod(\oo_{\widetilde{\mathscr{C}}})\ar[r]^{i} \ar[d]_{\widetilde{\alpha}} & \Mod(\oo_{\widehat{\mathscr{C}}})\ar[d]^{\widehat{\alpha}} \\
			\alpha\textnormal{-}\widetilde{\mathscr{C}}\ar[r]^{\iota}& \alpha\textnormal{-}\widehat{\mathscr{C}}.
		}
		\label{diagramme commutatif alpha associe}
	\end{equation}
\end{nothing}

\begin{nothing} \label{generalite o alpha module C}
Les foncteurs $\sigma_{*}$ \eqref{foncteur varsigma star} et $\sigma_{!}$ \eqref{foncteur varsigma !} (pour les $\oo$-modules) induisent des foncteurs que l'on note respectivement
\begin{eqnarray}
	\widehat{\sigma}_{*}:\alpha\textnormal{-}\widehat{\mathscr{C}} &\to & \Mod(\oo_{\widehat{\mathscr{C}}}), \label{fonteur varsigma prefaisceau star}\\
	\widehat{\sigma}_{!}:\alpha\textnormal{-}\widehat{\mathscr{C}} &\to & \Mod(\oo_{\widehat{\mathscr{C}}}). \label{fonteur varsigma prefaisceau !}
\end{eqnarray}
En vertu de (\cite{AG15} 1.4.9(i)), $\widehat{\sigma}_{*}$ transforme les faisceaux de $\alpha$-$\oo$-modules en des faisceaux de $\oo$-modules. Il définit donc un foncteur
\begin{eqnarray}
	\widetilde{\sigma}_{*}:\alpha\textnormal{-}\widetilde{\mathscr{C}} \to  \Mod(\oo_{\widetilde{\mathscr{C}}}). \label{fonteur varsigma faisceau star}
\end{eqnarray}
On désigne par $\widetilde{\sigma}_{!}$ le foncteur
\begin{equation}
	\widetilde{\sigma}_{!}:\alpha\textnormal{-}\widetilde{\mathscr{C}} \to \Mod(\oo_{\widetilde{\mathscr{C}}})\qquad M\mapsto \mm_{\widetilde{\mathscr{C}}}\otimes_{\oo_{\widetilde{\mathscr{C}}}}\widetilde{\sigma}_{*}(M).
	\label{fonteur varsigma faisceau !}
\end{equation}
On notera que $\widetilde{\sigma}_{!}$ s'identifie au composé du foncteur $\widehat{\sigma}_{!}$ \eqref{fonteur varsigma prefaisceau !} et du foncteur ``faisceau associé'' $\textsf{a}$. D'après \ref{propriete adjoint alpha}, le foncteur $\widehat{\sigma}_{*}$ (resp. $\widehat{\sigma}_{!}$) est un adjoint à droite (resp. gauche) de $\widehat{\alpha}$ (cf. \cite{AG15} 1.4.19); les morphismes d'adjonction $\widehat{\alpha}\circ \widehat{\sigma}_{*}\to \id$ et $\id \to \widehat{\alpha}\circ \widehat{\sigma}_{!}$ sont des isomorphismes. Le foncteur $\widetilde{\sigma}_{*}$ est un adjoint à droite de $\widetilde{\alpha}$ et le morphisme d'adjonction $\widetilde{\alpha}\circ\widetilde{\sigma}_{*}\to \id$ est un isomorphisme (\cite{AG15} 1.4.35). 
\end{nothing}

\begin{lemma}\label{lemma adjonction !}
	\textnormal{(i)} Le foncteur $\widetilde{\sigma}_{!}$ \eqref{fonteur varsigma faisceau !} est un adjoint à gauche de $\widetilde{\alpha}$ \eqref{alpha oo modules sheaf}.

	\textnormal{(ii)} Le morphisme d'adjonction $\id \to \widetilde{\alpha}\circ \widetilde{\sigma}_{!}$ est un isomorphisme.
\end{lemma}
\textit{Preuve}. Soient $M$ un faisceau de $\alpha$-$\oo$-modules sur $\mathscr{C}$ et $N$ un $\oo_{\widetilde{\mathscr{C}}}$-module. On a $\widetilde{\sigma}_{!}(M)=\textsf{a}(\widehat{\sigma}_{!}(\iota(M)))$. 

(i) Par adjonction et \eqref{diagramme commutatif alpha associe}, on a des isomorphismes
\begin{eqnarray*}
	\Hom_{\oo_{\widetilde{\mathscr{C}}}}(\textsf{a}(\widehat{\sigma}_{!}(\iota(M))),N)& \simeq & \Hom_{\oo_{\widehat{\mathscr{C}}}}(\widehat{\sigma}_{!}(\iota(M)),i(N)) \\
	&\simeq & \Hom_{\alpha\textnormal{-}\widehat{\mathscr{C}}}(\iota(M),\widehat{\alpha}(i(N))) \\
	&\simeq & \Hom_{\alpha\textnormal{-}\widetilde{\mathscr{C}}}(M,\widetilde{\alpha}(N)).
\end{eqnarray*}

(ii) Le morphisme d'adjonction $\iota(M)\to \widehat{\alpha}(\widehat{\sigma}_{!}(\iota(M)))$ est un isomorphisme. Compte tenu de \eqref{diagramme commutatif alpha associe}, on a des isomorphismes
\begin{equation}
	M\xrightarrow{\sim} \overline{\textsf{a}}(\widehat{\alpha}(\widehat{\sigma}_{!}(\iota(M))))\simeq \widetilde{\alpha}(\widetilde{\sigma}_{!}(M)).
\end{equation}

\begin{nothing}
	La donnée d'un monoïde commutatif unitaire de $\alpha$-$\widetilde{\mathscr{C}}$ est équivalente à la donnée d'un préfaisceau sur $\mathscr{C}$ à valeur dans $\alpha$-$\Alg(\oo)$ \eqref{alpha o algebra} dont le préfaisceau de $\alpha$-$\oo$-module sous-jacent est un faisceau (cf. \cite{AG15} 1.4.40). On note $\Alg(\alpha\textnormal{-}\widetilde{\mathscr{C}})$ la catégorie des monoïdes commutatifs unitaires de $\alpha$-$\widetilde{\mathscr{C}}$.

	De même, la donnée d'un monoïde commutatif unitaire de $\Mod(\oo_{\widetilde{\mathscr{C}}})$ est équivalente à la donnée d'une $(\oo_{\widetilde{\mathscr{C}}})$-algèbre de $\widetilde{\mathscr{C}}$. On note $\Alg(\oo_{\widetilde{\mathscr{C}}})$ la catégorie des monoïdes commutatifs unitaires de $\Mod(\oo_{\widetilde{\mathscr{C}}})$.

	Le foncteur $\widetilde{\alpha}$ \eqref{alpha oo modules sheaf} étant monoïdal, il induit un foncteur que l'on note encore
	\begin{equation}
		\widetilde{\alpha}:\Alg(\oo_{\widetilde{\mathscr{C}}})\to \Alg(\alpha\textnormal{-}\widetilde{\mathscr{C}}).		\label{alpha Alg sheaf}
	\end{equation}
	Compte tenu de \ref{Alg alpha sigma}, le foncteur $\widetilde{\sigma}_{*}$ \eqref{fonteur varsigma faisceau star} induit un foncteur que l'on note encore
	\begin{equation}
		\widetilde{\sigma}_{*}: \Alg(\alpha\textnormal{-}\widetilde{\mathscr{C}})\to \Alg(\oo_{\widetilde{\mathscr{C}}}).
		\label{sigma star Alg sheaf}
	\end{equation}
	D'après (\cite{AG15} 1.4.40), $\widetilde{\sigma}_{*}$ est un adjoint à droite de $\widetilde{\alpha}$; le morphisme d'adjonction $\widetilde{\alpha}\circ \widetilde{\sigma}_{*}\to \id$ est un isomorphisme et le morphisme d'adjonction $\id\to \widetilde{\sigma}_{*}\circ \widetilde{\alpha}$ induit un isomorphisme $\widetilde{\sigma}\xrightarrow{\sim}\widetilde{\alpha}\circ \widetilde{\sigma}_{*}\circ \widetilde{\alpha}$.
	\label{sigma alpha Alg sheaf}
\end{nothing}

\begin{nothing}
	Dans la suite de cette section, on fixe une $\oo_{\widetilde{\mathscr{C}}}$-algèbre $A$ de $\widetilde{\mathscr{C}}$. On désigne par $\Mod(A)$ la catégorie abélienne tensorielle des $A$-modules de $\widetilde{\mathscr{C}}$. Prenant pour $\varphi:\oo\to \Gamma(\widetilde{\mathscr{C}},A)$ l'homomorphisme canonique, on appelle \textit{catégorie des $\alpha$-$A$-modules} et l'on note $\alpha$-$\Mod(A)$ le quotient de la catégorie abélienne $\Mod(A)$ par la sous-catégorie épaisse des $A$-modules $\alpha$-nuls.
	
	On pose $A^{\alpha}=\widetilde{\alpha}(A)$ \eqref{alpha Alg sheaf} et on désigne par $\Mod(A^{\alpha})$ la catégorie des $A^{\alpha}$-modules unitaires de $\alpha$-$\widetilde{\mathscr{C}}$. La donnée d'une structure de $A^{\alpha}$-module unitaire sur un faisceau de $\alpha$-$\oo$-modules $M$ sur $\mathscr{C}$ est équivalente à la donnée pour tout $X\in \Ob(\mathscr{C})$ d'une structure de $\alpha(A(X))$-module unitaire sur $M(X)$ dans le sens de \ref{BMod unitaire} telle que pour tout morphisme $X\to Y$ de $\mathscr{C}$, le morphisme $M(Y)\to M(X)$ soit linéaire relativement au morphisme de $\oo^{\alpha}$-algèbres $\alpha(A(Y))\to \alpha(A(X))$, où $\alpha:\Alg(\oo)\to \alpha\textnormal{-}\Alg(\oo)$ est le foncteur \eqref{alpha o algebra} (cf. \cite{AG15} 1.4.41).

	Le foncteur $\widetilde{\alpha}$ \eqref{alpha oo modules sheaf} étant monoïdal, il définit un foncteur
	\begin{equation}
		\widetilde{\alpha}_{A}:\Mod(A)\to \Mod(A^{\alpha}).
		\label{alpha A sheaf}
	\end{equation}
	Celui-ci transforme les $\alpha$-isomorphismes en des isomorphismes (cf. \cite{AG15} 1.4.23(iii)). Il induit une équivalence de catégories (\cite{AG15} 1.4.41)
\begin{equation}
	\alpha\textnormal{-}\Mod(A) \xrightarrow{\sim} \Mod(A^{\alpha}).
\end{equation}
Dans la suite de cette section, pour tout $A$-module $M$, on pose abusivement $M^{\alpha}=\widetilde{\alpha}_{A}(M)$.
\label{alpha A Mod almost}
\end{nothing}
\begin{nothing}
	On pose $A'=\widetilde{\sigma}_{*}(A^{\alpha})$ \eqref{sigma star Alg sheaf}. D'après \ref{sigma alpha Alg sheaf}, on a un homomorphisme canonique de $\oo_{\widetilde{\mathscr{C}}}$-algèbres $\lambda:A\to A'$, qui induit un isomorphisme $\widetilde{\alpha}(A)\xrightarrow{\sim} \widetilde{\alpha}(A')$. Compte tenu de \ref{BMod unitaire}, le foncteur $\widetilde{\sigma}_{*}$ \eqref{fonteur varsigma faisceau star} induit un foncteur (cf. \cite{AG15} 1.4.41)
\begin{equation}
	\widetilde{\tau}_{*}': \Mod(A^{\alpha})\to \Mod(A').
\end{equation}
Composant avec le foncteur induit par $\lambda$, on obtient un foncteur
\begin{equation}
	\widetilde{\tau}_{*}: \Mod(A^{\alpha})\to \Mod(A).
	\label{tau * foncteur}
\end{equation}
On désigne par $\widetilde{\tau}_{!}$ le foncteur
\begin{equation}
	\widetilde{\tau}_{!}: \Mod(A^{\alpha})\to \Mod(A),\qquad M\mapsto \mm_{\widetilde{\mathscr{C}}}\otimes_{\oo_{\widetilde{\mathscr{C}}}}\widetilde{\tau}_{*}(M).
	\label{tau ! foncteur}
\end{equation}
	Pour tous $A^{\alpha}$-modules $M$ et $N$, on pose
	\begin{equation}
		\Hom_{A^{\alpha}}(M,N)=\Hom_{\Mod(A^{\alpha})}(M,N),
	\end{equation}
	qui est naturellement muni d'une structure de $\Gamma(\widetilde{\mathscr{C}},A)$-module, en particulier, d'une structure de $\oo$-module.
\end{nothing}

\begin{nothing}
	Soit $A\to B$ un morphisme de $\oo_{\widetilde{\mathscr{C}}}$-algèbres. On note $B^{\alpha}=\widetilde{\alpha}(B)$ \eqref{alpha Alg sheaf} et on désigne par $\Mod(B^{\alpha})$ la catégorie des $B^{\alpha}$-modules de $\alpha\textnormal{-}\widetilde{\mathscr{C}}$. On a des diagrammes commutatifs
	\begin{equation}
		\xymatrix{
			\Mod(B^{\alpha}) \ar[d] \ar[r]^{\widetilde{\tau}_{B!}}& \Mod(B) \ar[d]\\
			\Mod(A^{\alpha}) \ar[r]^{\widetilde{\tau}_{A!}}& \Mod(A), 			
		} \qquad
		\xymatrix{
			\Mod(B^{\alpha}) \ar[d] \ar[r]^{\widetilde{\tau}_{B*}}& \Mod(B) \ar[d]\\
			\Mod(A^{\alpha}) \ar[r]^{\widetilde{\tau}_{A*}}& \Mod(A), 			
		}
	\end{equation}
	où les flèches verticales sont des foncteurs d'oubli et $\widetilde{\tau}_{B*}$, $\widetilde{\tau}_{B!}$, $\widetilde{\tau}_{A*}$ et $\widetilde{\tau}_{A!}$ désignent les foncteurs \eqref{tau * foncteur} et \eqref{tau ! foncteur} relatifs à $B$ et $A$ respectivement. On peut donc omettre les indices $B$ et $A$ dans les notations de $\widetilde{\tau}_{B*}$, $\widetilde{\tau}_{B!}$, $\widetilde{\tau}_{A*}$ et $\widetilde{\tau}_{A!}$.
\end{nothing}

\begin{lemma}[\cite{AG15} 1.4.23] \label{lemma 1423}
	\textnormal{(i)} Pour qu'un morphisme de $A$-modules $f:M\to N$ soit un $\alpha$-isomorphisme, il faut et il suffit que le morphisme $\mm_{\widetilde{\mathscr{C}}}\otimes_{\oo_{\widetilde{\mathscr{C}}}}M\to \mm_{\widetilde{\mathscr{C}}}\otimes_{\oo_{\widetilde{\mathscr{C}}}}N$ induit par $f$ soit un isomorphisme.

	\textnormal{(ii)} Pour tous $A$-modules $M$ et $N$, on a un isomorphisme $\oo$-linéaire canonique
	\begin{equation}
		\Hom_{A^{\alpha}}(M^{\alpha},N^{\alpha})\xrightarrow{\sim} \Hom_{A}(\mm_{\widetilde{\mathscr{C}}}\otimes_{\oo_{\widetilde{\mathscr{C}}}}M,N).
		\label{almost hom vrai hom}
	\end{equation}
\end{lemma}

\begin{prop}
	\textnormal{(i)} Le foncteur $\widetilde{\tau}_{*}$ \eqref{tau * foncteur} est un adjoint à droite du foncteur $\widetilde{\alpha}_{A}$ \eqref{alpha A sheaf}.

	\textnormal{(ii)} Le morphisme d'adjonction $\widetilde{\alpha}_{A}\circ \widetilde{\tau}_{*}\to \id$ est un isomorphisme.

	\textnormal{(iii)} Le foncteur $\widetilde{\tau}_{!}$ \eqref{tau ! foncteur} est un adjoint à gauche du foncteur $\widetilde{\alpha}_{A}$.

	\textnormal{(iv)} Le morphisme d'adjonction $\id\to \widetilde{\alpha}_{A}\circ \widetilde{\tau}_{!}$ est un isomorphisme.
	\label{propriete adjoint alpha AMod}
\end{prop}

\textit{Preuve}. Les assertions (i) et (ii) sont démontrées dans (\cite{AG15} 1.4.41).

(iii) Soient $M$ un $A$-module et $N$ un $A^{\alpha}$-module. D'après (ii) et \eqref{almost hom vrai hom}, on a des isomorphismes
\begin{eqnarray*}
	\Hom_{A^{\alpha}}(N,M^{\alpha})& \xrightarrow{\sim} & \Hom_{A^{\alpha}}( (\widetilde{\tau}_{*}(N))^{\alpha},M^{\alpha} )\\
	&\xrightarrow{\sim} & \Hom_{A}( \mm_{\widetilde{\mathscr{C}}}\otimes_{\oo_{\widetilde{\mathscr{C}}}}\widetilde{\tau}_{*}(N),M)\\
	&\xrightarrow{\sim} & \Hom_{A}( \widetilde{\tau}_{!}(N),M).
\end{eqnarray*}

(iv) On notera que, pour tout $A^{\alpha}$-module $M$, le $\alpha$-$\oo$-module sous-jacent à $(\widetilde{\tau}_{!}(M))^{\alpha}$ s'identifie canoniquement à $\widetilde{\alpha}(\widetilde{\sigma}_{!}(M))$. L'assertion résulte du fait que le morphisme d'adjonction $\id\to \widetilde{\alpha}\circ \widetilde{\sigma}_{!}$ est un isomorphisme (cf. \ref{lemma adjonction !}(ii)).

\begin{coro}
	\textnormal{(i)} Les foncteurs $\widetilde{\alpha}_{A}$, $\widetilde{\tau}_{!}$ sont exacts et le foncteur $\widetilde{\tau}_{*}$ est exact à gauche.

	\textnormal{(ii)} Les foncteurs $\widetilde{\alpha}_{A}$ et $\widetilde{\tau}_{*}$ transforment les objets injectifs en des objets injectifs. En particulier, $\Mod(A^{\alpha})$ a suffisamment d'injectifs.
		
	\textnormal{(iii)} Le morphisme d'adjonction $\id \to \widetilde{\tau}_{*}\circ \widetilde{\alpha}_{A}$ (resp. $\widetilde{\tau}_{!}\circ\widetilde{\alpha}_{A}\to \id$) induit un isomorphisme $\widetilde{\alpha}\xrightarrow{\sim} \widetilde{\alpha}_{A}\circ \widetilde{\tau}_{*}\circ \widetilde{\alpha}_{A}$ (resp. $\widetilde{\alpha}_{A}\circ \widetilde{\tau}_{!}\circ \widetilde{\alpha}_{A}\xrightarrow{\sim} \widetilde{\alpha}_{A}$).
\label{coro inj exactitude}
\end{coro}

\textit{Preuve}. (i) Compte tenu de \ref{propriete adjoint alpha AMod}, $\widetilde{\alpha}_{A}$ est exact et $\widetilde{\tau}_{*}$ (resp. $\widetilde{\tau}_{!}$) est exact à gauche (resp. droit). Le $\oo_{\widetilde{\mathscr{C}}}$-module $\mm_{\widetilde{\mathscr{C}}}$ étant plat (\cite{SGAIV} V 1.7.1), l'exactitude de $\widetilde{\tau}_{!}$ s'ensuit compte tenu \eqref{tau ! foncteur}.

(ii) Le foncteur $\widetilde{\alpha}_{A}$ (resp. $\widetilde{\tau}_{*}$) est un adjoint à droite d'un foncteur exact, d'où l'assertion (\cite{SGAIV} V 0.2).

(iii) En effet, le morphisme composé
\begin{displaymath}
	\widetilde{\alpha}_{A}\to \widetilde{\alpha}_{A}\circ \widetilde{\tau}_{*}\circ \widetilde{\alpha}_{A}\to \widetilde{\alpha}_{A},
\end{displaymath}
où les flèches sont déduites des morphismes d'adjonction, s'identifie au morphisme identique. L'assertion pour $\widetilde{\tau}_{*}$ s'ensuit compte tenu de \ref{propriete adjoint alpha AMod}(ii). La démonstration de l'assertion pour $\widetilde{\tau}_{!}$ est similaire.

\begin{rem}
	Le foncteur $\widetilde{\alpha}_{A}$ ne transforme pas les objets projectifs en des objets projectifs. En effet, considérons le topos ponctuel et prenons pour $A$ l'anneau $\oo$. En vertu de \eqref{Hom oo mm M}, pour tout $\oo$-module $M$, on a un isomorphisme canonique
	\begin{displaymath}
		\Hom_{\oo^{\alpha}}(\oo^{\alpha},M^{\alpha})\xrightarrow{\sim}\Hom_{\oo}(\mm,M).
	\end{displaymath}
	Compte tenu de \ref{m pas projectif}, $\oo^{\alpha}$ n'est pas projectif dans la catégorie $\alpha\textnormal{-}\Mod(\oo)$.
\end{rem}

\begin{nothing}
	On désigne par $\mathbf{D}(\Mod(A))$ (resp. $\mathbf{D}(\Mod(A^{\alpha}))$) la catégorie dérivée de $\Mod(A)$ (resp. $\Mod(A^{\alpha})$). Les foncteurs exacts $\widetilde{\alpha}_{A}$ et $\widetilde{\tau}_{!}$ induisent des foncteurs exacts entre catégories dérivées
\begin{eqnarray}
	F: \mathbf{D}(\Mod(A))\to \mathbf{D}(\Mod(A^{\alpha})),\qquad G:\mathbf{D}(\Mod(A^{\alpha}))\to \mathbf{D}(\Mod(A)).
	\label{foncteurs F G}
\end{eqnarray}
\end{nothing}
\begin{prop}
	Le foncteur $G$ est un adjoint à gauche du foncteur $F$.
	\label{adjoint F and G}
\end{prop}
\textit{Preuve}. D'après (\cite{GR14} 1.1.12), on a des identités triangulaires
\begin{equation}
	\xymatrix{
		\widetilde{\alpha}_{A} \ar@{=}[rd]_{\id_{\widetilde{\alpha}_{A}}} \ar@{=>}[r]^-{\eta * \widetilde{\alpha}_{A}}& \widetilde{\alpha}_{A}\widetilde{\tau}_{!}\widetilde{\alpha}_{A} \ar@{=>}[d]^{\widetilde{\alpha}_{A} * \varepsilon} \\
		& \widetilde{\alpha}_{A}
	}
\qquad	\xymatrix{
	\widetilde{\tau}_{!} \ar@{=}[rd]_{\id_{\widetilde{\tau}_{!}}} \ar@{=>}[r]^-{\widetilde{\tau}_{!} * \eta}& \widetilde{\tau}_{!}\widetilde{\alpha}_{A}\widetilde{\tau}_{!} \ar@{=>}[d]^{\varepsilon * \widetilde{\tau}_{!}} \\
		& \widetilde{\tau}_{!}
	}
\end{equation}
où $\eta$ et $\varepsilon$ désignent les transformations naturelles d'adjonction. Comme $\widetilde{\alpha}_{A}$ et $\widetilde{\tau}_{!}$ sont exacts, on en déduit des identités triangulaires similaires pour la paire $(G,F)$. La proposition s'ensuit compte tenu de (\cite{GR14} 1.1.14).

\begin{nothing}
	Rappelons que la catégorie abélienne $\Mod(A^{\alpha})$ a suffisamment d'injectifs (\ref{coro inj exactitude}(ii)) et reprenons les définitions de \ref{definition Ext using Hom derive} pour cette catégorie. Pour tous $A^{\alpha}$-modules $M$, $N$ et tout entier $i\ge 0$, on pose $\Ext^{i}_{A^{\alpha}}(M,N)=\Ext^{i}_{\Mod(A^{\alpha})}(M,N)$ qui est muni d'une structure de $\oo$-module.

	D'après \eqref{def Ext general} et \ref{adjoint F and G}, pour tout $A^{\alpha}$-module $M$, tout $A$-module $N$ et tout entier $i\ge 0$, on a un isomorphisme canonique
	\begin{equation}
		\Ext_{A}^{i}(\widetilde{\tau}_{!}(M),N)\xrightarrow{\sim} \Ext_{A^{\alpha}}^{i}(M,N^{\alpha}).
		\label{Ext alpha Ext isomorphisme}
	\end{equation}
	
	Soient $M$ et $N$ deux $A$-modules. Le foncteur exact $\widetilde{\alpha}_{A}$ induit un morphisme canonique \eqref{exact foncteur induit Ext morphisme}
	\begin{equation}
		\Ext^{i}_{A}(M,N)\to \Ext^{i}_{A^{\alpha}}(M^{\alpha},N^{\alpha}).
		\label{morphisme ext to alpha ext}
	\end{equation}
	Composant avec l'isomorphisme $\Ext^{i}_{A^{\alpha}}(M^{\alpha},N^{\alpha})\xrightarrow{\sim}\Ext^{i}_{A}(\widetilde{\tau}_{!}(M^{\alpha}),N)$, on obtient un morphisme canonique
	\begin{equation}
		\Ext^{i}_{A}(M,N)\to \Ext^{i}_{A}(\widetilde{\tau}_{!}(M^{\alpha}),N).
		\label{Ext induit par adjonction}
	\end{equation}
	On notera que celui-ci est induit par le morphisme d'adjonction $G\circ F\to \id$. En particulier, il s'identifie au morphisme canonique induit par la flèche d'adjonction $\widetilde{\tau}_{!}(M^{\alpha})\to M$.
	\label{almost Ext}
\end{nothing}

\begin{lemma}
	Soient $M$ et $N$ deux $A$-modules de $\widetilde{\mathscr{C}}$. Le morphisme canonique \eqref{morphisme ext to alpha ext} est un $\alpha$-isomorphisme de $\oo$-modules.
	\label{coro alpha iso ext}
\end{lemma}
\textit{Preuve}. Il suffit de démontrer que le morphisme \eqref{Ext induit par adjonction} est un $\alpha$-isomorphisme. En vertu de \ref{coro inj exactitude}(iii), le morphisme canonique $f:\widetilde{\tau}_{!}(M^{\alpha})\to M$ est un $\alpha$-isomorphisme de $\Mod(A)$. D'après \ref{lemma 143}, pour tout $\gamma\in \mm$, il existe un morphisme $A$-linéaire $g_{\gamma}:M\to \widetilde{\tau}_{!}(M^{\alpha})$ tel que $g_{\gamma}\circ f=\gamma^{2}\id_{M}$ et $f\circ g_{\gamma}=\gamma^{2}\id_{\widetilde{\tau}_{!}(M^{\alpha})}$. On en déduit que le noyau et le conoyau de \eqref{Ext induit par adjonction} sont annulés par $\gamma^{2}$, d'où le lemme.

\begin{coro}
	Le foncteur $\widetilde{\alpha}_{A}$ \eqref{alpha A sheaf} induit une équivalence de catégories 
	\begin{equation}
		\Mod_{\mathbb{Q}}(A) \xrightarrow{\sim} \Mod_{\mathbb{Q}}(A^{\alpha}),
	\end{equation}
	où la source (resp. le but) désigne la catégorie des $A$-modules (resp. $A^{\alpha}$-modules) à isogénie près.
	\label{Mod A Q alpha Mod A Q}
\end{coro}
\textit{Preuve}. Par définition, le foncteur est essentiellement surjectif. La pleine fidélité résulte de \ref{coro alpha iso ext}.\\

$\hspace*{-1.2em}\bf{\arabic{section}.\stepcounter{theorem}\arabic{theorem}.}$	
Soient $M$ un $A^{\alpha}$-module et $N$ un $A$-module. En vertu de \ref{propriete adjoint alpha AMod}(iv), on a des isomorphismes canoniques et fonctoriels de $A^{\alpha}$-modules
	\begin{equation}
		(\widetilde{\tau}_{!}(M\otimes_{A^{\alpha}}N^{\alpha}))^{\alpha}\xrightarrow{\sim} M\otimes_{A^{\alpha}}N^{\alpha}, \qquad (\widetilde{\tau}_{!}(M))^{\alpha}\otimes_{A^{\alpha}}N^{\alpha}\xrightarrow{\sim} M\otimes_{A^{\alpha}}N^{\alpha}.
	\end{equation}
	On en déduit par \ref{lemma 1423}(ii) et \ref{propriete adjoint alpha AMod}(iii) un isomorphisme canonique et fonctoriel de $A$-modules
	\begin{equation}
		\mm_{\widetilde{\mathscr{C}}}\otimes_{\oo_{\widetilde{\mathscr{C}}}} \widetilde{\tau}_{!}(M\otimes_{A^{\alpha}}N^{\alpha})\xrightarrow{\sim} \mm_{\widetilde{\mathscr{C}}}\otimes_{\oo_{\widetilde{\mathscr{C}}}} \widetilde{\tau}_{!}(M)\otimes_{A}N
	\end{equation}
	Compte tenu de l'isomorphisme $\mm\otimes_{\oo}\mm\simeq \mm$ et de la définition de $\widetilde{\tau}_{!}$, celui-ci s'identifie à l'isomorphisme
	\begin{equation}
		\widetilde{\tau}_{!}(M\otimes_{A^{\alpha}}N^{\alpha})\xrightarrow{\sim} \widetilde{\tau}_{!}(M)\otimes_{A}N
		\label{isomorphisme tau ! tensor}
	\end{equation}

\begin{prop} \label{lemma suite spectrale Gabber}
	Soit $n$ un entier $\ge 1$ et supposons que $p^{n}A=0$. Soit $M$ un $A$-module plat sur $\oo_{n}$ \textnormal{(\cite{SGAIV} V 1.7)}.

	\textnormal{(i)} Pour tout $\oo_{n}$-module $P$ et tout entier $i\ge 0$, on a un isomorphisme canonique
	\begin{equation}
		\Ext_{\oo_{n}}^{i}(\mm/p^{n}\mm,P)\xrightarrow{\sim} \Ext_{\oo}^{i}(\mm,P).		
	\end{equation}

	\textnormal{(ii)} Le foncteur $\Hom_{A}(M,-):\Mod(A)\to \Mod(\oo_{n})$ transforme les objets injectifs en objets injectifs.

	\textnormal{(iii)} Soient $N_{1}$ et $N_{2}$ deux $A$-modules. On a un isomorphisme canonique et fonctoriel
		\begin{equation}
			\Hom_{A}(\mm_{\widetilde{\mathscr{C}}}\otimes_{\oo_{\widetilde{\mathscr{C}}}}N_{1},N_{2})\xrightarrow{\sim} \Hom_{\oo}(\mm,\Hom_{A}(N_{1},N_{2})).
		\label{Hom Aalpha m Hom}
	\end{equation}
	
	\textnormal{(iv)} Pour tout $A$-module $N$, on a une suite spectrale
	\begin{equation}
		\rE_{2}^{i,j}=\Ext^{i}_{\oo}(\mm, \Ext^{j}_{A}(M,N))\Rightarrow \rE^{i+j}=\Ext^{i+j}_{A}(\mm_{\widetilde{\mathscr{C}}}\otimes_{\oo_{\widetilde{\mathscr{C}}}}M,N).
		\label{suite spectrale Gabber}
	\end{equation}
\end{prop}

\textit{Preuve}. (i) Le foncteur $\Hom_{\oo}(-,P):\Mod(\oo)\to \Mod(\oo)$ s'identifie au foncteur composé
\begin{equation}
	\Mod(\oo)\xrightarrow{-\otimes_{\oo}\oo_{n}} \Mod(\oo_{n})\xrightarrow{\Hom_{\oo_{n}}(-,P)} \Mod(\oo).
\end{equation}
Comme $-\otimes_{\oo}\oo_{n}$ envoie les objets projectifs sur des objets projectifs, on en déduit une suite spectrale
	\begin{equation}
		\Ext_{\oo_{n}}^{i}(\Tor_{j}^{\oo}(\mm,\oo_{n}),P)\Rightarrow \Ext_{\oo}^{i+j}(\mm,P).
		\label{Ext oo Ext oon}
	\end{equation}
L'assertion résulte alors de la platitude de $\mm$.

(ii) Posons $\oo_{\widetilde{\mathscr{C}},n}=\oo_{\widetilde{\mathscr{C}}}/p^{n}\oo_{\widetilde{\mathscr{C}}}$. Soient $P$ un $\oo_{n}$-module et $P_{\widetilde{\mathscr{C}}}$ le faisceau associé au préfaisceau constant de valeur $P$ sur $\mathscr{C}$. Pour tout $A$-module $N$, par adjonction (\cite{SGAIV} IV 13.4.1), on a un isomorphisme canonique
\begin{equation}
	\Hom_{\oo_{\widetilde{\mathscr{C}},n}}(P_{\widetilde{\mathscr{C}}},N)\xrightarrow{\sim} \Hom_{\oo_{n}}(P,\Gamma(\widetilde{\mathscr{C}},N)).
	\label{Hom adjoint topos final}
\end{equation}
Soit $I$ un $A$-module injectif de $\widetilde{\mathscr{C}}$. Pour tout $\oo_{n}$-module $P$, on en déduit un isomorphisme canonique
\begin{equation}
	\Hom_{\oo_{\widetilde{\mathscr{C}},n}}(P_{\widetilde{\mathscr{C}}},\FHom_{A}(M,I))\xrightarrow{\sim}\Hom_{\oo_{n}}(P,\Hom_{A}(M,I))
\end{equation}
Par ailleurs, on a un isomorphisme canonique (\cite{SGAIV} IV 12.14)
\begin{equation}
	\Hom_{A}(P_{\widetilde{\mathscr{C}}}\otimes_{\oo_{\widetilde{\mathscr{C}},n}}M,I) \xrightarrow{\sim} \Hom_{\oo_{\widetilde{\mathscr{C}},n}}(P_{\widetilde{\mathscr{C}}},\FHom_{A}(M,I)).
\end{equation}
L'injectivité de $\Hom_{A}(M,I)$ résulte de celle de $I$ et de la platitude de $M$.

(iii) En vertu de (\cite{SGAIV} IV 12.14), on a un isomorphisme canonique et fonctoriel 
	\begin{eqnarray}
		\Hom_{A}(\mm_{\widetilde{\mathscr{C}}}\otimes_{\oo_{\widetilde{\mathscr{C}}}}N_{1},N_{2})\xrightarrow{\sim} \Hom_{\oo_{\widetilde{\mathscr{C}}}}(\mm_{\widetilde{\mathscr{C}}},\FHom_{A}(N_{1},N_{2})).
		\label{Hom o m Cartan}
	\end{eqnarray}
	Par adjonction (\cite{SGAIV} IV 13.4.1), on a un isomorphisme canonique et fonctoriel
	\begin{equation}
		\Hom_{\oo_{\widetilde{\mathscr{C}}}}(\mm_{\widetilde{\mathscr{C}}},\FHom_{A}(N_{1},N_{2}))\xrightarrow{\sim} \Hom_{\oo}(\mm,\Hom_{A}(N_{1},N_{2})).		
	\end{equation}

	(iv) D'après (iii), le foncteur
	\begin{equation}
		\Hom_{A}(\mm_{\widetilde{\mathscr{C}}}\otimes_{\oo_{\widetilde{\mathscr{C}}}}M,-): \Mod(A)\to \Mod(\oo_{n}),
	\end{equation}
	s'identifie au foncteur composé
	\begin{equation}
	\xymatrixcolsep{5pc}\xymatrix{
		\Mod(A) \ar[r]^-{\Hom_{A}(M,-)}& \Mod(\oo_{n}) \ar[r]^-{\Hom_{\oo_{n}}(\mm/p^{n}\mm,-)} & \Mod(\oo_{n}).
	}
	\end{equation}
	Comme $\Hom_{A}(M,-)$ envoie les objets injectifs sur des objets injectifs, on en déduit une suite spectrale
	\begin{equation}
		\Ext^{i}_{\oo_{n}}(\mm/p^{n}\mm, \Ext^{j}_{A}(M,N))\Rightarrow \Ext^{i+j}_{A}(\mm_{\widetilde{\mathscr{C}}}\otimes_{\oo_{\widetilde{\mathscr{C}}}}M,N).
		\label{suite spectrale Gabber faible}
	\end{equation}
	La suite spectrale \eqref{suite spectrale Gabber} s'ensuit compte tenu de (i).

	\begin{nothing} \label{morphisme de topos pullback alpha}
	Soit $\varphi:(\widetilde{\mathscr{C}}',A')\to (\widetilde{\mathscr{C}},A)$ un morphisme de topos annelés. On note $\alpha\textnormal{-}\widetilde{\mathscr{C}}'$ la catégorie des faisceaux de $\alpha$-$\oo$-modules de $\widetilde{\mathscr{C}}'$ \eqref{faisceaux de alpha oo modules} et on pose $A'^{\alpha}=\widetilde{\alpha}(A')$ \eqref{alpha Alg sheaf}. On désigne par $\Mod(A'^{\alpha})$ la catégorie des $A'^{\alpha}$-modules de $\alpha\textnormal{-}\widetilde{\mathscr{C}}'$. En vertu de \ref{lemma 1423}(i), le foncteur $\varphi^{*}:\Mod(A)\to \Mod(A')$ envoie les $\alpha$-isomorphismes sur des $\alpha$-isomorphismes. Il induit donc un foncteur 
		\begin{equation} \label{psi phi *}
		\psi: \Mod(A^{\alpha})\to \Mod(A'^{\alpha})
	\end{equation}
	qui s'insère dans un diagramme commutatif
	\begin{equation}
		\xymatrix{
			\Mod(A)\ar[r]^-{\varphi^{*}} \ar[d]_{\widetilde{\alpha}_{A}} & \Mod(A')\ar[d]^{\widetilde{\alpha}_{A'}}\\
			\Mod(A^{\alpha}) \ar[r]^-{\psi} & \Mod(A'^{\alpha}).
		}
		\label{diagramme commutatif alpha phi psi}
	\end{equation}
	Compte tenu de \ref{propriete adjoint alpha AMod}(iv), pour tout $A^{\alpha}$-module $M$, on a un isomorphisme canonique et fonctoriel
	\begin{equation}
		\psi(M)\xrightarrow{\sim} (\varphi^{*}(\widetilde{\tau}_{A!}(M)))^{\alpha}.
		\label{isomorphisme phi psi alpha}
	\end{equation}
	Le morphisme adjoint de cet isomorphisme
	\begin{equation}
		\widetilde{\tau}_{A'!}(\psi(M)) \rightarrow \varphi^{*}(\widetilde{\tau}_{A!}(M))
		\label{diagramme tau ! B B'}
	\end{equation}
	est un $\alpha$-isomorphisme en vertu de \ref{propriete adjoint alpha AMod}(iii). C'est même un isomorphisme en vertu de \ref{lemma 1423}(ii) et de l'isomorphisme $\mm\otimes_{\oo}\mm\simeq \mm$.
\end{nothing}

\begin{nothing} \label{prop morphisme plat topos}
	Conservons les notations de \ref{morphisme de topos pullback alpha} et supposons que le morphisme $\varphi$ soit plat, i.e. que le foncteur $\varphi^{*}:\Mod(A)\to \Mod(A')$ soit \textit{exact} (\cite{SGAIV} V 1.8). D'après (\cite{Ga62} III.1 Corollaire 3), le foncteur $\psi$ \eqref{psi phi *} est exact. En vertu de \eqref{diagramme commutatif alpha phi psi} et \eqref{diagramme tau ! B B'}, pour tout $A^{\alpha}$-module $M$, tout $A$-module $N$ et tout entier $i\ge 0$, les isomorphismes \eqref{Ext alpha Ext isomorphisme} induisent un diagramme commutatif	
	\begin{equation}
		\xymatrix{
			\Ext^{i}_{A^{\alpha}}(M,N^{\alpha}) \ar[r]^-{\psi} \ar[d]_{\wr} & \Ext^{i}_{A'^{\alpha}} (\psi(M),\psi(N^{\alpha})) \ar[d]^{\wr} \\
			\Ext^{i}_{A}(\widetilde{\tau}_{A!}(M),N) \ar[r]^-{\varphi^{*}} & \Ext^{i}_{A'}(\widetilde{\tau}_{A'!}(\psi(M)),\varphi^{*}(N)),
		}
		\label{diagramme commutatif Ext phi tau alpha}
	\end{equation}
	où les flèches horizontals sont les morphismes canoniques induits par les foncteurs exacts $\varphi^{*}$ et $\psi$  \eqref{exact foncteur induit Ext morphisme}.
\end{nothing}

\begin{nothing}	\label{fonctorialite suite spectrale Gabber}
	Conservons les notations et hypothèses de \ref{prop morphisme plat topos} et supposons de plus que $p^{n}A=0$. Soient $M$ un $A$-module qui est plat sur $\oo_{n}$ et $N$ un $A$-module. Comme $\varphi^{*}$ est exact, le $A'$-module $\varphi^{*}(M)$ est encore plat sur $\oo_{n}$. On désigne par $\rE$ (resp. $\rE'$) la suite spectrale \eqref{suite spectrale Gabber} associée aux $A$-modules $M$ et $N$ (resp. $A'$-modules $\varphi^{*}(M)$ et $\varphi^{*}(N)$). Le foncteur exact $\varphi^{*}$ induit un morphisme de suites spectrales (cf. \ref{explication ss Gabber}) 
	\begin{equation}
		u=(u_{r}^{i,j},u^{n}):\rE\to \rE'
		\label{morphisme suite spectrales Gabber}
	\end{equation}
	dont les morphismes
	\begin{eqnarray}
		&u_{2}^{i,j}&: \Ext^{i}_{\oo}(\mm, \Ext^{j}_{A}(M,N))\to \Ext^{i}_{\oo}(\mm, \Ext^{j}_{A'}(\varphi^{*}(M),\varphi^{*}(N)))\\
		&u^{n}&: \Ext^{n}_{A}(\mm_{\widetilde{\mathscr{C}}}\otimes_{\oo_{\widetilde{\mathscr{C}}}}M,N) \to \Ext^{n}_{A'}(\mm_{\widetilde{\mathscr{C}'}}\otimes_{\oo_{\widetilde{\mathscr{C}}'}}\varphi^{*}(M),\varphi^{*}(N))
	\end{eqnarray}
	s'identifient aux morphismes canoniques induits par $\varphi^{*}$ \eqref{exact foncteur induit Ext morphisme}. Compte tenu de l'isomorphisme canonique $\mm_{\widetilde{\mathscr{C}}}\otimes_{\oo_{\widetilde{\mathscr{C}}}}M\simeq \widetilde{\tau}_{!}(M^{\alpha})$ et \eqref{diagramme commutatif Ext phi tau alpha}, le morphisme $u^{n}$ s'identifie au morphisme canonique induit par $\psi$:
	\begin{equation}
		\Ext^{n}_{A^{\alpha}}(M^{\alpha},N^{\alpha})\to \Ext^{n}_{A'^{\alpha}}(\psi(M^{\alpha}),\psi(N^{\alpha})).
		\label{un morphisme Ext alpha}
	\end{equation}
\end{nothing}

\begin{definition}[\cite{AG15} 1.5.3] \label{type alpha fini topos}
	Soient $M$ un $A$-module de $\widetilde{\mathscr{C}}$ et $\gamma\in \mm$.

	\textnormal{(i)} On dit que $M$ est de \textit{type $\gamma$-fini} s'il existe un raffinement $(U_{i})_{i\in I}$ de l'objet final de $\widetilde{\mathscr{C}}$ tel que, pour tout $i\in I$, il existe un $(A|U_{i})$-module libre de type fini $N_{i}$ et un morphisme $(A|U_{i})$-linéaire
	\begin{equation}
		N_{i}\to M|U_{i}
	\end{equation}
	dont le conoyau est annulé par $\gamma$.

	\textnormal{(ii)} On dit que $M$ est de \textit{type $\alpha$-fini} s'il est de type $\gamma$-fini pour tout $\gamma\in \mm$.
\end{definition}

\begin{definition} \label{alpha plat plat alpha}
	\textnormal{(i)} On dit qu'un $A$-module $M$ est \textit{plat} (resp. \textit{$\alpha$-plat}) si, pour tout morphisme injectif de $A$-modules $f:N_{1}\to N_{2}$, le noyau de $\id_{M}\otimes f: M\otimes_{A}N_{1}\to M\otimes_{A}N_{2}$ est nul (resp. $\alpha$-nul).

	\textnormal{(ii)} On dit qu'un $A^{\alpha}$-module $M$ est \textit{plat} si, pour tout morphisme injectif de $A^{\alpha}$-modules $f:N_{1}\to N_{2}$, le noyau de $\id_{M}\otimes f: M\otimes_{A^{\alpha}}N_{1}\to M\otimes_{A^{\alpha}}N_{2}$ est nul.
\end{definition}

\begin{lemma} \label{lemma alpha plat}
	\textnormal{(i)} Pour qu'un $A$-module $M$ soit $\alpha$-plat, il faut et il suffit que $M^{\alpha}$ soit un $A^{\alpha}$-module plat.

	\textnormal{(ii)} Pour qu'un $A^{\alpha}$-module $N$ soit plat, il faut et il suffit que $\widetilde{\tau}_{!}(N)$ soit un $A$-module plat.
\end{lemma}

\textit{Preuve}. (i) Supposons que $M$ soit un $A$-module $\alpha$-plat et soit $f: N_{1}\to N_{2}$ une injection de $A^{\alpha}$-modules. D'après \ref{coro inj exactitude}(i), $f$ induit une injection de $A$-modules $\tau_{!}(N_{1})\to \tau_{!}(N_{2})$. On en déduit que le noyau de $\id_{M}\otimes \tau_{!}(f):M\otimes_{A}\tau_{!}(N_{1})\to M\otimes_{A}\tau_{!}(N_{2})$ est $\alpha$-nul. Comme le foncteur $\widetilde{\alpha}_{A}$ est exact et commute avec les produits tensoriels, le morphisme $\id_{M^{\alpha}}\otimes f: M^{\alpha}\otimes_{A^{\alpha}}N_{1}\to M^{\alpha}\otimes_{A^{\alpha}}N_{2}$ est une injection; d'où la platitude de $M^{\alpha}$.

D'autre part, supposons que $M^{\alpha}$ soit un $A^{\alpha}$-module plat et soit $f: N_{1}\to N_{2}$ une injection de $A$-modules. Posons $N=\Ker(M\otimes_{A}N_{1}\to M\otimes_{A}N_{2})$. Comme le foncteur $\widetilde{\alpha}_{A}$ est exact et $M^{\alpha}$ est plat, on a une injection de $A^{\alpha}$-modules $M^{\alpha}\otimes_{A^{\alpha}}N_{1}^{\alpha}\to M^{\alpha}\otimes_{A^{\alpha}}N_{2}^{\alpha}$. On en déduit que $N^{\alpha}=0$; d'où l'assertion cherchée.

(ii) La suffisance de la condition résulte de (i) et de l'isomorphisme $N\xrightarrow{\sim}(\widetilde{\tau}_{!}(N))^{\alpha}$ (\ref{propriete adjoint alpha AMod}(iv)). 

Supposons que $N$ soit un $A^{\alpha}$-module plat. En vertu de (i) et de l'isomorphisme $N\xrightarrow{\sim}(\widetilde{\tau}_{*}(N))^{\alpha}$ (\ref{propriete adjoint alpha AMod}(ii)), $\widetilde{\tau}_{*}(N)$ est $\alpha$-plat. Soit $f:N_{1}\to N_{2}$ un morphisme injectif de $A$-modules. Posons $L=\Ker(\widetilde{\tau}_{*}(N)\otimes_{A^{\alpha}}N_{1}\to \widetilde{\tau}_{*}(N)\otimes_{A^{\alpha}}N_{2})$ qui est $\alpha$-nul. Compte tenu de la platitude de $\mm_{\widetilde{\mathscr{C}}}$ sur $\oo_{\widetilde{\mathscr{C}}}$, on a $\mm_{\widetilde{\mathscr{C}}}\otimes_{\oo_{\widetilde{\mathscr{C}}}}L=\Ker(\widetilde{\tau}_{!}(N)\otimes_{A^{\alpha}}N_{1}\to \widetilde{\tau}_{!}(N)\otimes_{A^{\alpha}}N_{2})$. Ce dernier est nul; d'où la platitude de $\widetilde{\tau}_{!}(N)$. 

\begin{lemma} \label{pullback plat}
	Conservons les notations de \ref{morphisme de topos pullback alpha}. 

	\textnormal{(i)} Pour tout $A$-module plat $M$, $\varphi^{*}(M)$ est un $A'$-module plat. 
	
	\textnormal{(ii)} Pour tout $A^{\alpha}$-module plat $N$, $\psi(N)$ est un $A'^{\alpha}$-module plat.
\end{lemma}

\textit{Preuve}. L'assertion (i) est démontrée dans (\cite{SGAIV} V 1.7.1). D'après \ref{lemma alpha plat}(ii), le $A$-module $\widetilde{\tau}_{!}(N)$ est plat. En vertu de (i), le $A'$-module $\varphi^{*}(\widetilde{\tau}_{A!}(N))\simeq \widetilde{\tau}_{A'!}(\psi(N))$ \eqref{diagramme tau ! B B'} est plat. La $A^{\alpha}$-platitude de $\psi(N)$ s'ensuit compte tenu de \ref{lemma alpha plat}(ii).

\begin{nothing}
	On note $\widetilde{\mathscr{C}}^{\mathbb{N}^{\circ}}$ le topos des systèmes projectifs d'objets de $\widetilde{\mathscr{C}}$ \eqref{limite projective de topos}. On pose $\oo_{\widetilde{\mathscr{C}}^{\mathbb{N}^{\circ}}}=(\oo_{\widetilde{\mathscr{C}}})_{n\ge1}$ (resp. $\mm_{\widetilde{\mathscr{C}}^{\mathbb{N}^{\circ}}}=(\mm_{\widetilde{\mathscr{C}}})_{n\ge1}$) qui est isomorphe au faisceau associé au préfaisceau constant de valeur $\oo$ (resp. $\mm$) de $\widetilde{\mathscr{C}}^{\mathbb{N}^{\circ}}$. Soit $\breve{A}=(A_{n})_{n\ge 1}$ une $\oo_{\widetilde{\mathscr{C}}^{\mathbb{N}^{\circ}}}$-algèbre. On note $\Mod(\breve{A})$ la catégorie des $\breve{A}$-modules de $\widetilde{\mathscr{C}}^{\mathbb{N}^{\circ}}$.
	\label{notation breveA modules}
\end{nothing}
\begin{prop} \label{lemma FHom proj system}
	\textnormal{(i)} Soit $M$ un $A$-module de $\widetilde{\mathscr{C}}$. On a un isomorphisme canonique fonctoriel \eqref{tau * foncteur}
\begin{equation}
	\widetilde{\tau}_{*}(M^{\alpha}) \xrightarrow{\sim} \FHom_{\oo_{\widetilde{\mathscr{C}}}}(\mm_{\widetilde{\mathscr{C}}},M).
	\label{Hom oo mm M Amod}
\end{equation}

\textnormal{(ii)} Soit $M=(M_{n})_{n\ge 1}$ un $\breve{A}$-module de $\widetilde{\mathscr{C}}^{\mathbb{N}^{\circ}}$. On a un isomorphisme canonique fonctoriel
	\begin{equation}
		\FHom_{\oo_{\widetilde{\mathscr{C}}^{\mathbb{N}^{\circ}}}}(\mm_{\widetilde{\mathscr{C}}^{\mathbb{N}^{\circ}}},M)\simeq (\FHom_{\oo_{\widetilde{\mathscr{C}}}}(\mm_{\widetilde{\mathscr{C}}},M_{n}))_{n\ge 1}
	\end{equation}
\end{prop}
\textit{Preuve}. Soient $U$ un objet de $\widetilde{\mathscr{C}}$ et $j_{U}:\widetilde{\mathscr{C}}_{/U}\to \widetilde{\mathscr{C}}$ le morphisme de localisation de $\widetilde{\mathscr{C}}$ en $U$.

(i) On a un isomorphisme canonique fonctoriel $M(U)\xrightarrow{\sim} \Gamma(\widetilde{\mathscr{C}}_{/U},j_{U}^{*}(M))$. Comme $\widehat{\alpha}$ et $\widehat{\sigma}_{*}$ transforment les faisceaux en des faisceaux, d'après \eqref{Hom oo mm M}, on a un isomorphisme canonique fonctoriel
\begin{equation}
	\widetilde{\tau}_{*}(M^{\alpha})(U)\simeq \Hom_{\oo}(\mm,M(U)).
\end{equation}
Par adjonction (\cite{SGAIV} IV 13.4.1), on a un isomorphisme canonique fonctoriel
\begin{equation}
	\Hom_{\oo}(\mm,\Gamma(\widetilde{\mathscr{C}}_{/U},j_{U}^{*}(M)))\simeq \Hom_{\oo_{\widetilde{\mathscr{C}}}|_{U}}(j_{U}^{*}(\mm_{\widetilde{\mathscr{C}}}),j_{U}^{*}(M)).
\end{equation}
L'assertion s'ensuit.

(ii) Soit $n$ un entier $\ge 1$. On désigne par $\alpha_{n}:\widetilde{\mathscr{C}}\to \widetilde{\mathscr{C}}^{\mathbb{N}^{\circ}}$ le morphisme de topos canonique défini pour tout faisceau $F=(F_{n})_{n\ge 1}$ de $\widetilde{\mathscr{C}}^{\mathbb{N}^{\circ}}$ par $\alpha_{n}^{*}(F)=F_{n}$ (\cite{AGT} (III.7.1.2)). Le foncteur $\alpha_{n}^{*}$ admet un adjoint à droit $\alpha_{n!}$ (cf. \cite{AGT} (III.7.1.5)). On en déduit des isomorphismes:
\begin{eqnarray}
	\alpha_{n}^{*}(\FHom_{\oo_{\widetilde{\mathscr{C}}^{\mathbb{N}^{\circ}}}}(\mm_{\widetilde{\mathscr{C}}^{\mathbb{N}^{\circ}}},M))(U)&\simeq & \FHom_{\oo_{\widetilde{\mathscr{C}}^{\mathbb{N}^{\circ}}}}(\mm_{\widetilde{\mathscr{C}}^{\mathbb{N}^{\circ}}},M)(\alpha_{n!}(U))\\
	&\simeq & \Hom_{\oo_{\widetilde{\mathscr{C}}^{\mathbb{N}^{\circ}}}|_{\alpha_{n!}(U)}}(\mm_{\widetilde{\mathscr{C}}^{\mathbb{N}^{\circ}}}|_{\alpha_{n!}(U)},M|_{\alpha_{n!}(U)}).
	\label{isomoprhism aaa 1}
\end{eqnarray}

On note $[n]$ le sous-ensemble ordonné $\{1,2\dots,n\}$ de $\mathbb{N}$. On munit $\widetilde{\mathscr{C}}_{/U}\times [n]$ de la topologie totale relative au site fibré constant $\widetilde{\mathscr{C}}_{/U}\times [n]\to [n]$ de fibre $\widetilde{\mathscr{C}}$ (\cite{AGT} VI.7.1) et on désigne par $(\widetilde{\mathscr{C}}_{/U})^{[n]^{\circ}}$ le topos des faisceaux d'ensembles sur $\widetilde{\mathscr{C}}_{/U}\times [n]$ (cf. \cite{AGT} III.7.1). Le foncteur d'injection canonique $\widetilde{\mathscr{C}}_{/U}\times[n]\to \widetilde{\mathscr{C}}_{/U}\times \mathbb{N}$ induit un morphisme topos (\cite{AGT} III.7.8)
\begin{equation}
	\varphi_{n}:(\widetilde{\mathscr{C}}_{/U})^{[n]^{\circ}} \to (\widetilde{\mathscr{C}}_{/U})^{\mathbb{N}^{\circ}}.
\end{equation}
On note $j_{(U,n)}: (\widetilde{\mathscr{C}}^{\mathbb{N}^{\circ}})_{/\alpha_{n!}(U)}\to \widetilde{\mathscr{C}}^{\mathbb{N}^{\circ}}$ le morphisme de localisation de $\widetilde{\mathscr{C}}^{\mathbb{N}^{\circ}}$ en $\alpha_{n!}(U)$. D'après (\cite{AGT} III.7.9), on a une équivalence canonique de topos
\begin{equation}
	h:(\widetilde{\mathscr{C}}_{/U})^{[n]^{\circ}}\xrightarrow{\sim} (\widetilde{\mathscr{C}}^{\mathbb{N}^{\circ}})_{/\alpha_{n!}(U)}
\end{equation}
telle que $j_{(U,n)}\circ h$ soit le composé
\begin{equation}
	(\widetilde{\mathscr{C}}_{/U})^{[n]^{\circ}}\xrightarrow{\varphi_{n}}(\widetilde{\mathscr{C}}_{/U})^{\mathbb{N}^{\circ}}\xrightarrow{(j_{U})^{\mathbb{N}^{\circ}}}\widetilde{\mathscr{C}}^{\mathbb{N}^{\circ}}.
\end{equation}
On en déduit un isomorphisme canonique fonctoriel
\begin{equation}
	\Hom_{\oo_{\widetilde{\mathscr{C}}^{\mathbb{N}^{\circ}}}|_{\alpha_{n!}(U)}}(\mm_{\widetilde{\mathscr{C}}^{\mathbb{N}^{\circ}}}|_{\alpha_{n!}(U)},M|_{\alpha_{n!}(U)}) \xrightarrow{\sim} \Hom_{\oo_{(\widetilde{\mathscr{C}}_{/U})^{[n]^{\circ}}}}(\mm_{(\widetilde{\mathscr{C}}_{/U})^{[n]^{\circ}}},(j_{U}^{*}(M_{i}))_{1\le i \le n}).
	\label{isomoprhism aaa 2}
\end{equation}
D'après (\cite{AGT} VI.7.13), on a un isomorphisme canonique fonctoriel
\begin{equation}
	\Hom_{\oo_{(\widetilde{\mathscr{C}}_{/U})^{[n]^{\circ}}}}(\mm_{(\widetilde{\mathscr{C}}_{/U})^{[n]^{\circ}}},(j_{U}^{*}(M_{i}))_{1\le i \le n})\xrightarrow{\sim} \Hom_{\oo_{\widetilde{\mathscr{C}}_{/U}}}(\mm_{\widetilde{\mathscr{C}}_{/U}},j_{U}^{*}(M_{n})).
	\label{isomoprhism aaa 3}
\end{equation}
Compte tenu de \eqref{isomoprhism aaa 1}, \eqref{isomoprhism aaa 2} et \eqref{isomoprhism aaa 3}, on déduit un isomorphisme canonique fonctoriel
\begin{equation}
	\alpha_{n}^{*}(\FHom_{\oo_{\widetilde{\mathscr{C}}^{\mathbb{N}^{\circ}}}}(\mm_{\widetilde{\mathscr{C}}^{\mathbb{N}^{\circ}}},M))\xrightarrow{\sim} \FHom_{\oo_{\widetilde{\mathscr{C}}}}(\mm_{\widetilde{\mathscr{C}}},M_{n}),
\end{equation}
d'où la proposition.

\begin{nothing}
	Conservons les notations de \ref{notation breveA modules}. On désigne par $\alpha\textnormal{-}\widetilde{\mathscr{C}}^{\mathbb{N}^{\circ}}$ la catégorie des faisceaux de $\alpha$-$\oo$-modules de $\widetilde{\mathscr{C}}^{\mathbb{N}^{\circ}}$ \eqref{faisceaux de alpha oo modules}. Pour tout entier $n\ge 1$, on désigne par $\Mod(A_{n}^{\alpha})$ la catégorie des $A_{n}^{\alpha}$-modules de $\alpha\textnormal{-}\widetilde{\mathscr{C}}$ et par $\Mod(\breve{A}^{\alpha})$ la catégorie des $\breve{A}^{\alpha}$-modules de $\alpha\textnormal{-}\widetilde{\mathscr{C}}^{\mathbb{N}^{\circ}}$. On note $\mathbf{P}(\Mod(A_{\bullet}^{\alpha}))$ la catégorie des systèmes projectifs $(M_{n})_{n\ge 1}$, où $M_{n}$ est un $A_{n}^{\alpha}$-module et $M_{n+1}\to M_{n}$ est un morphisme $A_{n+1}^{\alpha}$-linéaire.

	Soient $M=(M_{n})_{n\ge 1}$ un objet de $\Mod(\breve{A})$ et $f=(f_{n})_{n\ge 1}$ un morphisme de $\Mod(\breve{A})$. Pour que $M$ soit $\alpha$-nul, il faut et il suffit que, pour tout $n\ge 1$, $M_{n}$ soit $\alpha$-nul dans la catégorie $\Mod(A_{n})$. Par suite, le morphisme $f$ est un $\alpha$-isomorphisme si et seulement si, pour tout entier $n\ge 1$, le morphisme $f_{n}$ est un $\alpha$-isomorphisme de $\Mod(A_{n})$. Par suite, le foncteur
	\begin{equation}
		a: \Mod(\breve{A})\to \mathbf{P}(\Mod(A_{\bullet}^{\alpha})) \qquad (M_{n})\mapsto (M_{n}^{\alpha})
		\label{proj foncteur can a}
	\end{equation}
	envoie les $\alpha$-isomorphismes sur des isomorphismes. Il induit donc un foncteur
	\begin{equation}
		b: \Mod(\breve{A}^{\alpha})\to \mathbf{P}(\Mod(A_{\bullet}^{\alpha})).
		\label{foncteur Mod Aalpha Proj systems}
	\end{equation}
	Les foncteurs $(\widetilde{\tau}_{A_{n}*}:\Mod(A_{n}^{\alpha})\to \Mod(A_{n}))_{n\ge 1}$ \eqref{tau * foncteur} induisent un foncteur
	\begin{equation}
		t: \mathbf{P}(\Mod(A_{\bullet}^{\alpha}))\to \Mod(\breve{A})\qquad (M_{n})\mapsto (\widetilde{\tau}_{*}(M_{n})).		\label{foncteur t alpha modules}
	\end{equation}
	On désigne par
	\begin{equation}
		s: \mathbf{P}(\Mod(A_{\bullet}^{\alpha})) \to \Mod(\breve{A}^{\alpha})
	\end{equation}
	le composé du foncteur $t$ et du foncteur canonique $\widetilde{\alpha}_{\breve{A}}:\Mod(\breve{A})\to \Mod(\breve{A}^{\alpha})$. Les isomorphismes $\widetilde{\alpha}_{A_{n}}\circ \widetilde{\tau}_{A_{n}*}\xrightarrow{\sim} \id$ (\ref{propriete adjoint alpha AMod}(ii)) induisent un isomorphisme
	\begin{equation}
		b\circ s\xrightarrow{\sim} \id.
	\end{equation}
	D'après \ref{lemma FHom proj system}, le foncteur composé
\begin{equation}
	t\circ a: \Mod(\breve{A})\to \Mod(\breve{A})
\end{equation}
s'identifie au foncteur composé $\widetilde{\tau}_{\breve{A}*}\circ \widetilde{\alpha}_{\breve{A}}$. En vertu de \ref{coro inj exactitude}(iii), le foncteur $\widetilde{\alpha}_{\breve{A}}:\Mod(\breve{A})\to \Mod(\breve{A}^{\alpha})$ est isomorphe au foncteur composé
	\begin{equation}
		\Mod(\breve{A})\xrightarrow{a} \mathbf{P}(\Mod(A_{\bullet}^{\alpha})) \xrightarrow{t} \Mod(\breve{A}) \xrightarrow{\widetilde{\alpha}_{\breve{A}}} \Mod(\breve{A}^{\alpha}).
	\end{equation}
	On en déduit un isomorphisme
	\begin{equation}
		\id \xrightarrow{\sim} s\circ b.
	\end{equation}
	Par suite, $b$ et $s$ sont des équivalences de catégories, quasi-inverses l'une de l'autre.
	\label{equivalence de cat Proj systems alpha}
\end{nothing}

\section{Déformations des faisceaux de $\alpha$-modules} \label{Deformation alpha}
Dans cette section, on se donne un topos $T$. On considère toujours $T$ comme muni de sa topologie canonique (\cite{SGAIV} II 2.5), qui en fait un site.

\begin{nothing}
	Rappelons d'abord la théorie des déformations pour les modules sur un topos annelé suivant (\cite{Il71} IV 3). Soit
	\begin{equation}
		p:A\to A_{0}
	\end{equation}
	une surjection de $\oo$-algèbres de $T$ dont le noyau $I$ est de carré-nul. Soient $M_{0}$, $J$ deux $A_{0}$-modules et
	\begin{equation}
		\widetilde{M}=(0\to J\to M\to M_{0}\to 0)
		\label{extension M0 J}
	\end{equation}
	une extension de $A$-modules. On a alors un diagramme commutatif
	\begin{equation}
		\xymatrix{
			& I\otimes_{A}M \ar[r] \ar[d] & M \ar[r] \ar@{=}[d] & M\otimes_{A} A_{0} \ar[r] \ar[d] & 0\\
			0 \ar[r] & J \ar[r] & M\ar[r] & M_{0}\ar[r] &0.
		}
		\label{diagramme commutative deformation tensor}
	\end{equation}
La flèche verticale de gauche induit un morphisme de $A_{0}$-modules
	\begin{equation}
		u(\widetilde{M}): I\otimes_{A_{0}}M_{0}\to J.
		\label{uM extension}
	\end{equation}
	En vertu de \eqref{diagramme commutative deformation tensor}, $u(\widetilde{M})$ est un épimorphisme si et seulement si le morphisme canonique $M\otimes_{A}A_{0}\to M_{0}$ est un isomorphisme. La correspondance $\widetilde{M}\mapsto u(\widetilde{M})$ définit un homomorphisme (cf. \cite{Il71} IV 3.1.1)
	\begin{equation}
		u:\Ext^{1}_{A}(M_{0},J) \to \Hom_{A_{0}}(I\otimes_{A_{0}}M_{0},J).
		\label{connection morphisme u}
	\end{equation}
On a une suite exacte (cf. \cite{Il71} IV 3.1.4):
	\begin{equation}
		0\to \Ext_{A_{0}}^{1}(M_{0},J)\to \Ext_{A}^{1}(M_{0},J)\xrightarrow{u} \Hom_{A_{0}}(I\otimes_{A_{0}}M_{0},J)\xrightarrow{\partial} \Ext_{A_{0}}^{2}(M_{0},J),
		\label{suite exacte deformation}
	\end{equation}
	où la première flèche est le morphisme défini par restriction des scalaires.
	\label{deformation generale}
\end{nothing}

\begin{theorem}[\cite{Il71} IV 3.1.5]
	Soient $M_{0}$, $J$ des $A_{0}$-modules et $u_{0}:I\otimes_{A_{0}}M_{0}\to J$ un morphisme $A_{0}$-linéaire. Alors:
	
	\textnormal{(i)} Il existe une obstruction
	\begin{equation}
		\partial(u_{0})\in \Ext^{2}_{A_{0}}(M_{0},J)
		\label{obstruction}
	\end{equation}
	dont l'annulation est nécessaire et suffisante pour l'existence d'une extension de $A$-modules $\widetilde{M}$ de $M_{0}$ par $J$ telle que $u(\widetilde{M})=u_{0}$ \eqref{connection morphisme u}.
	
	\textnormal{(ii)} Lorsque $\partial(u_{0})=0$, l'ensemble des classes d'isomorphismes de telles extensions $\widetilde{M}$ est un torseur sous $\Ext^{1}_{A_{0}}(M_{0},J)$.
	\label{theorem principal deformation modules}
\end{theorem}

\begin{lemma}[\cite{Il71} IV 3.1.1] \label{lemma plat deformation}
	Supposons que le topos $T$ ait suffisamment de points. Soit $\widetilde{M}=(0\to J\to M\to M_{0}\to 0)$ une extension de $A$-modules tel que le morphisme canonique $M\otimes_{A}A_{0}\to M_{0}$ soit un isomorphisme. Alors, pour que $M$ soit plat sur $A$, il faut et il suffit que $M_{0}$ soit plat sur $A_{0}$ et que $u(\widetilde{M})$ \eqref{uM extension} soit un isomorphisme.
\end{lemma}
\textit{Preuve}. Comme $T$ a suffisamment de points, on peut se ramener au cas où $T$ est le topos ponctuel. L'assertion résulte alors du critère de platitude (\cite{EGA III} 0.10.2.1).

\begin{rem}
	La suite exacte \eqref{suite exacte deformation} peut s'établir par un calcul direct (cf. \cite{Il71} IV 3.1.12). L'isomorphisme de Cartan
	\begin{equation}
		\RHom_{A}(M_{0},J)\xrightarrow{\sim} \RHom_{A_{0}}(M_{0}\otimes^{\rL}_{A}A_{0},J)
		\label{isomorphisme de Cartan derive}
	\end{equation}
	induit une suite spectrale
	\begin{equation}
		\rE^{i,j}_{2}=\Ext^{i}_{A_{0}}(\FTor_{j}^{A}(M_{0},A_{0}),J)\Rightarrow \rE^{i+j}=\Ext^{i+j}_{A}(M_{0},J).
		\label{suite spectrale Cartan Grothendieck}
	\end{equation}
	Comme $M_{0}$ est un $A_{0}$-module, on a alors un isomorphisme canonique $\FTor_{1}^{A}(M_{0},A_{0})\xrightarrow{\sim} I\otimes_{A_{0}}M$. La suite exacte des termes de bas degré de \eqref{suite spectrale Cartan Grothendieck} fournit la suite exacte \eqref{suite exacte deformation} (cf. \cite{Il71} IV 3.1.13).
	\label{rem suite spectrale Cartan Gro}
\end{rem}

\begin{nothing}
	Soit $\varphi:(T',A')\to (T,A)$ un morphisme plat de topos annelés (\cite{SGAIV} V 1.8). On pose $A_{0}'=\varphi^{*}(A_{0})$ et $I'=\varphi^{*}(I)$. On désigne par 
	\begin{equation}
		\rE'^{i,j}_{2}=\Ext^{i}_{A'_{0}}(\FTor_{j}^{A'}(\varphi^{*}(M_{0}),A'_{0}),\varphi^{*}(J))\Rightarrow \rE'^{i+j}=\Ext^{i+j}_{A'}(\varphi^{*}(M_{0}),\varphi^{*}(J))
		\label{}
	\end{equation}
	la suite spectrale \eqref{suite spectrale Cartan Grothendieck} associée aux $A'_{0}$-modules $\varphi^{*}(M_{0})$ et $\varphi^{*}(J)$. Le faisceau $\FTor_{A}^{j}(M_{0},A_{0})$ est calculé par une résolution plate de $M_{0}$. En vertu de \ref{pullback plat}(i) et l'exactitude du foncteur $\varphi^{*}$, on a un isomorphisme canonique
	\begin{equation}
		\varphi^{*}(\FTor_{A}^{j}(M_{0},A_{0}))\simeq \FTor_{A'}^{j}(\varphi^{*}(M_{0}),A_{0}').
	\end{equation}
	Le foncteur exact $\varphi^{*}$ induit un morphisme de suites spectrales (cf. \ref{explication ss Grothendieck})
	\begin{equation}
		u:\rE=(\rE_{2}^{i,j},\rE^{n})\to \rE'=(\rE'^{i,j}_{2},\rE'^{n})
		\label{morphisme de ss Grothendieck Ill}
	\end{equation}
	dont les morphismes
	\begin{eqnarray}
		&u_{2}^{i,j}&:\Ext^{i}_{A_{0}}(\FTor_{j}^{A}(M_{0},A_{0}),J)\to \Ext^{i}_{A'_{0}}(\FTor_{j}^{A'}(\varphi^{*}(M_{0}),A'_{0}),\varphi^{*}(J)),\\ 
		&u^{n}&:\Ext^{n}_{A}(M_{0},J)\to \Ext^{n}_{A'}(\varphi^{*}(M_{0}),\varphi^{*}(J))
	\end{eqnarray}
	s'identifient aux morphismes canoniques induits par $\varphi^{*}$ \eqref{exact foncteur induit Ext morphisme}. Les termes de bas degré du morphisme $u$ induisent un diagramme commutatif
	\begin{eqnarray}
		\xymatrix{
			0\ar[r] & \Ext_{A_{0}}^{1}(M_{0},J)\ar[r] \ar[d]_{\varphi^{*}}& \Ext_{A}^{1}(M_{0},J)\ar[r]^-{u} \ar[d]^{\varphi^{*}}& \Hom_{A_{0}}(I\otimes_{A_{0}}M_{0},J) \ar[d]^{\varphi^{*}} \\
			0\ar[r] & \Ext_{A_{0}^{'}}^{1}(\varphi^{*}(M_{0}),\varphi^{*}(J))\ar[r] & \Ext_{A^{'}}^{1}(\varphi^{*}(M_{0}),\varphi^{*}(J))\ar[r]^-{u'} & \Hom_{A_{0}^{'}}(I^{'}\otimes_{A_{0}^{'}}\varphi^{*}(M_{0}),\varphi^{*}(J))
		} \nonumber\\
	\xymatrix{
		\Hom_{A_{0}}(I\otimes_{A_{0}}M_{0},J)\ar[r]^-{\partial} \ar[d]_{\varphi^{*}}& \Ext_{A_{0}}^{2}(M_{0},J) \ar[d]^{\varphi^{*}} \\
		\Hom_{A_{0}^{'}}(I^{'}\otimes_{A_{0}^{'}}\varphi^{*}(M_{0}),\varphi^{*}(J))\ar[r]^-{\partial'} & \Ext_{A_{0}^{'}}^{2}(\varphi^{*}(M_{0}),\varphi^{*}(J)).} 
		\label{diagramme commutatif phi plat}
	\end{eqnarray}
	Cela signifie que la théorie de la déformation est fonctorielle par rapport au morphisme plat $\varphi$.
	\label{phi exact fonctoriel def}
\end{nothing}
\begin{nothing}
	On note $\alpha\textnormal{-}T$ la catégorie des $\alpha$-$\oo$-modules de $T$ \eqref{faisceaux de alpha oo modules} et on pose $A_{0}^{\alpha}=\widetilde{\alpha}(A_{0})$ et $A^{\alpha}=\widetilde{\alpha}(A)$ \eqref{alpha Alg sheaf}. Soient $M_{0}$ et $J$ deux $A_{0}^{\alpha}$-modules. D'après \ref{propriete adjoint alpha AMod}(iv), \eqref{Ext alpha Ext isomorphisme} et \eqref{isomorphisme tau ! tensor}, on a des isomorphismes canoniques 
	\begin{eqnarray}
		\Ext^{1}_{A}(\widetilde{\tau}_{!}(M_{0}),\widetilde{\tau}_{!}(J)) &\xrightarrow{\sim}& \Ext^{1}_{A^{\alpha}}(M_{0},J),\nonumber\\
		\Ext^{2}_{A_{0}}(\widetilde{\tau}_{!}(M_{0}),\widetilde{\tau}_{!}(J))&\xrightarrow{\sim}&\Ext^{2}_{A_{0}^{\alpha}}(M_{0},J), \nonumber\\
		\Hom_{A_{0}}(I\otimes_{A_{0}}\widetilde{\tau}_{!}(M_{0}),\widetilde{\tau}_{!}(J)) &\xrightarrow{\sim}& \Hom_{A_{0}^{\alpha}}(I^{\alpha}\otimes_{A_{0}^{\alpha}}M_{0}, J). \nonumber
	\end{eqnarray}
	On définit des morphismes
	\begin{eqnarray}
		u^{\alpha}:\Ext^{1}_{A^{\alpha}}(M_{0},J)\to \Hom_{A_{0}^{\alpha}}(I^{\alpha}\otimes_{A_{0}^{\alpha}}M_{0},J),\\
		\partial^{\alpha}: \Hom_{A_{0}^{\alpha}}(I^{\alpha}\otimes_{A^{\alpha}}M_{0},J) \to \Ext^{2}_{A_{0}^{\alpha}}(M_{0},J) ,
	\end{eqnarray}
	par les morphismes $u$ et $\partial$ relatifs aux $A_{0}$-modules $\widetilde{\tau}_{!}(M_{0})$ et $\widetilde{\tau}_{!}(J)$ \eqref{suite exacte deformation}. On en déduit par \eqref{suite exacte deformation} une suite exacte
	\begin{equation}
		0\to \Ext_{A_{0}^{\alpha}}^{1}(M_{0},J)\to \Ext_{A^{\alpha}}^{1}(M_{0},J)\xrightarrow{u^{\alpha}} \Hom_{A_{0}^{\alpha}}(I^{\alpha}\otimes_{A_{0}^{\alpha}}M_{0},J)\xrightarrow{\partial^{\alpha}} \Ext_{A_{0}^{\alpha}}^{2}(M_{0},J),
		\label{suite exacte deformation alpha}
	\end{equation}
	où la première flèche est le morphisme défini par restriction des scalaires.
	\label{deformation alpha generale}
\end{nothing}

\begin{theorem}
	Conservons les notations de \ref{deformation alpha generale}. Soient $M_{0}$, $J$ des $A_{0}^{\alpha}$-modules et $u_{0}:I^{\alpha}\otimes_{A_{0}^{\alpha}}M_{0}\to J$ un morphisme $A_{0}^{\alpha}$-linéaire. Alors:
	
	\textnormal{(i)} Il existe une obstruction
	\begin{equation}
		\partial^{\alpha}(u_{0})\in \Ext^{2}_{A_{0}^{\alpha}}(M_{0},J)
		\label{alpha obstruction}
	\end{equation}
	dont l'annulation est nécessaire et suffisante pour l'existence d'une extension de $A^{\alpha}$-modules $\widetilde{M}$ de $M_{0}$ par $J$ telle que $u^{\alpha}(\widetilde{M})=u_{0}$ \eqref{suite exacte deformation alpha}.
	
	\textnormal{(ii)} Lorsque $\partial^{\alpha}(u_{0})=0$, l'ensemble des classes d'isomorphismes de telles extensions $\widetilde{M}$ est un torseur sous $\Ext^{1}_{A_{0}^{\alpha}}(M_{0},J)$.
	\label{theorem principal deformation alpha modules}
\end{theorem}

\begin{nothing}
	Conservons les notations de \ref{phi exact fonctoriel def}. On note $\alpha\textnormal{-}T'$ la catégorie des faisceaux de $\alpha$-$\oo$-modules de $T'$ \eqref{faisceaux de alpha oo modules} et on pose $A'^{\alpha}=\widetilde{\alpha}(A')$ et $A'^{\alpha}_{0}=\widetilde{\alpha}(A_{0}')$ \eqref{alpha Alg sheaf}. Le foncteur exact $\varphi^{*}:\Mod(A)\to \Mod(A')$ s'étend en un foncteur exact $\psi:\Mod(A^{\alpha})\to \Mod(A'^{\alpha})$ (\ref{prop morphisme plat topos}). D'après \eqref{diagramme commutatif Ext phi tau alpha} et \eqref{diagramme commutatif phi plat}, on a un diagramme commutatif
	\begin{eqnarray}
		\xymatrix{
			0\ar[r] & \Ext_{A_{0}^{\alpha}}^{1}(M_{0},J)\ar[r] \ar[d]_{\psi}& \Ext_{A^{\alpha}}^{1}(M_{0},J)\ar[r]^-{u^{\alpha}} \ar[d]^{\psi}& \Hom_{A_{0}^{\alpha}}(I^{\alpha}\otimes_{A_{0}}M_{0},J) \ar[d]^{\psi} \\
			0\ar[r] & \Ext_{A'^{\alpha}_{0}}^{1}(\psi(M_{0}),\psi(J))\ar[r] & \Ext_{A'^{\alpha}}^{1}(\psi(M_{0}),\psi(J))\ar[r]^-{u'^{\alpha}} & \Hom_{A'^{\alpha}_{0}}(I'^{\alpha}\otimes_{A'^{\alpha}_{0}}\psi(M_{0}),\psi(J))
		} \nonumber\\
	\xymatrix{
		\Hom_{A_{0}^{\alpha}}(I^{\alpha}\otimes_{A_{0}^{\alpha}}M_{0},J)\ar[r]^-{\partial^{\alpha}} \ar[d]_{\psi}& \Ext_{A_{0}^{\alpha}}^{2}(M_{0},J) \ar[d]^{\psi} \\
		\Hom_{A'^{\alpha}_{0}}(I'^{\alpha}\otimes_{A'^{\alpha}_{0}}\psi(M_{0}),\psi(J))\ar[r]^-{\partial'^{\alpha}} & \Ext_{A'^{\alpha}_{0}}^{2}(\psi(M_{0}),\psi(J)),} 
		\label{diagramme commutatif phi plat alphat}
	\end{eqnarray}
	où les flèches horizontaux sont induites par le foncteur exact $\psi$ \eqref{exact foncteur induit Ext morphisme}. Cela signifie que la théorie de la déformation pour les $\alpha$-modules est fonctorielle par rapport au morphisme plat $\varphi$.
	\label{foncteur pullback alpha psi}
\end{nothing}

\begin{nothing}
	Soit $A$ une $\oo$-algèbre (ou $\mathcal{O}_{\overline{K}}$-algèbre) plate de $T$. Pour tout entier $n\ge 1$, on pose $A_{n}=A/p^{n}A$ et $I_{n}=p^{n}A/p^{2n}A$ l'idéal de $A_{2n}$. Le morphisme de multiplication par $p^{n}$ induit, pour tout entier $n\ge 1$, un isomorphisme canonique de $A_{n}$-modules $I_{n}\xrightarrow{\sim}A_{n}$. Pour tout entier $n\ge 1$, on pose $A_{n}^{\alpha}=\widetilde{\alpha}(A_{n})$ \eqref{alpha Alg sheaf}.
	
	Soient $n$ un entier $\ge 1$, $M_{n}$ un $A_{n}$-module et $N_{n}$ un $A_{n}^{\alpha}$-module. L'isomorphisme $I_{n}\xrightarrow{\sim} A_{n}$ induit un isomorphisme canonique de $A_{n}$-modules (resp. $A_{n}^{\alpha}$-modules)
	\begin{equation}
		v:I_{n}\otimes_{A_{n}}M_{n}\xrightarrow{\sim} M_{n}\qquad \textnormal{(resp.} \quad w:I_{n}^{\alpha}\otimes_{A_{n}^{\alpha}}N_{n}\xrightarrow{\sim} N_{n} \textnormal{)}.
		\label{u Mn canonique}
	\end{equation}
	On appelle \textit{déformation de $M_{n}$ sur $A_{2n}$} (resp. \textit{déformation de $N_{n}$ sur $A_{2n}^{\alpha}$}) toute extension de $A_{2n}$-modules (resp. $A_{2n}^{\alpha}$-modules)
	\begin{equation}
		\widetilde{M}=(0\to M_{n}\to M_{2n}\to M_{n}\to 0) \qquad \textnormal{(resp.} \quad \widetilde{N}=(0\to N_{n}\to N_{2n}\to N_{n}\to 0) \textnormal{)}
		\label{deformation of Mn}
	\end{equation}
	telle que $u(\widetilde{M})=v$ (resp. $u^{\alpha}(\widetilde{N})=w$) \eqref{u Mn canonique}.
	\label{notations deformation}
\end{nothing}

\begin{nothing}
	Pour toute déformation $0\to M_{n}\to M_{2n}\to M_{n}\to 0$ de $M_{n}$ sur $A_{2n}$, il est clair que $0\to M_{n}^{\alpha}\to M_{2n}^{\alpha}\to M_{n}^{\alpha}\to 0$ est une déformation de $M_{n}^{\alpha}$ sur $A_{2n}^{\alpha}$. D'autre part, si $0\to N_{n}\to N_{2n}\to N_{n}\to 0$ est une déformation de $N_{n}$ sur $A_{2n}^{\alpha}$, la suite exacte $0\to \widetilde{\tau}_{!}(N_{n})\to \widetilde{\tau}_{!}(N_{2n})\to \widetilde{\tau}_{!}(N_{n})\to 0$ est une déformation de $\widetilde{\tau}_{!}(N_{n})$ sur $A_{2n}$ en vertu de \ref{coro inj exactitude}(i) et \eqref{isomorphisme tau ! tensor}.
	\label{deformation alpha deformation}
\end{nothing}

\begin{lemma}	
	Supposons que le topos $T$ ait suffisamment de points et soit $n$ un entier $\ge 1$.
	
	\textnormal{(i)} Soit $M_{n}$ un $A_{n}$-module plat. Une déformation de $M_{n}$ sur $A_{2n}$ est équivalente à la donnée d'un $A_{2n}$-module plat $M_{2n}$ et d'un épimorphisme $A_{2n}$-linéaire $g:M_{2n}\to M_{n}$ tels que le morphisme canonique $M_{2n}\otimes_{A_{2n}}A_{n}\to M_{n}$ soit un isomorphisme.

	\textnormal{(ii)} Soit $N_{n}$ un $A_{n}^{\alpha}$-module plat \textnormal{(\ref{alpha plat plat alpha}(ii))}. Une déformation de $N_{n}$ sur $A_{2n}^{\alpha}$ est équivalente à la donnée d'un $A_{2n}^{\alpha}$-module plat $N_{2n}$ et d'un épimorphisme $A_{2n}^{\alpha}$-linéaire $g:N_{2n}\to N_{n}$ tels que le morphisme canonique $N_{2n}\otimes_{A_{2n}^{\alpha}}A_{n}^{\alpha}\to N_{n}$ soit un isomorphisme.
	\label{lemma deformation equivalence}
\end{lemma}

\textit{Preuve}. (i) Soit $0\to M_{n}\to M_{2n}\to M_{n}\to 0$ une déformation de $M_{n}$ sur $A_{2n}$. Comme $v$ est un isomorphisme, le morphisme $M_{2n}\otimes_{A_{2n}}A_{n}\to M_{n}$ est un isomorphisme (cf. \ref{deformation generale}). D'après \ref{lemma plat deformation}, $M_{2n}$ est plat sur $A_{2n}$.

D'autre part, soient $M_{2n}$ un $A_{2n}$-module plat et $g:M_{2n}\to M_{n}$ un morphisme tel que $M_{2n}\otimes_{A_{2n}}A_{n}\xrightarrow{\sim} M_{n}$. Par platitude de $M_{2n}$, on en déduit un isomorphisme $I_{n}\otimes_{A_{2n}}M_{2n}\xrightarrow{\sim} \Ker(g)$, d'où l'assertion.

(ii) Soit $0\to N_{n}\to N_{2n}\to N_{n}\to 0$ une déformation de $N_{n}$ sur $A_{2n}^{\alpha}$. D'après \ref{lemma alpha plat}(ii), $\widetilde{\tau}_{!}(N_{n})$ est $A_{n}$-plat. En vertu de (i) et \ref{deformation alpha deformation}, $\widetilde{\tau}_{!}(N_{2n})$ est $A_{2n}$-plat et on a un isomorphisme $\widetilde{\tau}_{!}(N_{2n})\otimes_{A_{2n}}A_{n}\simeq \widetilde{\tau}_{!}(N_{n})$. On en déduit par \ref{propriete adjoint alpha AMod}(iv) un isomorphisme $N_{2n}\otimes_{A_{2n}^{\alpha}}A_{n}^{\alpha}\simeq N_{n}$. La platitude de $N_{2n}$ résulte de \ref{lemma alpha plat}(i).

D'autre part, soient $N_{2n}$ un $A_{2n}^{\alpha}$-module plat et $g:N_{2n}\to N_{n}$ un morphisme tel que $N_{2n}\otimes_{A_{2n}^{\alpha}}A_{n}^{\alpha}\xrightarrow{\sim} N_{n}$. Par platitude de $N_{2n}$, on en déduit un isomorphisme $I_{n}^{\alpha}\otimes_{A_{2n}^{\alpha}}N_{2n}\xrightarrow{\sim} \Ker(g)$, d'où l'assertion.

\begin{nothing}	\label{Mod ptf oon}
	Dans la suite de cette section, on se donne un schéma connexe $X$ et un point géométrique $\overline{x}$ de $X$. Soit $n$ un entier $\ge 1$. On note $X_{\fet}$ le topos fini étale de $X$ et $\oo_{n}$ le faisceau constant de valeur $\oo_{n}$ de $X_{\fet}$. Le foncteur fibre en $\overline{x}$
	\begin{equation}
		\Mod(X_{\fet},\oo_{n})\to \Mod(\oo_{n}),\qquad \mathbb{L}\mapsto \mathbb{L}_{\overline{x}}
		\label{foncteur fibre x Mod}
	\end{equation}
	induit une équivalence de catégories \eqref{foncteur nu}
	\begin{equation}
		\Mod(X_{\fet},\oo_{n}) \xrightarrow{\sim}\Rep_{\oo_{n}}(\pi_{1}(X,\overline{x})),
		\label{equi Rep oon Mod oon}
	\end{equation}
	où le but désigne la catégorie des $\oo_{n}$-représentations de $\pi_{1}(X,\overline{x})$ \eqref{oo oon reprentations general}. On désigne par $\Mod^{\aptf}(X_{\fet},\oo_{n})$ la catégorie des $\oo_{n}$-modules $\alpha$-plats de type $\alpha$-fini de $X_{\fet}$ (\ref{type alpha fini topos}, \ref{alpha plat plat alpha}).
\end{nothing}

\begin{lemma} \label{lemma ptf oon module alphat deformation}
	\textnormal{(i)} Pour qu'un $\oo_{n}$-module $\mathbb{L}$ de $X_{\fet}$ soit $\alpha$-plat de type $\alpha$-fini, il faut et il suffit que sa fibre $\mathbb{L}_{\overline{x}}$ \eqref{foncteur fibre x Mod} soit un $\oo_{n}$-module $\alpha$-plat de type $\alpha$-fini.
	
	\textnormal{(ii)} Le foncteur \eqref{equi Rep oon Mod oon} induit une équivalence de catégories \eqref{oo oon rep aptf}
	\begin{equation} \label{equivalence de categories Mod aptf Rep aptf}		
		\Mod^{\aptf}(X_{\fet},\oo_{n})\xrightarrow{\sim} \Rep_{\oo_{n}}^{\aptf}(\pi_{1}(X,\overline{x}))
	\end{equation}
\end{lemma}
\textit{Preuve}. (i) Comme le foncteur \eqref{foncteur fibre x Mod} est conservatif, la $\alpha$-platitude d'un $\oo_{n}$-module de $X_{\fet}$ est équivalente à la $\alpha$-platitude de sa fibre en $\overline{x}$. 

Si $\mathbb{L}$ est de type $\alpha$-fini, sa fibre $\mathbb{L}_{\overline{x}}$ est évidement de type $\alpha$-fini. Supposons que $\mathbb{L}_{\overline{x}}$ soit de type $\alpha$-fini et montrons que $\mathbb{L}$ est de type $\alpha$-fini. Pour tout $\gamma\in\mm$, il existe un $\oo_{n}$-module libre de type fini $M$ et un morphisme $\oo_{n}$-linéaire $u:M\to \mathbb{L}_{\overline{x}}$ dont le noyau est annulé par $\gamma$. Soit $e_{1},\cdots,e_{m}$ une base de $M$. L'assertion recherchée étant locale pour $X_{\fet}$, quitte à remplacer $X$ par un revêtement étale, on peut supposer que $\pi_{1}(X,\overline{x})$ fixe les éléments $u(e_{1}),\cdots,u(e_{m})$ de $\mathbb{L}_{\overline{x}}$. Munissant $M$ de l'unique $\oo_{n}$-représentation de $\pi_{1}(X,\overline{x})$ telle que $e_{1},\cdots,e_{m}$ soient fixes, l'homomorphisme $u$ est alors $\pi_{1}(X,\overline{x})$-équivariant. Par suite, $\mathbb{L}$ est un $\oo_{n}$-module de type $\alpha$-fini de $X_{\fet}$.

(ii) Cela résulte de (i).


\begin{nothing} \label{notation Mod ptf oon alpha}
	On désigne par $\alpha\textnormal{-}X_{\fet}$ la catégorie des $\alpha$-$\oo$-modules de $X_{\fet}$ \eqref{faisceaux de alpha oo modules}. Soit $n$ un entier $\ge 1$. On pose $\oo_{n}^{\alpha}=\widetilde{\alpha}(\oo_{n})$ \eqref{alpha Alg sheaf} et on désigne par $\Mod(X_{\fet},\oo_{n}^{\alpha})$ la catégorie des $\oo_{n}^{\alpha}$-modules de $\alpha\textnormal{-}X_{\fet}$ \eqref{alpha A Mod almost}. On a un foncteur canonique \eqref{alpha oo modules sheaf} 
	\begin{equation}
		\widetilde{\alpha}:\Mod(X_{\fet},\oo_{n})\to \Mod(X_{\fet},\oo_{n}^{\alpha}).
	\end{equation}
On désigne par $\Mod^{\ptf}(X_{\fet},\oo_{n}^{\alpha})$ l'image essentielle de la catégorie $\Mod^{\aptf}(X_{\fet},\oo_{n})$ dans $\Mod(X_{\fet},\oo_{n}^{\alpha})$.
\end{nothing}

\begin{lemma} \label{coro deformation alpha plat}
	Soit $\mathscr{L}$ un objet de $\Mod^{\ptf}(X_{\fet},\oo_{n}^{\alpha})$. 
	 
	\textnormal{(i)} Le $\oo_{n}^{\alpha}$-module $\mathscr{L}$ est plat \eqref{alpha plat plat alpha}.

	\textnormal{(ii)} Toute déformation de $\mathscr{L}$ sur $\oo_{2n}^{\alpha}$ \textnormal{(\ref{lemma deformation equivalence}(ii))} est un objet de $\Mod^{\ptf}(X_{\fet},\oo_{2n}^{\alpha})$. 
\end{lemma}
\textit{Preuve}. (i) Cela résulte de \ref{lemma alpha plat}(i).

(ii) Soit $\mathbb{L}$ un $\oo_{n}$-module $\alpha$-plat de type $\alpha$-fini de $X_{\fet}$ tel que $\mathbb{L}^{\alpha}\simeq \mathscr{L}$. On en déduit par adjonction (\ref{propriete adjoint alpha AMod}(iii)) un $\alpha$-isomorphisme $\widetilde{\tau}_{!}(\mathscr{L})\to \mathbb{L}$; d'où la $\alpha$-finitude de $\widetilde{\tau}_{!}(\mathscr{L})$. D'après \ref{lemma alpha plat}(ii), $\widetilde{\tau}_{!}(\mathscr{L})$ est $\oo_{n}$-plat. Soit $\mathscr{L}'$ une déformation de $\mathscr{L}$ sur $\oo_{2n}^{\alpha}$. En vertu de \ref{deformation alpha deformation} et \ref{lemma deformation equivalence}(i), $\widetilde{\tau}_{!}(\mathscr{L}')$ est une déformation de $\widetilde{\tau}_{!}(\mathscr{L})$ sur $\oo_{2n}$ et est donc $\oo_{2n}$-plat. Par ailleurs, on a une suite exacte de $\oo_{2n}$-modules
\begin{equation}
	0\to \widetilde{\tau}_{!}(\mathscr{L})_{\overline{x}} \to \widetilde{\tau}_{!}(\mathscr{L}')_{\overline{x}} \to \widetilde{\tau}_{!}(\mathscr{L})_{\overline{x}}\to 0.
\end{equation}
On en déduit par \ref{lemma Gabber exactes tf} et \ref{lemma ptf oon module alphat deformation}(i) la $\alpha$-finitude de $\widetilde{\tau}_{!}(\mathscr{L}')$. Par suite, $\widetilde{\tau}_{!}(\mathscr{L}')$ appartient à $\Mod^{\aptf}(X_{\fet},\oo_{2n})$. L'assertion résulte alors de l'isomorphisme $\mathscr{L}'\xrightarrow{\sim} (\widetilde{\tau}_{!}(\mathscr{L}'))^{\alpha}$ (\ref{propriete adjoint alpha AMod}(iv)).

\begin{nothing}\label{notations Mod aptf oon oo}
	On note $X_{\fet}^{\mathbb{N}^{\circ}}$ le topos des systèmes projectifs d'objets de $X_{\fet}$ \eqref{limite projective de topos} et $\breve{\oo}$ l'anneau $(\oo_{n})_{n\ge 1}$ de $X_{\fet}^{\mathbb{N}^{\circ}}$. On désigne par $\Mod^{\aptf}(X_{\fet}^{\mathbb{N}^{\circ}},\breve{\oo})$ la sous-catégorie pleine de $\Mod(X_{\fet}^{\mathbb{N}^{\circ}},\breve{\oo})$ formée des $\breve{\oo}$-modules $(\mathbb{L}_{n})_{n\ge 1}$ tels que, pour tout entier $n\ge 1$, $\mathbb{L}_{n}$ soit un objet de $\Mod^{\aptf}(X_{\fet},\oo_{n})$ et que le morphisme canonique $\mathbb{L}_{n+1}\otimes_{\oo_{n+1}}\oo_{n}\to \mathbb{L}_{n}$ soit un \textit{$\alpha$-isomorphisme}. Compte tenu de \ref{oo oon rep aptf}, les foncteurs \eqref{equivalence de categories Mod aptf Rep aptf} induisent une équivalence de catégories
	\begin{equation} \label{systeme proj Mod Rep}
	\Mod^{\aptf}(X_{\fet}^{\mathbb{N}^{\circ}},\breve{\oo})\xrightarrow{\sim} \Rep_{\breve{\oo}}^{\aptf}(\pi_{1}(X,\overline{x})).
	\end{equation}
	
	On désigne par $\mathbf{P}(\Mod(X_{\fet},\oo_{\bullet}^{\alpha}))$ la catégorie des systèmes projectifs $(\mathscr{L}_{n})_{n\ge 1}$, où $\mathscr{L}_{n}$ est un $\oo_{n}^{\alpha}$-module de $\alpha\textnormal{-}X_{\fet}$ et $\mathscr{L}_{n+1}\to \mathscr{L}_{n}$ est un morphisme $\oo_{n+1}^{\alpha}$-linéaire \eqref{equivalence de cat Proj systems alpha}. Rappelons que le foncteur canonique
	\begin{equation}
		\Mod(X_{\fet}^{\mathbb{N}^{\circ}},\breve{\oo})\to \mathbf{P}(\Mod(X_{\fet},\oo_{\bullet}^{\alpha}))\qquad (\mathbb{L}_{n})_{n\ge 1}\mapsto (\mathbb{L}_{n}^{\alpha})_{n\ge 1},
		\label{Mod X oo to P oon module}
	\end{equation}
	induit une équivalence de catégories \eqref{foncteur Mod Aalpha Proj systems}
	\begin{equation}
		\Mod(X_{\fet}^{\mathbb{N}^{\circ}},\breve{\oo}^{\alpha})\xrightarrow{\sim}\mathbf{P}(\Mod(X_{\fet},\oo_{\bullet}^{\alpha})),
		\label{Mod X oo alpha P oon module}
	\end{equation}

	On désigne par $\mathbf{P}(\Mod^{\aptf}(X_{\fet},\oo_{\bullet}^{\alpha}))$ la sous-catégorie pleine de $\mathbf{P}(\Mod(X_{\fet},\oo_{\bullet}^{\alpha}))$ formée des systèmes projectifs $(\mathscr{L}_{n})_{n\ge 1}$ tels que, pour tout entier $n\ge 1$, $\mathscr{L}_{n}$ soit un objet de $\Mod^{\ptf}(X_{\fet},\oo_{n}^{\alpha})$ \eqref{notation Mod ptf oon alpha} et que le morphisme canonique $\mathscr{L}_{n+1}\otimes_{\oo_{n+1}^{\alpha}}\oo_{n}^{\alpha}\xrightarrow{\sim} \mathscr{L}_{n}$ soit un \textit{isomorphisme}. On note $\mathbf{P}(\Mod^{\ptf}(X_{\fet},\oo_{\bullet}^{\alpha}))_{\mathbb{Q}}$ la catégorie des objets de $\mathbf{P}(\Mod^{\ptf}(X_{\fet},\oo_{\bullet}^{\alpha}))$ à isogénie près. On désigne par $\Sysl(X_{\fet}^{\mathbb{N}^{\circ}},\breve{\oo})$ (resp. $\Sysl_{\mathbb{Q}}(X_{\fet}^{\mathbb{N}^{\circ}},\breve{\oo})$) la catégorie des $\breve{\oo}$-modules localement libres de type fini de $X_{\fet}^{\mathbb{N}^{\circ}}$ (resp. $\breve{\oo}$-modules localement libres de type fini de $X_{\fet}^{\mathbb{N}^{\circ}}$ à isogénie près).
\end{nothing}

\begin{lemma} \label{lemma PMod essentielle}
	La catégorie $\mathbf{P}(\Mod^{\aptf}(X_{\fet},\oo_{\bullet}^{\alpha}))$ est l'image essentielle de la catégorie $\Mod^{\aptf}(X_{\fet}^{\mathbb{N}},\breve{\oo})$ par le foncteur \eqref{Mod X oo to P oon module}.
\end{lemma}
\textit{Preuve}. Il est clair que l'image essentielle de $\Mod^{\aptf}(X_{\fet}^{\mathbb{N}},\breve{\oo})$ est contenu dans $\mathbf{P}(\Mod^{\aptf}(X_{\fet},\oo_{\bullet}^{\alpha}))$. Soit $\mathscr{L}=(\mathscr{L}_{n})_{n\ge 1}$ un objet de $\mathbf{P}(\Mod^{\aptf}(X_{\fet},\oo_{\bullet}^{\alpha}))$. En vertu de la preuve de \ref{coro deformation alpha plat}(ii), $\widetilde{\tau}_{!}(\mathscr{L}_{n})$ est $\alpha$-plat de type $\alpha$-fini. L'isomorphisme canonique $\mathscr{L}_{n+1}\otimes_{\oo_{n+1}^{\alpha}}\oo_{n}^{\alpha}\xrightarrow{\sim} \mathscr{L}_{n}$ induit un isomorphisme $\widetilde{\tau}_{!}(\mathscr{L})_{n+1}\otimes_{\oo_{n+1}}\oo_{n}\xrightarrow{\sim} \widetilde{\tau}_{!}(\mathscr{L}_{n})$ \eqref{isomorphisme tau ! tensor}. L'assertion résulte alors de l'isomorphisme $(\mathscr{L}_{n})_{n\ge 1}\xrightarrow{\sim}( (\widetilde{\tau}_{!}(\mathscr{L}_{n}))^{\alpha})_{n\ge 1}$ (\ref{propriete adjoint alpha AMod}(iv)).


\begin{lemma} \label{equivalence final LLft aptf Q}
	Le foncteur canonique 
	\begin{equation}
		\Sysl(X_{\fet}^{\mathbb{N}^{\circ}},\breve{\oo})\to \Mod^{\aptf}(X_{\fet}^{\mathbb{N}^{\circ}},\breve{\oo})
		\label{foncteur canonique Sysl Modaptf}
	\end{equation}
	induit une équivalence de catégories 
	\begin{equation}
		\Sysl_{\mathbb{Q}}(X_{\fet}^{\mathbb{N}^{\circ}},\breve{\oo}) \xrightarrow{\sim} \mathbf{P}(\Mod^{\ptf}(X_{\fet},\oo_{\bullet}^{\alpha}))_{\mathbb{Q}}.
		\label{equivalence final}
	\end{equation}
\end{lemma}
\textit{Preuve}. D'après \ref{coro Rep aptf oo aptf oo Q}, \ref{prop oo aptf C rep}, \eqref{Rep C to Mod C} et \eqref{systeme proj Mod Rep}, le foncteur \eqref{foncteur canonique Sysl Modaptf}, qui correspond au foncteur canonique $\Rep_{\breve{\oo}}^{\ltf}(\pi_{1}(X,\overline{x}))\to \Rep_{\breve{\oo}}^{\aptf}(\pi_{1}(X,\overline{x}))$, induit une équivalence de catégories $\Sysl_{\mathbb{Q}}(X_{\fet}^{\mathbb{N}^{\circ}},\breve{\oo})$ $\xrightarrow{\sim} \Mod^{\aptf}_{\mathbb{Q}}(X_{\fet}^{\mathbb{N}^{\circ}},\breve{\oo})$, où le but désigne la catégorie des objets de $\Mod^{\aptf}(X_{\fet}^{\mathbb{N}^{\circ}},\breve{\oo})$ à isogénie près. L'équivalence de catégories $\Mod^{\aptf}_{\mathbb{Q}}(X_{\fet}^{\mathbb{N}^{\circ}},\breve{\oo})\simeq \mathbf{P}(\Mod^{\ptf}(X_{\fet},\oo_{\bullet}^{\alpha}))_{\mathbb{Q}}$ résulte de \ref{Mod A Q alpha Mod A Q} et \ref{lemma PMod essentielle}.

\section{Déformations des représentations et déformations des $\alpha$-modules de Faltings} \label{Dem finale}
\begin{nothing}	\label{basic setting alpha Faltings}
	Soient $X$ un $S$-schéma propre à réduction semi-stable de fibre générique géométrique connexe, $n$ un entier $\ge 1$. On note $(\widetilde{E},\overline{\mathscr{B}})$ (resp. $(\widetilde{E}_{s},\overline{\mathscr{B}}_{n})$) le topos annelé de Faltings (resp. la fibre spéciale du topos annelé de Faltings) associé au $S$-schéma $X$ \eqref{Es Bn topos Faltings}. On désigne par $\alpha\textnormal{-}X_{\overline{\eta},\fet}$ (resp. $\alpha\textnormal{-}\widetilde{E}_{s}$) la catégorie des $\alpha$-$\oo$-modules de $X_{\overline{\eta},\fet}$ (resp. $\widetilde{E}_{s}$) \eqref{faisceaux de alpha oo modules}. On pose $\oo_{n}^{\alpha}=\widetilde{\alpha}(\oo_{n})$, $\overline{\mathscr{B}}_{n}^{\alpha}=\widetilde{\alpha}(\overline{\mathscr{B}}_{n})$ \eqref{alpha Alg sheaf} et on désigne par $\Mod(X_{\overline{\eta},\fet},\oo_{n}^{\alpha})$ la catégorie des $\oo_{n}^{\alpha}$-modules de $\alpha\textnormal{-}X_{\overline{\eta},\fet}$, et par $\Mod(\overline{\mathscr{B}}_{n}^{\alpha})$ la catégorie des $\overline{\mathscr{B}}_{n}^{\alpha}$-modules de $\alpha$-$\widetilde{E}_{s}$ \eqref{alpha A Mod almost}.
	
	Comme $\overline{\mathscr{B}}_{n}$ est plat sur $\oo_{n}$ (\cite{AGT} III.9.2), le foncteur $\beta_{n}^{*}:\Mod(X_{\overline{\eta},\fet},\oo_{n})\to \Mod(\overline{\mathscr{B}}_{n})$ \eqref{topos beta n 2} est exact et il induit un foncteur exact (\ref{prop morphisme plat topos}) que l'on note
	\begin{equation}
		b_{n}:\Mod(X_{\overline{\eta},\fet},\oo_{n}^{\alpha})\to \Mod(\overline{\mathscr{B}}_{n}^{\alpha}).
		\label{foncteur bn}
	\end{equation}
	Pour tout entier $i \ge 0$ et tous $\oo_{n}^{\alpha}$-modules $\mathscr{L}$, $\mathscr{L}'$ de $X_{\overline{\eta},\fet}$, celui-ci induit un morphisme canonique \eqref{exact foncteur induit Ext morphisme}
	\begin{equation}
		\Ext^{i}_{\oo_{n}^{\alpha}}(\mathscr{L},\mathscr{L}')\to \Ext^{i}_{\overline{\mathscr{B}}_{n}^{\alpha}}(b_{n}(\mathscr{L}),b_{n}(\mathscr{L}')).
		\label{Ext fet to Ext B}
	\end{equation}
\end{nothing}
\begin{prop} \label{alpha iso to iso Ext}
	Soient $\mathbb{L}$ et $\mathbb{L}'$ deux $\oo_{n}$-modules localement libres de type fini de $X_{\overline{\eta},\fet}$. Posons $\mathscr{L}=\mathbb{L}^{\alpha}$ et $\mathscr{L}'=\mathbb{L}'^{\alpha}$. Alors le morphisme \eqref{Ext fet to Ext B} est un isomorphisme si $i=0,1$ et est un monomorphisme si $i=2$. 
\end{prop}
	\textit{Preuve}. Le foncteur exact $\beta_{n}^{*}$ induit, pour tout entier $i\ge 0$, un morphisme canonique \eqref{exact foncteur induit Ext morphisme}
	\begin{equation}
		\Ext^{i}_{\oo_{n}}(\mathbb{L},\mathbb{L}')\to \Ext^{i}_{\overline{\mathscr{B}}_{n}}(\beta_{n}^{*}(\mathbb{L}),\beta_{n}^{*}(\mathbb{L}')).
		\label{Ext fet to Ext B2}
	\end{equation}
	Comme $\mathbb{L}$ est un $\oo_{n}$-module localement libre de type fini de $X_{\fet}$, on a un isomorphisme canonique
	\begin{equation}
		\beta_{n}^{*}(\FHom_{\oo_{n}}(\mathbb{L},\mathbb{L}'))\xrightarrow{\sim} \FHom_{\overline{\mathscr{B}}_{n}}(\beta_{n}^{*}(\mathbb{L}),\beta_{n}^{*}(\mathbb{L}')).
		\label{beta n FHom pullback}
	\end{equation}
	D'après \eqref{Ext 1 iso to H1} et \eqref{beta n FHom pullback}, le morphisme \eqref{Ext fet to Ext B2} s'identifie au morphisme canonique induit par $\beta_{n}^{*}$
	\begin{equation}
		\rH^{i}(X_{\overline{\eta},\fet},\FHom_{\oo_{n}}(\mathbb{L},\mathbb{L}'))\to \rH^{i}(\widetilde{E}_{s},\beta_{n}^{*}(\FHom_{\oo_{n}}(\mathbb{L},\mathbb{L}'))).
		\label{morphisme Faltings FHom}
	\end{equation}
	On notera que $\FHom_{\oo_{n}}(\mathbb{L},\mathbb{L}')$ est un $\oo_{n}$-module localement libre de type fini de $X_{\overline{\eta},\fet}$. Le morphisme \eqref{morphisme Faltings FHom} est donc un $\alpha$-isomorphisme pour $i=0,1$, et est un $\alpha$-monomorphisme pour $i=2$ en vertu de \ref{presque iso faltings modifie n}. Par suite, il en est de même de \eqref{Ext fet to Ext B2}.

	Comme $\mathbb{L}$ est localement libre de type fini,  on a deux suites spectrales (cf. \ref{fonctorialite suite spectrale Gabber})
	\begin{eqnarray}
		\rE_{2}^{i,j}=\Ext^{i}_{\oo}(\mm,\Ext^{j}_{\oo_{n}}(\mathbb{L},\mathbb{L}')) &\Rightarrow& \rE^{i+j}=\Ext^{i+j}_{\oo_{n}^{\alpha}}(\mathscr{L},\mathscr{L}') \label{ss Ext oon}\\
		\rE_{2}'^{i,j}=\Ext^{i}_{\oo}(\mm,\Ext^{j}_{\overline{\mathscr{B}}_{n}}(\beta_{n}^{*}(\mathbb{L}),\beta_{n}^{*}(\mathbb{L}')))&\Rightarrow &\rE'^{i+j}=\Ext^{i+j}_{\overline{\mathscr{B}}_{n}^{\alpha}}(b_{n}(\mathscr{L}),b_{n}(\mathscr{L}')) \label{ss Ext Bn}
	\end{eqnarray}
	et les foncteurs exacts $\beta_{n}^{*}$ et $b_{n}$ induisent un morphisme de suites spectrales 
	\begin{equation}
		u=(u_{r}^{i,j},u^{m}):\rE\to \rE'.
	\end{equation}
	D'après \ref{prop Tsuji alpha iso}(iii)-(iv), le morphisme canonique induit par \eqref{Ext fet to Ext B2}
	\begin{equation}
		u_{2}^{i,j}: \Ext^{i}_{\oo}(\mm,\Ext^{j}_{\oo_{n}}(\mathbb{L},\mathbb{L}')) \to \Ext^{i}_{\oo}(\mm,\Ext^{j}_{\overline{\mathscr{B}}_{n}}(\beta_{n}^{*}(\mathbb{L}),\beta_{n}^{*}(\mathbb{L}'))),
		\label{u2ij ss oon Bn}
	\end{equation}
	est un isomorphisme si $j=0,1$ et est un monomorphisme si $i=0$ et $j=2$. D'après \ref{lemma suite spectrale}, le morphisme canonique $u^{m}$ \eqref{Ext fet to Ext B} est un isomorphisme si $m=0,1$; et est un monomorphisme si $m=2$.
\begin{prop}
	Soient $m$ un entier $\ge 1$ et $\mathscr{L}_{m}$ un objet de $\Mod^{\ptf}(X_{\overline{\eta},\fet},\oo_{m}^{\alpha})$ \eqref{notation Mod ptf oon alpha} tels que le morphisme \eqref{Ext fet to Ext B} induit par $b_{m}$
	\begin{equation}
		\Ext^{i}_{\oo_{m}^{\alpha}}(\mathscr{L}_{m},\mathscr{L}_{m}) \to \Ext^{i}_{\overline{\mathscr{B}}_{m}^{\alpha}}(b_{m}(\mathscr{L}_{m}),b_{m}(\mathscr{L}_{m}))
		\label{Exti oom to Bm alpha}
	\end{equation}
	soit un isomorphisme pour $i=0,1$ et soit un monomorphisme pour $i=2$. Supposons qu'il existe une déformation $M$ de $b_{m}(\mathscr{L}_{m})$ sur $\overline{\mathscr{B}}_{2m}^{\alpha}$ \textnormal{(\ref{pullback plat}(ii), \ref{notations deformation})}. Alors:

	\textnormal{(i)} Il existe une déformation $\mathscr{L}_{2m}$ de $\mathscr{L}_{m}$ sur $\oo_{2m}^{\alpha}$ telle que $b_{2m}(\mathscr{L}_{2m})$ soit isomorphe à $M$ en tant que déformation.
	
	\textnormal{(ii)} Le morphisme \eqref{Ext fet to Ext B} induit par $b_{2m}$
	\begin{equation}
		\Ext^{i}_{\oo_{2m}^{\alpha}}(\mathscr{L}_{2m},\mathscr{L}_{2m}) \to \Ext^{i}_{\overline{\mathscr{B}}_{2m}^{\alpha}}(b_{2m}(\mathscr{L}_{2m}),b_{2m}(\mathscr{L}_{2m}))
		\label{Exti oo2m to B2m alpha}
	\end{equation}
	est un isomorphisme pour $i=0,1$ et est un monomorphisme pour $i=2$.
	\label{def de rep et alpha def de modules}
\end{prop}

\textit{Preuve}. (i) On pose $I=p^{m}\oo/p^{2m}\oo$ et $I'=p^{m}\overline{\mathscr{B}}/p^{2m}\overline{\mathscr{B}}$. En vertu de \ref{foncteur pullback alpha psi}, le foncteur $b_{2m}$ induit un diagramme commutatif
\begin{eqnarray}
	&\xymatrix{
		0 \ar[r] & \Ext^{1}_{\oo_{m}^{\alpha}}(\mathscr{L}_{m},\mathscr{L}_{m}) \ar[r] \ar[d]_{b_{m}}& \Ext^{1}_{\oo_{2m}^{\alpha}}(\mathscr{L}_{m},\mathscr{L}_{m})\ar[d]^{b_{2m}} \\
		0 \ar[r] & \Ext^{1}_{\overline{\mathscr{B}}_{m}^{\alpha}}( b_{m}(\mathscr{L}_{m}),b_{m}(\mathscr{L}_{m})) \ar[r] & \Ext^{1}_{\overline{\mathscr{B}}_{2m}^{\alpha}}( b_{m}(\mathscr{L}_{m}),b_{m}(\mathscr{L}_{m}))}
	\label{diag com fondamental}\\
	&	\xymatrix{
		\Ext^{1}_{\oo_{2m}^{\alpha}}(\mathscr{L}_{m},\mathscr{L}_{m})\ar[r]^-{u_{1}} \ar[d]_{b_{2m}} & \Hom_{\oo_{m}^{\alpha}}(I^{\alpha}\otimes_{\oo_{m}^{\alpha}}\mathscr{L}_{m},\mathscr{L}_{m}) \ar[r]^-{\partial_{1}} \ar[d]^{b_{m}} & \Ext^{2}_{\oo_{m}^{\alpha}}(\mathscr{L}_{m},\mathscr{L}_{m}) \ar[d]^{b_{m}} \\
		\Ext^{1}_{\overline{\mathscr{B}}_{2m}^{\alpha}}(b_{m}(\mathscr{L}_{m}),b_{m}(\mathscr{L}_{m})) \ar[r]^-{u_{2}} & \Hom_{\overline{\mathscr{B}}_{m}^{\alpha}}(I'^{\alpha}\otimes_{\overline{\mathscr{B}}_{m}^{\alpha}}b_{m}(\mathscr{L}_{m}),b_{m}(\mathscr{L}_{m})) \ar[r]^-{\partial_{2}} & \Ext^{2}_{\overline{\mathscr{B}}_{m}^{\alpha}}(b_{m}(\mathscr{L}_{m}),b_{m}(\mathscr{L}_{m}))
		}\nonumber
\end{eqnarray}
où la première (resp. seconde) ligne est la suite exacte \eqref{suite exacte deformation alpha} associée aux $\oo_{m}^{\alpha}$-modules $\mathscr{L}_{m}$ et $\mathscr{L}_{m}$ (resp. $\overline{\mathscr{B}}_{m}^{\alpha}$-modules $b_{m}(\mathscr{L}_{m})$ et $b_{m}(\mathscr{L}_{m})$). Notons $w_{1}$ (resp. $w_{2}$) l'isomorphisme \eqref{u Mn canonique} associé au $\oo_{m}^{\alpha}$-module plat $\mathscr{L}_{m}$ de $X_{\overline{\eta},\fet}$ (resp. $\overline{\mathscr{B}}_{m}^{\alpha}$-module plat $b_{m}(\mathscr{L}_{m})$). On a alors $b_{m}(w_{1})=w_{2}$.

Comme il existe une déformation $M$ de $b_{m}(\mathscr{L}_{m})$ sur $\overline{\mathscr{B}}_{2m}^{\alpha}$, $\partial_{2}(w_{2})$ est nul d'après \ref{theorem principal deformation alpha modules}(i). La dernière flèche verticale de \eqref{diag com fondamental} est injective. Par suite $\partial_{1}(w_{1})=0$. Compte tenu \ref{theorem principal deformation alpha modules}(i), il existe une déformation de $\mathscr{L}_{m}$ sur $\oo_{2m}^{\alpha}$. En vertu de \ref{theorem principal deformation alpha modules}(ii), l'assertion résulte du fait que la première flèche verticale du diagramme \eqref{diag com fondamental} est un isomorphisme.

(ii) Posons $\mathbb{L}_{m}=\tau_{!}(\mathscr{L}_{m})$ et $\mathbb{L}_{2m}=\tau_{!}(\mathscr{L}_{2m})$. D'après \eqref{diagramme tau ! B B'} et \eqref{diagramme commutatif Ext phi tau alpha}, l'hypothèse revient à dire que le morphisme induit par $\beta_{m}^{*}$
\begin{equation}
	\Ext^{i}_{\oo_{m}}(\mathbb{L}_{m},\mathbb{L}_{m})\to \Ext^{i}_{\overline{\mathscr{B}}_{m}}(\beta_{m}^{*}(\mathbb{L}_{m}),\beta_{m}^{*}(\mathbb{L}_{m}))
	\label{hypothese Ext i oom}
\end{equation}
est un isomorphisme pour $i=0,1$ et est un monomorphisme pour $i=2$. De même, il suffit de démontrer la même propriété pour le morphisme induit par $\beta_{2m}^{*}$
\begin{equation}
	\Ext^{i}_{\oo_{2m}}(\mathbb{L}_{2m},\mathbb{L}_{2m})\to \Ext^{i}_{\overline{\mathscr{B}}_{2m}}(\beta_{2m}^{*}(\mathbb{L}_{2m}),\beta_{2m}^{*}(\mathbb{L}_{2m})).
	\label{Ext oo2m B2m}
\end{equation}
Comme $\tau_{!}$ est exact, $\mathbb{L}_{2m}$ est une extension de $\mathbb{L}_{m}$ par $\mathbb{L}_{m}$. Par dévissage, il suffit de démontrer que le morphisme induit par $\beta_{2m}^{*}$ 
\begin{equation}
	\Ext^{i}_{\oo_{2m}}(\mathbb{L}_{m},\mathbb{L}_{m})\to \Ext^{i}_{\overline{\mathscr{B}}_{2m}}(\beta_{m}^{*}(\mathbb{L}_{m}),\beta_{m}^{*}(\mathbb{L}_{m}))
	\label{Ext oo2m B2m Lm}
\end{equation}
est un isomorphisme pour $i=0,1$, et est un monomorphisme pour $i=2$. Rappelons qu'on a deux suites spectrales \eqref{phi exact fonctoriel def}
\begin{eqnarray}
	\rE^{i,j}_{2}=\Ext^{i}_{\oo_{m}}(\FTor_{j}^{\oo_{2m}}(\mathbb{L}_{m},\oo_{m}),\mathbb{L}_{m})&\Rightarrow& \rE^{i+j}=\Ext^{i+j}_{\oo_{2m}}(\mathbb{L}_{m},\mathbb{L}_{m}),
	\label{ss Ext Tor oom}\\
	\rE'^{i,j}_{2}=\Ext^{i}_{\overline{\mathscr{B}}_{m}}(\FTor_{j}^{\overline{\mathscr{B}}_{2m}}(\beta_{m}^{*}(\mathbb{L}_{m}),\overline{\mathscr{B}}_{m}),\beta_{m}^{*}(\mathbb{L}_{m}))&\Rightarrow& \rE'^{i+j}=\Ext^{i+j}_{\overline{\mathscr{B}}_{2m}}(\beta_{m}^{*}(\mathbb{L}_{m}),\beta_{m}^{*}(\mathbb{L}_{m})).
	\label{ss Ext Tor Bm}
\end{eqnarray}
et que le foncteur $\beta_{2m}^{*}$ induit un morphisme de suites spectrales
\begin{equation}
	u=(u_{r}^{i,j},u^{n}):\rE\to \rE'.
\end{equation}

Comme $I\simeq \oo_{m}$, $I'\simeq \overline{\mathscr{B}}_{m}$, on en déduit, pour tout entier $j\ge 2$, des isomorphismes 
\begin{equation}
	\FTor_{j}^{\oo_{2m}}(\mathbb{L}_{m},\oo_{m})\simeq \FTor_{j-1}^{\oo_{2m}}(\mathbb{L}_{m},\oo_{m}),\qquad \FTor_{j}^{\overline{\mathscr{B}}_{2m}}(\beta_{m}^{*}(\mathbb{L}_{m}),\overline{\mathscr{B}}_{m})\simeq \FTor_{j-1}^{\overline{\mathscr{B}}_{2m}}(\beta_{m}^{*}(\mathbb{L}_{m}),\overline{\mathscr{B}}_{m})
\end{equation}
Pour $j=1$, on a $\FTor_{1}^{\oo_{2m}}(\mathbb{L}_{m},\oo_{m})\simeq \mathbb{L}_{m}$ et $\FTor_{1}^{\overline{\mathscr{B}}_{2m}}(\beta_{m}^{*}(\mathbb{L}_{m}),\overline{\mathscr{B}}_{m})\simeq \beta_{m}^{*}(\mathbb{L}_{m})$. Le morphisme $u^{i,j}_{2}:\rE^{i,j}_{2}\to \rE'^{i,j}_{2}$ s'identifie alors au morphisme \eqref{hypothese Ext i oom}. D'après \ref{lemma suite spectrale} et l'hypothèse, le morphisme $u^{i}$ \eqref{Ext oo2m B2m Lm} est un isomorphisme si $i=0,1$ et est un monomorphisme si $i=2$; d'où l'assertion.

	\begin{theorem}
		Soient $C$ une courbe propre et lisse sur $\overline{K}$, $\check{C}=C\otimes_{\overline{K}}\mathfrak{C}$ et $F$ un fibré vectoriel de Deninger-Werner sur $\check{C}$ \textnormal{(\ref{def of categorie DW rat}(i))}. Alors, il existe un trait $S'$ fini sur $S$, un $S'$-modèle semi-stable et régulier $X$ de $C$ \eqref{Notations C X LPft Vect} et un fibré vectoriel de Deninger-Werner et de Weil-Tate $\mathcal{F}$ sur $\check{\overline{X}}=X\times_{S'}\check{\overline{S}}$ tel que $\mathcal{F}_{\check{\overline{\eta}}}\simeq F$ \textnormal{(\ref{def of categorie DW}, \ref{FV de WT sur X})}.
		\label{DW implique WT pf}
	\end{theorem}
	\textit{Preuve}. On fixe un entier $n\ge 1$. D'après \ref{DW module pn admissible}, quitte à remplacer $S$ par une extension finie, il existe un $S$-modèle semi-stable et régulier $X$ de $C$ et un fibré vectoriel de Deninger-Werner $\mathcal{F}$ sur $\check{\overline{X}}$ de fibre générique $F$ tel qu'on ait un isomorphisme
	\begin{equation}
		\gamma_{n}:\sigma_{n}^{*}(\mathcal{F}_{n})\xrightarrow{\sim} \beta_{n}^{*}(\mathbb{V}_{n}(\mathcal{F})),
		\label{Bn associe proof}
	\end{equation}
	où $\mathbb{V}_{n}(\mathcal{F})$ est la $\oo_{n}$-représentation de $\pi_{1}(C,\overline{x})$ associée à $\mathcal{F}$ (\ref{l'action du groupe fondamentale}). On note $\mathbb{L}_{n}$ le $\oo_{n}$-module localement libre de type fini de $C_{\fet}$ associé à $\mathbb{V}_{n}(\mathcal{F})$ (cf. \eqref{Rep on to Mod on}) et on pose $\mathscr{L}_{n}=\mathbb{L}_{n}^{\alpha}$.

	On démontre par récurrence que pour tout entier $m\in \{2^{l}n|l\ge 0\}$ il existe un objet $\mathscr{L}_{m}$ de $\Mod^{\ptf}(C_{\fet},\oo_{m}^{\alpha})$ vérifiant l'hypothèse de \ref{def de rep et alpha def de modules} et un isomorphisme de $\overline{\mathscr{B}}_{m}^{\alpha}$-modules
	\begin{equation}
		\gamma_{m}: (\sigma_{m}^{*}(\mathcal{F}_{m}))^{\alpha}\xrightarrow{\sim} b_{m}(\mathscr{L}_{m})
		\label{alpha iso gamma m}
	\end{equation}
	tels que les isomorphismes $(\gamma_{m})$ soient compatibles. L'assertion pour $m=n$ résulte de \ref{alpha iso to iso Ext} et \eqref{Bn associe proof}. Supposons que l'assertion soit vraie pour $m\ge n$ et montrons-là pour $2m$.

	Le $\overline{\mathscr{B}}_{2m}$-module $\sigma_{2m}^{*}(\mathcal{F}_{2m})$ est une déformation de $\sigma_{m}^{*}(\mathcal{F}_{m})$ sur $\overline{\mathscr{B}}_{2m}$ \ref{lemma deformation equivalence}(i). On l'identifie à une déformation de $b_{m}(\mathscr{L}_{m})$ sur $\overline{\mathscr{B}}_{2m}^{\alpha}$ via le $\alpha$-isomorphisme $\gamma_{m}$ \eqref{alpha iso gamma m}. D'après \ref{def de rep et alpha def de modules}(i), il existe une déformation $\mathscr{L}_{2m}$ de $\mathscr{L}_{m}$ sur $\oo_{2m}^{\alpha}$ et un isomorphisme de $\overline{\mathscr{B}}_{2m}^{\alpha}$-modules
	\begin{equation}
		\gamma_{2m}: (\sigma_{2m}^{*}(\mathcal{F}_{2m}))^{\alpha}\xrightarrow{\sim} b_{2m}(\mathscr{L}_{2m})
	\end{equation}
	compatible avec $\gamma_{m}$. En vertu de \ref{coro deformation alpha plat}(ii) et \ref{def de rep et alpha def de modules}(ii), $\mathscr{L}_{2m}$ est un objet de $\Mod^{\ptf}(C_{\fet},\oo_{2m}^{\alpha})$ et il satisfait l'hypothèse de \ref{def de rep et alpha def de modules}.

	Pour tout entier $m\ge 1$, choisissons un entier $l\ge 0$ tel que $2^{l-1}n<m\le 2^{l}n$ et posons $\mathscr{L}_{m}=\mathscr{L}_{2^{l}n}\otimes_{\oo_{2^{l}n}^{\alpha}}\oo_{m}^{\alpha}$. On en déduit des isomorphismes compatibles
	\begin{equation}
		\gamma_{m}:(\sigma_{m}^{*}(\mathcal{F}_{m}))^{\alpha}\xrightarrow{\sim} b_{m}(\mathscr{L}_{m}).
		\label{alpha iso gamma m all}
	\end{equation}

	On désigne par $\alpha\textnormal{-}C_{\fet}^{\mathbb{N}^{\circ}}$ (resp. $\alpha\textnormal{-}\widetilde{E}_{s}^{\mathbb{N}^{\circ}}$) la catégorie des $\alpha$-$\oo$-modules de $C_{\fet}^{\mathbb{N}^{\circ}}$ (resp. $\widetilde{E}_{s}^{\mathbb{N}^{\circ}}$) \eqref{faisceaux de alpha oo modules}. On pose $\breve{\oo}^{\alpha}=\widetilde{\alpha}(\breve{\oo})$ (resp. $\breve{\overline{\mathscr{B}}}^{\alpha}=\widetilde{\alpha}(\breve{\overline{\mathscr{B}}})$) \eqref{alpha Alg sheaf} et on désigne par $\Mod(C_{\fet}^{\mathbb{N}^{\circ}},\breve{\oo}^{\alpha})$ la catégorie des $\breve{\oo}^{\alpha}$-modules de $\alpha\textnormal{-}C_{\fet}^{\mathbb{N}^{\circ}}$, et par $\Mod(\breve{\overline{\mathscr{B}}}^{\alpha})$ la catégorie des $\breve{\overline{\mathscr{B}}}^{\alpha}$-modules de $\alpha$-$\widetilde{E}_{s}^{\mathbb{N}^{\circ}}$ \eqref{alpha A Mod almost}. Posons $\breve{\mathcal{F}}=(\mathcal{F}_{m})_{m\ge 1}$ et $\mathscr{L}=(\mathscr{L}_{m})_{m\ge 1}\in \mathbf{P}(\Mod^{\ptf}(C_{\fet},\oo_{\bullet}^{\alpha}))$ \eqref{notations Mod aptf oon oo} que l'on considère aussi comme un objet de $\Mod(C_{\fet}^{\mathbb{N}^{\circ}},\breve{\oo}^{\alpha})$ (cf. \ref{equivalence de cat Proj systems alpha}). Le foncteur $\breve{\beta}^{*}$ induit un foncteur que l'on note \eqref{morphisme de topos pullback alpha} 
	\begin{equation}
		b: \Mod(C_{\fet}^{\mathbb{N}^{\circ}},\breve{\oo}^{\alpha})\to \Mod(\breve{\overline{\mathscr{B}}}^{\alpha}).
	\end{equation}
	Sous l'équivalence de catégories $\Mod(\breve{\overline{\mathscr{B}}}^{\alpha})\xrightarrow{\sim} \mathbf{P}(\Mod(\overline{\mathscr{B}}_{\bullet}))$ \eqref{foncteur Mod Aalpha Proj systems}, le $\breve{\overline{\mathscr{B}}}^{\alpha}$-module $b(\mathscr{L})$ correspond au système projectif $(b_{m}(\mathscr{L}_{m}))_{m\ge 1}$. En vertu de \ref{equivalence de cat Proj systems alpha}, les isomorphismes compatibles \eqref{alpha iso gamma m all} induisent un isomorphisme de $\breve{\overline{\mathscr{B}}}^{\alpha}$-modules
	\begin{equation} \label{gamma BQ associe alpha}
		\gamma: (\breve{\sigma}^{*}(\breve{\mathcal{F}}))^{\alpha}\xrightarrow{\sim} b(\mathscr{L}),
	\end{equation}
	où $\breve{\sigma}$ est le morphisme de topos annelés \eqref{morphisme de topos sigma limit}. En vertu de \ref{equivalence final LLft aptf Q}, il existe un $\breve{\oo}$-module localement libre de type fini $\mathbb{L}$ de $C_{\fet}^{\mathbb{N}^{\circ}}$ tel que $\mathbb{L}_{\mathbb{Q}}\simeq \mathscr{L}_{\mathbb{Q}}$ dans la catégorie $\Mod_{\mathbb{Q}}(C_{\fet}^{\mathbb{N}^{\circ}},\breve{\oo})$. En vertu de \ref{Mod A Q alpha Mod A Q}, on en déduit par \eqref{gamma BQ associe alpha} un isomorphisme de $\breve{\overline{\mathscr{B}}}_{\mathbb{Q}}$-modules
	\begin{equation}
		(\breve{\sigma}^{*}(\breve{\mathcal{F}}))_{\mathbb{Q}}\simeq (\breve{\beta}^{*}(\mathbb{L}))_{\mathbb{Q}},
	\end{equation}
d'où le théorème.


	\begin{coro}
	Tout fibré vectoriel de Deninger-Werner sur $\check{C}$ est de Weil-Tate \eqref{Weil-Tate}.
		\label{DW implique WT coro}
	\end{coro}

	\begin{prop}
		La restriction du foncteur $\mathscr{V}_{\check{C}}:\BB_{\check{C}}^{\WT}\to \Rep_{\mathfrak{C}}^{\cont}(\pi_{1}(C,\overline{x}))$ \eqref{varrho check C WT} à la sous-catégorie $\BB_{\check{C}}^{\DW}$ s'identifie au foncteur de Deninger-Werner $\mathbb{V}_{\check{C}}$ \eqref{rho C check}.
	\label{foncteur DW WT pf}
	\end{prop}
	\textit{Preuve}. Soient $F$ un fibré vectoriel de Deninger-Werner sur $\check{C}$. D'après \ref{DW implique WT pf}, quitte à remplacer $S$ par une extension finie, il existe un $S$-modèle semi-stable et régulier $X$ de $C$, un fibré vectoriel de Deninger-Werner et de Weil-Tate $\mathcal{F}$ sur $\check{\overline{X}}$ de fibre générique $F$. En vertu de \ref{compatible of 2 rho}, on en déduit un isomorphisme de $\Rep_{\mathfrak{C}}^{\cont}(\pi_{1}(C,\overline{x}))$:
	\begin{equation}
		\mathbb{V}_{\check{C}}(F)\xrightarrow{\sim} \mathscr{V}_{\check{C}}(F).
		\label{isomorphisme VDW VWT}
	\end{equation}

	Il reste à démontrer que celui-ci est fonctoriel en $F$. Soit $f:F\to F'$ un morphisme $\BB_{\check{C}}^{\DW}$. En vertu de \ref{semistable rcm coro} et \ref{DW implique WT pf}, quitte à remplacer $S$ par une extension finie, il existe un $S$-modèle semi-stable et régulier $X$ de $C$ et deux fibrés vectoriels de Deninger-Werner et de Weil-Tate $\mathcal{F}$ et $\mathcal{F}'$ sur $\check{\overline{X}}$ tels que $\mathcal{F}_{\check{\overline{\eta}}}\simeq F$ et $\mathcal{F}'_{\check{\overline{\eta}}}\simeq F'$. D'après \eqref{Hom OX Hom OC pf}, $f$ s'étend en un morphisme de $\mathcal{F}$ à $\mathcal{F}'$. La fonctorialité de \eqref{isomorphisme VDW VWT} résulte alors de \ref{compatible of 2 rho}.

\section{Appendice}
\begin{nothing}
	On se donne une catégorie abélienne $\mathscr{A}$ ayant suffisamment d'injectifs, un complexe borné inférieurement $K^{\bullet}$ de $\mathscr{A}$ et un objet $X$ de $\mathscr{A}$. Choisissons $L^{\bullet,\bullet}$ une résolution injective de Cartan-Eilenberg de $K^{\bullet}$. La suite spectrale d'hypercohomologie (\cite{EGA III} (0.11.4.3.2))
	\begin{equation}
		\rE^{i,j}_{2}=\Ext^{i}_{\mathscr{A}}(X,\rH^{j}(K^{\bullet}))\Rightarrow \rE^{i+j}=\rH^{i+j}(\RHom_{\mathbf{D}(\mathscr{A})}(X,K^{\bullet})), \label{ss hypercohomologie inj}
	\end{equation}
	est la suite spectrale du bi-complexe $\Hom_{\mathscr{A}}(X,L^{\bullet,\bullet})$ associée à la filtration canonique $F^{l}_{\rII}=(\sum_{i+j=n,l\ge i}$ $\Hom_{\mathscr{A}}(X,L^{\bullet,\bullet}))_{n\in \mathbb{Z}}$ (\cite{EGA III} 0.11.3.2). Rappelons cette dernière. Pour la page $0$, on a (\cite{Weib} 5.6.2)
	\begin{equation} \label{E0 app}
		\rE_{0}^{i,j}=\Hom_{\mathscr{A}}(X,L^{j,i}),
	\end{equation}
	et le morphisme $d_{0}^{i,j}:\rE_{0}^{i,j}\to \rE_{0}^{i,j+1}$ est induit par la différentielle $L^{j,i}\to L^{j+1,i}$. 

	Notons, pour tous entiers $i,l$, $\rZ^{l}_{\rI}(L^{\bullet,i})$ (resp. $\rB^{l}_{\rI}(L^{\bullet,i})$, resp. $\rH^{l}_{\rI}(L^{\bullet,i})$) le $l$-ème sous-groupe de cocycles (resp. cobords, resp. groupe de cohomologie) du complexe $(L^{j,i})_{j\in \mathbb{Z}}$. Le complexe $(\rZ^{l}_{\rI}(L^{\bullet,i}))_{i\ge 0}$ (resp. $(\rB^{l}_{\rI}(L^{\bullet,i}))_{i\ge 0}$, resp. $(\rH^{l}_{\rI}(L^{\bullet,i}))_{i\ge 0}$) est une résolution injective de $\rZ^{l}(K^{\bullet})$ (resp. $\rB^{l}(K^{\bullet})$, resp. $\rH^{l}(K^{\bullet})$) (cf. \cite{EGA III} 0.11.4.2). On a alors (\cite{Weib} 5.6.2)
	\begin{equation} \label{E1 app}
	\rE_{1}^{i,j}=\Hom_{\mathscr{A}}(X,\rH^{j}_{\rI}(L^{\bullet,i})).
\end{equation}
Le morphisme $d_{1}^{i,j}:\rE^{i,j}_{1}\to \rE^{i+1,j}_{1}$ est induite par la différentielle $\rH^{j}(L^{\bullet,i})\to \rH^{j}(L^{\bullet,i+1})$. On en déduit
\begin{equation} \label{E2 app}
	\rE_{2}^{i,j}=\Ext^{i}_{\mathscr{A}}(X,\rH^{j}(K^{\bullet})).
\end{equation}
\end{nothing}

\begin{nothing} \label{notations generales App}
	On se donne de plus une catégorie abélienne $\mathscr{A}'$ ayant suffisamment d'injectifs et un foncteur \textit{exact} $F:\mathscr{A}\to \mathscr{A}'$. Soient $K'^{\bullet}$ un complexe borné inférieurement de $\mathscr{A}'$ et $g:F(K^{\bullet})\to K'^{\bullet}$ un morphisme de complexes. On désigne par $\rE'$ la suite spectrale d'hypercohomologie
	\begin{equation}
		\rE'^{i,j}_{2}=\Ext^{i}_{\mathscr{A}'}(F(X),\rH^{j}(K'^{\bullet}))\Rightarrow \rE'^{i+j}=\rH^{i+j}(\RHom_{\mathbf{D}(\mathscr{A}')}(F(X),K'^{\bullet})).
	\end{equation}
	Choisissons une résolution injective de Cartan-Eilenberg $L'^{\bullet,\bullet}$ de $K'^{\bullet}$. On notera que $F(L^{\bullet,\bullet})$ est encore une résolution de Cartan-Eilenberg de $F(K^{\bullet})$ (mais pas nécessairement injective). Le morphisme $g$ s'étend en un morphisme de bi-complexes (\cite{EGA III} 0.11.4.2, cf. aussi \cite{CE56} XVII 1.2)
	\begin{equation}
		G:F(L^{\bullet,\bullet})\to L'^{\bullet,\bullet}.
		\label{morphisme de bicomplexes G}
	\end{equation}
	On en déduit par le foncteur $F$ un morphisme de bi-complexes 
	\begin{equation}
		\Hom_{\mathscr{A}}(X,L^{\bullet,\bullet})\to \Hom_{\mathscr{A}'}(F(X),L'^{\bullet,\bullet})
		\label{morphisme bicomplexes}
	\end{equation}
	et par suite un morphisme de suites spectrales (\cite{EGA III} 0.11.2.3)
\begin{equation}
	u=(u_{r}^{i,j},u^{n}):\rE\to \rE'.
	\label{ap u suites spectrales}
\end{equation}
\end{nothing}
\begin{nothing}
	Le morphisme $G$ \eqref{morphisme de bicomplexes G} induit au morphisme $(u_{0}^{i,j}:\rE_{0}^{i,j}\to \rE'^{i,j}_{0})_{i,j}$ de $u$ \eqref{E0 app} et il induit, pour tout entier $j$, un morphisme de complexes
	\begin{equation}
		(F(\rH_{\rI}^{j}(L^{\bullet,i})))_{i\ge 0}\to (\rH_{\rI}^{j}(L'^{\bullet,i}))_{i\ge 0}.
		\label{G induit KK}
	\end{equation}
	Le morphisme composé, induit par celui-ci et le foncteur exact $F$,
	\begin{equation}
		\Hom_{\mathscr{A}}(X,\rH_{\rI}^{j}(L^{\bullet,i}))\to \Hom_{\mathscr{A}'}(F(X),F(\rH_{\rI}^{j}(L^{\bullet,i}))) \to \Hom_{\mathscr{A}'}(F(X),\rH_{\rI}^{j}(L'^{\bullet,i}))
	\end{equation}
	s'identifie au morphisme $u_{1}^{i,j}$ de $u$ \eqref{E1 app}. Les morphismes $(u_{1}^{i,j})_{i,j}$ induisent le morphisme $u_{2}^{i,j}$ \eqref{E2 app}
	\begin{equation}
		\Ext^{i}_{\mathscr{A}}(X,\rH^{j}(K^{\bullet}))\to \Ext^{i}_{\mathscr{A}'}(F(X),\rH^{j}(K'^{\bullet}))
		\label{morphisme u2}
	\end{equation}
	En outre, le morphisme de complexes $g$ induit, pour tout entier $j$, un morphisme de groupes de cohomologie
	\begin{equation} \label{g induit KK}
		 F(\rH^{j}(K^{\bullet}))\simeq \rH^{j}(F(K^{\bullet}))\to \rH^{j}(K'^{\bullet}).
	\end{equation}
	Rappelons que le complexe $(\rH_{\rI}^{j}(L^{\bullet,i}))_{i\ge 0}$ (resp. $(\rH_{\rI}^{j}(L'^{\bullet,i}))_{i\ge 0}$) est une résolution injective de $\rH^{j}(K^{\bullet})$ (resp. $\rH^{j}(K'^{\bullet})$). Le morphisme de complexes \eqref{G induit KK} étend le morphisme \eqref{g induit KK}. Le morphisme $u_{2}^{i,j}$ \eqref{morphisme u2} s'identifie alors au morphisme composé
	\begin{equation}
		\Ext^{i}_{\mathscr{A}}(X,\rH^{j}(K^{\bullet}))\to \Ext^{i}_{\mathscr{A}'}(F(X),F(\rH^{j}(K^{\bullet})))\to \Ext^{i}_{\mathscr{A}'}(F(X),\rH^{j}(K'^{\bullet})),
	\end{equation}
	où la première flèche est le morphisme canonique \eqref{exact foncteur induit Ext morphisme} induit par le foncteur exact $F$ et la seconde flèche est induite par \eqref{g induit KK}.
\end{nothing}
\begin{nothing}
	Le morphisme $g$ et le foncteur exact $F$ induisent, pour tout entier $n$, un morphisme \eqref{exact foncteur induit Ext morphisme}
	\begin{equation}
		\rE^{n}=\Ext^{n}_{\mathbf{D}(\mathscr{A})}(X,K^{\bullet})\to \rE'^{n}=\Ext^{n}_{\mathbf{D}(\mathscr{A}')}(F(X),K'^{\bullet}),
	\end{equation}
	qui n'est autre que le morphisme $u^{n}$ de $u$.
\end{nothing}
\begin{nothing}
	Reprenons les notations de \ref{fonctorialite suite spectrale Gabber}. On considère le cas où $\mathscr{A}=\mathscr{A}'=\Mod(\oo)$ et on prend pour $F$ le foncteur identique et pour $X$ le $\oo$-module $\mm$. Choisissons $I^{\bullet}$ (resp. $I'^{\bullet}$) une résolution injective de $N$ (resp. $\varphi^{*}(N)$) et posons $K^{\bullet}=\Hom_{A}(M,I^{\bullet})$ et $K'^{\bullet}=\Hom_{A'}(\varphi^{*}(M),I'^{\bullet})$. On a alors un morphisme de complexes $\varphi^{*}(I^{\bullet})\to I'^{\bullet}$ et par suite un morphisme de complexes $g:K^{\bullet}\to K'^{\bullet}$. On notera que, dans ce cas, la suite spectrale \eqref{suite spectrale Gabber} s'identifie à la suite spectrale \eqref{ss hypercohomologie inj}. Appliquant ce qui précède, on en déduit le morphisme de suites spectrales $u$ \eqref{morphisme suite spectrales Gabber} et la description des morphismes $u_{2}^{i,j}$ et $u^{n}$.
	\label{explication ss Gabber}
\end{nothing}
\begin{nothing}
	On considère le cas dual de \ref{notations generales App}. On se donne deux catégories abéliennes $\mathscr{A}$ et $\mathscr{A}'$ ayant suffisamment d'objets projectifs, un foncteur \textit{exact} $F:\mathscr{A}\to \mathscr{A}'$ et $X$ un objet de $\mathscr{A}$. Soient $K^{\bullet}$ un complexe borné supérieurement de $\mathscr{A}$, $K'^{\bullet}$ un complexe borné supérieurement de $\mathscr{A}'$ et $g:K'^{\bullet}\to F(K^{\bullet})$ un morphisme de complexes. En remplaçant les résolutions injectives par les résolutions projectives, on construit un morphisme $u:\rE\to \rE'$ entre les deux suites spectrales d'hypercohomologie
	\begin{eqnarray}
		\rE^{i,j}_{2}=\Ext^{i}_{\mathscr{A}}(\rH^{j}(K^{\bullet}),X)&\Rightarrow& \rE^{i+j}=\rH^{i+j}(\RHom_{\mathbf{D}(\mathscr{A})}(K^{\bullet},X)),\label{ss hypercohomologie proj}\\
		\rE'^{i,j}_{2}=\Ext^{i}_{\mathscr{A}'}(\rH^{j}(K'^{\bullet}),F(X))&\Rightarrow& \rE'^{i+j}=\rH^{i+j}(\RHom_{\mathbf{D}(\mathscr{A}')}(K'^{\bullet},F(X))).
	\end{eqnarray}
	De même, les morphismes $u^{i,j}_{2}$ et $u^{n}$ s'identifient aux morphismes canoniques induits par le morphisme $g$ et le foncteur exact $F$.
\end{nothing}
\begin{nothing} \label{explication ss Grothendieck}
	Reprenons les notations de \ref{phi exact fonctoriel def}. On considère le cas où $\mathscr{A}=\Mod(A)$ et $\mathscr{A}'=\Mod(A')$ et on prend pour $F$ le foncteur exact $\varphi^{*}$ et pour $X$ le $A$-module $J$. Choisissons $P^{\bullet}\to A_{0}$ (resp. $P'^{\bullet}\to A_{0}'$) une résolution projective de $A$-modules (resp. $A'$-modules) et posons $K^{\bullet}=M_{0}\otimes_{A}P^{\bullet}$ et $K'^{\bullet}=\varphi^{*}(M_{0})\otimes_{A'}P'^{\bullet}$. Comme $A_{0}'=\varphi^{*}(A_{0})$, on a alors un morphisme de complexes $P'^{\bullet}\to \varphi^{*}(P^{\bullet})$. Il induit un morphisme de complexes $g:K'^{\bullet}\to \varphi^{*}(K^{\bullet})$. On notera que, dans ce cas, la suite spectrale \eqref{suite spectrale Cartan Grothendieck} s'identifie à la suite spectrale \eqref{ss hypercohomologie proj}. Appliquant ce qui précède, on en déduit le morphisme de suites spectrales $u$ \eqref{morphisme de ss Grothendieck Ill} et la description des morphismes $u_{2}^{i,j}$ et $u^{n}$.
\end{nothing}


\begin{thebibliography}{50}
	\bibitem{Ab00} A. Abbes, \textit{Réduction semi-stable des courbes d'après Artin, Deligne, Grothendieck, Mumford, Saito, Winters...}, dans "Courbes semi-stables et groupe fondamental en géométrie algébrique", J.-B. Bost, F. Loeser et M. Raynaud (éditeurs), Progress in Mathematics 187, Birkhäuser (2000), 59-110.
	\bibitem{Ab10} A. Abbes, \textit{Éléments de Géométrie Rigide: Volume I. Construction et étude géométrique des espaces rigides}. Springer, 2010.	
	\bibitem{AG15} A. Abbes, M. Gros, \textit{La suite spectrale de Hodge-Tate}, prépublication (2015), \href{http://arxiv.org/abs/1509.03617}{arXiv:1509.03617}.
	\bibitem{AGT} A. Abbes, M. Gros, T. Tsuji, \textit{The $p$-adic Simpson correspondence}, to appear in Annals of Mathematics Studies, Princeton University Press.
	\bibitem{Ach14} P. Achinger. \textit{$K (\pi, 1)$-neighborhoods and comparison theorems}; to appear in Compositio Mathematica.
	\bibitem{SGAIV} M. Artin, A. Grothendieck, J. L. Verdier, \textit{Théorie des topos et cohomologie étale des schémas}, SGA 4, Springer-Verlag, Tome 1, LNM 269 (1972) ; Tome 2, LNM 270 (1972) ; Tome 3, LNM 305 (1973).
	\bibitem{SGAVI} P. Berthelot, A. Grothendieck, L. Illusie, \textit{Théorie des intersections et théorème de Riemann-Roch}, SGA 6, LNM 225, Springer-Verlag (1971).
	\bibitem{BLR90} S. Bosch, W. Lütkebohmert, M. Raynaud. \textit{Néron models}. Springer, (1990).
	\bibitem{CE56} H. Cartan et S. Eilenberg. \textit{Homological algebra}. Princeton University (1956).
	\bibitem{BouAlgcom} N. Bourbaki, \textit{Algèbre commutative}, Chapitres 1-9, Hermann (1985).
	\bibitem{SGAIII} M. Demazure, A. Grothendieck. \textit{Schémas en groupes}, SGA3. LNM 151-153, Springer-Verlag, (1970).
	\bibitem{DM82} P. Deligne P, J. S. Milne. \textit{Tannakian categories} Hodge cycles, motives, and Shimura varieties. Springer Berlin Heidelberg, 1981: 101-228.
	\bibitem{DW03} C. Deninger, A. Werner. \textit{Vector bundles and $p$-adic representations I}; preprint (2003), \href{http://arxiv.org/abs/math/0309273}{arXiv:0309.273}.
	\bibitem{DW05} C. Deninger, A. Werner. \textit{Vector bundles on $p$-adic curves and parallel transport}; Annales scientifiques de l'Ecole normale supérieure. Vol. 38. No. 4 (2005): 553-597.
	\bibitem{SGAI} A. Grothendieck, \textit{Revêtements étales et groupe fondamental}, SGA 1, Lecture Notes in Mathematics 224, Springer-Verlag (1971).
	\bibitem{EGAInew} A. Grothendieck, J.A. Dieudonné, \textit{Éléments de Géométrie Algébrique I}, Seconde édition, Springer-Verlag (1971).	
	\bibitem{EGA I} A. Grothendieck, J. Dieudonné: \textit{\'{E}léments de Géométrie Algébrique I}. Publ. Math. IHES 4 (1960).
	\bibitem{EGA III} A. Grothendieck et J. Dieudonné, \textit{\'{E}léments de Géométrie Algébrique, III \'{E}tude cohomologique des faisceaux cohérents}. Publ. Math. IHES 11 (1961), 17 (1963).
	\bibitem{EGA IV} A. Grothendieck et J. Dieudonné, \textit{Éléments de Géométrie Algébrique, IV Étude locale des schémas et des morphismes de schémas}. Pub. Math. IHES 20 (1964), 24 (1965), 28 (1966), 32 (1967).
	\bibitem{Epp73} H.P. Epp. \textit{Eliminating wild ramification}; Inventiones mathematicae 19.3 (1973): 235-249.
	\bibitem{Fal87} G. Faltings, \textit{$p$-adic Hodge theory}, J. Amer. Math. Soc. 1 (1988), 255-299.
	\bibitem{Fal02} G. Faltings, \textit{Almost étale extensions}, dans Cohomologies p-adiques et applications arithmétiques. II, Astérisque 279 (2002), 185-270.
	\bibitem{Fal05} G. Faltings, \textit{A $p$-adic Simpson correspondence}, Adv. Math. 198 (2005), 847-862.
	\bibitem{Ga62} P. Gabriel, \textit{Des catégories abéliennes}, Bulletin de la Société Mathématique de France 90 (1962): 323-448.
	\bibitem{Gi71} J. Giraud, \textit{Cohomologie non abélienne}, Springer-Verlag (1971).
	\bibitem{Gi73} D. Gieseker, \textit{Stable vector bundles and the frobenius morphism}; Annales scientifiques de l'École Normale Supérieure 6.1 (1973): 95-101.
	\bibitem{GR03} O. Gabber, L. Ramero, \textit{Almost Ring Theory}; Lecture Notes in Mathematics 1800, Springer (2003).
	\bibitem{GR14} O. Gabber, L. Ramero, \textit{Foundations for Almost Ring Theory} Release 6.8; preprint, \href{http://arxiv.org/abs/math/0409584v9}{arXiv:math/0409584}, (2004).
	\bibitem{Il71} L. Illusie, \textit{Complexe cotangente et déformations I}. Springer Lect. Notes Math. 239 (1971).
	\bibitem{Il14} L. Illusie, \textit{Produits orientés}, dans \textit{Travaux de Gabber sur l'uniformisation locale et la cohomologie étale des schémas quasi-excellents}, Séminaire à l’École Polytechnique dirigé par L. Illusie, Y. Laszlo et F. Orgogozo; Astérisque 363-364 (2014).
	\bibitem{Ka89} K. Kato, \textit{Logarithmic structures of Fontaine-Illusie}; Algebraic analysis, geometry, and number theory (Baltimore, MD, 1988), (1989), 191-224.
	\bibitem{Lich68} S. Lichtenbaum, \textit{Curves over discrete valuation rings}; Am. J. Math. 90 (1968), 380–405.
	\bibitem{Liu} Q. Liu. \textit{Algebraic Geometry and Arithmetic Curves}; Oxford: Oxford university press, (2002).
	\bibitem{Liu06} Q. Liu. \textit{Stable reduction of finite covers of curves}; Compositio Mathematica 142.01 (2006): 101-118.
	\bibitem{MB85} L. Moret-Bailly, \textit{Métriques permises}; Séminaire sur les pinceaux arithmétiques: la conjecture de Mordell, Astérisque, 127 (1985), 29-88.
	\bibitem{Milne} J.S Milne. \textit{Etale Cohomology} (PMS-33). No. 33. Princeton university press, 1980.
	\bibitem{Mu70} D. Mumford, Abelian variety, Oxford University Press (1970).
	\bibitem{NS65} M.S. Narasimhan, C.S. Seshadri, \textit{Stable and unitary vector bundles on a compact Riemann surface}, Ann. of Math. 82 (1965), 540-567.	
	\bibitem{Ray70} M. Raynaud, \textit{Spécialisation du foncteur de Picard}; Publications Mathematiques de l'IHES, 1970, 38(1): 27-76.
	\bibitem{Ray90} M. Raynaud, \textit{p-groupes et réduction semi-stable des courbes}, in {The Grothendieck Festschrift, vol. III}, eds P. Cartier et al., Progress in Mathematics, vol. 88 (Birkh\"{a}user, Basel, 1990), 179-197.
	\bibitem{RG71} M. Raynaud, L. Gruson, \textit{Critères de platitude et de projectivité}; Inventiones mathematicae 13.1 (1971): 1-89.
	\bibitem{Ro06} J.-E. Roos, \textit{Derived functors of inverse limits revisited}, J. London Math. Soc. (2) 73 (2006), 65-83.
	\bibitem{Sa04} T. Saito, \textit{Log smooth extension of a family of curves and semi-stable reduction}; Journal of Algebraic Geometry, 2004, 13(2): 287-322.
	\bibitem{Sch13} P. Scholze, \textit{$p$-adic Hodge theory for rigid-analytic varieties}; Forum of Mathematics, Pi. Cambridge University Press, 2013, 1: e1.
	\bibitem{Se75} J-P. Serre, \textit{Groupes algébriques et corps de classes}. Hermann, 1975.
	\bibitem{Sim92} C. Simpson, \textit{Higgs bundles and local systems}, Pub. Math. IHÉS, 75 (1992), 5-95.
	\bibitem{Ta67} J. Tate, \textit{$p$-divisible groups}; Proc. Conf. Local Fields (Driebergen, 1966), Springer, Berlin, 158-183 (1967).
	\bibitem{To07} J. Tong, \textit{Application d'Albanese pour les courbes et contractions}; Mathematische Annalen 338 (2007): 405–420.
	
	\bibitem{Weib} C.A. Weibel. \textit{An introduction to homological algebra}; Cambridge university press, 1995.
	
	\bibitem{Stacks} The Stacks Project Authors. \textit{Stacks Project}. \href{http://stacks.math.columbia.edu}{http://stacks.math.columbia.edu}, 2015.
\end{thebibliography}
\end{document}